\documentclass[a4paper,11pt]{book} 

\usepackage{vmargin}
\usepackage[Lenny]{fncychap}
\usepackage{fancyhdr}
\usepackage[pagebackref = true, hyperindex = true, frenchlinks = true,
colorlinks = true, citecolor = NavyBlue, linkcolor = Mulberry,
urlcolor = LimeGreen, linktocpage]{hyperref}
\usepackage[nottoc,notlof,notlot]{tocbibind} 

\newcommand{\Auteur}{Robin Langer}
\newcommand{\Titre}{Cylindric plane partitions, Lambda determinants, Commutants in semicircular systems}
\newcommand{\Date}{6 Decembre 2013}
\author{\Auteur}
\title{\Titre}
\date{\today}

\usepackage{amssymb}
\usepackage{latexsym}
\usepackage{amsfonts}
\usepackage[all]{xy}

\usepackage{amsmath}
\usepackage{amsthm}
\usepackage[usenames,dvipsnames]{xcolor}
\usepackage{tikz}

\usepackage{makeidx}

\setmarginsrb{3.2cm}{2cm}{3.2cm}{2.5cm}{1cm}{1cm}{1cm}{1cm}

\makeindex

\newtheorem{lemma}{Lemma}[section]
\newtheorem{theorem}{Theorem}[section]
\newtheorem{definition}{Definition}[section]
\newtheorem{proposition}{Proposition}[section]
\newtheorem{remark}{Remark}[section]
\newtheorem{corollary}{Corollary}[section]

\DeclareMathOperator{\Poly}{Poly}
\DeclareMathOperator{\Inv}{Inv}
\DeclareMathOperator{\Dinv}{Dinv}

\DeclareMathOperator{\GL}{GL}
\DeclareMathOperator{\img}{Image}
\DeclareMathOperator{\alg}{Alg}
\DeclareMathOperator{\comm}{Comm}

\DeclareMathOperator{\diag}{diag}

\DeclareMathOperator{\cols}{columns}

\DeclareMathOperator{\RPP}{RPP}
\DeclareMathOperator{\CPP}{CPP}
\DeclareMathOperator{\PP}{PP}

\DeclareMathOperator{\bin}{Bin}
\DeclareMathOperator{\inv}{inv}
\DeclareMathOperator{\dinv}{dinv}

\DeclareMathOperator{\lab}{label}

\DeclareMathOperator{\ALCD}{ALCD}
\DeclareMathOperator{\IP}{P}

\DeclareMathOperator{\hook}{hook}
\DeclareMathOperator{\cohook}{cohook}

\DeclareMathOperator{\tr}{tr}

\DeclareMathOperator{\sign}{sign}

\DeclareMathOperator{\lhs}{lhs}

    \makeatletter
    \renewcommand{\@makechapterhead}[1]{%
    \vspace*{50 pt}%
    {\setlength{\parindent}{0pt} \raggedright \normalfont
    \bfseries\Huge
    \ifnum \value{secnumdepth}>1
       \if@mainmatter\thechapter.\ \fi%
    \fi
    #1\par\nobreak\vspace{40 pt}}}
    \makeatother

\begin{document}

\title{Cylindric plane partitions, Lambda determinants, Commutants in semicircular systems}
\author{Robin Langer}

\frontmatter
\lhead[\oldstylenums \thepage]{\rightmark}
\rhead[\leftmark]{\oldstylenums \thepage}


\thispagestyle{empty}

\vspace*{\fill}

\begin{center}
    \huge
    \textsc{Th\`ese de doctorat}
\end{center}

\begin{center}
    \large
    pour l'obtention du grade de
\end{center}

\begin{center}
    \Large
    Docteur de l'universit\'e Paris-Est
\end{center}

\begin{center}
    \large
    \textbf{Sp\'ecialit\'e Informatique}
\end{center}

\begin{center}
    \large
    \textit{au titre de l'\'Ecole Doctorale de Math\'ematiques \\
    et des Sciences et Techniques de l'Information et de la Communication}
\end{center}

\vspace{3em}

\begin{center}
    \large
    Pr\'esent\'ee et soutenue publiquement par
\end{center}

\begin{center}
    \Large
    \Auteur
\end{center}

\begin{center}
    \large
    le \Date
\end{center}

\vspace{3em}

\begin{center}
    \Large
    \textbf{\Titre}
\end{center}

\vspace{3em}

\begin{center}
    \large
    \textbf{Devant le jury compos\'e par}
\end{center}

\vspace{1em}

\begin{center}
    \Large
    \begin{tabular}{ll}
        Piotr \'Sniady        \qquad & \qquad Rapporteur \\
        Jean-Christophe Aval           \qquad & \qquad Rapporteur \\
        Philippe Biane          \qquad & \qquad Directeur de th\'ese \\
        Jean-Yves Thibon            \qquad & \qquad Examinateur \\
        Florent Hivert  \qquad & \qquad Examinateur \\
        Marc Van Leeuwen           \qquad & \qquad Examinateur \\
        Jiang Zeng        \qquad & \qquad Examinateur \\
    \end{tabular}
\end{center}

\vspace{2.1em}

\begin{center}
    \footnotesize
    Laboratoire d'informatique Gaspard-Monge \\
    UMR 8049 LIGM \\
    5, bd Descartes, Champs-sur-Marne, 77454 Marne-la-Vall\'ee Cedex 2, France
\end{center}

\vspace*{\fill}

\cleardoublepage

\lhead[\oldstylenums \thepage]{}
\rhead[Acknowledgments]{\oldstylenums \thepage}
\section*{Acknowledgments}

The person to whom I owe the greatest debt is {\color{blue}Philippe Biane}. Firstly for accepting me as his student and secondly for displaying unending patience in the face of my sometimes erratic work style. Philippe has always had interesting ideas for problems to work on and great skill at explaining the background context of such problems. He always had plenty of time for me and was available for all my questions. Having such a wide range of knowledge in many diverse fields of mathematics, being able to talk to Philippe was like having access to a walking talking library. 

The second person to whom I owe a great mathematical debt is {\color{blue}Paul Zinn-Justin}. His high speed crash courses on integrability, representation theory and geometry were not always easy to follow - but I learned a lot all the same. In particular he helped me to fill in many gaps that I had in basic linear algebra from my undergraduate years.
 
I could not have asked for a better office mate than {\color{blue} Mathieu Josuat-Verges}. He would often let me think aloud when I was ``stuck'' and just by asking a simple question to clarify what I was trying to explain, it would suddenly become clear to me what the problem was.

I would like to thank EVERYBODY in the {\color{blue}combinatorics group at Marne-la-Vallee} for having such a team spirit and for creating an environment in which maths is cool, and solving problems is fun and doing research is not some horrible chore which you're only doing for the money [which is crap anyway] and would really rather be sitting on the beach drinking martinis [been there, done that, it gets boring]. Seriously folks, if you're willing to rip each other's throats out over tiny little promotions and grants then you should go into finance. Its much more profitable.

I would especially like to thank {\color{blue}Jean-Christophe Novelli} and {\color{blue}Jean-Yves Thibon} for organizing the weekly ``group de travail'' as well as the great course on ``combinatorial Hopf algebras''. {\color{blue}Alain Lascoux} will be sorely missed for his regular cryptic remarks which often turn out to be quite profound and just the ``trick'' you were looking for to re-frame a difficult problem.

To other thesis students in my laboratory (some now graduated): Samuele Giraudo, Remi Maurice, Vincent Vong, Gregory Chatel, Olivier Bouillot, Hayat Cheballah and Vivian Pons - its great to have such a large group of students all working on closely related subjects.

The Marne-la-Vallee secretaries: Corrine Palescandolo, Angelique Crombez and Sylvie Cach have always been very friendly and helpful. The same goes for Noëlle Delgado, the secretary of LIAFA.

I would like to thank {\color{blue} Guillaume Chapuy} for organizing two very fruitful ``group de travail'', the first on Macdonald polynomials and the second on the KP hierarchy. I would also like to thank my former office mates Adeline Pierrot and Ir\'ene Marcovici for defending my people against a modern day blood libel which is sadly still alive today in France even amongst supposed ``enlightened progressives'', and even within the French mathematics community.

I must thank {\color{blue} Sylvie Corteel} for suggesting the subject of cylindric plane partitions and for passing me a pre-print of {\color{blue} Shingo Adachi} which motivated my work in Chapter 4. Thanks also to {\color{blue} Xavier Viennot} who explained to me the essentials of Chapter 2 at the Calanques de Marseille.

I would very much like to thank {\color{blue} Ekaterina Vassilieva} for adding some much needed class to the French combinatorial community, and for showing me around to the best sites in Paris which I might not otherwise have discovered alone.

I am grateful to {\color{blue} Olga Azenhas} for organizing a great summer school in Portugal on ``Algebraic and Enumerative combinatorics'' where I was given the opportunity to present my work to an international audience. Likewise I am grateful to {\color{blue} Masao Ishikawa} for organizing the workshop ``Algebraic Combinatorics related to Young diagrams and Statistical Physics'' in Kyoto Japan, where again I was able to present my work to an international audience.

Both the rapporteurs, {\color{blue}Piotr Sniady} and {\color{blue}Jean-Christophe Aval} were very generous in taking the time to carefully read my manuscript. I am especially grateful to Piotr Sniady for pointing out an error in the proof of Theorem 6.5.1 which fortunately I was able to fix.
The remainder of the Jury: {\color{blue}Jiang Zeng}, {\color{blue} Marc van Leeuwen}, {\color{blue} Jean-Yves Thibon} and {\color{blue} Florent Hivert} also deserve my appreciation for taking time out of their own research in order to attend my defense.

Finally I'd like to thank my housemate Perrine Dubertrand for putting up with me not doing any housework while I was writing up, Deborah and Bernard Camain for warmly welcoming me into their home for shabbat dinners, my (ex)-husband Denchai Rotsapha for supporting my ambitions to go back to school (instead of having a baby), my mother Kerry Miller for always placing a high value on education, my father Albert Langer for teaching me about the halting problem and Godel's theorem on the back of a napkin when I was still a kid, my sister Wendy Langer for sharing with me her passion for physics and finally Mrs McLaren who was the awesomest high school teacher ever, and probably the only reason I ever managed to finish high school at all [behavioural problems, meh].

This thesis would not have been possible without the financial assistance of CNRS, University of Marne-la-Vallee and most importantly, Smorgan Steeel. Nor would this thesis have been possible were it not for the Her Majesty The Queen of England's personal request that the ``governors of the province beyond the river'' grant me permission to pass freely without let or hindrance and to afford me such assistance and protection as may be necessary. Baruch Hashem.

\cleardoublepage

\lhead[\oldstylenums \thepage]{}
\rhead[Résumé]{\oldstylenums \thepage}
\begin{center}
    {\bf \Titre}
\end{center}

\section*{\normalsize Resum\'e}

Cette th\`ese se compose de trois parties. La premi\`ere partie est consacr\'ee aux partitions planes cylindriques, la deuxi\`eme aux lambda-d\'eterminants et enfin la troisi\`eme aux commutateurs dans les syst\`eemes semi-circulaires.

La classe des  partitions planes cylindriques est une g\'en\'eralisation naturelle de celle des partitions planes inverses. Borodin a donn\'e r\'ecemment une s\'erie g\'en\'eratrice pour les partitions planes cylindriques. Notre premier r\'esultat est une preuve bijective de cette identit\'e utilisant  les diagrammes de croissance de Fomin for la correspondance RSK g\'en\'eralis\'ee. Le deuxi\`eme r\'esultat
est un $(q,t)$-analogue de la formule de Borodin, qui g\'en\'eralise un r\'eultat d'Okada. Enfin le trois\`eme r\'esultat de la premi\`ere partie est une description combinatoire explicite du poids de Macdonald intervenant dans cette formule, qui utilise un mod\`ele de chemins non-intersectant pour les partitions planes cylindriques.

Les matrices \`a signes alternants ont \'et\'e d\'ecouvertes par Robbins et Rumsey alors qu'ils \'etudiaient les  
 \emph{$\lambda$-d\'eterminants}.  Dans la deuxi\`eme partie de cette th\`ese nous d\'emontrons une g\'en\'eralisation \`a plusieurs param\`etres de ce $\lambda$-d\'eterminant, generalisant un r\'esultat r\'ecent de   di Francesco.
Comme le $\lambda$-d\'eterminant,  notre formule  est un exemple du \emph{ph\'enom\`ene de Laurent}.

Les syst\`emes semi-circulaires ont \'et\'e introduits par 
 Voiculescu afin d'\'etudier les alg\`ebres de  von Neumann des groupes libres. Dans la trois\`eme partie de la th\`ese, nous \'etudions les commutateurs dans l'alg\`ebre engendr\'ee par un syst\`eme semi-circuliare. Nous avons mis en \'evidence une matrice poss\'edant une structure auto-silimaire int\'eressante, qui nous permet de donner une formule explicite pour la projection sur l'epace des commutateurs de degr\'e donn\'e. En utilisant cette expression, nous donnons une preuve simple du fait que les syst\`eme semi-circulaires engendrent des facteurs.



\subsection*{\small Keywords:}
Cylindric partitions, Borodin's identity, growth diagrams, local rules, Schur functions, Pieri rules, Cauchy identity, Macdonald polynomials, commutation relations, RSK correspondence, non-intersecting lattice paths on the cylinder, alternating sign matrices, domino tilings of Aztec diamond, Bruhat order, lambda determinants, Laurent phenomenon, semicircular systems, von Neumann algebras, free probability theory, Chebyshev polynomials.

\newpage

\section*{\normalsize Abstract}

This thesis is divided into three parts. The first part deals with cylindric plane partitions. The second with lambda-determinants and the third with commutators in semicircular systems.

Cylindric plane partitions may be thought of as a natural generalization of reverse plane partitions. A generating series for the enumeration of cylindric plane partitions was recently given by Borodin. The first result of section one is a new bijective proof of Borodin's identity which makes use of Fomin's growth diagram framework for generalized RSK correspondences. The second result is a $(q,t)$-analog of Borodin's identity which extends previous work by Okada in the reverse plane partition case. The third result is an explicit combinatorial interpretation of the Macdonald weight occurring in the $(q,t)$-analog using the non-intersecting lattice path model for cylindric plane partitions.

Alternating sign matrices were discovered by Robbins and Rumsey whilst studying $\lambda$-determinants.  In the second part of this thesis we prove a multi-parameter generalization of the $\lambda$-determinant, generalizing a recent result by di Francesco.
Like the original $\lambda$-determinant, our formula exhibits the Laurent phenomenon.

Semicircular systems were first introduced by Voiculescu as a part of his study of von Neumann algebras. In the third part of this thesis we study certain commutator subalgebras of the semicircular system. We find a projection matrix with an interesting self-similar structure. Making use of our projection formula we given an alternative, elementary proof that the semicircular system is a factor.

\subsection*{\small Keywords:}
Cylindric partitions, Borodin's identity, growth diagrams, local rules, Schur functions, Pieri rules, Cauchy identity, Macdonald polynomials, commutation relations, RSK correspondence, non-intersecting lattice paths on the cylinder, alternating sign matrices, domino tilings of Aztec diamond, Bruhat order, lambda determinants, Laurent phenomenon, semicircular systems, von Neumann algebras, free probability theory, Chebyshev polynomials.

\cleardoublepage

\lhead[\oldstylenums \thepage]{Contents}
\rhead[Contents]{\oldstylenums \thepage}
\tableofcontents
\cleardoublepage

\newdimen{\cellsize}
\newcommand\Bigboxes{\setlength{\cellsize}{24pt}\def\boxformat{}}%
\newcommand\bigboxes{\setlength{\cellsize}{18pt}\def\boxformat{}}
\newcommand\medboxes{\setlength{\cellsize}{14pt}\def\boxformat{}}
\newcommand\smallboxes{\setlength{\cellsize}{8pt}\def\boxformat{\scriptstyle}}
\medboxes
\newsavebox{\cellcontent}
\def\hidehrule#1#2{\kern-#1
  \hrule height#1 depth#2 \kern-#2 }%
\def\hidevrule#1#2{\kern-#1{\dimen\cellcontent=#1%
    \advance\dimen\cellcontent by#2\vrule width\dimen\cellcontent}\kern-#2 }%
\def\makeblankbox#1#2{\hbox{\lower\dp\cellcontent\vbox{\hidehrule{#1}{#2}%
    \kern-#1 
    \hbox to \wd\cellcontent{\hidevrule{#1}{#2}%
      \raise\ht\cellcontent\vbox to #1{}
      \lower\dp\cellcontent\vtop to #1{}
      \hfil\hidevrule{#2}{#1}}%
    \kern-#1\hidehrule{#2}{#1}}}}
\newcommand\cellify[1]{\defaultcell%
\sbox{\cellcontent}{\vbox to \cellsize{%
\vfill%
\hbox to \cellsize{\hfill$\boxformat #1$\hfill}
\vfill}}%
\rlap{\drawnbox}
\usebox{\cellcontent}}
\newcommand\tableau[1]{\vtop{\let\\\cr
\baselineskip -16000pt \lineskiplimit 16000pt \lineskip 0pt
\ialign{&\cellify{##}\cr#1\crcr}}}
\newcommand\defaultcell{\gdef\drawnbox{
\makeblankbox{0.2pt}{0.2pt}
}}
\newcommand\graycell{\gdef\drawnbox{%
\rlap{\color{gray}\vrule width \cellsize height \cellsize}%
\makeblankbox{0.2pt}{0.2pt}
}}
\newcommand\bluecell{\gdef\drawnbox{%
\rlap{\color{Blue}\vrule width \cellsize height \cellsize}%
\makeblankbox{0.2pt}{0.2pt}
}}
\newcommand\redcell{\gdef\drawnbox{%
\rlap{\color{Red}\vrule width \cellsize height \cellsize}%
\makeblankbox{0.2pt}{0.2pt}
}}
\newcommand\greencell{\gdef\drawnbox{%
\rlap{\color{Green}\vrule width \cellsize height \cellsize}%
\makeblankbox{0.2pt}{0.2pt}
}}
\newcommand\orangecell{\gdef\drawnbox{%
\rlap{\color{Orange}\vrule width \cellsize height \cellsize}%
\makeblankbox{0.2pt}{0.2pt}
}}
\newcommand\yellowcell{\gdef\drawnbox{%
\rlap{\color{Yellow}\vrule width \cellsize height \cellsize}%
\makeblankbox{0.2pt}{0.2pt}
}}
\newcommand\thickyellowcell{\gdef\drawnbox{%
\rlap{\color{Yellow}\vrule width \cellsize height \cellsize}%
\makeblankbox{0.2pt}{0.1\cellsize}
}}
\newcommand\thickbluecell{\gdef\drawnbox{%
\rlap{\color{Blue}\vrule width \cellsize height \cellsize}%
\makeblankbox{0.2pt}{0.1\cellsize}
}}
\newcommand\thickredcell{\gdef\drawnbox{%
\rlap{\color{Red}\vrule width \cellsize height \cellsize}%
\makeblankbox{0.2pt}{0.1\cellsize}
}}
\newcommand\thickgreencell{\gdef\drawnbox{%
\rlap{\color{Green}\vrule width \cellsize height \cellsize}%
\makeblankbox{0.2pt}{0.1\cellsize}
}}
\newcommand\thickcell{\gdef\drawnbox{
\makeblankbox{0.2pt}{0.1\cellsize}%
}}
\newcommand\missingcell{\gdef\drawnbox{}}
\newcommand\vdotscell{\gdef\drawnbox{\kern-1.6pt\vbox{\baselineskip=4pt\lineskiplimit=0pt\hbox{}\hbox{.}\hbox{.}\hbox{.}\hbox{}}}}
\newcommand\hdotscell{\gdef\drawnbox{\vbox to \cellsize{\hbox{\kern1pt$\ldotp\ldotp\ldotp$}}}}
\newcommand\vhdotscell{\gdef\drawnbox{\rlap{\kern-1.6pt\vbox{\baselineskip=4pt\lineskiplimit=0pt\hbox{}\hbox{.}\hbox{.}\hbox{.}\hbox{}}}\vbox to \cellsize{\hbox{\kern1pt$\ldotp\ldotp\ldotp$}}}}
%
\newcommand\vertlinecell{\gdef\drawnbox{\unitlength=\cellsize%
\begin{picture}(1,1)
\put(0,0){\line(0,1){1}}
\end{picture}}}
\newcommand\horizlinecell{\gdef\drawnbox{\unitlength=\cellsize%
\begin{picture}(1,1)
\put(0,1){\line(1,0){1}}
\end{picture}}}
\newcommand\defaultcella{\gdef\drawnbox{
\unitlength=\cellsize%
\begin{picture}(1,1)
\put(0,0){\line(1,0){1}}
\put(0,0){\line(0,1){1}}
\put(1,0){\line(0,1){1}}
\put(0,1){\line(1,0){1}}
\end{picture}}}
\newcommand\thickcella{
\gdef\drawnbox{%
\unitlength=\cellsize%
\begin{picture}(1,1)
\linethickness{0.1\cellsize}
\put(0.0,0.05){\line(1,0){1}}
\put(0.05,0){\line(0,1){1}}
\put(0.95,0){\line(0,1){1}}
\put(0,0.96){\line(1,0){1}}
\end{picture}}}
\newcommand\graycella{\gdef\drawnbox{
\rlap{\color{Gray}\vrule width \cellsize height \cellsize}%
\unitlength=\cellsize%
\begin{picture}(1,1)
\put(0,0){\line(1,0){1}}
\put(0,0){\line(0,1){1}}
\put(1,0){\line(0,1){1}}
\put(0,1){\line(1,0){1}}
\end{picture}}%
}

\mainmatter

\lhead[\oldstylenums \thepage]{Introduction}
\rhead[Introduction]{\oldstylenums \thepage}
\addcontentsline{toc}{chapter}{Introduction}
\chapter*{Introduction}

This thesis is divided into three parts. The three parts are entirely independent and may be read in any order. The first part is significantly longer than the other two.

Part I deals with \emph{cylindric plane partitions}. The main tools used are the theory of symmetric functions and Macdonald polynomials \cite{macdonald} and Fomin's theory of \emph{growth diagrams} \cite{fomin-growth,fomin-comm}. The results of this section are entirely my own. They were presented as a poster at FPSAC 2012 and will appear in the Electronic Journal of combinatorics.

Part II  deals with a multi-parameter generalization of the \emph{$\lambda$-determinant} of Robbins and Rumsey \cite{robbins} which was originally conjectured by Alain Lascoux. Our formula exhibits the \emph{Laurent phenomenon} \cite{laurent}. It also generalizes a recent result by di Franceso \cite{cluster}. Our approach is completely different from that of di Francesco. We follow the original proof of Robbins and Rumsey closely, analyzing carefully the Bruhat order structure on pairs of interlacing alternating sign matrices (or equivalently, domino tilings of the Aztec Diamond \cite{aztec}). The idea for the proof was suggested by my advisor, Philippe Biane.
This work is not yet published.

Part III is a study of commutators in semicircular systems \cite{semicircular}. This section is somewhat of a work in progress. Although we have some preliminary results, we have not yet had the chance to apply them seriously.
This is joint work with Philippe Biane.

\section*{Cylindric Plane Partitions}
\label{sec:in}

\subsection*{Summary of results}

There are three main results in this section. The first is a bijective proof of Borodin's identity. The second is a Macdonald polynomial analog. The third is a combinatorial interpretation of the weight function which appears in the Macdonald polynomial analog of Borodin's identity.

Our bijection actually proves a refined version of Borodin's identity. The refined version of the reverse plane partition case is due to Gasner \cite{gasner-rpp}. A bijective proof using Fomin's growth diagram framework was previously given in the reverse plane partition case by Krattenthaler \cite{krattenthaler}. I was inspired to attempt the cylindric case after reading a well-written paper by Adachi \cite{adachi}.

The Macdonald polynomial analog is also proved in the full generality of the refined case. 
The Macdonald analog of the reverse plane partition case is due to Okada \cite{okada}.
The Hall--Littlewood case of cylindric identity is due to Corteel and Savelief, Cyrille and Vuleti{\'c} \cite{vuletic}.
The commutation relations which are key to the whole approach are due to Haiman, Garcia and Tesler \cite{garsia}.
When $q=0$ our combinatorial formula for the weight function reduces to the Hall--Littlewood version given in \cite{vuletic}.

\subsection*{Bijective proof of Borodin's identity}
Cylindric plane partitions were first introduced by Gessel and Krattenthaler \cite{gessel-1997}. 
For any binary string $\pi$ of length $T$,
a \emph{cylindric plane partition} with profile $\pi$ may be defined as a sequence of integer partitions:
\begin{equation*} 
(\mu^0, \mu^1, \ldots \mu^T) \qquad \qquad \mu^0 = \mu^T 
\end{equation*}
such that if $\pi_k = 1$ then $\mu^k / \mu^{k-1}$ is a \emph{horizontal strip}.
Otherwise if $\pi_k = 0$ then $\mu^{k-1} / \mu^k$ is a horizontal strip (see Section \ref{sec:strip}).

The \emph{weight} of a cylindric partition is given by
$  |\mathfrak{c}| = |\mu_1| + |\mu_2| + \cdots |\mu_T|$.
In the special case where $\mu_0 = \mu_T = \emptyset$ we recover the usual definition of a \emph{reverse plane partition}.
If, in addition to this there are no inversions in the profile, we have a \emph{regular plane partition} (see Section \ref{sec:pp}).

For those readers who are more familiar with the definition of a plane partition as an array of integers which is weakly decreasing along both rows and columns, the bijection with the ``interlacing sequence'' model is obtained by reading from right to left along the $NW \to SE$ diagonals. For example, the plane partition:
\[ \tableau{4 & 3 & 2 & 2 & 0 \\ 3 & 2 & 1 & 1 & 0 \\ 1 & 1 & 1 & 0 & 0\\ 1 & 0 & 0 & 0 & 0 \\ 0 & 0 & 0 & 0 & 0} \]
corresponds to the following sequence of partitions:
\[ \mathfrak{c} = (\emptyset, (2), (2,1), (3,1), (4,2,1), (3,1), (1), (1), \emptyset ) \]

A regular plane partition may also be thought of as a pair of \emph{semi-standard Young tableaux} of the same shape. In the case of our example, the two tableaux are:
\[\tableau{\missingcell & \missingcell & \missingcell & 4  \\ \missingcell & \missingcell & 4 & 2 \\ 4& 3 & 1 & 1} \qquad
\tableau{\missingcell & \missingcell & \missingcell & 4  \\ \missingcell & \missingcell & 4 & 3 \\ 4& 3 & 3 & 1}
\] 
We are using neither the French nor the English notation for Young diagrams.
The first tableau in the pair corresponds to the first half of the sequence read off the $NW \to SE$ diagonals, from right to left:
\[ \emptyset \rightharpoonup (2) \rightharpoonup (2,1) \rightharpoonup (3,1) \rightharpoonup (4,2,1)\]
while the second tableau in the pair corresponds to the second half of the sequence read off the $NW \to SE$ diagonals, from right to left:
\[\emptyset \rightharpoonup (1) \rightharpoonup (1) \rightharpoonup (3,1) \rightharpoonup (4,2,1)) \]

The \emph{RSK} correspondence gives a bijection between pairs of standard Young tableau of the same shape, and integer matrices with non-negative integers. This is essentially the \emph{Cauchy identity} for symmetric functions:
\[ \prod_{i,j \geq 1} \frac{1}{1-x_i,y_j} = \sum_\lambda S_\lambda(X) S_\lambda(Y) \]

In 2007 Borodin \cite{borodin}  gave a symmetric function theoretic proof of the following hook-product formula for the enumeration of cylindric plane partitions of given profile $\pi$ of length $T$ . In 2008, a very different proof involving the representation theory of $\widehat{sl}(n)$ was given by Tingley \cite{tingley}: 
\begin{equation*}\label{intro:borodin} 
\sum_{\mathfrak{c} \in \CPP(\pi)} z^{|\mathfrak{c}|} = 
\prod_{n \geq 0} \left ( \frac{1}{1-z^{(n+1)nT}} \prod_{\substack{i < j \\ \pi_i > \pi_j}} \frac{1}{1-z^{j-i + nT}} 
\prod_{\substack{i > j \\ \pi_i > \pi_j }} \frac{1}{1 - z^{j-i + (n+1)T}} \right ) 
\end{equation*}

This identity generalizes, not only MacMahon's identity for regular plane partitions: 
\begin{equation*}
\sum_{\mathfrak{c} \in \PP} z^{|\mathfrak{c}|} = \left ( \frac{1}{1-z^n} \right )^n 
\end{equation*}
but also Stanley's identity for reverse plane partitions:

\begin{equation*} 
\sum_{\mathfrak{c} \in \RPP(\pi)} z^{|\mathfrak{c}|} = 
\prod_{\substack{i < j \\ \pi_i > \pi_j}} \frac{1}{1-z^{j-i}} 
\end{equation*}

The \emph{Robinson correspondence} \cite{Robinson} gives a bijection between permutations and pairs of \emph{standard tableaux} of the same shape. Fomin's growth diagram framework \cite{fomin-growth} gives a particularly elegant way of understanding this bijection. 
Fomin's growth diagram framework is strictly equivalent to Viennot's geometric construction \cite{shadow} (see Section \ref{sec:robinson}).

Underlying Fomin's approach is an action of the \emph{Weyl algebra} on integer partitions (see Section \ref{weyl}). The creation operator $c$ adds a box to a partition in every way possible. The annihilation operator $c^*$ removes a box from a partition in every way possible. The \emph{canonical commutation relations}:
\[ [c^*,c] = 1 \]
can be understood as saying that each integer partition always has one more outside corner than inside corner.

Fomin also showed \cite{fomin-comm} that his abstract framework can be applied in the context of a wide class of commutation relations, including those commutation relations between vertex operators used to study plane partitions by Okounkov Reshetikhin \cite{okounkov}. 
The specific local rules which are needed in this case were described explicitly by vanLeeuwen \cite{vanLeeuwen}.
Interestingly, the local rules which are needed for the full RSK correspondence can be derived directly from those which apply in the special case of the Robinson correspondence (see Section \ref{sec:RSK}).

The commutation relations which we are interested in are those between the \emph{Pieri operator} $\Omega[Xz]$ and the \emph{dual Pieri operator} $\Omega^*[Xz]$. These operators act on Schur functions as follows:

\begin{equation*}
\Omega[Xz] S_\mu[X]
= \sum_r S_\mu[X] h_r[X] z^r = \sum_{\lambda \in U_r(\mu)} S_\lambda[X] z^r \end{equation*}
\begin{equation*} \Omega^*[Xz] S_\lambda = S_\lambda[X+z] = \sum_{\mu \in D_r(\mu)} S_\mu[X] z^r \end{equation*}
Here $U_r(\mu)$ denotes the set of all partitions which can be obtained from $\mu$ by adding a horizontal $r$-strip and $D_r(\lambda$) denotes the set of all partitions which can be obtained from $\lambda$ by removing a horizontal $r$-strip.

The Pieri operators satisfy the following commutation relations (see Section \ref{omega-commutation}): 
\begin{equation*} \label{intro:commute}
\Omega^*[Xu] \Omega[Xv] = \frac{1}{1-uv} \Omega[Xv] \Omega^*[Xu]
\end{equation*}

In section \ref{sec:algebraic} we give a slightly simplified algebraic proof of Borodin's identity using only the above commutation relation and a certain ``traciality property''. The definitions in Section \ref{sec:cylindric-defs} are necessary to understand the framework of the bijective proof which is contained in Section \ref{higher-order}. The bijection itself is described in Section \ref{sec:bijection}. The proof that the bijection is weight preserving is given in Section \ref{sec:weight}.

\subsection*{Macdonald Polynomial analog}

\emph{Macdonald polynomials} $\{P_\lambda(X)\}$ \cite{macdonald} are a family of symmetric polynomials over the ring $\mathbb{Q}(q,t)$ of rational functions in $q$ and $t$. The Macdonald polynomials bear a number of remarkable similarities with the Schur functions. In particular the Macdonald polynomials satisfy the following $(q,t)$-analog of the Cauchy identity (see Section \ref{sec:macdonald}):
\begin{equation*} \prod_{i,j \geq 1} \frac{(tx_i y_j;q)_\infty}{(x_iy_j;q)_\infty} =
\sum_\lambda \prod_{s \in \lambda} \frac{(1 - q^{a_\lambda(s)} t^{\ell_\lambda(s)+1})}{(1 - q^{a_\lambda(s)+1} t^{\ell_\lambda(s)}) } P_\lambda(X;q,t) P_\lambda(Y;q,t)
\end{equation*}
The Macdonald polynomials also satisfy a $(q,t)$-analog of the Pieri rules:
\begin{equation*}
 \Omega[Xz]_{q,t} \, P_\mu(X;q,t) = \sum_{\lambda \in U(\mu)}\psi_{\lambda / \mu}(q,t) \, P_\lambda(X;q,t) z^{|\lambda|-|\mu|}
\end{equation*}
\begin{equation*}
 \Omega^*[Xz]_{q,t} \, P_\lambda(X;q,t) = \sum_{\mu \in D(\lambda)}\varphi_{\lambda / \mu}(q,t) \, P_\mu(X;
q,t) z^{|\lambda|-|\mu|}
\end{equation*}
The Pieri coefficients (\cite{macdonald} page 341) are given by:
\begin{equation*}
 \varphi_{\lambda / \mu}(q,t) = \prod_{s \in C_{\lambda / \mu}}\frac{1 - q^{a_\lambda(s)} t^{\ell_\lambda(s)+1}}{1 -q^{a_\lambda(s)+1} t^{\ell_\lambda(s)}} 
\prod_{s \in {C}_{\lambda / \mu}}\frac{1 - q^{a_\mu(s)+1} t^{\ell_\mu(s)}}{1 -q^{a_\mu(s)} t^{\ell_\mu(s)+1}}
\end{equation*}
\begin{equation*}
 \psi_{\lambda / \mu}(q,t) = \prod_{s \not\in C_{\lambda / \mu}}\frac{1 - q^{a_\lambda(s)+1} t^{\ell_\lambda(s)}}{1 - q^{a_\lambda(s)} t^{\ell_\lambda(s)+1}} 
\prod_{s \not\in {C}_{\lambda / \mu}}\frac{1 -q^{a_\mu(s)} t^{\ell_\mu(s)+1}}{1 -q^{a_\mu(s)+1} t^{\ell_\mu(s)}}
\end{equation*}
Here $C_{\lambda / \mu}$ denotes the set of columns of $\lambda$ which are longer than the corresponding columns of $\mu$.
In Section \ref{sec:qtborodin} we prove the following $(q,t)$-analog of Borodin's identity:

\begin{theorem}\label{intro:qt-borodin}
\begin{equation} 
\sum_{\mathfrak{c} \in CPP(\pi)} W_\mathfrak{c}(q,t) z^{|\mathfrak{c}|} = 
\prod_{n \geq 0} \left ( \frac{1}{1 - z^{(n+1)T}}  
\prod_{\substack{i < j \\ \pi_i > \pi_j}} \frac{(tz^{j - i + nT};q)_\infty}{(z^{j - i + nT};q)_\infty} 
\prod_{\substack{i > j \\ \pi_i > \pi_j}} \frac{(tz^{j-i + (n+1)T};q)_\infty}{(z^{j-i + (n+1)T};q)_\infty} \right )
\end{equation}
\end{theorem}
where the weight function is given by:
\begin{equation}  
W_{\mathfrak{c}}(q,t) = \prod_{\substack{k =0 \\ \pi_k = 1}}^T \varphi_{\mu^k / \mu^{k-1}}(q,t) \prod_{\substack{k = 0\\ \pi_k = 0}}^T \psi_{\mu^{k-1} / \mu^k}(q,t) 
\end{equation}
Our proof relies on the following $(q,t)$-analog of the commutation relation which is due to Haiman, Garcia and Tesler \cite{garsia}:
\begin{equation} \label{intro:qt-commute}
\Omega_{q,t}^*[Xu] \, \Omega_{q,t}[Xv] = \frac{(tuv;q)_\infty}{(uv;q)_\infty} \, \Omega_{q,t}[Xv] \, \Omega_{q,t}^*[Xu]
\end{equation}

\subsection*{Simplification of weight function}

Although we have defined cylindric plane partitions as certain sequences of integer partitions which differ by a horizontal strip, it is also possible to define them families of non-intersecting lattice paths on a cylinder (see Section \ref{sec:lattice}). Using this latter definition, the weight function 
$W_{\mathfrak{c}}(q,t)$ may be greatly simplified.

Recall that in the plethystic notation \cite{garsia} if: 
\[ a(q,t) = \sum_{n,m} a_{n,m} \, q^n t^m \]
with $a_{n,m} \in \mathbb{Z} $ and $a_{0,0} = 0$, then we have:
\begin{equation*} \Omega \left [ a(q,t) \right ] = \prod_{n,m} \frac{1 }{(1 - q^n t^m)^{a_{n,m}}} \end{equation*}
In section \ref{sec:weight-proof} we make use of the plethystic notation to give the cylindric weight function the following explicit combinatorial description:
\begin{theorem}\label{intro:weight}
\begin{equation}
W_\mathfrak{c}(q,t) = \Omega \left [ (q-t) \mathcal{D}_\mathfrak{c}(q,t) \right ] 
\end{equation}
where the alphabet $\mathcal{D}_\mathfrak{c}(q,t)$ is given by:
\begin{equation}\mathcal{D}_\mathfrak{c}(q,t) =  \sum_{s \in peak(\mathfrak{c})} \, q^{a_\mathfrak{c}(s)} t^{\ell_\mathfrak{c}(s)} - \sum_{s \in valley(\mathfrak{c})} \, q^{a_\mathfrak{c}(s)} t^{\ell_\mathfrak{c}(s)}
\end{equation}
\end{theorem}

The precise definition of ``valley'' and ``peak'' cubes depend on the lattice path picture. They are defined in section \ref{sec:cube-lattice}.

\subsection*{Outline}

For the impatient reader, here is a quick road map. Page links to important definitions can be found in the index.

The most important part of Chapter $1$ is the section pertaining to cylindric diagrams, cylindric inversion co-ordinates and cylindric hook lengths, as well as the rotation operator. The purpose of section \ref{sec:pp} is to act as a motivator for the definition or arbitrary cylindric diagrams in section \ref{alcd}. The definition of arms and legs in \ref{sec:arms-legs} will be needed for the definition of the Macdonald polynomial weight defined in section \ref{sec:macdonald-pieri}.

Chapter $2$ does not contain original work, and with the exception of Section \ref{sec:local-explicit} is not strictly speaking needed in order to construct the cylindric bijection. It does nevertheless provide motivation and intuition for understanding Section \ref{higher-order}.

The introduction to the theory of symmetric functions in chapter $3$ is very brief. The key point is the Pieri formula and the commutation relations in both the Schur and Macdonald cases. Sections \ref{sec:sym-local} and \ref{sec:rep} attempt to clarify the relation between the algebra and the combinatorics.

Section \ref{sec:algebraic} contains all the algebra that is needed to understand section \ref{sec:qtborodin}. It is not strictly needed to understand the bijection, since the bijection can be formulated in a purely combinatorial manner. Section \ref{higher-order} sets up notation required for working with local rules and encoding the recursive structure of the bijection. Section \ref{sec:bijection} defines cylindric growth diagrams and verifies that they act as an ``interpolation'' object between the two sides of the identity to be proved. Finally in section \ref{sec:weight} we verify that the bijection which has been constructed in the previous two sections is weight preserving.

In section \ref{sec:qtborodin} we prove Theorem \ref{intro:qt-borodin}. In Section \ref{sec:lattice} we describe how cylindric plane partitions may be interpreted as non-intersecting lattice paths on a cylinder.
It is here that we define peak and valley cubes, amongst other things. The diagonal reading of a family of non-intersecting lattice paths together with Proposition \ref{prop:crucial} is crucial for understanding the combinatorial reformulation of the weight function. The rest can be safely ignored.

In section \ref{sec:weight-proof} we prove Theorem \ref{intro:weight}.
The key point to understand is that in the ``interlacing sequence'' model the cubes are grouped according to which partition in the sequence they belong to. In the lattice path model, cubes from the same column number but different partitions are grouped together.
To get from one definition of the weight function to the other, we simply switch between these two models.

\section*{Lambda determinants}

\subsection*{Main result}
An \emph{alternating sign matrix} is a square matrix of $0$'s $1$'s and $-1$'s such that the sum of each row and column is $1$ and the non-zero entries in each row and column alternate in sign. For example:
\[A = 
\left ( \begin{matrix}
0 & 0 & 1 & 0 \\
0 & 1 & -1 & 1 \\
1 & 0 & 0 & 0 \\
0 & 0 & 1 & 0
\end{matrix} \right )
\]
A permutation matrix is an special case of an alternating sign matrix with no $(-1)$'s.
The total number of alternating sign matrices of size $n$ is given by:

\[A_n = \prod_{k=0}^{n-1} \frac{(3k+1)!}{(n+k)!} \qquad \qquad  1, 1, 2, 7, 42, 429, 7436, \ldots\]

The first proof of this result was given by Zeilberger \cite{zeil}. A much simplified proof was given by Kuperberg \cite{kuperberg}. Kuperberg's proof made use of ideas from the theory of integrable systems, the Yang-Baxter equation and the  six vertex model with domain wall boundary conditions. It also made use of a recurrence relation due to Izergin and Korepin \cite{korepin}.

The theory of alternating sign matrices is currently an active area of research. One long standing open problem is to find an explicit bijection between the set of alternating sign matrices and \emph{totally symmetric self-complementary plane partitions} \cite{tsscp}. Alternating sign matrices also played a key role in the recently proven \emph{Razumov-Stroganoff conjecture} \cite{razumov}.

The first time which alternating sign matrices appeared in the literature was in the famous paper by Robbins and Rumsey \cite{robbins} on the \emph{lambda-determinant}. The lambda determinant may be defined as follows:

For each $k=0 \ldots n$ let us denote by $x_n[k]$ the doubly indexed collection of variables $x_n[k]_{i,j}$ with indices running from $i,j = 1..(n-k+1)$. 
One should think of the variables as forming a square pyramid with base $n+1$ by $n+1$. The index $k$ determines the ``height'' 
of the variable in the pyramid.

The variables are initialized as follows:
\begin{align*} 
x_n[0]_{i,j} & = 1 \mbox{ for all } i,j = 1..(n+1) \\
x_n[1]_{i,j} & = M_{i,j} \mbox{ for all } i,j = 1..n \\
\end{align*}
The value of the remaining variables is calculated via the following recurrence:
\begin{align} \label{intro:lambda-rec}
x_n[k+1]_{i,j} =  \frac{ x_n[k]_{i,j} x_n[k]_{i+1,j+1} + \lambda \, x_n[k]_{i,j+1} x_n[k]_{i+1,j}}{x_n[k-1]_{i+1,j+1}}
\end{align}
The end result \cite{robbins} is that:

\begin{equation*}  x_n[n]_{1,1} = \sum_{B \,\in \, \mathfrak{A}_n} \lambda^{\inv(B)}(1+\lambda)^{N(B)} \prod_{i,j=1}^n M_{i,j}^{B_{i,j}}  \end{equation*}
Here $\mathfrak{A}_n$ denotes the set of all alternating sign matrices of size $n$, $\inv(B)$ denotes the \emph{inversion number} of $B$ and $N(B)$ denotes the number of negative ones in $B$.

Note that the $\lambda$-determinant exhibits the \emph{Laurent phenomenon} \cite{laurent}. From the recursive definition we expect the value of $x_n[n]_{1,1}$ to be a rational function. The fact that it is a Laurent polynomial is very surprising.

When $\lambda = -1$ the $\lambda$-determinant reduces to the regular determinant, and the recursive method for calculating the determinant above reduces to the algorithm known as \emph{Dodgson condensation} \cite{dodgson}.

Our main result is to replace the recurrence given in equation \ref{intro:lambda-rec} by the following recurrence:
\begin{align} 
x_n[k+1]_{i,j} =  \frac{ \mu_{i,n-k+1-j} x_n[k]_{i,j} x_n[k]_{i+1,j+1} + \lambda_{i,j} x_n[k]_{i,j+1} x_n[k]_{i+1,j}}{x_n[k-1]_{i+1,j+1}}
\end{align}
This is a more general case of the recurrence considered by di Francesco \cite{cluster}. The closed form expression which we fine for $x_n[n]_{1,1}$ is the following:
\begin{equation*}
x_n[n]_{1,1} = \sum_{|B| = n} M^B \left ( 
\prod_{(i,j) \in \inv(B)} \lambda_{i,j}
\prod_{(i,j) \in \dinv(B)} \mu_{i,n+1-j}
\prod_{B_{i,j} = -1} 
\left (\mu_{i,n+1-j} + \lambda_{i,j} \right )
\right )
\end{equation*}
where $\inv(B)$ denotes the set of \emph{inversions} of $B$ and $\dinv(B)$ denotes the set of \emph{dual inversions} of $B$.
Note that our formula also exhibits the Laurent phenomenon.

It is also possible to consider more general initial conditions:
\begin{align*} 
x_n[0]_{i,j} & = N_{i,j} \mbox{ for all } i,j = 1..(n+1) \\
x_n[1]_{i,j} & = M_{i,j} \mbox{ for all } i,j = 1..n \\
\end{align*}
and give a closed form expression for $x_n[k+1]_{1,1}$. In order to do this we must first introduce the idea of \emph{interlacing matrices}, which are closely related to \emph{domino tilings of the Aztec diamond} \cite{aztec}. We may state the formula here, however the reader will have to wait until sections \ref{sec:left-corner} and \ref{sec:right-cornersum} for the definitions of $F_\lambda(B)$ and $G^n_\mu(B)$.

\begin{equation} 
x_n[k+1]_{1,1} = \sum_{\substack{(A,B) \\ |B| = k, |A| = k-1}} 
\frac{F_\lambda(B)}{s(F_\lambda(A))} \frac{G^n_\mu(B)}{t(G^n_\mu(A))} M^B s(N)^{-A} 
\end{equation}
The sum is over all pairs of interlacing matrices.

\subsection*{Outline}

In Section \ref{lambda:perm} we define the Bruhat order. We show how a permutation can be represented by a monotone triangle, and look at the inversions and \emph{dual inversion} of a permutation.

In Section \ref{sec:asm} we define alternating sign matrices, and show that they complete the Bruhat order as a lattice. We extend the definition of monotone triangle to alternating sign matrices, as well as the definition of inversion and dual inversion.

In Section \ref{sec:interlacing} we define left corner sum matrices and left interlacing matrices, We show that pairs of left interlacing matrices are in bijection with domino tilings of the Aztec diamond.
It is here that we define the notation for $F(B)$.

In Section \ref{sec:duality2} we define right corner sum matrices and right interlacing matrices We study the duality between left and right interlacing matrices.
It is here that we define the notation for $G(B)$.

In Section \ref{sec:main-theorem} we prove our main theorem. The proof is by recurrence and makes use of results established in Section \ref{sec:duality2}.

\section*{Semicircular systems}

For any Hilbert space $H$ we may define its Fock space $\mathfrak{F}$ to be the metric closure of the tensor algebra of $H$. That is:
\[ \mathfrak{F} = \overline{T(H)}\]
where:
\[ T(H) = \oplus_{n \geq 0} H^{\otimes n} \]
Let us fix $\Omega$ to be some element of $H^{\otimes 0}$ of norm $1$.
For any $v \in H$ one can construct the \emph{creation operator}:
\[ c_v[x] = v \otimes x \]
as well as the \emph{annihilation operator}
\begin{align*}
 c^*_v[\Omega] & = 0 \\
 c^*_v[x \otimes y] & = \langle v | x \rangle y
\end{align*}
An operator of the form $A_v = c_v + c^*_v$ may be thought of as a \emph{semi-circular random variable} \cite{semicircular}.
Let $\mathfrak{A}$ denote the von Neumann algebra generated by the semi-circular random variables of the form $A_v$.
The map:
\[ A \mapsto A[\Omega] \]
gives an embedding of $\mathfrak{A}$ into $\mathfrak{F}$ (as a vector space). We are interested in subspaces of $\mathfrak{F}$ of the form:
\[  V_A = \{[A,y], y \in \mathfrak{A} \}\]
where $A$ is some fixed element of $\mathfrak{A}$.

In Section \ref{semicircular:matrix} we study the \emph{Gramm-matrix} of a natural, non-orthogonal basis of $V_A$.
We find that this matrix has a curious self-similar structure.
In Section \ref{sec:semicircular:main} we find an explicit projection formula for the projection of any $B$ which is not in the subalgebra of $\mathfrak{A}$ generated by $A$ onto the subspace $V_A$. This allows us to prove, in particular, that the center of $\mathfrak{A}$ is trivial. Although this is already well known \cite{semicircular}, our proof is particularly simple and elementary.

\cleardoublepage

\lhead[\oldstylenums \thepage]{\S~\thesection \; --- \; \rightmark}
\rhead[Chapitre~\thechapter \; --- \; \leftmark]{\oldstylenums \thepage}

\part{Cylindric Plane Partitions}
\chapter{Partitions}

\section{Integer Partitions}\label{sec:int}

\index{Young diagram}

An \emph{integer partition} is simply a weakly decreasing list of non-negative integers. 
It is often convenient to represent an integer partition visually as a \emph{Young diagram}, which is a collection of ``boxes'' in the Cartesian plane which are ``stacked up'' in the bottom right hand corner. 
\[
\tableau{
\missingcell & \missingcell & \missingcell & \missingcell  & \missingcell  \\
\missingcell & \missingcell & \missingcell  &  & \\
\missingcell & \missingcell &  & &  \\
\missingcell  & \missingcell &  & &  \\
 & &&&\\
} \\
\]
\[ \lambda = (5,3,3,2) \]

Note that our convention differs from both the standard french and English conventions. If the sum of the parts of $\lambda$ is equal to $n$, then we say that $\lambda$ is a partition of $n$ and write $|\lambda| = n$. If $\lambda$ has exactly $k$ distinct non-zero parts, then we say that $\lambda$ has length $k$ and write $\ell(\lambda) = k$.
The generating series for integer partitions is given by:

\begin{equation}
\sum_{\lambda} z^{|\lambda|} t^{\ell(\lambda)} = \prod_{n \geq 1} \frac{1}{1-t z^n} 
\end{equation}

The \emph{conjugate} of the integer partition $\lambda = (\lambda_1, \lambda_2, \ldots, \lambda_k)$ is defined to be $\lambda' = (\lambda'_1, \lambda'_2, \ldots \lambda'_r)$ where $\lambda'_j = \#\{i \,|\, \lambda_i \geq j \}$.
For example the conjugate of the partition $\lambda = (5,3,3,2)$ is $\lambda' = (4,4,3,1,1)$. 
In terms of Young diagrams, conjugating a partition is equivalent to reflecting about the main diagonal.
\[
\tableau{
\missingcell & \missingcell & \missingcell & \\
\missingcell & \missingcell & \missingcell & \\
\missingcell &&& \\
&&&\\
&&&
}
\]
\[ \lambda' = (4,4,3,1,1) \]

\subsection{Inversions}\label{sec:inv}

\index{Partition!profile, minimum}

The \emph{minimum profile} of an integer partition is the binary string which traces out the ``jagged boundary'' of the associated Young diagram. Reading from the top right hand corner to the bottom left hand corner, a zero is recorded for every vertical step and a one for every horizontal step. For example the minimum profile of our example partition $\lambda = (5,3,3,2)$ is $110100110$:
\[
\tableau{
\missingcell & \missingcell & \missingcell & \missingcell 1 & \missingcell 1 \\
\missingcell & \missingcell & \missingcell 1 & 0 & \\
\missingcell & \missingcell & 0 & &  \\
\missingcell 1 & \missingcell1 & 0 & &  \\
0 & &&&\\
} 
\]

The minimum profile of an integer partition necessarily starts with a one and ends with a zero.  An integer partition is uniquely determined by its minimum profile.

\index{Binary String!inversion}

An \emph{inversion} in a binary string $\pi$ is a pair of indices $(i,j)$ such that $i < j$ and $\pi_i > \pi_j$. 
There is a natural bijection between the ``boxes'' of an integer partition $\lambda$ and the inversions in the minimum profile $\pi$ of the partition $\lambda$. The box marked with a star below has inversion co-ordinates $(2,6)$ because the one lying above it is position $2$ in the profile, while the zero lying to the left of it is position $6$ in the profile.

\[
\tableau{
\missingcell & \missingcell & \missingcell & \missingcell {\color{blue} 1} & \missingcell 1 \\
\missingcell & \missingcell & \missingcell 1 & 0 & \\
\missingcell & \missingcell & 0 &  &  \\
\missingcell 1 & \missingcell1 & {\color{blue} 0} & * &  \\
0 & &&&\\
} 
\]
\[ \pi = 1{\color{blue} 1} 010{\color{blue}0}110\]

\index{Box!arm length}
\index{Box!leg length}
\index{Box!hook length}
\index{Box!inversion coordinates}

\subsection{Arms, legs and hooks}\label{sec:arms-legs}

Let $s$ be a box of the partition $\lambda$ with profile $\pi$. Suppose that $s$ has ``inversion coordinates'' $(i,j)$.
The \emph{arm length} of $s$ is given by: 
\begin{equation}
a_\lambda(s) = \#\{ i < k < j \, | \, \pi_k = 1\}
\end{equation}
The \emph{leg length} of $s$
is given by:
\begin{equation}
\ell_\lambda(s) = \#\{ i < k < j \, | \, \pi_k = 0\}
\end{equation} 
The \emph{hook length} of $s$ is given by: 
\begin{equation}
h_\lambda(s) = a_\lambda(s) + b_\lambda(s) + 1 = j-i
\end{equation}

The arm length of the box $s$ counts the number of boxes in the same row as $s$ lying to the left, while the leg length of $s$ counts the number of boxes in the same column as $s$ but above. 
The arm length of our example box, marked with a star in the diagram above, is equal to $1$ while the leg length is equal to $2$.

\index{Box!inside corner}
\index{Box!outside corner}

A box with arm length zero and leg length zero is said to be an \emph{outside corner}. Equivalently an outside corner corresponds to a subword of the profile of the form $10$. An \emph{inside corner} is defined to be a subword of the profile of the form $01$. 
The outside corners of our example profile are the following:

\[ 1{\color{blue}10}100110 \qquad 110{\color{blue}10}0110 \qquad 1101001{\color{blue}10}\]
The inside corners of our example profile are the following:
\[ 11{\color{red}01}00110 \qquad 11010{\color{red}01}10 \]
The number of outside corners is always equal to one more than the number of inside corners.
An integer partition always has exactly one more outside corner than inside corner.

\subsection{Partial orders}\label{sec:partial}

\subsubsection{Generalized profiles}
\index{Partition!profile, generalized}

A \emph{generalized profile} is an arbitrary string of zeros and ones.
A generalized profile may be thought of as an integer partition pre-conceived as sitting inside some larger rectangle.
For example the generalized profile of our example partition $\lambda = (5,3,3,2)$ thought of as sitting inside an $8$ by $8$ box is $0000110100110111$:

\[
\tableau{
\missingcell & \missingcell & \missingcell &\missingcell & \missingcell & \missingcell & \missingcell & \missingcell & \missingcell 0 \\
\missingcell & \missingcell & \missingcell &\missingcell & \missingcell & \missingcell & \missingcell & \missingcell & \missingcell 0 \\
\missingcell & \missingcell & \missingcell &\missingcell & \missingcell & \missingcell & \missingcell & \missingcell & \missingcell 0 \\
\missingcell & \missingcell & \missingcell &\missingcell & \missingcell & \missingcell & \missingcell 1 & \missingcell 1 & \missingcell 0\\
\missingcell & \missingcell & \missingcell &\missingcell & \missingcell & \missingcell 1 & 0 & \\
\missingcell & \missingcell & \missingcell &\missingcell & \missingcell & 0 & &  \\
\missingcell & \missingcell & \missingcell &\missingcell 1 & \missingcell1 & 0 & &  \\
\missingcell 1 & \missingcell 1 & \missingcell 1 & 0 & &&&\\
} 
\]

\hfill \break

Let $\bin(n,m)$ denote the set of all binary strings with $n$ zeros and $m$ ones. This set is in bijection with the set of all integer partitions whose Young diagrams fit inside an $n$ by $m$ box.
For any pair of integer partitions $\mu$ and $\lambda$ we may find some $(n,m)$ such that both $\mu$ and $\lambda$ admit generalized profiles lying in $\bin(n,m)$.

\subsubsection{Young lattice}
\index{Young lattice}

It is possible to define a \emph{partial order} on the set of all integer partitions. For any two partitions $\mu$ and $\lambda$ we say that $\mu \subseteq \lambda$ if the Young diagram of $\mu$ fits inside the Young diagram of $\lambda$. 
For any pair of partitions $\alpha$ and $\beta$ there is a unique smallest partition containing both $\alpha$ and $\beta$ which we denote by $\alpha \cup \beta$. Similarly there is a unique largest partition contained in both $\alpha$ and $\beta$ which we denote by $\alpha \cap \beta$. In other words, our partial order forms a \emph{lattice} which is known as the \emph{Young lattice}

\[
\tableau{
\missingcell & \missingcell & \missingcell & \missingcell  & \missingcell  \\
\missingcell & \missingcell & \missingcell  &  & \\
\missingcell & \missingcell &  & &  \graycell\\
\missingcell  & \missingcell &  \graycell & \graycell & \graycell \\
 & & \graycell & \graycell & \graycell \\
} \\
\]
\[ (3,3,1) \subset (5,3,3,2) \]

\subsubsection{Partial order on binary strings}\label{sec:bin-order}
\index{Binary String!partial order}

It is also possible to define a partial order on $\bin(n,m)$
whose covering relations are given by
$\pi \prec \pi'$ if and only if there is some $i$ such that $\pi_i = 0 = \pi'_{i+1}$ and $\pi_{i+1} = 1 = \pi'_i$
and for all other $k$ we have $\pi'_k = \pi_k$. 
In other words $\pi'$ is obtained from $\pi$ by adding an inversion.

If $\lambda$ is the partition with generalized profile $\pi$ and $\lambda'$ is the partition with generalized profile $\pi'$ with $\pi,\pi' \in \bin(n,m)$ then
$\pi \prec \pi'$ if and only if $\lambda'$ can be obtained from $\lambda$ by adding a single box.

We shall denote by $\pi_{\min}$ the binary string with $n$ zeros follows by $m$ ones
and $\pi_{\max}$ the binary string with $m$ ones follows by $n$ zeros.
Note that $\pi_{\min}$ corresponds to the empty partition, while $\pi_{\max} = (m,m,\ldots,m)$.

\subsubsection{Dominance and lexicographic order}\label{sec:dominance}
\index{Partition!dominance order}
\index{Partition!lexicographic order}

Finally, the \emph{dominance order} on integer partitions is defined by $\mu \trianglelefteq \lambda$ if and only if:
\begin{align*}
\lambda_1 & \geq \mu_1 \\
\lambda_1 + \lambda_2 & \geq \mu_1 + \mu_2 \\
& \cdots \\
\lambda_1 + \lambda_2 + \cdots \lambda_k & \geq \mu_1 + \mu_2 + \cdots + \mu_k
\end{align*} 

The {\emph{lexicographical order} is a total order defined on integer partitions by $\lambda \geq \mu$ if and only if there exists some $m$
such that $\lambda_m > \mu_m$ and $\lambda_i = \mu_i$ for all $i \leq m$.
Note that $\mu \trianglelefteq \lambda$ implies $\mu \leq \lambda$ but the converse is false. 

\subsection{Horizontal and vertical strips}\label{sec:strip}

\index{Horizontal strip!regular}

For any pair of partitions $\lambda$ and $\mu$ satisfying $\mu \subseteq \lambda$ we say that $\lambda / \mu$ is a {\emph{horizontal strip} and write $\mu \rightharpoonup \lambda$ if and only if 
 \[ \lambda_1 \geq \mu_1 \geq \lambda_2 \geq \mu_2 \geq \ldots \]

Equivalently, $\mu \rightharpoonup \lambda$ if and only if each column of $\lambda$ contains at most one more box than the corresponding column of $\mu$.

\[
\tableau{
\missingcell & \missingcell & \missingcell & \missingcell  & \missingcell  \\
\missingcell & \missingcell & \missingcell  &  & \\
\missingcell & \missingcell & \graycell & \graycell &  \graycell\\
\missingcell  & \missingcell &  \graycell & \graycell & \graycell \\
 & & \graycell & \graycell & \graycell \\
} \\
\]
\[ (3,3,3) \rightharpoonup (5,3,3,2) \]

In terms of profiles, 
$\lambda / \mu$ is a horizontal strip if and only if the generalized profile of $\lambda$ can be obtained from the generalized profile of $\mu$ by ``hopping'' some of the ones, which may be thought of as ``particles'', a single step to the left.

\begin{center}
\begin{tikzpicture}[scale=0.8]

\begin{scope}

\path (0,0) node[shape=circle,inner sep=0.3mm](C){$0$};
\path (1,0) node[shape=circle,inner sep=0.3mm](C){$0$};
\path (2,0) node[shape=circle,inner sep=0.3mm](C){$0$};
\path (3,0) node[shape=circle,inner sep=0.3mm](C1){$1$};
\path (4,0) node[shape=circle,inner sep=0.3mm](C2){$1$};
\path (5,0) node[shape=circle,inner sep=0.3mm](C3){$1$};
\path (6,0) node[shape=circle,inner sep=0.3mm](C){$0$};
\path (7,0) node[shape=circle,inner sep=0.3mm](C){$0$};
\path (8,0) node[shape=circle,inner sep=0.3mm](C){$0$};
\path (9,0) node[shape=circle,inner sep=0.3mm](C4){$1$};
\path (10,0) node[shape=circle,inner sep=0.3mm](C5){$1$};
\path (11,0) node[shape=circle,inner sep=0.3mm](C6){$1$};

\path (0,1) node[shape=circle,inner sep=0.3mm](D){$0$};
\path (1,1) node[shape=circle,inner sep=0.3mm](D){$0$};
\path (2,1) node[shape=circle,inner sep=0.3mm](D1){$1$};
\path (3,1) node[shape=circle,inner sep=0.3mm](D2){$1$};
\path (4,1) node[shape=circle,inner sep=0.3mm](D){$0$};
\path (5,1) node[shape=circle,inner sep=0.3mm](D3){$1$};
\path (6,1) node[shape=circle,inner sep=0.3mm](D){$0$};
\path (7,1) node[shape=circle,inner sep=0.3mm](D){$0$};
\path (8,1) node[shape=circle,inner sep=0.3mm](D4){$1$};
\path (9,1) node[shape=circle,inner sep=0.3mm](D5){$1$};
\path (10,1) node[shape=circle,inner sep=0.3mm](D){$0$};
\path (11,1) node[shape=circle,inner sep=0.3mm](D6){$1$};

\path (5,2) node[shape=circle,inner sep=0.3mm](C0){};

\draw[thick] (D1) -- (C1);
\draw[thick] (D2) -- (C2);
\draw[thick] (D3) -- (C3);
\draw[thick] (D4) -- (C4);
\draw[thick] (D5) -- (C5);
\draw[thick] (D6) -- (C6);

\end{scope}

\end{tikzpicture}
\end{center}

\index{Vertical strip!regular}

Similarly, for any pair of partitions $\lambda$ and $\mu$ satisfying $\mu \subseteq \lambda$ we say that $\lambda / \mu$ is a {\emph{vertical strip} and write $\mu \rightharpoondown \lambda$ if and only if 
 \[ \lambda'_1 \geq \mu'_1 \geq \lambda'_2 \geq \mu'_2 \geq \ldots \]

Equivalently, $\mu \rightharpoondown \lambda$ if and only if each column of $\lambda$ contains at most one more box than the corresponding \emph{row} of $\mu$.

\[
\tableau{
\missingcell & \missingcell & \missingcell & \missingcell  & \missingcell  \\
\missingcell & \missingcell & \missingcell  &  \graycell & \graycell \\
\missingcell & \missingcell &  & \graycell &  \graycell\\
\missingcell  & \missingcell &   & \graycell & \graycell \\
 & \graycell & \graycell & \graycell & \graycell \\
} \\
\]
\[ (4,2,2,2) \rightharpoondown (5,3,3,2) \]

In terms of profiles, 
$\lambda / \mu$ is a vertical strip if and only if the generalized profile of $\lambda$ can be obtained from the generalized profile of $\mu$ by ``hopping'' some of the zeros, which may be thought of as ``holes'', a single step to the right.

\begin{center}
\begin{tikzpicture}[scale=0.8]

\begin{scope}

\path (11,0) node[shape=circle,inner sep=0.3mm](C1){$1$};
\path (10,0) node[shape=circle,inner sep=0.3mm](C2){$1$};
\path (9,0) node[shape=circle,inner sep=0.3mm](E1){$0$};
\path (8,0) node[shape=circle,inner sep=0.3mm](C0){$1$};
\path (7,0) node[shape=circle,inner sep=0.3mm](C0){$1$};
\path (6,0) node[shape=circle,inner sep=0.3mm](E2){$0$};
\path (5,0) node[shape=circle,inner sep=0.3mm](E3){$0$};
\path (4,0) node[shape=circle,inner sep=0.3mm](E4){$0$};
\path (3,0) node[shape=circle,inner sep=0.3mm](C6){$1$};
\path (2,0) node[shape=circle,inner sep=0.3mm](C0){$1$};
\path (1,0) node[shape=circle,inner sep=0.3mm](E5){$0$};
\path (0,0) node[shape=circle,inner sep=0.3mm](E6){$0$};

\path (11,1) node[shape=circle,inner sep=0.3mm](D1){$1$};
\path (10,1) node[shape=circle,inner sep=0.3mm](F1){$0$};
\path (9,1) node[shape=circle,inner sep=0.3mm](D2){$1$};
\path (8,1) node[shape=circle,inner sep=0.3mm](D3){$1$};
\path (7,1) node[shape=circle,inner sep=0.3mm](F2){$0$};
\path (6,1) node[shape=circle,inner sep=0.3mm](F3){$0$};
\path (5,1) node[shape=circle,inner sep=0.3mm](D4){$1$};
\path (4,1) node[shape=circle,inner sep=0.3mm](F4){$0$};
\path (3,1) node[shape=circle,inner sep=0.3mm](D5){$1$};
\path (2,1) node[shape=circle,inner sep=0.3mm](D6){$1$};
\path (1,1) node[shape=circle,inner sep=0.3mm](F5){$0$};
\path (0,1) node[shape=circle,inner sep=0.3mm](F6){$0$};

\path (5,2) node[shape=circle,inner sep=0.3mm](C0){};

\draw[thick] (E1) -- (F1);
\draw[thick] (E2) -- (F2);
\draw[thick] (E3) -- (F3);
\draw[thick] (E4) -- (F4);
\draw[thick] (E5) -- (F5);
\draw[thick] (E6) -- (F6);

\end{scope}

\end{tikzpicture}
\end{center}

If $\lambda / \mu$ is a horizontal strip, then the conjugate $\lambda' / \mu'$ is a vertical strip.
The profile of $\lambda'$ is obtained from the profile of $\lambda$ by reversing the string and interchanging the role of zeros and ones.

\section{Plane partitions}\label{sec:pp}
Throughout this section, all labels are assumed to be non-negative integers.

\subsection{Regular plane partitions}



\index{Plane partition!regular}

A \emph{regular plane partition} is a labelled rectangle whose labels are weakly decreasing from north to south, and from east to west. 
The \emph{weight} of a regular partition is the sum of the labels. For example, the following regular plane partition has weight $22$:
\[ \tableau{4 & 3 & 2 & 2 & 0 \\ 3 & 2 & 1 & 1 & 0 \\ 1 & 1 & 1 & 0 & 0\\ 1 & 0 & 0 & 0 & 0 \\ 0 & 0 & 0 & 0 & 0} \]

In general we may think of the rectangle as extending infinitely in both the eastern and southern directions. We require that only a finite number of labels are non-zero.

Reading from right to left along the $NW \to SE$ diagonals we obtain a sequence of integer partitions which differ each by a horizontal strip. The sequence increases first, and then decreases:
\[ \mathfrak{c} = (\emptyset, (2), (2,1), (3,1), (4,2,1), (3,1), (1), (1), \emptyset ) \]


\index{MacMahon identity}

The generating series for regular plane partitions is given by MacMahon's famous formula:
\begin{equation} \label{macmahon}
\sum_{\mathfrak{c} \in \PP} z^{|\mathfrak{c}|} = \prod_{n \geq 1} \left ( \frac{1}{1-z^n} \right )^n 
\end{equation}

The right hand side of MacMahon's identity can be expressed in terms of \emph{hook lengths} as follows:
\[ \prod_{n \geq 1} \left ( \frac{1}{1-z^n} \right )^n = \prod_{i,j \geq 1} \frac{1}{1-z^{i+j-1}} = \prod_{s \in \square} \frac{1}{1-z^{h_\square(s)}} \]
The third product is over all boxes of the infinite rectangle which we denote by the symbol $\square$. 

\index{Arbitrarily labelled!rectangle}

Let us define an \emph{arbitrarily labelled rectangle} to be a labelled rectangle with no conditions whatsoever on the labels. For example:
\[\tableau{0 & 1 & 0 & 0 \\
0 & 2 & 2 & 0 \\
3 & 0 & 0 & 0 \\
0 & 0 & 0 & 1 \\}
\]

The \emph{weight} of an arbitrarily labelled rectangle is given by a sum over the boxes of the diagram of [the label of the box] times [the hook length of the box]. In our example above the weight is given by:
\[
1 *2 + 2 *3 + 2*4 + 3*3 + 1 * 7 = 32
\]

The right hand side of MacMahon's identity may be interpreted as a weighted sum over all arbitrarily labelled rectangles.

\subsection{Reverse plane partitions}

\index{Plane partition!reverse}

A \emph{reverse plane partition} is a labelled Young diagram with the property that the labels are weakly decreasing in both the eastern and southern directions. For example:

\[
\tableau{
\missingcell & \missingcell & \missingcell & \missingcell  & \missingcell  \\
\missingcell & \missingcell & \missingcell  & 3  & 3 \\
\missingcell & \missingcell & 4 & 3 & 3 \\
\missingcell  & \missingcell & 4 & 2 & 1  \\
 4 & 4 &3 &1 &1 \\
} 
\]

Reading from right to left along the $NW \to SE$ diagonals we have the following sequence of integer partitions:
\[ \emptyset \rightharpoonup (3) \rightharpoonup (3,3) \leftharpoonup (3,1) \rightharpoonup (4,2,1) \leftharpoonup (4,1) \leftharpoonup (3) \rightharpoonup (4) \rightharpoonup (4) \leftharpoonup \emptyset \]

Let $\lambda$ be the shape of the Young diagram, and let $\pi$ be its profile. An equivalent definition for a reverse plane partition of shape $\lambda$ is the following:
\begin{definition} \label{def:rpp}
For any binary string $\pi$ of length $T$,
a \emph{reverse plane partition} with profile $\pi$ is a sequence of integer partitions:
\begin{equation*} 
(\emptyset, \mu^1, \mu^2, \ldots \mu^T,\emptyset) 
\end{equation*}
such that if $\pi_k = 1$ then $\mu^k / \mu^{k-1}$ is a \emph{horizontal strip}, otherwise if $\pi_k = 0$ then $\mu^{k-1} / \mu^k$ is a horizontal strip.
\end{definition}

\index{Stanley identity}

The following formula for the enumeration of reverse plane partitions with arbitrary profile $\pi$ is due to Stanley \cite{stanley-rpp}:
\begin{equation} \label{stanley}
\sum_{\mathfrak{c} \in \RPP(\pi)} z^{|\mathfrak{c}|} = 
\prod_{\substack{i < j \\ \pi_i > \pi_j}} \frac{1}{1-z^{j-i}} 
\end{equation}

Note that the right hand side could have been written in the form:
\[
\prod_{\substack{i < j \\ \pi_i > \pi_j}} \frac{1}{1-z^{j-i}} = \prod_{s \in \lambda} \frac{1}{1-z^{h_\lambda(s)}}
\]

\index{Arbitrarily labelled!diagram}

Let us define an \emph{arbitrarily labelled diagrams of shape $\lambda$} to be a labelled Young diagram of shape $\lambda$ with no conditions whatsoever on the labels. For example:
\[
\tableau{
\missingcell & \missingcell & \missingcell & \missingcell  & \missingcell  \\
\missingcell & \missingcell & \missingcell  & 0  & 3 \\
\missingcell & \missingcell & 0 & 1 & 0 \\
\missingcell  & \missingcell & 2 & 0 & 1  \\
 0 & 0 &0 &1 &0 \\
} 
\]

The \emph{weight} of an arbitrarily labelled diagram of shape $\lambda$ is given by a sum over the boxes of the diagram of [the label of the box] times [the hook length of the box]. In our example above the weight is given by:
\[
3 *2 + 1*3 + 2*2 + 1*5 + 1 * 7 = 25
\]

The right hand side of Stanley's identity may be interpreted as a weighted sum over all arbitrarily labelled diagrams of shape $\lambda$.

\section{Cylindric Plane Partitions} \label{sec:cylindric-defs}
\subsection{Definition}\label{cpp}

\index{Plane partition!cylindric}

Cylindric plane partitions were first introduced by Gessel and Krattenthaler \cite{gessel-1997}. 
We shall work with a modified, though equivalent, definition.

\begin{definition} \label{cylindric-def}
For any binary string $\pi$ of length $T$,
a \emph{cylindric plane partition} with profile $\pi$ may be defined as a sequence of integer partitions:
\begin{equation} 
(\mu^0, \mu^1, \ldots \mu^T) \qquad \qquad \mu^0 = \mu^T 
\end{equation}
such that if $\pi_k = 1$ then $\mu^k / \mu^{k-1}$ is a \emph{horizontal strip}, otherwise if $\pi_k = 0$ then $\mu^{k-1} / \mu^k$ is a horizontal strip.
\end{definition}

In the special case where $\mu^0 = \mu^T = \emptyset$ we recover the usual definition of a reverse plane partition.

\index{Cube!interlacing sequence definition}

\begin{definition}\label{def:cube}
A \emph{cube} of a cylindric plane partition is defined to be a \emph{box} of one of the underlying integer partitions.
\end{definition}

\begin{definition}
The \emph{weight} of the cylindric partition $\mathfrak{c} = (\mu^0, \mu^1, \ldots \mu^T)$ is given by
$  |\mathfrak{c}| = |\mu^1| + |\mu^2| + \cdots |\mu^T|$.
\end{definition}

In other words, the weight of a cylindric plane partition is the number of cubes. Note that to avoid double counting, we do not include the boxes of the partition $\mu^0$ in the definition of the weight of $\mathfrak{c}$.

\subsection{Cylindric diagrams}\label{sec:cylindric-diagram}

\index{Cylindric diagram!definition}

A \emph{cylindric diagram} may be thought of as an infinite partition with periodic profile, which has been wrapped around a cylinder. Here is an example of a cylindric diagram with profile $\pi = 10100$ and period $T = 5$. The ``fundamental domain'' is coloured in yellow. 

\[
 \tableau{
\missingcell & \missingcell & \missingcell & \missingcell & \missingcell & \missingcell & \missingcell & \missingcell 1  \\
\missingcell & \missingcell & \missingcell & \missingcell & \missingcell & \missingcell & \missingcell 1&0  \\
\missingcell & \missingcell & \missingcell & \missingcell & \missingcell & \missingcell &0& \\
\missingcell & \missingcell & \missingcell & \missingcell & \missingcell & \missingcell 1&0& \\
\missingcell & \missingcell & \missingcell & \missingcell & \missingcell 1&0& \thickcell&  \\
\missingcell & \missingcell & \missingcell & \missingcell &0&&  & \thickcell\\
\missingcell & \missingcell & \missingcell & \missingcell 1&0&&  & \\
\missingcell & \missingcell & \missingcell 1&\yellowcell 0& \thickcell &  && \\
\missingcell & \missingcell &\yellowcell0&\yellowcell&\yellowcell & \thickcell&& \\
\missingcell & \missingcell 1&\yellowcell0&\yellowcell& \yellowcell &\yellowcell&\thickcell& \\
\missingcell1 &0 & \thickyellowcell & \yellowcell &\yellowcell&\yellowcell&\yellowcell& \thickcell\\
0&&& \thickyellowcell & \yellowcell &\yellowcell&\yellowcell& \yellowcell\\
0&&&& \thickyellowcell &\yellowcell&\yellowcell& \yellowcell\\
}
\]
Note that the profile is read from top to bottom, right to left. The $1$'s represent horizontal steps while the $0$'s represent vertical steps. Although we have not drawn all of it, the diagram above is to be understood to extend infinitely in the Eastern and Southern directions.

Cylindric plane partitions are often represented as labelled cylindric diagrams. For example, the cylindric plane partition
\[ \mathfrak{c} = ((3,2,2), (4,3,2,1), (4,3,2), (6,4,3,2), (5,3,2), (3,2,2))\] 
with profile $10100$
may be represented as:
\[  \tableau{
\missingcell & 4 & \thickcell{3} \\
6 & 4 & 3 & \thickcell{2} \\
5 & 4 & 3 & 2 & \thickcell{2} \\
\thickcell{3} & 3 & 3 & 2 & 1 \\
\missingcell & \thickcell{2} & 2 & 2 & 0 \\
\missingcell & \missingcell & \thickcell{2} & 0 & 0 
}
\]
The individual partitions in the interlacing sequence picture are read off the $NW \to SE$ diagonals in order from right to left.

\begin{lemma}
The labels of the cylindric diagram associated to a cylindric plane partition are weakly decreasing along both ``cylindric rows'' and ``cylindric columns''.
\begin{proof}
This is an immediate consequence of the horizontal strip condition on diagonals
\end{proof}
\end{lemma}

\index{Borodin identity}

The following hook-product formula for the enumeration of cylindric plane partitions of given profile was first given by Borodin \cite{borodin}. A very different proof involving the representation theory of $\widehat{sl}(n)$ was later given by Tingley \cite{tingley}: 
\begin{equation} \label{eq:borodin}
\sum_{\mathfrak{c} \in \CPP(\pi)} z^{|\mathfrak{c}|} = 
\prod_{n \geq 0} \left ( \frac{1}{1-z^{(n+1)T}} \prod_{\substack{i < j \\ \pi_i > \pi_j}} \frac{1}{1-z^{j-i + nT}} 
\prod_{\substack{i > j \\ \pi_i > \pi_j }} \frac{1}{1 - z^{j-i + (n+1)T}} \right ) 
\end{equation}

Here $T$ denotes the length of the profile $\pi$. 

\subsection{Cylindric inversion coordinates}\label{cylindric-inversion}

It is natural to index the boxes of the cylindric diagram via ``cylindric inversion coordinates'' $(i,j,k)$ where $\pi_i = 1$, $\pi_j = 0$ and if $j < i$ then $k \geq 1$ otherwise $k \geq 0$.
 Consider the box labelled $*$ below:
\index{Cylindric diagram!cylindric inversion coordinates}
\[
 \tableau{
\missingcell & \missingcell & \missingcell & \missingcell & \missingcell & \missingcell & \missingcell & \missingcell 1  \\
\missingcell & \missingcell & \missingcell & \missingcell & \missingcell & \missingcell & \missingcell 1&0  \\
\missingcell & \missingcell & \missingcell & \missingcell & \missingcell & \missingcell &0& \\
\missingcell & \missingcell & \missingcell & \missingcell & \missingcell & \missingcell 1&0& \\
\missingcell & \missingcell & \missingcell & \missingcell & \missingcell 1&0& \thickcell&  \\
\missingcell & \missingcell & \missingcell & \missingcell &0&&  & \thickcell\\
\missingcell & \missingcell & \missingcell & \missingcell 1&0&&  & \\
\missingcell & \missingcell & \missingcell 1&\yellowcell 0& \thickcell &  && \\
\missingcell & \missingcell &\yellowcell0&\yellowcell&\yellowcell & \thickcell&& \\
\missingcell & \missingcell 1&\yellowcell0&\yellowcell& \yellowcell &\yellowcell&\thickcell& \\
\missingcell1 &0 & \thickyellowcell & \yellowcell &\yellowcell&\yellowcell * &\yellowcell& \thickcell\\
0&&& \thickyellowcell & \yellowcell &\yellowcell&\yellowcell& \yellowcell\\
0&&&& \thickyellowcell &\yellowcell&\yellowcell& \yellowcell\\
}
\]

It corresponds to the following inversion in the infinite profile:
\[
\cdots |10100 |{\color{red}\underline{1}}0100 |10100 |1{\color{red}\underline{0}}100 | \cdots
\]

The inversion coordinates of the box marked $*$ are $(1,2,2)$. The $1$ of the inversion occurs in position $1$ of the profile. The $0$ of the inversion occurs in position $2$ of the profiles. There are two ``bars'' $|$ between the $1$ and the $0$, hence the $k$ coordinate is $2$.

Here are the cylindric inversion coordinates of each box of our example cylindric diagram:

\begin{align*}
i \mbox{ coordinate} \quad & \quad & j \mbox{ coordinate} \quad & \quad & k \mbox{ coordinate}\quad \\
\tableau{
\missingcell & 2 & \thickcell{2} \\
4 & 4 & 4 & \thickcell{4} \\
5 & 5 & 5 & 5& \thickcell{5} \\
\thickcell{2} & 2 & 2 & 2 & 2 \\
\missingcell & \thickcell{4} & 4 & 4 & 4 \\
\missingcell & \missingcell & \thickcell{5} & 5 & 5 
}
& \quad &
\tableau{
\missingcell & 1 & \thickcell{3} \\
3 & 1 & 3 & \thickcell{1} \\
3 & 1 & 3 & 1 & \thickcell{3} \\
\thickcell{3} & 1 & 3 & 1 & 3 \\
\missingcell & \thickcell{1} & 3 & 1 & 3 \\
\missingcell & \missingcell & \thickcell{3} & 1 & 3 
}
& \quad &
\tableau{
\missingcell & 0 & \thickcell{1} \\
0 & 0 & 1 & \thickcell{1} \\
0 & 0 & 1 & 1 & \thickcell{2} \\
\thickcell{1} & 1 & 2 & 2 & 3 \\
\missingcell & \thickcell{1} & 2 & 2 & 3 \\
\missingcell & \missingcell & \thickcell{2} & 2 & 3 
} 
\end{align*}

Two boxes lie in the same ``cylindric row'' if they have the same $i$-coordinate, and in the same ``cylindric column'' if they have the same $j$-coordinate. 

\[
 \tableau{
\missingcell & \missingcell & \missingcell & \missingcell & \missingcell & \missingcell & \missingcell & \missingcell 1  \\
\missingcell & \missingcell & \missingcell & \missingcell & \missingcell & \missingcell & \missingcell 1&0  \\
\missingcell & \missingcell & \missingcell & \missingcell & \missingcell & \missingcell &0& \\
\missingcell & \missingcell & \missingcell & \missingcell & \missingcell & \missingcell 1&0& \\
\missingcell & \missingcell & \missingcell & \missingcell & \missingcell 1&0& \thickcell&  \\
\missingcell & \missingcell & \missingcell & \missingcell &0&&  & \thickcell\\
\missingcell & \missingcell & \missingcell & \missingcell 1&0&&  & \\
\missingcell & \missingcell & \missingcell 1&\redcell0& \thickcell &  && \\
\missingcell & \missingcell &\bluecell0&\bluecell&\bluecell & \thickcell&& \\
\missingcell & \missingcell 1&\greencell0&\greencell& \greencell &\greencell&\thickcell& \\
\missingcell1 &0 & \thickredcell & \redcell &\redcell&\redcell&\redcell& \thickcell\\
0&&& \thickbluecell & \bluecell &\bluecell&\bluecell& \bluecell\\
0&&&& \thickgreencell &\greencell&\greencell&\greencell \\
}
\qquad \qquad
 \tableau{
\missingcell & \missingcell & \missingcell & \missingcell & \missingcell & \missingcell & \missingcell & \missingcell 1  \\
\missingcell & \missingcell & \missingcell & \missingcell & \missingcell & \missingcell & \missingcell 1&0  \\
\missingcell & \missingcell & \missingcell & \missingcell & \missingcell & \missingcell &0& \\
\missingcell & \missingcell & \missingcell & \missingcell & \missingcell & \missingcell 1&0& \\
\missingcell & \missingcell & \missingcell & \missingcell & \missingcell 1&0& \thickcell&  \\
\missingcell & \missingcell & \missingcell & \missingcell &0&&  & \thickcell\\
\missingcell & \missingcell & \missingcell & \missingcell 1&0&&  & \\
\missingcell & \missingcell & \missingcell 1& \bluecell 0& \thickcell &  && \\
\missingcell & \missingcell &\redcell0&\bluecell& \redcell & \thickcell&& \\
\missingcell & \missingcell 1&\redcell0&\bluecell&\redcell  &\bluecell&\thickcell& \\
\missingcell1 &0 & \thickredcell & \bluecell &\redcell&\bluecell&\redcell& \thickcell\\
0&&& \thickbluecell & \redcell & \bluecell & \redcell & \bluecell \\
0&&&& \thickredcell & \bluecell & \redcell & \bluecell \\
}
\]

The $k$-coordinate may be thought of as a sort of ``depth'' or ``winding number''.

\[
 \tableau{
\missingcell & \missingcell & \missingcell & \missingcell & \missingcell & \missingcell & \missingcell & \missingcell 1  \\
\missingcell & \missingcell & \missingcell & \missingcell & \missingcell & \missingcell & \missingcell 1&0  \\
\missingcell & \missingcell & \missingcell & \missingcell & \missingcell & \missingcell &0& \\
\missingcell & \missingcell & \missingcell & \missingcell & \missingcell & \missingcell 1&0& \\
\missingcell & \missingcell & \missingcell & \missingcell & \missingcell 1&0& \thickcell&  \\
\missingcell & \missingcell & \missingcell & \missingcell &0&&  & \thickcell\\
\missingcell & \missingcell & \missingcell & \missingcell 1&0&&  & \\
\missingcell & \missingcell & \missingcell 1& \redcell 0& \thickcell &  && \\
\missingcell & \missingcell &\redcell0&\redcell& \bluecell & \thickcell&& \\
\missingcell & \missingcell 1&\redcell0&\redcell&\bluecell  &\bluecell&\thickcell& \\
\missingcell1 &0 & \thickbluecell & \bluecell &\greencell&\greencell&\yellowcell& \thickcell\\
0&&& \thickbluecell & \greencell & \greencell & \yellowcell & \yellowcell \\
0&&&& \thickgreencell & \greencell & \yellowcell & \yellowcell \\
} \\ \\
\]

\subsection{Cylindric hook length}

\index{Cylindric diagram!cylindric hook length}

The \emph{cylindric hook length} of a box is the hook length of the box relative to the larger partition.

\[
 \tableau{
\missingcell & \missingcell & \missingcell & \missingcell & \missingcell & \missingcell & \missingcell & \missingcell 1  \\
\missingcell & \missingcell & \missingcell & \missingcell & \missingcell & \missingcell & \missingcell 1&0  \\
\missingcell & \missingcell & \missingcell & \missingcell & \missingcell & \missingcell &0& \\
\missingcell & \missingcell & \missingcell & \missingcell & \missingcell & \missingcell 1&0& \\
\missingcell & \missingcell & \missingcell & \missingcell & \missingcell 1&\yellowcell0& \thickcell&  \\
\missingcell & \missingcell & \missingcell & \missingcell &0&\yellowcell&  & \thickcell\\
\missingcell & \missingcell & \missingcell & \missingcell 1&0&\yellowcell&  & \\
\missingcell & \missingcell & \missingcell 1&0& \thickcell &\yellowcell  && \\
\missingcell & \missingcell &0&& & \thickyellowcell&& \\
\missingcell & \missingcell 1&0&&  &\yellowcell&\thickcell& \\
\missingcell1 & \yellowcell 0 &  \thickyellowcell & \yellowcell &\yellowcell&\redcell 11&& \thickcell\\
0&&& \thickcell &  &&& \\
0&&&& \thickcell &&& \\
}
\]

We shall use the notation $h_{\hat{\lambda}(\pi)}(b)$ to denote the cylindric hook length of the box $b$ relative to the cylindric diagram $\hat{\lambda}(\pi)$

\begin{lemma}\label{cylindrichook}
The cylindric hook length of a box with cylindric inversion coordinates $(i,j,k)$ is given by $j - i + kT$.
\end{lemma}

Here are the hook lengths of the boxes in our example partition:

\[  \tableau{
\missingcell & 1 & \thickcell{4} \\
1 & 3 & 6 & \thickcell{8} \\
2 & 4 & 7 & 9 & \thickcell{12} \\
\thickcell{4} & 6 & 9 & 11 & 14 \\
\missingcell & \thickcell{8} & 11 & 13 & 16 \\
\missingcell & \missingcell & \thickcell{12} & 14 & 17 
}
\]

As a consequence of Lemma \ref{cylindrichook} Borodin's identity (equation \ref{eq:borodin}) may be rewritten in the form:
\begin{align}
\sum_{\mathfrak{c} \in \CPP(\pi)} z^{|\mathfrak{c}|} & = \left ( \sum_{\gamma} z^{T\, |\gamma|} \right ) 
\left ( \sum_{s \in \widehat{\lambda}(\pi)} \frac{1}{1 - z^{h_{\widehat{\lambda}(\pi)}(s)}} \right ) 
\end{align}

\subsection{Arbitrarily labelled cylindric diagrams}\label{alcd}

\index{Arbitrarily labelled!cylindric diagram}

An \emph{arbitrarily labelled cylindric diagram} $\mathfrak{d}$ with profile $\pi$ is simply an assignment of non-negative integers to the boxes of the associated cylindrical diagram
in such a way that only finitely many of the labels are non-zero. We shall use the notation $\ALCD(\pi)$ to denote the set of all arbitrarily labelled cylindric diagrams with profile $\pi$.

\begin{definition} \label{depth}
The \emph{depth} of an arbitrarily labelled cylindric diagram $\mathfrak{d}$ is the smallest $k$ such that all boxes with cylindric inversion coordinates $(i,j,k')$ with $k' \geq k$ have label zero.
\end{definition}

The weight of an arbitrarily labelled cylindric diagram is given by the sum over boxes in the cylindric diagram of the label of the box times the hook length of the box:
\begin{equation} |\mathfrak{d}| = \sum_{b \in \hat{\lambda}(\pi)} \lab(b) \,\, |\hook(b)| \end{equation}

For example, the following arbitrarily labelled cylindric diagram has depth 2 and weight 26.

\[ \tableau{
\missingcell & 1 & \thickcell{0} \\
5 & 0 & 0 & \thickcell{0} \\
0 & 1 & 1 & 1 & \thickcell{0} \\
\thickcell{0} & 0 & 0 & 0 & 0 \\
\missingcell & \thickcell{0} & 0 & 0 & 0 \\
\missingcell & \missingcell & \thickcell{0} & 0 & 0 
}
\]

The right hand side of Borodin's identity may be interpreted combinatorially as a weighted sum over pairs $(\gamma, \mathfrak{d})$
where $\gamma$ is an integer partition and $\mathfrak{d}$ is an arbitrarily labelled cylindric diagram with profile $\pi$.
The weight of the pair $(\gamma, \mathfrak{d})$ is given by $T|\gamma| + |\mathfrak{d}|$.
where $\gamma$ is an integer partition and $\mathfrak{d}$ is an \emph{arbitrarily labelled cylindric diagrams}

\subsection{Rotation operator}

\index{Cylindric diagram!rotation operator}

\begin{definition}\label{shift}
The rotation operator on binary strings is defined by:
\begin{equation} \sigma(\pi)_i = \pi_{(i + 1) \mod T} \end{equation}
where $T$ is the length of $\pi$.
\end{definition}

There is a natural bijection between the cylindric diagram with profile $\pi$ and the 
cylindric diagram with profile $\sigma(\pi)$. Nevertheless since the same box will have
different cylindric inversion coordinates, depending on the choice of rotation of the profile, 
we prefer to consider two cylindric diagrams which differ by a rotation to be two distinct objects.
Likewise for arbitrarily labelled cylindric diagrams.

\begin{lemma}
The rotation operator $\sigma$ naturally induces well-defined weight preserving maps:
\begin{align*} 
\sigma : \ALCD(\pi) & \to \ALCD(\sigma(\pi)) \\
\sigma : \CPP(\pi) & \to \CPP(\sigma(\pi))
\end{align*}
\end{lemma}

\chapter{Correspondences}

Throughout this section, all labels are assumed to be non-negative integers.

\section{Robinson correspondence}\label{sec:robinson}
\subsection{Standard Young Tableau}\label{sec:SSYT}

\tikzstyle{gray}=[circle,fill=black,draw,inner sep = 0.5pt]
\tikzstyle{black}=[circle,fill=black,draw,inner sep = 3pt]
\tikzstyle{blue}=[circle,fill=blue,draw,inner sep = 2pt]
\tikzstyle{red}=[circle,fill=red,draw,inner sep = 2pt]
\tikzstyle{green}=[circle,fill=green,draw,inner sep = 2pt]
\tikzstyle{purple}=[circle,fill=purple,draw,inner sep = 2pt]
\tikzstyle{yellow}=[circle,fill=yellow,draw,inner sep = 2pt]
\tikzstyle{orange}=[circle,fill=orange,draw,inner sep = 2pt]
\tikzstyle{pink}=[circle,fill=pink,draw,inner sep = 2pt]

\index{Tableau!standard}

A \emph{standard tableau} is a labelled Young diagram with $n$ boxes for which the labels are strictly decreasing along both columns and rows,
and for which each number from $1$ to $n$ occurs exactly once.
For example:

\[
\tableau{
\missingcell & \missingcell & \missingcell & \missingcell  & \missingcell  \\
\missingcell & \missingcell & \missingcell  & 13 & 10 \\
\missingcell & \missingcell & 11 & 9& 7 \\
\missingcell  & \missingcell & 8 & 6 & 3 \\
12 & 5&4&2&1\\
} \\
\]

\index{Yamanouchi word}

A \emph{Yamanouchi word} is a word with the property that for any initial subword, and for any $k$, the number of times that $k$ appears in that initial subword is greater than or equal to the number of times that $k+1$ appears in that subword.

Standard tableau are naturally in bijection with Yamanouchi words. For example, the tableau above corresponds to the word:
\[ 1121123234314 \]
The $k$-th number in the Yamanouchi word tells you on which row you will find the label $k$. Note that we are counting the rows from bottom to top.

\subsection{Viennot's shadow method}\label{sec:viennot}

\index{Robinson correspondence!regular}

A \emph{permutation} is a $\{0,1\}$ matrix with exactly one $1$ in each row and each column.
A \emph{partial permutation} is a $\{0,1\}$ matrix with at most one $1$ in each row and each column.
Note that a partial permutation matrix need not be square.
The \emph{Robinson correspondence} \cite{Robinson} gives a bijection between permutations and pairs of standard Young tableau of the same shape.

Instead of defining the algorithm rigorously, we shall just give an example using \emph{Viennot's shadow technique} \cite{shadow}. We begin with the matrix of our permutation, but transformed into the grid.
In our case, we have chosen:
\[
\sigma = \left [ \begin{matrix} 1 & 2 & 3 & 4 & 5 & 6 & 7 \\ 5 & 3 & 6 & 1 & 4 & 7 & 2 \end{matrix} \right ]
\]
Here is the grid. The large black dots correspond to ones in the permutation matrix, while the small black dots correspond to zeros.
\begin{center}
\begin{tikzpicture}[scale=0.8]
\foreach \x in {1,...,7}
\foreach \y in {1,...,7}
\node  at (\x,\y) [gray]{};

\node (Black1) at (1,3) [black]{};
\node (Black2) at (2,5) [black]{};
\node (Black3) at (4,7) [black]{};

\node (Blue1) at (2,3) []{};
\node (Blue2) at (4,5) []{};

\node (Black4) at (3,2) [black]{};
\node (Black5) at (5,4) [black]{};
\node (Black6) at (7,6) [black]{};

\node (Black7) at (6,1) [black]{};

\end{tikzpicture}
\end{center}

Next imagine that the sun is shining down from the top left hand corner, and that it casts a shadow:

\begin{center}
\begin{tikzpicture}[scale=0.8]
\foreach \x in {1,...,7}
\foreach \y in {1,...,7}
\node  at (\x,\y) [gray]{};

\node (Black1) at (1,3) [black]{};
\node (Black2) at (2,5) [black]{};
\node (Black3) at (4,7) [black]{};

\node (Blue1) at (2,3) [blue]{};
\node (Blue2) at (4,5) [blue]{};

\draw[draw=black,thick] (1,0) -- (Black1)-- (Blue1) -- (Black2) -- (Blue2) -- (Black3) -- (8,7);

\node (Black4) at (3,2) [black]{};
\node (Black5) at (5,4) [black]{};
\node (Black6) at (7,6) [black]{};

\node (Blue3) at (5,2) [blue]{};
\node (Blue4) at (7,4) [blue]{};

\draw[draw=black,thick] (3,0) -- (Black4)-- (Blue3) -- (Black5) -- (Blue4) -- (Black6) -- (8,6);

\node (Black7) at (6,1) [black]{};

\draw[draw=black,thick] (6,0) -- (Black7) -- (8,1);

\end{tikzpicture}
\end{center}

We have essentially regrouped the permutation by decreasing sequences:
\[
\left [ \begin{matrix} 1 & 2 & 4 \\ 5 & 3 & 1 \end{matrix} \right ] \left [ \begin{matrix} 3 & 5 & 7 \\ 6 & 4 & 2 \end{matrix} \right ] \left [ \begin{matrix} 6 \\ 7 \end{matrix} \right ]
\]

Remove the first number from the top row and the last letter from the bottom row of each block then shift everything across to obtain a partial permutation:
\[
\left [ \begin{matrix} . & 2 & 4 \\ 5 & 3 & . \end{matrix} \right ] \left [ \begin{matrix} . & 5 & 7 \\ 6 & 4 & . \end{matrix} \right ] \left [ \begin{matrix} . \\ . \end{matrix} \right ]
\mapsto
\left [ \begin{matrix} 2 & 4 \\ 5 & 3 \end{matrix} \right ] \left [ \begin{matrix} 5 & 7 \\ 6 & 4 \end{matrix} \right ] 
\mapsto
\left [ \begin{matrix} 1 & 2 & 3 & 4 & 5 & 6 & 7 \\ * & 5 & * & 3 & 6 & * & 4\end{matrix} \right ]
\]

Now we repeat the process, beginning by drawing the partial permutation onto the grid:

\begin{center}
\begin{tikzpicture}[scale=0.8]
\foreach \x in {1,...,7}
\foreach \y in {1,...,7}
\node  at (\x,\y) [gray]{};

\node (Blue1) at (2,3) [blue]{};
\node (Blue2) at (4,5) [blue]{};

\node (Blue3) at (5,2) [blue]{};
\node (Blue4) at (7,4) [blue]{};

\end{tikzpicture}
\end{center}

The sun shines down from the top left hand corner and causes shadows:

\begin{center}
\begin{tikzpicture}[scale=0.8]
\foreach \x in {1,...,7}
\foreach \y in {1,...,7}
\node  at (\x,\y) [gray]{};

\node (Blue1) at (2,3) [blue]{};
\node (Blue2) at (4,5) [blue]{};

\node (Red1) at (4,3) [red]{};

\draw[draw=blue,thick] (2,0) -- (Blue1)-- (Red1) -- (Blue2) --  (8,5);

\node (Blue3) at (5,2) [blue]{};
\node (Blue4) at (7,4) [blue]{};

\node (Red2) at (7,2) [red]{};

\draw[draw=blue,thick] (5,0) -- (Blue3)-- (Red2) -- (Blue4)  -- (8,4);

\end{tikzpicture}
\end{center}

Regrouping the partial permutation into decreasing sequences:
\[
\left [ \begin{matrix} 2 & 4 \\ 5 & 3 \end{matrix} \right ] \left [ \begin{matrix}5 & 7 \\ 6 & 4 \end{matrix} \right ]
\]

Extracting the partial permutation:
\[
\left [ \begin{matrix} . & 4 \\ 5 & . \end{matrix} \right ] \left [ \begin{matrix}. & 7 \\ 6 & . \end{matrix} \right ]
\mapsto
\left [ \begin{matrix} 4 \\ 5 \end{matrix} \right ] \left [ \begin{matrix} 7 \\ 6  \end{matrix} \right ]
\mapsto
\left [ \begin{matrix} 1 & 2 & 3 & 4 & 5 & 6 & 7 \\ * & * & * & 5 & * & * & 6\end{matrix}\right ]
\]

And once more we draw the partial permutation:

\begin{center}
\begin{tikzpicture}[scale=0.8]
\foreach \x in {1,...,7}
\foreach \y in {1,...,7}
\node  at (\x,\y) [gray]{};

\node (Red1) at (4,3) [red]{};

\node (Red2) at (7,2) [red]{};

\end{tikzpicture}
\end{center}

The sun shines down from the top left hand corner, casting shadows:

\begin{center}
\begin{tikzpicture}[scale=0.8]
\foreach \x in {1,...,7}
\foreach \y in {1,...,7}
\node  at (\x,\y) [gray]{};

\node (Red1) at (4,3) [red]{};

\draw[draw=red,thick] (4,0) -- (Red1) --  (8,3);

\node (Red2) at (7,2) [red]{};

\draw[draw=red,thick] (7,0) -- (Red2) --  (8,2);

\end{tikzpicture}
\end{center}

Here is where the algorithm terminates, since the next partial permutation in the sequence is empty:
\[
\left [ \begin{matrix} 4 \\ 5 \end{matrix} \right ] \left [ \begin{matrix} 7 \\ 6 \end{matrix}\right ]
\mapsto
\left [ \begin{matrix} . \\ . \end{matrix} \right ] \left [ \begin{matrix} . \\ . \end{matrix}\right ]
\mapsto
\left [ \begin{matrix} 1 & 2 & 3 & 4 & 5 & 6 & 7 \\ * & * & * & * & * & * & * \end{matrix} \right ]
\]

Putting all this information together, we obtain the following \emph{Viennot shadow diagram}:

\index{Viennot shadow}

\begin{center}
\begin{tikzpicture}[scale=0.8]
\foreach \x in {1,...,7}
\foreach \y in {1,...,7}
\node  at (\x,\y) [gray]{};

\node at (1,-1) {1};
\node at (2,-1) {\color{blue}2};
\node at (3,-1) {1};
\node at (4,-1) {\color{red}3};
\node at (5,-1) {\color{blue}2};
\node at (6,-1) {1};
\node at (7,-1) {\color{red}3};

\node at (9,7) {1};
\node at (9,6) {1};
\node at (9,5) {\color{blue}2};
\node at (9,4) {\color{blue}2};
\node at (9,3) {\color{red}3};
\node at (9,2) {\color{red}3};
\node at (9,1) {1};

\node (Black1) at (1,3) [black]{};
\node (Black2) at (2,5) [black]{};
\node (Black3) at (4,7) [black]{};

\node (Blue1) at (2,3) [blue]{};
\node (Blue2) at (4,5) [blue]{};

\node (Red1) at (4,3) [red]{};

\draw[draw=black,thick] (1,0) -- (Black1)-- (Blue1) -- (Black2) -- (Blue2) -- (Black3) -- (8,7);

\draw[draw=blue,thick] (2,0) -- (Blue1)-- (Red1) -- (Blue2) --  (8,5);

\draw[draw=red,thick] (4,0) -- (Red1) --  (8,3);

\node (Black4) at (3,2) [black]{};
\node (Black5) at (5,4) [black]{};
\node (Black6) at (7,6) [black]{};

\node (Blue3) at (5,2) [blue]{};
\node (Blue4) at (7,4) [blue]{};

\node (Red2) at (7,2) [red]{};

\draw[draw=black,thick] (3,0) -- (Black4)-- (Blue3) -- (Black5) -- (Blue4) -- (Black6) -- (8,6);

\draw[draw=blue,thick] (5,0) -- (Blue3)-- (Red2) -- (Blue4)  -- (8,4);

\draw[draw=red,thick] (7,0) -- (Red2) --  (8,2);

\node (Black7) at (6,1) [black]{};

\draw[draw=black,thick] (6,0) -- (Black7) -- (8,1);

\end{tikzpicture}
\end{center}

The pair of standard tableau which we are searching for can now be read, as Yamanouchi words from the right most column
and the bottom most row.
The black path is labelled $1$, the blue path is labelled $2$ and the red path is labelled $3$.
Reading down the final column we have:
\[
1  1  2  2  3  3 1
\]
while reading along the bottom row we have:
\[
1 2 1 3 2 1 3
\]
Thus, under the Robinson correspondence, we have:
\[
\left [ \begin{matrix} 1 & 2 & 3 & 4 & 5 & 6 & 7 \\ 5 & 3 & 6 & 1 & 4 & 7 & 2 \end{matrix} \right ]
\mapsto
\left ( \quad
\tableau{\missingcell & 6 & 5 \\ \missingcell & 4 & 3 \\ 7 & 2 & 1 \\}
\quad , \quad 
\tableau{\missingcell & 7 & 4 \\ \missingcell & 5 & 2 \\ 6 & 3 & 1} \quad
\right )
\]

\subsection{Reverse Viennot shadow method}\label{sec:viennot-reverse}

To reverse the algorithm, we begin with the two Yamanouchi words, and the empty permutation:
\[ ((1,1,2,2,3,3,1), (1,2,1,3,2,1,3))  \]
The largest number, the threes, occur in positions ${\color{blue}5}$ and $\color{blue}6$ on the left, and positions $\color{red}4$ and $\color{red}7$ on the right.
This gives us:
\[ ((1,1,2,2,.,.,1), (1,2,1,.,2,1,.))  \qquad \left [ \begin{matrix} {\color{red}4} \\ {\color{blue}5} \end{matrix} \right ] \left [ \begin{matrix} {\color{red}7} \\ {\color{blue}6} \end{matrix}\right ]\]
Next we note that the twos occur at positions $\color{blue}3$ and $\color{blue}4$ on the left, and positions $\color{red}2$ and $\color{red}5$ on the right:
\[
\left [ \begin{matrix} 4 \\ 5 \end{matrix} \right ] \left [ \begin{matrix} 7 \\ 6  \end{matrix} \right ]
\mapsto
\left [ \begin{matrix} . & 4 \\ 5 & . \end{matrix} \right ] \left [ \begin{matrix}. & 7 \\ 6 & . \end{matrix} \right ]
\mapsto
\left [ \begin{matrix} {\color{red}2} & 4 \\ 5 & {\color{blue}3} \end{matrix} \right ] \left [ \begin{matrix}{\color{red}5} & 7 \\ 6 & {\color{blue}4} \end{matrix} \right ]
\]
We are left with:
\[ ((1,1,.,.,.,.,1), (1,.,1,.,.,1,.)) \]
The ones on the left occur at position $\color{blue}1$,$\color{blue}2$ and $\color{blue}7$ while the ones on the right occur at positions $\color{red}1$, $\color{red}3$ and $\color{red}6$.
This gives us:

\[
\left [ \begin{matrix} 2 & 4 \\ 5 & 3 \end{matrix} \right ] \left [ \begin{matrix}5 & 7 \\ 6 & 4 \end{matrix} \right ]
\mapsto
\left [ \begin{matrix} . & 2 & 4 \\ 5 & 3 & . \end{matrix} \right ] \left [ \begin{matrix} . & 5 & 7 \\ 6 & 4 & . \end{matrix} \right ] \left [ \begin{matrix} . \\ . \end{matrix} \right ]
\mapsto
\left [ \begin{matrix} {\color{red}1} & 2 & 4 \\ 5 & 3 & {\color{blue}1} \end{matrix} \right ] \left [ \begin{matrix} {\color{red}3} & 5 & 7 \\ 6 & 4 & {\color{blue}2} \end{matrix} \right ] \left [ \begin{matrix} {\color{red}6} \\ {\color{blue}7} \end{matrix} \right ]
\]
The final result is the permutation:
\[
\sigma = \left [ \begin{matrix} 1 & 2 & 3 & 4 & 5 & 6 & 7 \\ 5 & 3 & 6 & 1 & 4 & 7 & 2 \end{matrix} \right ]
\]

\subsection{Fomin growth diagrams}\label{sec:fomin}

Viennot's shadow method is strictly equivalent to \emph{Fomin growth diagram} technique \cite{fomin-growth}.

We begin by labelling the faces of the first and last rows of our shadow diagram with the empty partition:

\begin{center}
\begin{tikzpicture}[scale=1.2]
\foreach \x in {1,...,7}
\foreach \y in {1,...,7}
\node  at (\x,\y) [gray]{};

\node at (0.5,0.5) []{$()$};
\node at (0.5,1.5) []{$()$};
\node at (0.5,2.5) []{$()$};
\node at (0.5,3.5) []{$()$};
\node at (0.5,4.5) []{$()$};
\node at (0.5,5.5) []{$()$};
\node at (0.5,6.5) []{$()$};
\node at (0.5,7.5) []{$()$};
\node at (1.5,7.5) []{$()$};
\node at (2.5,7.5) []{$()$};
\node at (3.5,7.5) []{$()$};
\node at (4.5,7.5) []{$()$};
\node at (5.5,7.5) []{$()$};
\node at (6.5,7.5) []{$()$};
\node at (7.5,7.5) []{$()$};

\node (Black1) at (1,3) [black]{};
\node (Black2) at (2,5) [black]{};
\node (Black3) at (4,7) [black]{};

\node (Blue1) at (2,3) [blue]{};
\node (Blue2) at (4,5) [blue]{};

\node (Red1) at (4,3) [red]{};

\draw[draw=black,thick] (1,0) -- (Black1)-- (Blue1) -- (Black2) -- (Blue2) -- (Black3) -- (8,7);

\draw[draw=blue,thick] (2,0) -- (Blue1)-- (Red1) -- (Blue2) --  (8,5);

\draw[draw=red,thick] (4,0) -- (Red1) --  (8,3);

\node (Black4) at (3,2) [black]{};
\node (Black5) at (5,4) [black]{};
\node (Black6) at (7,6) [black]{};

\node (Blue3) at (5,2) [blue]{};
\node (Blue4) at (7,4) [blue]{};

\node (Red2) at (7,2) [red]{};

\draw[draw=black,thick] (3,0) -- (Black4)-- (Blue3) -- (Black5) -- (Blue4) -- (Black6) -- (8,6);

\draw[draw=blue,thick] (5,0) -- (Blue3)-- (Red2) -- (Blue4)  -- (8,4);

\draw[draw=red,thick] (7,0) -- (Red2) --  (8,2);

\node (Black7) at (6,1) [black]{};

\draw[draw=black,thick] (6,0) -- (Black7) -- (8,1);

\end{tikzpicture}
\end{center}

Next we continue filling in the faces of the diagram in such a way that a new box is added to the first row every time that one crosses one of the black shadows:

\begin{center}
\begin{tikzpicture}[scale=1.2]
\foreach \x in {1,...,7}
\foreach \y in {1,...,7}
\node  at (\x,\y) [gray]{};

\node at (0.5,0.5) []{$()$};
\node at (0.5,1.5) []{$()$};
\node at (0.5,2.5) []{$()$};
\node at (0.5,3.5) []{$()$};
\node at (0.5,4.5) []{$()$};
\node at (0.5,5.5) []{$()$};
\node at (0.5,6.5) []{$()$};
\node at (0.5,7.5) []{$()$};
\node at (1.5,7.5) []{$()$};
\node at (2.5,7.5) []{$()$};
\node at (3.5,7.5) []{$()$};
\node at (4.5,7.5) []{$()$};
\node at (5.5,7.5) []{$()$};
\node at (6.5,7.5) []{$()$};
\node at (7.5,7.5) []{$()$};

\node at (1.5,6.5) []{$()$};
\node at (2.5,6.5) []{$()$};
\node at (3.5,6.5) []{$()$};
\node at (4.5,6.5) []{$(1)$};
\node at (5.5,6.5) []{$(1)$};
\node at (6.5,6.5) []{$(1)$};
\node at (7.5,6.5) []{$(1)$};

\node at (1.5,5.5) []{$()$};
\node at (2.5,5.5) []{$()$};
\node at (3.5,5.5) []{$()$};
\node at (4.5,5.5) []{$(1)$};
\node at (5.5,5.5) []{$(1)$};
\node at (6.5,5.5) []{$(1)$};
\node at (7.5,5.5) []{$(2)$};

\node at (1.5,4.5) []{$()$};
\node at (2.5,4.5) []{$(1)$};
\node at (3.5,4.5) []{$(1)$};

\node at (1.5,3.5) []{$()$};
\node at (2.5,3.5) []{$(1)$};
\node at (3.5,3.5) []{$(1)$};

\node at (1.5,2.5) []{$(1)$};

\node at (1.5,1.5) []{$(1)$};

\node at (1.5,0.5) []{$(1)$};

\node (Black1) at (1,3) [black]{};
\node (Black2) at (2,5) [black]{};
\node (Black3) at (4,7) [black]{};

\node (Blue1) at (2,3) [blue]{};
\node (Blue2) at (4,5) [blue]{};

\node (Red1) at (4,3) [red]{};

\draw[draw=black,thick] (1,0) -- (Black1)-- (Blue1) -- (Black2) -- (Blue2) -- (Black3) -- (8,7);

\draw[draw=blue,thick] (2,0) -- (Blue1)-- (Red1) -- (Blue2) --  (8,5);

\draw[draw=red,thick] (4,0) -- (Red1) --  (8,3);

\node (Black4) at (3,2) [black]{};
\node (Black5) at (5,4) [black]{};
\node (Black6) at (7,6) [black]{};

\node (Blue3) at (5,2) [blue]{};
\node (Blue4) at (7,4) [blue]{};

\node (Red2) at (7,2) [red]{};

\draw[draw=black,thick] (3,0) -- (Black4)-- (Blue3) -- (Black5) -- (Blue4) -- (Black6) -- (8,6);

\draw[draw=blue,thick] (5,0) -- (Blue3)-- (Red2) -- (Blue4)  -- (8,4);

\draw[draw=red,thick] (7,0) -- (Red2) --  (8,2);

\node (Black7) at (6,1) [black]{};

\draw[draw=black,thick] (6,0) -- (Black7) -- (8,1);

\end{tikzpicture}
\end{center}

Keep ``growing'' the diagram in such a way that every time you cross a blue line, you add a box to the second row:

\begin{center}
\begin{tikzpicture}[scale=1.2]
\foreach \x in {1,...,7}
\foreach \y in {1,...,7}
\node  at (\x,\y) [gray]{};

\node at (0.5,0.5) []{$()$};
\node at (0.5,1.5) []{$()$};
\node at (0.5,2.5) []{$()$};
\node at (0.5,3.5) []{$()$};
\node at (0.5,4.5) []{$()$};
\node at (0.5,5.5) []{$()$};
\node at (0.5,6.5) []{$()$};
\node at (0.5,7.5) []{$()$};

\node at (1.5,7.5) []{$()$};
\node at (2.5,7.5) []{$()$};
\node at (3.5,7.5) []{$()$};
\node at (4.5,7.5) []{$()$};
\node at (5.5,7.5) []{$()$};
\node at (6.5,7.5) []{$()$};
\node at (7.5,7.5) []{$()$};

\node at (1.5,6.5) []{$()$};
\node at (2.5,6.5) []{$()$};
\node at (3.5,6.5) []{$()$};
\node at (4.5,6.5) []{$(1)$};
\node at (5.5,6.5) []{$(1)$};
\node at (6.5,6.5) []{$(1)$};
\node at (7.5,6.5) []{$(1)$};

\node at (1.5,5.5) []{$()$};
\node at (2.5,5.5) []{$()$};
\node at (3.5,5.5) []{$()$};
\node at (4.5,5.5) []{$(1)$};
\node at (5.5,5.5) []{$(1)$};
\node at (6.5,5.5) []{$(1)$};
\node at (7.5,5.5) []{$(2)$};

\node at (1.5,4.5) []{$()$};
\node at (2.5,4.5) []{$(1)$};
\node at (3.5,4.5) []{$(1)$};
\node at (4.5,4.5) []{$(1,1)$};
\node at (5.5,4.5) []{$(1,1)$};
\node at (6.5,4.5) []{$(1,1)$};
\node at (7.5,4.5) []{$(2,1)$};

\node at (1.5,3.5) []{$()$};
\node at (2.5,3.5) []{$(1)$};
\node at (3.5,3.5) []{$(1)$};
\node at (4.5,3.5) []{$(1,1)$};
\node at (5.5,3.5) []{$(2,1)$};
\node at (6.5,3.5) []{$(2,1)$};
\node at (7.5,3.5) []{$(2,2)$};

\node at (1.5,2.5) []{$(1)$};
\node at (2.5,2.5) []{$(1,1)$};
\node at (3.5,2.5) []{$(1,1)$};

\node at (1.5,1.5) []{$(1)$};
\node at (2.5,1.5) []{$(1,1)$};
\node at (3.5,1.5) []{$(2,1)$};

\node at (1.5,0.5) []{$(1)$};
\node at (2.5,0.5) []{$(1,1)$};
\node at (3.5,0.5) []{$(2,1)$};

\node (Black1) at (1,3) [black]{};
\node (Black2) at (2,5) [black]{};
\node (Black3) at (4,7) [black]{};

\node (Blue1) at (2,3) [blue]{};
\node (Blue2) at (4,5) [blue]{};

\node (Red1) at (4,3) [red]{};

\draw[draw=black,thick] (1,0) -- (Black1)-- (Blue1) -- (Black2) -- (Blue2) -- (Black3) -- (8,7);

\draw[draw=blue,thick] (2,0) -- (Blue1)-- (Red1) -- (Blue2) --  (8,5);

\draw[draw=red,thick] (4,0) -- (Red1) --  (8,3);

\node (Black4) at (3,2) [black]{};
\node (Black5) at (5,4) [black]{};
\node (Black6) at (7,6) [black]{};

\node (Blue3) at (5,2) [blue]{};
\node (Blue4) at (7,4) [blue]{};

\node (Red2) at (7,2) [red]{};

\draw[draw=black,thick] (3,0) -- (Black4)-- (Blue3) -- (Black5) -- (Blue4) -- (Black6) -- (8,6);

\draw[draw=blue,thick] (5,0) -- (Blue3)-- (Red2) -- (Blue4)  -- (8,4);

\draw[draw=red,thick] (7,0) -- (Red2) --  (8,2);

\node (Black7) at (6,1) [black]{};

\draw[draw=black,thick] (6,0) -- (Black7) -- (8,1);

\end{tikzpicture}
\end{center}

When you cross a red line, you add a box to the third row:

\index{Fomin growth diagram}

\begin{center}
\begin{tikzpicture}[scale=1.4]
\foreach \x in {1,...,7}
\foreach \y in {1,...,7}
\node  at (\x,\y) [gray]{};

\node at (0.5,0.5) []{$()$};
\node at (0.5,1.5) []{$()$};
\node at (0.5,2.5) []{$()$};
\node at (0.5,3.5) []{$()$};
\node at (0.5,4.5) []{$()$};
\node at (0.5,5.5) []{$()$};
\node at (0.5,6.5) []{$()$};
\node at (0.5,7.5) []{$()$};

\node at (1.5,7.5) []{$()$};
\node at (2.5,7.5) []{$()$};
\node at (3.5,7.5) []{$()$};
\node at (4.5,7.5) []{$()$};
\node at (5.5,7.5) []{$()$};
\node at (6.5,7.5) []{$()$};
\node at (7.5,7.5) []{$()$};

\node at (1.5,6.5) []{$()$};
\node at (2.5,6.5) []{$()$};
\node at (3.5,6.5) []{$()$};
\node at (4.5,6.5) []{$(1)$};
\node at (5.5,6.5) []{$(1)$};
\node at (6.5,6.5) []{$(1)$};
\node at (7.5,6.5) []{$(1)$};

\node at (1.5,5.5) []{$()$};
\node at (2.5,5.5) []{$()$};
\node at (3.5,5.5) []{$()$};
\node at (4.5,5.5) []{$(1)$};
\node at (5.5,5.5) []{$(1)$};
\node at (6.5,5.5) []{$(1)$};
\node at (7.5,5.5) []{$(2)$};

\node at (1.5,4.5) []{$()$};
\node at (2.5,4.5) []{$(1)$};
\node at (3.5,4.5) []{$(1)$};
\node at (4.5,4.5) []{$(1,1)$};
\node at (5.5,4.5) []{$(1,1)$};
\node at (6.5,4.5) []{$(1,1)$};
\node at (7.5,4.5) []{$(2,1)$};

\node at (1.5,3.5) []{$()$};
\node at (2.5,3.5) []{$(1)$};
\node at (3.5,3.5) []{$(1)$};
\node at (4.5,3.5) []{$(1,1)$};
\node at (5.5,3.5) []{$(2,1)$};
\node at (6.5,3.5) []{$(2,1)$};
\node at (7.5,3.5) []{$(2,2)$};

\node at (1.5,2.5) []{$(1)$};
\node at (2.5,2.5) []{$(1,1)$};
\node at (3.5,2.5) []{$(1,1)$};
\node at (4.5,2.5) []{$(1,1,1)$};
\node at (5.5,2.5) []{$(2,1,1)$};
\node at (6.5,2.5) []{$(2,1,1)$};
\node at (7.5,2.5) []{$(2,2,1)$};

\node at (1.5,1.5) []{$(1)$};
\node at (2.5,1.5) []{$(1,1)$};
\node at (3.5,1.5) []{$(2,1)$};
\node at (4.5,1.5) []{$(2,1,1)$};
\node at (5.5,1.5) []{$(2,2,1)$};
\node at (6.5,1.5) []{$(2,2,1)$};
\node at (7.5,1.5) []{$(2,2,2)$};

\node at (1.5,0.5) []{$(1)$};
\node at (2.5,0.5) []{$(1,1)$};
\node at (3.5,0.5) []{$(2,1)$};
\node at (4.5,0.5) []{$(2,1,1)$};
\node at (5.5,0.5) []{$(2,2,1)$};
\node at (6.5,0.5) []{$(3,2,1)$};
\node at (7.5,0.5) []{$(3,2,2)$};

\node (Black1) at (1,3) [black]{};
\node (Black2) at (2,5) [black]{};
\node (Black3) at (4,7) [black]{};

\node (Blue1) at (2,3) [blue]{};
\node (Blue2) at (4,5) [blue]{};

\node (Red1) at (4,3) [red]{};

\draw[draw=black,thick] (1,0) -- (Black1)-- (Blue1) -- (Black2) -- (Blue2) -- (Black3) -- (8,7);

\draw[draw=blue,thick] (2,0) -- (Blue1)-- (Red1) -- (Blue2) --  (8,5);

\draw[draw=red,thick] (4,0) -- (Red1) --  (8,3);

\node (Black4) at (3,2) [black]{};
\node (Black5) at (5,4) [black]{};
\node (Black6) at (7,6) [black]{};

\node (Blue3) at (5,2) [blue]{};
\node (Blue4) at (7,4) [blue]{};

\node (Red2) at (7,2) [red]{};

\draw[draw=black,thick] (3,0) -- (Black4)-- (Blue3) -- (Black5) -- (Blue4) -- (Black6) -- (8,6);

\draw[draw=blue,thick] (5,0) -- (Blue3)-- (Red2) -- (Blue4)  -- (8,4);

\draw[draw=red,thick] (7,0) -- (Red2) --  (8,2);

\node (Black7) at (6,1) [black]{};

\draw[draw=black,thick] (6,0) -- (Black7) -- (8,1);

\end{tikzpicture}
\end{center}

The pair of standard partitions may be read of the bottom row and the final column as a sequence of integer partitions, each of which may be obtained from the previous by adding a single box.

\subsection{Skew standard tableau}\label{sec:skew-standard}

\index{Tableau!skew standard}

A \emph{standard skew tableau} may be defined as a sequence of partitions:
\[ S = (\lambda_1,\lambda_2, \ldots \lambda_n) \]
such that if $i > j$ then $\lambda_j \subseteq \lambda_i$ and such that each partition in the sequence differs by at most one box.
We say that $S$ is of \emph{shape} $\lambda_n / \lambda_1$.

It is often convenient to represent a skew standard tableau as a \emph{labelled skew diagram}. For example,
the skew standard tableau:
\[ S = ((2,1),(2,2),(2,2),(2,2,1),(2,2,2),(3,2,2))\]
may be represented by the diagram:
\[
\tableau{\missingcell & 4 & 2 \\ \missingcell & 1 & \graycell \\ 5 & \graycell & \graycell}
\]
The box $s$ is labelled $k$ if it first appears in the $k$-th partition of $S$.

\index{Tableau!standard}

A regular standard Young tableau is a special case of a skew standard tableau in which $\lambda_1 = ()$ and for which each partition in the sequence differs by exactly one box.

\subsection{Fomin's local rules}\label{sec:fomin-local}

\index{Local rule!Fomin, forward}

\emph{Fomin's local rules} \cite{fomin-growth} tell us how to find the partition $\lambda$ given the partitions $\rho$, $\mu$ and $\nu$ and a value of $x$ which may be either $0$ or $1$, corresponding to a small dot and a large dot respectively.

\begin{center}
\begin{tikzpicture}[scale=1.4]
\node  at (2,2) []{$\color{blue}x$};

\node at (2.5,2.5) []{$\color{blue}\nu$};
\node at (2.5,1.5) []{$\color{red}\lambda$};
\node at (1.5,2.5) []{$\color{blue}\rho$};
\node at (1.5,1.5) []{$\color{blue}\mu$};

\end{tikzpicture}
\end{center}

The blue data is input, and the red data is output. The rules are as follows:
\begin{enumerate}
\item If $x=1$ and $\mu = \nu = \rho$ then $\lambda$ is the partition obtained from $\rho$ by adding a box to the first row.
\item If $x=0$ and $\mu = \nu = \rho$ then $\lambda = \rho$.
\item If $x=0$ and $\mu = \rho \neq \nu $ then $\lambda = \nu$ (this corresponds to the case of a shadow passing through the face vertically).
\item If $x=0$ and $\nu = \rho \neq \mu$ then $\lambda = \mu$ (this corresponds to the case of a shadow passing through the face horizontally).
\item If $x=0$ and $\mu$, $\nu$ and $\rho$ are pairwise distinct then $\lambda=\mu\cup\nu$ (this corresponds to the case of two shadows passing through each other).
\item If $x=0$ and $\mu = \nu \neq \rho$ and $\mu$ differs from $\rho$ on the $k$-th row, then $\lambda$ is the partition which is obtained from $\mu$ by adding a box to row $k+1$ (this corresponds to the case of two shadows ``bouncing off'' each other).
\end{enumerate}

\index{Fomin growth diagram!local rules}

In the diagram below we have highlighted an example of rule $1$ in yellow, rule $2$ in blue, rule $3$ in red, rule $4$ in pink, rule $5$ in orange and rule $6$ in purple.

\begin{center}
\begin{tikzpicture}[scale=1.4]
\foreach \x in {1,...,7}
\foreach \y in {1,...,7}
\node  at (\x,\y) [gray]{};

\node at (0.5,0.5) []{$()$};
\node at (0.5,1.5) []{$()$};
\node at (0.5,2.5) []{$()$};
\node at (0.5,3.5) []{$()$};
\node at (0.5,4.5) []{$()$};
\node at (0.5,5.5) []{$()$};
\node at (0.5,6.5) []{$()$};
\node at (0.5,7.5) []{$()$};

\node at (1.5,7.5) []{$()$};
\node at (2.5,7.5) []{$()$};
\node at (3.5,7.5) []{$()$};
\node at (4.5,7.5) []{$()$};
\node at (5.5,7.5) [pink]{$()$};
\node at (6.5,7.5) [pink]{$()$};
\node at (7.5,7.5) []{$()$};

\node at (1.5,6.5) []{$()$};
\node at (2.5,6.5) []{$()$};
\node at (3.5,6.5) [red]{$()$};
\node at (4.5,6.5) [red]{$(1)$};
\node at (5.5,6.5) [pink]{$(1)$};
\node at (6.5,6.5) [pink]{$(1)$};
\node at (7.5,6.5) []{$(1)$};

\node at (1.5,5.5) []{$()$};
\node at (2.5,5.5) []{$()$};
\node at (3.5,5.5) [red]{$()$};
\node at (4.5,5.5) [red]{$(1)$};
\node at (5.5,5.5) []{$(1)$};
\node at (6.5,5.5) []{$(1)$};
\node at (7.5,5.5) []{$(2)$};

\node at (1.5,4.5) []{$()$};
\node at (2.5,4.5) [blue]{$(1)$};
\node at (3.5,4.5) [blue]{$(1)$};
\node at (4.5,4.5) []{$(1,1)$};
\node at (5.5,4.5) []{$(1,1)$};
\node at (6.5,4.5) [purple]{$(1,1)$};
\node at (7.5,4.5) [purple]{$(2,1)$};

\node at (1.5,3.5) []{$()$};
\node at (2.5,3.5) [blue]{$(1)$};
\node at (3.5,3.5) [blue]{$(1)$};
\node at (4.5,3.5) []{$(1,1)$};
\node at (5.5,3.5) []{$(2,1)$};
\node at (6.5,3.5) [purple]{$(2,1)$};
\node at (7.5,3.5) [purple]{$(2,2)$};

\node at (1.5,2.5) []{$(1)$};
\node at (2.5,2.5) [yellow]{$(1,1)$};
\node at (3.5,2.5) [yellow]{$(1,1)$};
\node at (4.5,2.5) []{$(1,1,1)$};
\node at (5.5,2.5) []{$(2,1,1)$};
\node at (6.5,2.5) []{$(2,1,1)$};
\node at (7.5,2.5) []{$(2,2,1)$};

\node at (1.5,1.5) []{$(1)$};
\node at (2.5,1.5) [yellow]{$(1,1)$};
\node at (3.5,1.5) [yellow]{$(2,1)$};
\node at (4.5,1.5) []{$(2,1,1)$};
\node at (5.5,1.5) []{$(2,2,1)$};
\node at (6.5,1.5) [orange]{$(2,2,1)$};
\node at (7.5,1.5) [orange]{$(2,2,2)$};

\node at (1.5,0.5) []{$(1)$};
\node at (2.5,0.5) []{$(1,1)$};
\node at (3.5,0.5) []{$(2,1)$};
\node at (4.5,0.5) []{$(2,1,1)$};
\node at (5.5,0.5) []{$(2,2,1)$};
\node at (6.5,0.5) [orange]{$(3,2,1)$};
\node at (7.5,0.5) [orange]{$(3,2,2)$};

\node (Black1) at (1,3) [black]{};
\node (Black2) at (2,5) [black]{};
\node (Black3) at (4,7) [black]{};

\node (Blue1) at (2,3) [blue]{};
\node (Blue2) at (4,5) [blue]{};

\node (Red1) at (4,3) [red]{};

\draw[draw=black,very thick] (1,0) -- (Black1)-- (Blue1) -- (Black2) -- (Blue2) -- (Black3) -- (8,7);

\draw[draw=blue,very thick] (2,0) -- (Blue1)-- (Red1) -- (Blue2) --  (8,5);

\draw[draw=red,very thick] (4,0) -- (Red1) --  (8,3);

\node (Black4) at (3,2) [black]{};
\node (Black5) at (5,4) [black]{};
\node (Black6) at (7,6) [black]{};

\node (Blue3) at (5,2) [blue]{};
\node (Blue4) at (7,4) [blue]{};

\node (Red2) at (7,2) [red]{};

\draw[draw=black,very thick] (3,0) -- (Black4)-- (Blue3) -- (Black5) -- (Blue4) -- (Black6) -- (8,6);

\draw[draw=blue,very thick] (5,0) -- (Blue3)-- (Red2) -- (Blue4)  -- (8,4);

\draw[draw=red,very thick] (7,0) -- (Red2) --  (8,2);

\node (Black7) at (6,1) [black]{};

\draw[draw=black,very thick] (6,0) -- (Black7) -- (8,1);

\end{tikzpicture}
\end{center}

Although it is not immediately obvious that Fomin's local rules are exhaustive as given, one can prove by recurrence that they cover every possible case which can arise, given suitable initial conditions.

\subsection{Skew Robinson correspondence}\label{sec:fomin-growth}

It is now possible to forget about Viennot's shadows, and perform the Robinson correspondence using only Fomin's local rules.
In fact, we shall do something slightly more general.

\index{Robinson correspondence!skew}
\index{Fomin growth diagram!input}

The input to algorithm is a pair $(A,B)$ of skew standard tableau of shape $\alpha / \mu$ and $\beta / \mu$ respectively, together
with a partial permutation $\sigma$.
The partial permutation must be compatible with the pair $(A,B)$ in the sense that
if there is a $1$ in column $k$ then $A_k = A_{k+1}$, similarly, if there is a $1$ in row $k$ then $B_k = B_{k+1}$.
For example:

\begin{center}
\begin{tikzpicture}[scale=1.4]
\foreach \x in {4,...,8}
\foreach \y in {0,...,4}
\node  at (\x,\y) [gray]{};

\node at (3.5,4.5) []{$(1)$};
\node at (4.5,4.5) []{$(1)$};
\node at (5.5,4.5) []{$(1,1)$};
\node at (6.5,4.5) []{$(1,1)$};
\node at (7.5,4.5) []{$(1,1)$};
\node at (8.5,4.5) []{$(2,1)$};

\node at (3.5,3.5) []{$(1)$};

\node at (3.5,2.5) []{$(1,1)$};

\node at (3.5,1.5) []{$(1,1)$};

\node at (3.5,0.5) []{$(2,1)$};

\node at (3.5,-0.5) []{$(2,1)$};

\node (Black5) at (6,4) [black]{};

\node (Black7) at (7,0) [black]{};

\end{tikzpicture}
\end{center}

Another way of expressing the compatibility condition between $(A,B)$ and $\sigma$ is that if the partial permutation matrix $\sigma$
is $n$ by $m$ then each of the labels from $1$ to $m$ must occur at most once once either in the labelled diagram for $A$
or in the top row of the two line notation for $\sigma$. Similarly each of the labels from $1$ to $n$ must occur
at most once either in the labelled diagram for $B$ or in the bottom row of the two line notation for $\sigma$.
\[
A = \tableau{\missingcell & 2 \\ 5 & \graycell} \qquad \qquad B = \tableau{\missingcell & 2 \\ 4 & \graycell}
\qquad \qquad \sigma = \left [ \begin{matrix} . & . & 3 & 4 & . \\ . & . & 1 & 5 & . \end{matrix} \right ]
\]

Starting from these initial conditions, we may make use of Fomin's local rules to fill out the rest of the data in the diagram:

\begin{center}
\begin{tikzpicture}[scale=1.4]
\foreach \x in {4,...,8}
\foreach \y in {0,...,4}
\node  at (\x,\y) [gray]{};

\node at (3.5,4.5) []{$(1)$};
\node at (4.5,4.5) []{$(1)$};
\node at (5.5,4.5) []{$(1,1)$};
\node at (6.5,4.5) []{$(1,1)$};
\node at (7.5,4.5) []{$(1,1)$};
\node at (8.5,4.5) []{$(2,1)$};

\node at (3.5,3.5) []{$(1)$};
\node at (4.5,3.5) []{$(1,1)$};
\node at (5.5,3.5) []{$(1,1)$};
\node at (6.5,3.5) []{$(2,1)$};
\node at (7.5,3.5) []{$(2,1)$};
\node at (8.5,3.5) []{$(2,2)$};

\node at (3.5,2.5) []{$(1,1)$};
\node at (4.5,2.5) []{$(1,1)$};
\node at (5.5,2.5) []{$(1,1,1)$};
\node at (6.5,2.5) []{$(2,1,1)$};
\node at (7.5,2.5) []{$(2,2,1)$};
\node at (8.5,2.5) []{$(2,2,1)$};

\node at (3.5,1.5) []{$(1,1)$};
\node at (4.5,1.5) []{$(1,1)$};
\node at (5.5,1.5) []{$(1,1,1)$};
\node at (6.5,1.5) []{$(2,1,1)$};
\node at (7.5,1.5) []{$(2,2,1)$};
\node at (8.5,1.5) []{$(2,2,1)$};

\node at (3.5,0.5) []{$(2,1)$};
\node at (4.5,0.5) []{$(2,1)$};
\node at (5.5,0.5) []{$(2,1,1)$};
\node at (6.5,0.5) []{$(2,2,1)$};
\node at (7.5,0.5) []{$(2,2,1)$};
\node at (8.5,0.5) []{$(2,2,2)$};

\node at (3.5,-0.5) []{$(2,1)$};
\node at (4.5,-0.5) []{$(2,1)$};
\node at (5.5,-0.5) []{$(2,1,1)$};
\node at (6.5,-0.5) []{$(2,2,1)$};
\node at (7.5,-0.5) []{$(3,2,1)$};
\node at (8.5,-0.5) []{$(3,2,2)$};

\node (Black5) at (6,4) [black]{};

\node (Black7) at (7,0) [black]{};

\end{tikzpicture}
\end{center}

The output of the algorithm, may be read off the right-most column and the bottom-most row:

\begin{center}
\begin{tikzpicture}[scale=1.4]
\foreach \x in {4,...,8}
\foreach \y in {0,...,4}
\node  at (\x,\y) [gray]{};

\node at (8.5,4.5) []{$(2,1)$};

\node at (8.5,3.5) []{$(2,2)$};

\node at (8.5,2.5) []{$(2,2,1)$};
\node at (8.5,1.5) []{$(2,2,1)$};

\node at (8.5,0.5) []{$(2,2,2)$};

\node at (3.5,-0.5) []{$(2,1)$};
\node at (4.5,-0.5) []{$(2,1)$};
\node at (5.5,-0.5) []{$(2,1,1)$};
\node at (6.5,-0.5) []{$(2,2,1)$};
\node at (7.5,-0.5) []{$(3,2,1)$};
\node at (8.5,-0.5) []{$(3,2,2)$};

\end{tikzpicture}
\end{center}

\index{Fomin growth diagram!output}

The output is a pair of standard skew tableaux $(A',B')$ of shape $\lambda / \alpha$ and $\lambda / \beta$ respectively.

\[ A' = \tableau{\missingcell & 5 & 2 \\ \missingcell & 3 & \graycell \\ 4 & \graycell & \graycell}
\qquad \qquad 
B'= \tableau{\missingcell & 4 & 2 \\ \missingcell & 1 & \graycell \\ 5 & \graycell & \graycell} \]

Note that the label $k$ is missing from $A'$ if and only if it is missing from both $A$ and the top row of $\sigma$.
Similarly the label $k$ is missing from $B'$ if and only if it is missing from both $B$ and the bottom row of $\sigma$.

\index{Fomin growth diagram!standardization}
 
We shall say that a growth diagram is \emph{standard} if it has no repeated rows or columns. The output pair $(A',B')$ of a growth diagram will have the property that each and every label occurs exactly once if and only if the corresponding growth diagram is standard. The \emph{standardization} of a growth diagram is the diagram obtained from the original by removing any repeated rows or columns. 

\subsection{Fomin's reverse local rules}\label{sec:fomin-local-reverse}

\index{Local rule! Fomin, reverse}

\emph{Fomin's reverse local rules} tell us how to find the pair $(\rho,x)$ given the partitions $\lambda$, $\mu$ and $\nu$.

\begin{center}
\begin{tikzpicture}[scale=1.4]
\node  at (2,2) []{$\color{red}x$};

\node at (2.5,2.5) []{$\color{blue}\nu$};
\node at (2.5,1.5) []{$\color{blue}\lambda$};
\node at (1.5,2.5) []{$\color{red}\rho$};
\node at (1.5,1.5) []{$\color{blue}\mu$};

\end{tikzpicture}
\end{center}

The value of $x$ may be either $0$ or $1$, corresponding to a small dot and a large dot respectively. The blue data is input, and the red data is output. Fomin's reverse local rule may be described as follows:
\begin{enumerate}
\item If $\mu = \nu \neq \lambda$ and if $\lambda$ differs from $\mu$ on the first row then $\rho = \mu$ and $x=1$.
\item If $\mu = \lambda = \nu$ then $\rho = \lambda$ and $x = 0$.
\item If $\mu = \lambda \neq \nu$ then $\rho = \nu$ and $x=0$.
\item If $\nu = \lambda \neq \mu$ then $\rho = \mu$ and $x=0$.
\item If $\lambda$, $\mu$ and $\nu$ are pairwise distinct, then $\rho=\mu\cap\nu$ and $x=0$.
\item If $\mu = \nu \neq \lambda$ and if $\lambda$ differs from $\mu$ on the $k$-th row for $k > 1$ then $\rho$ is the partition obtained from $\mu$ by removing a box on the $(k-1)$-th row and $x=0$.
\end{enumerate}

Staring with the output configuration from the previous section, one may use Fomin's reverse local rules to ``grow'' the diagram backwards and recover the initial conditions.

\index{Robinson correspondence!skew}

This algorithm is slightly more general than we need. It corresponds, in fact, to the \emph{skew Robinson correspondence} \cite{sagan}.
We recover the special case of the Robinson correspondence when $(A',B')$ are regular standard Young tableau.

\section{RSK and Burge correspondences}\label{sec:RSK}

The RSK correspondence and the Burge correspondence both give a bijection between non-negative integer matrices and pairs of semi-standard Young tableaux of the same shape.

\subsection{Semi-standard tableau}\label{sec:semi-standard}

\index{Tableau!semi-standard}

A \emph{semi-standard Young tableau} is a labelled Young diagram which is weakly decreasing along rows and strictly decreasing along columns. For example:

\[
\tableau{
\missingcell & \missingcell & \missingcell & \missingcell  & \missingcell  \\
\missingcell & \missingcell & \missingcell  & 5  & 4 \\
\missingcell & \missingcell & 4 & 4 & 3 \\
\missingcell  & \missingcell & 3 & 2 & 2  \\
 2 & 2 &1 &1 &1 \\
} 
\]

A semistandard Young tableau may also be represented by a sequence of integer partitions which differ successively by a horizontal strip. The $k$-th integer partition in this sequence is given by the subtableau covered by the first $k$ labels:
\[
\emptyset \rightharpoonup (3) \rightharpoonup (5,2) \rightharpoonup (5,3,1) \rightharpoonup (5,3,3,1) \rightharpoonup (5,3,3,2)
\]

\index{Tableau!semi-standard!content}

The \emph{content} $c(T)$ of a semi-standard Young tableau is the vector: 
\[(c_1(T), c_2(T), c_3(T), \ldots )\]
where $c_k(T)$ denotes the number of times the label $k$ occurs in the tableau $T$.
The content of our example tableau is $(3,4,2,3,1)$.

We remark that a regular plane partition may be thought of as a pair of semi-standard tableau of the same shape. For example,
the regular plane partition from section \ref{sec:pp}:
\[ \tableau{4 & 3 & 2 & 2 & 0 \\ 3 & 2 & 1 & 1 & 0 \\ 1 & 1 & 1 & 0 & 0\\ 1 & 0 & 0 & 0 & 0 \\ 0 & 0 & 0 & 0 & 0} \]
is in bijection with the following two tableaux:
\[\tableau{\missingcell & \missingcell & \missingcell & 4  \\ \missingcell & \missingcell & 4 & 2 \\ 4& 3 & 1 & 1} \qquad
\tableau{\missingcell & \missingcell & \missingcell & 4  \\ \missingcell & \missingcell & 4 & 3 \\ 4& 3 & 3 & 1}
\] 
The first tableau in the pair corresponds to the first half of the sequence read off the $NW \to SE$ diagonals, from right to left:
\[ \emptyset \rightharpoonup (2) \rightharpoonup (2,1) \rightharpoonup (3,1) \rightharpoonup (4,2,1)\]
while the second tableau in the pair corresponds to the second half of the sequence read off the $NW \to SE$ diagonals, from right to left:
\[\emptyset \rightharpoonup (1) \rightharpoonup (1) \rightharpoonup (3,1) \rightharpoonup (4,2,1)) \]

\index{Tableau!standardization}
\index{Tableau!de-standardization}

The \emph{standardization} of a semi-standard Young tableau is the unique standard Young tableau with the following properties: 
\begin{enumerate}
\item The labels $1,2,\ldots c_1(T)$ lie on the first row (read from the bottom).
\item The set of labels $i$ satisfying: $ c_k(T) < i \leq c_{k+1}(T)$ forms a horizontal strip.
\item If $ c_k(T) < i,j \leq c_{k+1}(T)$ and $i<j$ then $i$ lies to the right of $j$.  
\item when you replace the first $c_1(T)$ labels with $1$, the next $c_2(T)$ labels with $2$, the next $c_3(T)$ labels with $3$ and so on,
you recover the original semi-standard Young tableau.
\end{enumerate}
The standardization of our example semi-standard Young tableau is the following:
\[
\tableau{
\missingcell & \missingcell & \missingcell & \missingcell  & \missingcell  \\
\missingcell & \missingcell & \missingcell  & 13  & 10 \\
\missingcell & \missingcell & 12 & 11 & 8 \\
\missingcell  & \missingcell & 9 & 5 & 4  \\
 7 & 6&3 &2 &1 \\
} 
\]

\subsection{Horizontal and vertical strips}\label{sec:remark}

\index{Permutation!partial!increasing chain}
\index{Permutation!partial!decreasing chain}

A partial permutation matrix $\sigma$ is said to form an \emph{increasing chain} if there exists $k,\ell \geq 0$ such that after removing the first $k$ rows and the first $\ell$ columns, the remaining matrix is that of the identity permutation. Similarly, a partial permutation matrix $\sigma$ is said to form an \emph{decreasing chain} if there exists $k, \ell \geq 0$ such that after removing the \emph{final} $k$ rows and the final $\ell$ columns, the remaining permutation is that of the \emph{maximal permutation}.

\index{Horizontal strip!properly labelled}

Let us define a \emph{properly labelled horizontal strip}
to be a standard skew tableau with the property that if $i > j$ then the label $i$ always occurs to the left of label $j$.
The following important lemma is due to vanLeeuwen \cite{vanLeeuwen}.

\index{Fomin growth diagram!patch}

\begin{lemma}\label{van-horizontal}
Suppose that the skew Robinson correspondence sends $(A,B,\sigma)$ to $(A',B')$.
The pair of skew standard tableau $A'$ and $B'$ are properly labelled horizontal strips
if and only if the pair of skew tableau $A$ and $B$ are also properly labelled horizontal strips and if after the growth diagram has been standardized the partial permutation $\sigma$ form a \emph{decreasing chain}.

\end{lemma}

An example:

\begin{center}
\begin{tikzpicture}[scale=1.4] 
\foreach \x in {4,...,7}
\foreach \y in {1,...,4}
\node  at (\x,\y) [gray]{};

\node at (3.5,4.5) []{$2$};
\node at (4.5,4.5) []{$21$};
\node at (5.5,4.5) []{$22$};
\node at (6.5,4.5) []{$22$};
\node at (7.5,4.5) []{$22$};

\node at (3.5,3.5) []{$21$};
\node at (4.5,3.5) []{$211$};
\node at (5.5,3.5) []{$221$};
\node at (6.5,3.5) []{$221$};
\node at (7.5,3.5) []{$221$};

\node at (3.5,2.5) []{$22$};
\node at (4.5,2.5) []{$221$};
\node at (5.5,2.5) []{$222$};
\node at (6.5,2.5) []{$222$};
\node at (7.5,2.5) []{$222$};

\node at (3.5,1.5) []{$22$};
\node at (4.5,1.5) []{$221$};
\node at (5.5,1.5) []{$222$};
\node at (6.5,1.5) []{$322$};
\node at (7.5,1.5) []{$322$};

\node at (3.5,0.5) []{$22$};
\node at (4.5,0.5) []{$221$};
\node at (5.5,0.5) []{$222$};
\node at (6.5,0.5) []{$322$};
\node at (7.5,0.5) []{$422$};

\node (Black1) at (6,2) [black]{};
\node (Black2) at (7,1) [black]{};

\draw[draw=black,thick] (6,0) -- (Black1) --  (8,2);
\draw[draw=black,thick] (7,0) -- (Black2) --  (8,1);

\node (Red1) at (4,4) [red]{};
\node (Red2) at (5,3) [red]{};

\draw[draw=blue,thick] (3,4) -- (Red1) --  (4,5);
\draw[draw=blue,thick] (3,3) -- (Red2) --  (5,5);

\draw[draw=red,thick] (8,4) -- (Red1) --  (4,0);
\draw[draw=red,thick] (8,3) -- (Red2) --  (5,0);





\end{tikzpicture}
\end{center}

\[
\sigma = \left [ \begin{matrix} 1 & 2 & 3 & 4 \\ * & * & 3 & 4\end{matrix}\right]
\qquad \qquad A = B = \tableau{2 & 1 \\ \graycell & \graycell}
\qquad \qquad
A' = B' = \tableau{\missingcell & \missingcell & 2 & 1 \\
\missingcell & \missingcell & \graycell & \graycell \\
4 & 3 & \graycell & \graycell}
\]

\index{Fomin growth diagram!horizontal type}

We shall say that a growth diagram is of \emph{horizontal type} if the output tableaux $(A',B')$ are both properly labelled horizontal strips.

\index{Horizontal strip!properly labelled}
\index{Fomin growth diagram!vertical type}

Similarly, let us define a \emph{properly labelled vertical strip}
to be a skew standard tableau with the property that if $i > j$ then the label $i$ always occurs to the \emph{right} of label $j$.
The following lemma is also due to vanLeeuwen \cite{vanLeeuwen}:

\begin{lemma}\label{van-vertical}
Suppose that the skew Robinson correspondence sends $(A,B,\sigma)$ to $(A',B')$.
The pair of skew standard tableau $A'$ and $B'$ are properly labelled vertical strips
if and only if the pair of skew tableau $A$ and $B$ are also properly labelled vertical strips and, after the growth diagram has been standardized, the partial permutation $\sigma$ form an \emph{increasing chain}.
\end{lemma}
For example:

\begin{center}
\begin{tikzpicture}[scale=1.4]
\foreach \x in {4,...,7}
\foreach \y in {1,...,4}
\node  at (\x,\y) [gray]{};

\node at (3.5,4.5) []{$11$};
\node at (4.5,4.5) []{$11$};
\node at (5.5,4.5) []{$11$};
\node at (6.5,4.5) []{$21$};
\node at (7.5,4.5) []{$22$};

\node at (3.5,3.5) []{$11$};
\node at (4.5,3.5) []{$11$};
\node at (5.5,3.5) []{$21$};
\node at (6.5,3.5) []{$22$};
\node at (7.5,3.5) []{$221$};

\node at (3.5,2.5) []{$11$};
\node at (4.5,2.5) []{$21$};
\node at (5.5,2.5) []{$22$};
\node at (6.5,2.5) []{$221$};
\node at (7.5,2.5) []{$2211$};

\node at (3.5,1.5) []{$21$};
\node at (4.5,1.5) []{$22$};
\node at (5.5,1.5) []{$221$};
\node at (6.5,1.5) []{$2211$};
\node at (7.5,1.5) []{$22111$};

\node at (3.5,0.5) []{$22$};
\node at (4.5,0.5) []{$221$};
\node at (5.5,0.5) []{$2211$};
\node at (6.5,0.5) []{$22111$};
\node at (7.5,0.5) []{$221111$};

\node (Black1) at (4,3) [black]{};
\node (Black2) at (5,4) [black]{};

\node (Blue1) at (4,2) [blue]{};
\node (Blue2) at (5,3) [blue]{};
\node (Blue3) at (6,4) [blue]{};

\node (Red1) at (4,1) [red]{};
\node (Red2) at (5,2) [red]{};
\node (Red3) at (6,3) [red]{};
\node (Red4) at (7,4) [red]{};

\node (Green2) at (5,1) [green]{};
\node (Green3) at (6,2) [green]{};
\node (Green4) at (7,3) [green]{};

\node (Purple3) at (6,1) [purple]{};
\node (Purple4) at (7,2) [purple]{};

\node (Yellow4) at (7,1) [yellow]{};

\draw[draw=black,thick] (3,2) -- (Blue1)-- (Black1) -- (Blue2) -- (Black2) -- (Blue3) -- (6,5);
\draw[draw=blue,thick] (3,1) -- (Red1)-- (Blue1)-- (Red2) -- (Blue2) -- (Red3) -- (Blue3) -- (Red4) -- (7,5);
\draw[draw=red,thick] (4,0) -- (Red1)-- (Green2)-- (Red2) -- (Green3) -- (Red3) -- (Green4) -- (Red4) -- (8,4);
\draw[draw=green,thick] (5,0) -- (Green2)-- (Purple3) -- (Green3) -- (Purple4) -- (Green4)-- (8,3);
\draw[draw=purple,thick] (6,0) -- (Purple3) -- (Yellow4) -- (Purple4) -- (8,2);
\draw[draw=yellow,thick] (7,0) -- (Yellow4) -- (8,1);





\end{tikzpicture}
\end{center}

\[
\sigma = \left [ \begin{matrix} 1 & 2 & 3 & 4 \\ 2 & 1 & * & *\end{matrix} \right ] 
\qquad \qquad
A = B = \tableau{2 & \graycell \\ 1 & \graycell}
\quad \qquad
A' = B' = \tableau{ \missingcell & 4 \\ \missingcell & 3 \\ \missingcell & 2 \\ \missingcell & 1 \\ \graycell & \graycell \\ \graycell & \graycell}
\]

\index{Fomin growth diagram!vertical type}
\index{Vertical strip!properly labelled}

We shall say that a growth diagram is of \emph{vertical type} if the output tableau $(A',B')$ are both properly labelled vertical strips.

\index{Local rule!weight condition}

In both the horizontal strip and the vertical strip case case the following weight condition is satisfied. Suppose that $A'$ is of shape $\lambda / \alpha$ and $B'$ is of shape $\lambda / \beta$ while $A$ is of shape $\alpha / \mu$ and $B$ is of shape $\beta / \mu$. If the partial permutation $\sigma$ contains exactly $m$ ones then:
\begin{equation}
|\lambda| + |\mu| = |\alpha| + |\beta| + m
\end{equation}

\subsection{Block permutation matrices}\label{sec:block-perm}

\index{Permutation!block}

Suppose that $M$ is a non-negative integer matrix. Let $c_k$ denote the sum of the entries on the $k$-th row of $M$ and let $c'_k$ denote the sum of the entries in the $k$-th column of $M$.
If the total sum of all the entries of $M$ is equal to $n$, then there are two canonical ways in which we may associate an $n$ by $n$ permutation matrix $P$ to the matrix $M$.

In both cases we begin by dividing up our matrix $P$ into ``block rows'' and ``block columns''. The $k$-th ``block row'' of $P$ contains $c_k$ normal rows of $P$, while the $k$-th ``block column'' of $P$ contains $c'_k$ columns of $P$. A ``block'' of $P$ is an intersection of a ``block row'' and a ``block column''.

\index{RSK permutation}

The first way of associating a permutation matrix $P$ to a non-negative integer matrix $M$ is to place ones into the block matrix $P$ in such a way that:
\begin{itemize} 
\item The number of ones in the block $(i,j)$ of $P$ is equal to $M_{i,j}$.
\item The ones are strictly decreasing from left to right as you move along any block row.
\item The ones are strictly decreasing from top to bottom as you move along any block column.
\end{itemize}
For example:

\[
\left [ \begin{matrix} 1 & 3 \\ 2 & 1\end{matrix} \right ]
\mapsto
\left [ \begin{array}{ccc|cccc} 
1 & 0 & 0 & 0 & 0 & 0 & 0 \\
0 & 0 & 0 & 1 & 0 & 0 & 0 \\
0 & 0 & 0 & 0 & 1 & 0 & 0 \\
0 & 0 & 0 & 0 & 0 & 1 & 0 \\
\hline
0 & 1 & 0 & 0 & 0 & 0 & 0 \\
0 & 0 & 1 & 0 & 0 & 0 & 0 \\
0 & 0 & 0 & 0 & 0 & 0 & 1 \\
\end{array}
\right ]
\]
In this case we shall refer to $P$ as the \emph{RSK permutation associated to $M$}.
If the row sums and column sums are known, then this operation is invertible

\index{Burge permutation}

The second way of associating a permutation matrix $P$ to a non-negative integer matrix $M$ is to place ones into the block matrix $P$ in such a way that:
\begin{itemize} 
\item The number of ones in the block $(i,j)$ of $P$ is equal to $M_{i,j}$.
\item The ones are strictly \emph{increasing} from left to right as you move along any block row.
\item The ones are strictly \emph{increasing} from bottom to top as you move along any block column.
\end{itemize}
For example:
\[
\left [ \begin{matrix} 1 & 3 \\ 2 & 1\end{matrix} \right ]
\mapsto
\left [ \begin{array}{ccc|cccc} 
0 & 0 & 0 & 0 & 0 & 0 & 1 \\
0 & 0 & 0 & 0 & 0 & 1 & 0 \\
0 & 0 & 0 & 0 & 1 & 0 & 0 \\
0 & 0 & 1 & 0 & 0 & 0 & 0 \\
\hline
0 & 0 & 0 & 1 & 0 & 0 & 0 \\
0 & 1 & 0 & 0 & 0 & 0 & 0 \\
1 & 0 & 0 & 0 & 0 & 0 & 0 \\
\end{array}
\right ]
\]
In this case we shall refer to $P$ as the \emph{Burge permutation associated to $M$}.
Again, if the row sums and column sums are known, then this operation is invertible.

\subsection{Growth diagram patches}\label{sec:remark2}

\index{Fomin growth diagram!patch}

Its not hard to see that any subdiagram of a growth diagram is again a growth diagram. In particular, given a block permutation matrix $P$ we may apply the Robinson correspondence, and then consider each of the blocks as separate, individual growth diagrams.

\begin{lemma}[vanLeeuwen]\label{block-growth}
If $P$ is an RSK permutation, then each block of $P$ is of horizontal type.
Similarly, if $P$ is a Burge permutation, then each block of $P$ is of vertical type.
\end{lemma}

\subsection{Reverse algorithm (RSK)}\label{sec:RSK-reverse}

\index{RSK correspondence!reverse}

The RSK correspondence may be realized with the aid of the Robinson correspondence. The reverse algorithm is as follows:

Suppose that your two semi-standard Young tableau are $(T,T')$ with content vectors $c = (c_1, \ldots, c_r)$ and $c' = (c'_1, \ldots, c'_s)$ respectively. For example:
\[
T = \tableau{\missingcell & 3 & 3 \\ \missingcell &2 & 2 \\ 3 & 1 & 1 } \qquad \qquad T' = \tableau{\missingcell & 4 & 3 \\ \missingcell & 3 & 2 \\ 3 & 2 & 1}
\]
We have $c = (2,2,3)$ and $c' = (1,2,3,1)$.

The first step is to standardize both $T$ and $T'$ to obtain a pair of standard tableau $(S,S')$.
In our case we have:
\[ \left ( \quad
\tableau{\missingcell & 6 & 5 \\ \missingcell & 4 & 3 \\ 7 & 2 & 1 \\}
\quad , \quad 
\tableau{\missingcell & 7 & 4 \\ \missingcell & 5 & 2 \\ 6 & 3 & 1} \quad
\right )
\]
Note that this is the same pair of tableau which appeared in our example in section \ref{sec:viennot}

The next step is to apply the reverse Robinson correspondence to get some permutation matrix. In our case we have:
\[
P = \left [
\begin{matrix}
0 & 0 & 0 & 1 & 0 & 0 & 0 \\
0 & 0 & 0 & 0 & 0 & 0 & 1 \\
0 & 1 & 0 & 0 & 0 & 0 & 0 \\
0 & 0 & 0 & 0 & 1 & 0 & 0 \\
1 & 0 & 0 & 0 & 0 & 0 & 0 \\
0 & 0 & 1 & 0 & 0 & 0 & 0 \\
0 & 0 & 0 & 0 & 0 & 1 & 0 \\
\end{matrix}
\right ]
\]

Finally we use the content vectors of $T$ and $T'$ to find the matrix $M$ for which $P$ is the associated RSK permutation:
 \[
\left [
\begin{array}{c|cc|ccc|c}
0 & 0 & 0 & 1 & 0 & 0 & 0 \\
0 & 0 & 0 & 0 & 0 & 0 & 1 \\
\hline
0 & 1 & 0 & 0 & 0 & 0 & 0 \\
0 & 0 & 0 & 0 & 1 & 0 & 0 \\
\hline
1 & 0 & 0 & 0 & 0 & 0 & 0 \\
0 & 0 & 1 & 0 & 0 & 0 & 0 \\
0 & 0 & 0 & 0 & 0 & 1 & 0 \\
\end{array}
\right ]
\mapsto
\left [
\begin{matrix}
0 & 0 & 1 & 1 \\
0 & 1 & 1 & 0 \\
1 & 1 & 1 & 0
\end{matrix}
\right ]
\]

Although in our example the resulting matrix only contains zeros and ones, in general the entries of $M$ may be any non-negative integer. 

\subsection{Forward algorithm (RSK)}

\index{RSK correspondence!forward}

Given a non-negative integer matrix $M$, begin by recording the row and column sums of $M$ into the vectors $c$ and $c'$.
Next, let $P$ be the RSK permutation associated to $M$, and apply the Robinson correspondence to $P$.
Suppose that $(S,S')$ is the resulting pair of standard tableau. The content vectors $(c,c')$ may now be used to \emph{de-standardize} $S$ and $S'$ into a pair of semi-standard Young tableau $(T,T')$

\subsection{Forward algorithm (Burge)}

\index{Burge correspondence!forward}

Given a non-negative integer matrix $M$, the forward algorithm is as follows:
\begin{enumerate}
\item Record the content vectors $c$ and $c'$ by reading of the row sums and column sums of $M$ respectively.
\item Find the Burge permutation $P$ associated to $M$.
\item Apply the Robinson correspondence to the permutation matrix $P$ to obtain a pair of standard Young tableau $(S,S')$.
\item Conjugate the pair of standard Young tableau $(S,S')$ to obtain a new pair of standard Young tableau $(R,R')$
\item De-standardize the tableau $(R,R')$ using the content vectors $(c,c')$ to obtain the final pair of standard Young tableau $(T,T')$.
\end{enumerate}

\subsection{Reverse algorithm (Burge)}

\index{Burge correspondence!reverse}

Given a pair $(T,T')$ of semi-standard Young tableau of the same shape, the reverse algorithm is as follows:
\begin{enumerate}
\item Record the content vectors $(c,c')$ of $(T,T')$.
\item Standardize the semi-standard Young tableaux $(T,T')$ to obtain a pair of standard Young tableau $(R,R')$.
\item Conjugate the pair of standard Young tableau $(R,R')$ to obtain a new pair of standard Young tableau $(S,S')$
\item Apply the reverse Robinson correspondence to $(S,S')$ to obtain the permutation $P$.
\item Use the content vectors $(c,c')$ to find the matrix $M$ for which $P$ is the associated Burge permutation.
\end{enumerate}

\subsection{Local rules}\label{sec:local}
\index{Local rule!RSK}
\index{Local rule!Burge}

As a consequence of lemma \ref{block-growth}, both the RSK and the Burge correspondence can be performed directly,
without passing through the standardization and de-standardization procedure. 

The growth diagram for the RSK correspondence looks very similar to that for the Robinson correspondence, the main difference being that neighboring partitions differ by a horizontal strip, rather than a single box. In the case of the Burge correspondence, neighboring partitions differ by a vertical strip - but the final result is conjugated.

For any partition $\lambda$ let $U(\lambda)$ denote the set of integer partitions which can be obtained from $\lambda$ by adding a horizontal strip, and let $D(\lambda)$ denote the set of integer partitions which can be obtained from $\lambda$ by removing a horizontal strip.
A forward local rule is a map with type signature:

\begin{equation} 
\mathfrak{U}_{\alpha,\beta} : (\mathbb{Z}_{\geq 0}, D(\alpha) \cap D(\beta)) \to U(\alpha) \cap U(\beta)  
\end{equation}

\index{Local rule!weight condition}

such that if $\lambda = \mathfrak{U}_{\alpha,\beta}(m,\mu)$ then the following weight conditions are satisfied:

\begin{equation} 
|\lambda| + |\mu| = |\alpha| + |\beta| + m
\end{equation}

A reverse local rules is a map with type signature:

\begin{equation} 
\mathfrak{D}_{\alpha,\beta} : U(\alpha) \cap U(\beta) \to (\mathbb{Z}_{\geq 0}, D(\alpha) \cap D(\beta))
\end{equation}

As a consequence of lemma \ref{van-horizontal} and \ref{van-vertical}, the local rules for both the RSK correspondence and the Burge correspondence may be derived directly from Fomin's local rules.

We shall describe here only the local rule associated to the Burge correspondence. The reader is referred to \cite{vanLeeuwen} for the RSK version of the local rule.

\subsection{Burge local rule}\label{sec:local-explicit}

In the case of the Burge correspondence, the operator $\mathfrak{D}_{\alpha,\beta}$ is defined as follows.
Suppose that $(m,\mu) = \mathfrak{D}_{\alpha,\beta}(\lambda)$. Let $\overline{A}$ denote the set of columns of $\lambda$ which are longer than the corresponding columns of $\alpha$
and let $\overline{B}$ denote the set of columns of $\lambda$ which are longer than the corresponding columns of $\beta$.

Next, for each $i \in \overline{A} \cap \overline{B} $ let $\delta(i) \not\in \overline{A} \cup \overline{B}$ denote the largest integer such that $\delta(i) < i$
and $\delta(i) \neq \delta(j)$ for any $j \in \overline{A} \cap \overline{B}$ such that $j > i$.
Let:
\[ \overline{C} = \{\delta(i) > 0 \, | \, i \in \overline{A} \cap \overline{B} \} \]
Finally let $\mu$ be the partition obtained from $\lambda$ by removing a box from the end of each of the columns indexed by $\overline{A} \cup \overline{B} \cup \overline{C}$ and let
\[m = \#\{ \delta(i) \leq 0 \, | \, i \in \overline{A} \cap \overline{B} \} \]

Here is an example:
\[ \boxed{\mathfrak{D}_{(6,5,5,3),(6,6,5,2)}(7,6,5,3,1) = (1,(6,5,4,2))}\]

The calculation proceeds as follows:

\begin{align*}
\lambda' & = (5,4,4,3,3,2,1) \\
\overline{A} = \cols(\lambda / \alpha) & = \{1,6,7\} \\
\overline{B} = \cols(\lambda / \beta) & = \{1,3,7\} \\
\end{align*}

\tikzstyle{blacknode}=[shape=circle,fill,inner sep=0.5mm]
\tikzstyle{rednode}=[shape=circle,red, fill, inner sep=0.5mm]
\tikzstyle{graynode}=[shape=rectangle, gray, fill, inner sep=0.5mm]
\tikzstyle{bluenode}=[shape=circle, blue, fill, inner sep=0.5mm]

\begin{center}\begin{tikzpicture}[scale=0.8]

\path (-1,0) node[bluenode](n0){};

\path (0,0) node[blacknode](n1){};
\path (1,0) node[graynode](n2){};
\path (2,0) node[rednode](n3){};
\path (3,0) node[graynode](n4){};
\path (4,0) node[graynode](n5){};
\path (5,0) node[rednode](n6){};
\path (6,0) node[blacknode](n7){};
\path (7,0) node[graynode](n8){};
\path (8,0) node[graynode](n9){};

\path (0,-0.5) node[](){1};
\path (1,-0.5) node[](){2};
\path (2,-0.5) node[](){3};
\path (3,-0.5) node[](){4};
\path (4,-0.5) node[](){5};
\path (5,-0.5) node[](){6};
\path (6,-0.5) node[](){7};
\path (7,-0.5) node[](){8};
\path (8,-0.5) node[](){9};

\path[<-] (n5) edge [bend left] (n7);
\path[<-] (n0) edge [bend left] (n1);
\end{tikzpicture}\end{center}

\begin{align*}
\overline{C} & = \{5\} \\
\mu' & = (4,4,3,3,2,1) \\
m & = 1 
\end{align*}

The inverse operator $\mathfrak{U}_{\alpha,\beta}$ is defined similarly.
Suppose that $\lambda = \mathfrak{U}_{\alpha,\beta}(m,\mu)$. Let $A$ denote the set of columns of $\alpha$ which are longer than the corresponding columns of $\mu$
and let $B$ denote the set of columns of $\beta$ which are longer than the corresponding columns of $\mu$.

Next, for each $i \in A \cap B$ let $\epsilon(i) \not \in A \cup B$ be the smallest integer such that $\epsilon(i) > i$ and $\epsilon(i) \neq \epsilon(j)$ for any $j \in A \cap B$ with $j < i$.
Let 
\[ C = \{ \epsilon(i) \,|\, i \in A \cap B \}\] 
Finally let $D$ denote the first $m$ elements of the complement of the set $A \cup B \cup C$ 
and let $\lambda$ be the partition obtained from $\mu$ by adding a box to the end of each of the columns in $A \cup B \cup C \cup D$.

Here is the inverse of our example:

\[ \boxed{\mathfrak{U}_{(6,5,5,3),(6,6,5,2)}(1,(6,5,4,2)) = (7,6,5,3,1) }\]

The calculation is straightforward:

\begin{align*}
\mu' & = (4,4,3,3,2,1) \\
A = \cols(\alpha / \mu) & = \{3,5\} \\
B = \cols(\beta / \mu) & = \{5, 6\} \\
m & = 1
\end{align*}

\tikzstyle{blacknode}=[shape=circle,fill,inner sep=0.5mm]
\tikzstyle{rednode}=[shape=circle,red, fill, inner sep=0.5mm]
\tikzstyle{graynode}=[shape=rectangle, gray, fill, inner sep=0.5mm]
\tikzstyle{bluenode}=[shape=circle, blue, fill, inner sep=0.5mm]

\begin{center}\begin{tikzpicture}[scale=0.8]

\path (-1,0) node[bluenode](n0){};
\path (0,0) node[graynode](n1){};
\path (1,0) node[graynode](n2){};
\path (2,0) node[rednode](n3){};
\path (3,0) node[graynode](n4){};
\path (4,0) node[blacknode](n5){};
\path (5,0) node[rednode](n6){};
\path (6,0) node[graynode](n7){};
\path (7,0) node[graynode](n8){};
\path (8,0) node[graynode](n9){};

\path (0,-0.5) node[](){1};
\path (1,-0.5) node[](){2};
\path (2,-0.5) node[](){3};
\path (3,-0.5) node[](){4};
\path (4,-0.5) node[](){5};
\path (5,-0.5) node[](){6};
\path (6,-0.5) node[](){7};
\path (7,-0.5) node[](){8};
\path (8,-0.5) node[](){9};

\path[->] (n5) edge [bend left] (n7);
\path[->] (n0) edge [bend left] (n1);
\end{tikzpicture}\end{center}

\begin{align*}
C & = \{7\} \\
\lambda' & = (5,4,4,3,3,2,1) \\ 
\end{align*}

\chapter{Symmetric functions and Macdonald polynomials}

In this chapter we recall some of the theory of symmetric functions. 

\section{Symmetric functions}

\subsection{Compositions}

A \emph{composition} is simply a list of non-negative integers. For example:
\[ c = (2,0,1,2) \]

The \emph{symmetric group} $\mathcal{S}_m$ acts naturally on the set of compositions with $m$ parts. Each orbit of $\mathcal{S}_m$ contains a unique partition. 
For example, the partition corresponding to our example composition above is:
\[ \lambda = (2,2,1,0) \]

Note that we allow for the possibility of trailing zeros.
For any composition $\lambda$ let $r(\lambda,i)$ denote the number of parts of $\lambda$ equal to $i$. 
and let:
$$r_\lambda = \prod_{i \geq 0} r(\lambda,i)!$$

Then $r_\lambda$ is the order of the subgroup of $\mathcal{S}_m$ that stabilizes the composition ${\lambda}$. 
The number of distinct elements in the order of $\lambda$ is equal to $m!/r_\lambda$.
In our example $r_\lambda = 2$ and the number of distinct compositions is equal to $12$. 

\subsection{Multivariable polynomials}

Consider the multivariate polynomial ring $\mathbb{Q}[x_1, \ldots, x_m]$. 
Each monomial in $\mathbb{Q}[x_1, \ldots, x_m]$ corresponds to a composition. For example, the monomial $x_1^2 x_3 x_4^2$ in $\mathbb{Q}[x_1,x_2,x_3,x_4]$ corresponds to our example composition $(2,0,1,2)$.

For notational convenience,
if $\eta = (\eta_1, \eta_2, \ldots, \eta_m)$ is a composition, then by $X^\eta$ we mean the monomial $x_1^{\eta_1} x_2^{\eta_2} \ldots x_m^{\eta_m}$. There is a natural addition on the space of compositions which corresponds to multiplication in the polynomial ring. if $\eta = (\eta_1, \ldots, \eta_m)$ and $\gamma = (\gamma_1,  \ldots, \gamma_m)$ then $\eta + \gamma = (\eta_1 + \gamma_1,  \ldots, \eta_m + \gamma_m)$, and $X^{(\eta + \gamma)} = X^\eta X^\gamma$.

\subsection{Monomial symmetric functions}

For each partition $\lambda$ we define the {\it monomial symmetric function} to be:
$$ m_\lambda(X) = \frac{1}{r_\lambda} \sum_{\sigma \in S_n} X^{\sigma(\lambda)} $$

In other words the monomial symmetric
function $m_\lambda(X)$ is the sum of all {\it distinct} permutations of the monomial $X^\lambda$.
An example:

\begin{eqnarray*}
m_{(2,2,1,0)}(X) 
=  x_1^{2}x_2^{2}x_3 + x_1^{2}x_2x_3^{2} + x_1x_2^{2}x_3^{2}
 + x_1^{2}x_2^{2}x_4 + x_1^{2}x_2x_4^{2} + x_1x_2^{2}x_4^{2} \\
 + x_1^{2}x_3^{2}x_4 + x_1^{2}x_3x_4^{2} + x_1x_3^{2}x_4^{2} 
 + x_2^{2}x_3^{2}x_4 + x_2^{2}x_3x_4^{2} + x_2x_3^{2}x_4^{2} 
\end{eqnarray*}

The monomial symmetric functions have the property that they are invariant under the action of the symmetric group on the set of variables.
Furthermore, any function which is symmetric under the action of the symmetric group on the variables may be expressed as a linear combination of monomial symmetric functions. 

\subsection{Infinitely many variables}

If there are $m$ variables, then we only need to consider integer partitions with at most $m$ parts.
Sometimes it will be necessary to work with an infinite number of variables $\{x_1,x_2,x_3, \ldots \}$. 
In this case we take any partition which is padded at the end with an infinite number of zeros.

\section{Plethystic notation}\label{sec:plethystic}

\subsection{Alphabets}

We shall use the notation $\Lambda$ to denote the ring of symmetric functions, in infinitely many variables, over the field of rational numbers \cite{macdonald}.
In what follows we shall make extensive use of the \emph{plethystic notation} \cite{garsia,lascoux}.

\index{Plethystic notation}

In the plethystic notation addition corresponds to the union of two sets and multiplication corresponds to the Cartesian product. 
For example, we write:
\begin{equation} X = x_1 + x_2 + \cdots \end{equation}
to denote the set of variables $\{x_1, x_2, \ldots \}$. We also write:
\begin{equation} XY = (x_1 + x_2 + \cdots)(y_1 + y_2 + \ldots) \end{equation}
to denote the set of variables $\{x_1 y_1, x_1 y_2, \ldots, x_2 y_1, \ldots x_2 y_2 \ldots \}$.

\subsection{Complete and elementary symmetric functions}

The \emph{complete symmetric functions} may be defined by their generating series:
\begin{equation} \label{eq:complete}
 \Omega[Xz] = \prod_{i\geq1} \frac{1}{1-x_iz} = \sum_{n\geq0} h_n(X) z^n 
\end{equation}

The generating series for the \emph{elementary symmetric functions} may be expressed using the \emph{plethystic negation} of an alphabet \cite{lascoux,garsia}:
\begin{equation}
  \Omega[-Xz] = \prod_{i\geq1} (1-x_iz) = \sum_{n\geq0} (-1)^n e_n(X) z^n 
\end{equation}

It helps to think of the generating function for the complete symmetric functions as a sort of exponential function:
\begin{align}
\Omega[X+Y] &= \Omega[X] \Omega[Y] \\
\Omega[-X] &= \Omega^{-1}[X]
\end{align}

For any integer partition $\lambda$ we define:
\begin{align}
h_\lambda(X) & = h_{\lambda_1}(X) h_{\lambda_2}(X) h_{\lambda_3}(X) \cdots \\
e_\lambda(X) & = e_{\lambda_1}(X) e_{\lambda_2}(X) e_{\lambda_3}(X) \cdots 
\end{align}

One can show that, just like the monomial symmetric functions $\{m_\lambda(X)\}$ both $\{h_\lambda(X)\}$ and $\{e_\lambda(X)\}$ are bases for the ring of symmetric functions. 

\subsection{Hall inner production}

\index{Symmetric function!$\Omega(Xz)$ definition}

It is a straightforward exercise to verify that:
\begin{equation}
  \Omega[XY] = \prod_{i,j\geq1} \frac{1}{1-x_i y_j} = \sum_\lambda h_\lambda(X) m_\lambda(Y) 
\end{equation}

\index{Schur function!Hall inner product}

The \emph{Hall inner product} $\langle - | - \rangle$ may be characterized by th property that:
\begin{equation}\label{eq:hall}
\langle m_\lambda(X), h_\mu(X) \rangle = \delta_{\lambda, \mu}
\end{equation}

Equivalently, the Hall inner product is characterized by the fact that for any pair of bases $\{f_\lambda(X)\}$ and $\{g_\lambda(X)\}$ such that:
\[ \langle f_\lambda(X), g_\mu(X) \rangle = \delta_{\lambda, \mu} \]
we have:
\begin{equation}
  \Omega[XY] = \prod_{i,j\geq1} \frac{1}{1-x_i y_j} = \sum_\lambda f_\lambda(X) g_\lambda(Y) 
\end{equation}

\section{Schur functions}\label{sec:schur}

\subsection{Definition in terms of semistandard Young tableau}

The most important basis for the ring of symmetric functions is the \emph{Schur basis}.
There are many different ways in which the Schur functions may be defined. Each has its advantages and its disadvantages.
For now we shall make use of the following definition. 

\index{Schur function!definition}
For any integer partition $\lambda$, the \emph{Schur function} is given by:
\begin{equation}\label{schur}
 S_\lambda(X) = \sum_{T \in sh(\lambda)} x_1^{c_1(T)} x_2^{c_2(T)} x_3^{c_3(T)} \cdots 
\end{equation}
where the sum is over all semistandard Young tableau of shape $\lambda$ and: 
\[ (c_1(T),c_2(T), \ldots )\] 
denotes the \emph{content vector} of $T$.
Note that it is not immediately clear that the Schur function as defined above is symmetric, See \cite{schutz} for a proof. 

\subsection{Cauchy identity and RSK}

\index{Schur function!Cauchy identity}

The \emph{Cauchy identity} is an immediate consequence of the definition of Schur functions in terms of semi-standard Young tableau:

\begin{equation}\label{eq:cauchy}
  \Omega[XY] = \prod_{i,j\geq1} \frac{1}{1-x_i y_j} = \sum_\lambda S_\lambda(X) S_\lambda(Y) 
\end{equation}

The left hand side may be interpreted as a weighted sum over non-negative integer matrices.
The \emph{refined weight} of a matrix $M$ is given by: 
\[ (x_1^{d_1} x_2^{d_2} \cdots )(y_1^{d'_1} y_2^{d'_2} \cdots) \]
where $d_k$ is the sum of all the entries in the $k$-th row of $M$ and $d'_k$ is the sum of all the entries in
the $k$-th column of $M$.
 
The right hand side may be interpreted as a weighted sum over pairs of semi-standard Young tableau of the same shape.
The refined weight of a pair of semi-standard Young tableau $(T,T)$ is given by:
\[ (x_1^{c_1} x_2^{c_2} \cdots) (y_1^{c'_1} y_2^{c'_2} \cdots) \]
where $(c_1,c_2, \ldots )$ is the content vector of $T$ while $(c'_1,c'_2, \ldots )$ is the content vector of $T'$.
This is none other than the RSK correspondence (see Section \ref{sec:RSK-reverse}).
We remark that MacMahon's identity may be obtained from Cauchy's identity by specializing the variables to:
\[ x_k = y_k = z^{k-\frac{1}{2}} \]

\subsection{Pieri rules}

The following recurrence is an immediate consequence of the definition of the Schur function given in equation \ref{schur}: 
\begin{equation}
S_\lambda[X+z] = \sum_{\mu \in D(\mu)} S_\mu[X] z^{|\lambda|- |\mu|}
\end{equation}
Here $D(\lambda)$ denotes the set of partitions which can be obtained from $\lambda$ by removing a horizontal strip.
This is known as the dual Pieri rule.

Assuming that we have proved the symmetry of the Schur function, we may write:
\begin{equation}\label{fish1}
S_\lambda(X) = \sum_{\mu} K_{\lambda,\mu}\, m_\mu(X)
\end{equation}
where $K_{\lambda,\mu}$ is the \emph{Kostka number} which counts the number of semi-standard Young tableau of shape $\lambda$ and \emph{content} $\mu$ (see Section{\ref{sec:semi-standard}).
Note that $K_{\lambda,\mu} = 0$ unless $\lambda$ is greater than $\mu$ in the dominance order (see Section \ref{sec:dominance}).

Since Cauchy's identity tells us that the Schur-functions form an orthonormal basis with respect to the Hall inner product, we have by duality that:
\begin{equation}\label{fish2}
S_\mu(X) = \sum_{\lambda} K_{\lambda,\mu} \,h_\lambda(X)
\end{equation}

\index{Schur function!dual Pieri formula}
The following recurrence is an immediate consequence of equation \ref{fish2}:
\begin{equation}\label{p1} \Omega[Xz] S_\mu[X] 
 = \sum_{\lambda \in U(\mu)} S_\lambda[X] z^{|\lambda|- |\mu|}\end{equation}
where  $U(\mu)$ denotes the set of partitions which can be obtained from $\mu$ by adding a horizontal strip.
This is known as the Pieri rule.

\index{Schur function!Pieri formula}

Combining equation \ref{p1} and equation \ref{fish1} and using another duality argument, we discover that:
\begin{equation}\label{p2} \Omega^*[Xz]S_\lambda[X] = S_\lambda[X+z]  \end{equation}
where $\Omega^*[Xz]$ is defined to be adjoint to the operator $\Omega[Xz]$ with respect to the Hall inner product.
\begin{equation}
 \langle f(X) \, | \, \Omega^*[Xz] g(X) \rangle = \langle \Omega[Xz] f(X) \, | \, g(X) \rangle  
\end{equation}

\index{Symmetric function!$\Omega^*(X)$ definition}

\section{Local rules}\label{sec:sym-local}
To show how these ideas are connected, let us assume that the Pieri rule and the dual Pieri rule are true and attempt to deduce Cauchy's identity as a consequence. We will then prove an important commutation relation between the Pieri operator and the dual Pieri operator

\subsection{Alternative proof of Cauchy identity}

\begin{proposition}\label{cauchy}
If $X_n$ is the alphabet $X_n = x_1 + x_2 + \cdots + x_n$
and $Y_n$ is the alphabet $Y_n = y_1 + y_2 + \cdots + y_n$ then:
\[ \sum_\lambda S_\lambda[X_n] S_\lambda[Y_n] = \prod_{i,j=1}^n \frac{1}{1-x_iy_j} \]
\begin{proof}
The case $n=1$ is trivial.
We assume that the result is true for all $k < n$.
By the Pieri formula (equation \ref{p1}), and the induction assumption, we have on the right hand side:

\begin{align*}
RHS(n) & = \left (\sum_r h_r[X_{n-1}] x_n^r \right ) 
\left (\sum_s h_s[Y_{n-1}] y_n^s \right ) \left (  \sum_\mu S_\mu[X_{n-1}] S_\mu[Y_{n-1}] \right )\frac{1}{1-x_n y_n} \\
& = \sum_{\alpha, \beta} S_{\alpha}[X_{n-1}] S_{\beta}[Y_{n-1}] \left ( \frac{1}{1-x_n y_n} \sum_{\mu \in D(\alpha) \cap D(\beta)} x_n^{|\alpha / \mu|} y_n^{|\beta / \mu|} \right )
\end{align*}

while by the dual Pieri formula (equation \ref{p2}), we have on the left hand side:

\begin{align*}
LHS(n) & = \sum_\lambda S_\lambda[X_{n-1} + x_n] S_\lambda[X_{n-1} + y_n] \\
& = \sum_{\alpha, \beta} S_\alpha[X_{n-1}] S_\beta[Y_{n-1}] \left ( \sum_{\lambda \in U(\alpha) \cap U(\beta)} x_n^{|\lambda / \alpha|} y_n^{|\lambda / \beta|} \right )
\end{align*}

In order to complete this proof, we need, for any given $\alpha$ and $\beta$ a bijection between the set of pairs $(m,\mu)$ with $\mu \in D(\alpha) \cup D(\beta)$ and $m$ a non-negative integer and the set of partitions $\lambda \in U(\alpha) \cup U(\beta)$ in such a way that:
\[ |\lambda|  = |\alpha | + |\beta| - |\mu| + m \]
A \emph{local rule}  is just such a bijection (see section \ref{sec:local}). 
\begin{equation}
\mathfrak{U}_{\alpha,\beta} : (\mathbb{Z}_{\geq 0}, D(\alpha) \cup D(\beta)) \to U(\alpha) \cup U(\beta)
\end{equation}

\end{proof}
\end{proposition}

\index{Local rule!symmetric functions}

\subsection{Commutation Relations}\label{omega-commutation}

\begin{proposition}\label{commute}
The operators $\Omega[Xv]$ and $\Omega^*[Xu]$ satisfy the following commutation relations:

\begin{equation}\label{commutation-pieri}
\boxed{ \Omega^*[Xu] \, \Omega[Xv] = \frac{1}{1-uv} \, \Omega[Xv] \,\Omega^*[Xu]}
\end{equation}
\begin{proof}
For any integer partition $\lambda$ we have:
\begin{align*}
 \Omega^*[Xu] \, \Omega[Xv] S_\lambda(X)
& = \Omega[(X+u)v] S_\lambda(X + u) \\
& = \Omega[uv] \Omega[Xv] \Omega^*[Xu] S_\lambda(X) \\
& = \frac{1}{1-uv} \Omega[Xv] \Omega^*[Xu] S_\lambda(X)
\end{align*}
\end{proof}
\end{proposition}

A local rule may be thought of as a bijective proof of proposition \ref{commute}, which may be expressed in the form:

\begin{equation}
\langle \, S_\alpha \, |  \,\Omega^*[Xu] \, \Omega[Xv] \, S_\beta  \,\rangle
= \frac{1}{1-uv} \, \langle  \,S_\alpha  \,|  \,\Omega[Xv] \,\Omega^*[Xu]  \,S_\beta  \,\rangle 
\end{equation}
for all integer partitions $\alpha$ and $\beta$.
In other words, a local rule is a map with type signature:

\begin{equation} 
\mathfrak{D}_{\alpha,\beta} : U(\alpha) \cap U(\beta) \to (\mathbb{Z}_{\geq 0}, D(\alpha) \cap D(\beta))
\end{equation}

Such that if $(m, \mu) = \mathfrak{D}_{\alpha,\beta}(\lambda)$ then the following weight conditions are satisfied:

\begin{align} 
|\lambda / \alpha| & = |\beta / \mu| + m \\
|\lambda / \beta| & = |\alpha / \mu| + m 
\end{align}

By adding these two weight conditions together and dividing by $2$ we obtain:
\begin{equation} 
 |\lambda|  = |\alpha | + |\beta| - |\mu| + m 
\end{equation}

This is precisely the inverse of the local rule described at the end of proposition \ref{cauchy}. Again see section \ref{sec:local}.

\section{Robinson correspondence revisited}\label{sec:rep}

\subsection{Representation Theory}
\subsubsection{Permutations}

Any permutation may be expressed as the product of disjoint cycles. For example, the permutation: 
\[ \sigma = \left [ \begin{matrix} 1 & 2 & 3 & 4 & 5 & 6 & 7 & 8 & 9 \\ 4&8&3&6&5&7&9&2&1 \end{matrix} \right ] \]
may be expressed in \emph{disjoint cycle notation} as  $\sigma = (1,4,6,7,9)(2,8)(3)(5)$.
Since disjoint cycles commute with one another, they can always be ordered from longest to shortest.
The \emph{cycle type} of a permutation is the integer partition which determines the length of the cycles.
The cycle type of our example partition is $(5,2,1,1)$.

A pair of permutations $\sigma_1$ and $\sigma_2$ are said to be \emph{conjugate} if there exists some third permutation $\tau$ such that
$\sigma_1 = \tau \sigma_2 \tau^{-1}$. Two permutations are conjugate if and only if they have the same cycle type.

\subsubsection{Characters}

A \emph{representation} of the symmetric group $\mathfrak{S}_n$ is a homomorphism from the $\mathfrak{S}_n$ to $\GL(V)$ for some complex vector space $V$. A representation is said to be \emph{irreducible} if there is no non-trivial subspace $U \subseteq V$ with the property that $\sigma.u \in U$ for all $u \in U$ and for all $\sigma \in \mathfrak{S}_n$.

One can show that the irreducible representations of the symmetric group are naturally indexed by integer partitions \cite{yellow}. Furthermore the \emph{dimension} of the irreducible representation indexed by $\lambda$, which we denote by $f^\lambda$, is given by the number of standard tableau of shape $\lambda$. \cite{yellow}.

Suppose that $\rho : \mathfrak{S}_n \to \GL(V)$ is a representation of the symmetric group. The \emph{character}
associated to $\rho$ is the map $\chi : \mathfrak{S}_n \to \mathbb{C}$ given by $\chi(\sigma) = \tr(\sigma)$.
Since the trace is invariant on conjugacy classes, it makes sense to write $\chi(\mu)$ where $\mu$ is a partition of $n$.

Note that for the identity permutation we have:
\begin{equation}\label{eq:magic}
\chi^\lambda(1^n) = f^\lambda
\end{equation}

\subsection{Newton power sums}
\subsubsection{Definition}

The Newton power sums are defined by:
\begin{equation} 
p_n(X) = x_1^n + x_2^n + \cdots 
\end{equation}
As in the case of the complete symmetric functions and the elementary symmetric functions, if we define:
\begin{equation} 
p_\lambda(X) = p_{\lambda_1}(X) p_{\lambda_2}(X) \cdots
\end{equation}
then $\{p_\lambda(X)\}$ is a basis for the ring of symmetric functions. 

\subsubsection{Relation to complete symmetric functions}
For any integer partition $\lambda$ let us define:
\begin{equation}
z_\lambda = \prod_i m_i(\lambda)! i^{m_i(\lambda)}
\end{equation}
where $m_i(\lambda)$ denotes the number of parts of $\lambda$ of length $i$.
As a consequence of the fact that:
\begin{equation}
\log \left ( \frac{1}{1-x} \right ) = \sum_n \frac{x^n}{n}
\end{equation}
we have that:
\begin{equation} \Omega(Xz) = \exp \left ( \sum_{n \geq 1} \frac{p_n(X)}{n} z^n \right )\end{equation}
In particular:
\begin{equation}
h_n(X) = \sum_{\lambda \vdash n} \frac{p_\lambda(X)}{z_\lambda}
\end{equation}

\subsubsection{Dual basis}

Using the notation from Section \ref{sec:plethystic} we have:
\begin{equation}
p_n(X)p_n(Y) = p_n(XY)
\end{equation}
Now:
\begin{align}
\Omega(XY) & = \exp \left ( \sum_{n \geq 1} \frac{p_n(XY)}{n} z^n \right ) \\
& = \exp \left ( \sum_{n \geq 1} \frac{p_n(X)p_n(Y)}{n} z^n \right ) \\
& = \sum_{\lambda} \frac{p_\lambda(X)}{z_\lambda}P_\lambda(Y)
\end{align}
It follows that the dual basis to the $\{P_\lambda(X)\}$ is given by:
\[\left \{ \frac{p_\lambda(X)}{z_\lambda} \right \}\]
\index{Schur function!definition}

\subsection{Algebraic proof of Robinson correspondence}

An alternative definition for the Schur functions is the following:
\begin{equation}
S_\lambda(X) = \sum_\mu \frac{\chi^\lambda(\mu)}{z_\lambda} p_\mu(X)
\end{equation}
where $\chi^\lambda(\mu)$ denotes the \emph{character} of the representation of $\mathfrak{S}_n$ 
indexed by $\lambda$ evaluated on the conjugacy class of type $\mu$.

The Robinson correspondence (section \ref{sec:robinson}) tells us that:
\begin{equation} \label{magic}
\boxed{n! = \sum_\lambda (f_\lambda)^2}
\end{equation}
This identity may be derived from Cauchy's formula (equation \ref{eq:cauchy}) as follows.
On the left hand side we have:
\begin{equation}
\Omega(XY) = \sum_n h_n(XY) = \sum_\lambda \frac{p_\lambda(XY)}{z_\lambda} 
\end{equation}
and so:
\begin{equation}
\langle \Omega(XY) | p_{(1^n)}(XY) \rangle_{XY} = \frac{1}{z_{(1)^n}} = \frac{1}{n!}
\end{equation}
While on the right hand side we have:
\begin{align}
& \left \langle \sum_\lambda S_\lambda(X) S_\lambda(Y) \,|\, p_{(1^n)}(X) p_{(1^n)}(Y) \right \rangle_{XY} \\
& = \sum_\lambda \left \langle S_\lambda(X) | p_{(1^n)}(X) \rangle_{X}\,\langle S_\lambda(Y) | p_{(1^n)}(Y) \right \rangle_{Y} \\
& = \sum_\lambda \left ( \frac{\chi^\lambda((1^n))}{n!}\right )^2 \\
& = \sum_\lambda \left ( \frac{f^\lambda}{n!}\right )^2
\end{align}
Putting the two sides together we recover equation \ref{magic}.

\subsection{Canonical commutation relations}\label{weyl}
\index{Canonical commutation relations}

Let $p_1^*$ denote the adjoint of the operator ``multiplication by $p_1$'' with respect to the Hall inner product. The pair $(p_1, p_1^*)$ satisfy the \emph{canonical commutation relations}:
\begin{equation}
[ p_1^*, p_1 ] = 1
\end{equation}
Noting that $p_1(X) = h_1(X)$ this can be understood by the fact that an integer partition always has one more outside corner box than inside corner box. Fomin's original local rule can be understood as a bijective proof of the identity:
\begin{equation}\label{fomin-com}
\langle \, S_\alpha \, |  \, p_1^* p_1 S_\beta  \,\rangle
= 1 + \, \langle  \,S_\alpha  \,|  \,p_1 p_1^*  \,S_\beta  \,\rangle 
\end{equation}
Equation \ref{magic} can be proved algebraically by remarking that:
\begin{align}
\sum_\lambda f_\lambda f_\lambda &= 
\sum_\lambda \langle 1 p_1^n | S_\lambda \rangle \langle S_\lambda | p_1^n 1 \rangle \\
& = \langle 1 p_1^n | p_1^n 1 \rangle \\
 & =\langle 1 | (p_1^*)^n(p_1)^n 1 \rangle \\
& = n!
\end{align}
The last step makes use of the commutation relations in equation \ref{fomin-com}
together with the fact that $p_1^* 1 = 0$. Each term of the form $p_1^*$ must ``annihilate'' with some term of the form $p_1$.
There are exactly $n!$ ways in which this can happen.

\section{Macdonald Polynomials}

\subsection{Plethystic notation}
\index{Plethystic notation}

In this section we recall some of the theory of Macdonald polynomials.
Let $\Lambda_{q,t}$ denote the ring of symmetric functions over the field $\mathbb{Q}(q,t)$ of rational functions in the indeterminants $q$ and $t$.
Making use of the plethystic notation, the expression: 
\[ \frac{1-t}{1-q}\]
may be interpreted as the alphabet:
\begin{equation*}
 \{ 1, q, q^2, q^3, \cdots, -t, -tq, -tq^2, -tq^3. \cdots \} 
\end{equation*}
So that, for example:
\begin{equation} 
\Omega \left [ \frac{1-t}{1-q}z\right ] = \frac{(tz;q)_\infty}{(z;q)_\infty} 
\end{equation}
\index{Pochhammer symbol}

where the \emph{Pochhammer symbol} is defined by:
\begin{equation} 
(a;q)_\infty = \prod_{n \geq 0} (1-aq^n) 
\end{equation}

\subsection{$(q,t)$-Pieri operators}
\index{Macdonald polynomial!Pieri operator}
We shall define:
\begin{equation} 
 \Omega_{q,t}[Xz] = \Omega \left [ \frac{1-t}{1-q}Xz \right ] = \prod_{i \geq 1} \frac{(tx_iz;q)_\infty}{(x_iz;q)_\infty} 
\end{equation}

Note that when $q=t$ this reduces to the generating series for complete symmetric function (see equation \ref{eq:complete}).
\index{Macdonald polynomial!inner product}
\index{Macdonald polynomial!Cauchy identity}
Let $\langle -\,|\,-\rangle_{q,t}$ denote the \emph{Macdonald inner product} associated to the $(q,t)$-deformed Cauchy Kernel: 
\begin{equation}
 \Omega_{q,t}[XY] = \Omega \left [ XY \frac{1-t}{1-q} \right ] = \prod_{i,j \geq 1} \frac{(tx_i y_j;q)_\infty}{(x_i y_j;q)_\infty} 
\end{equation}
When $q=t$ this reduces to the usual Cauchy kernel (see equation \ref{eq:cauchy}).
\index{Macdonald polynomial!dual Pieri operator} 
The operator $\Omega^*_{q,t}[Xz]$ is defined to be adjoint to the operator $\Omega_{q,t}[Xz]$ with respect to the Macdonald inner product.

\begin{equation*}
 \langle f(X) \, | \, \Omega^*_{q,t}[Xz] g(X) \rangle_{q,t} = \langle \Omega_{q,t}[Xz] f(X) \, | \, g(X) \rangle_{q,t}  
\end{equation*}

When $q=t$ the Macdonald inner product reduces to the Hall inner product (see equation \ref{eq:hall}).

\subsection{Definition}\label{sec:macdonald}
\index{Macdonald polynomial!definition}

The Macdonald polynomials $\{P_\lambda(X;q,t)\}$ may be defined as the unique basis for $\Lambda_{q,t}$ which is both orthogonal
with respect to the Macdonald inner product:
\begin{equation}
\langle P_\lambda(X;q,t), P_\mu(X;q,t) \rangle_{q,t} = 0 \mbox{ if } \lambda \neq \mu
\end{equation}
and which admits a triangular change of basis with respect to the monomial symmetric functions:
\begin{equation}
P_\lambda(X) = m_\lambda(X) + \sum_{\mu < \lambda} \Box_{\lambda \mu} \, m_\mu(X)
\end{equation}
Here we are using the \emph{lexicographical order} on integer partitions. The dual basis is denoted by $\{Q_\lambda(X;q,t)\}$.
When $q=t$ the Macdonald polynomials reduce to the Schur functions.

\subsubsection{Curious properties}

It is a surprising fact \cite{macdonald} (page 338 -- 340) that, just like the Schur functions, the Macdonald polynomials admit an expansion in terms of the monomial symmetric functions of the form:
\begin{equation}
P_\lambda(X) = m_\lambda(X) + \sum_{\mu \triangleleft \lambda} \Box_{\lambda \mu} \, m_\mu(X)
\end{equation}
Here we are using the \emph{dominance order} on integer partition.

Another surprising fact about the Macdonald polynomials is that, when expanded in terms of the Schur functions, the coefficients $K_{\lambda,\mu}(q,t)$ are polynomials in $q$ and $t$ with non-negative integer co-efficients.
\begin{equation}
P_\lambda(X) = \sum_\mu K_{\lambda,\mu}(q,t) S_\mu(X)
\end{equation}

This was known as \emph{Macdonald's positivity conjecture} until it was eventually proved by Haiman using difficult ideas from algebraic geometry \cite{Haiman}. It remains an open problem to find a combinatorial proof.

There is a simple closed form expression for the inner product of a Macdonald polynomial with itself:

\begin{equation}
\langle P_\lambda(X;q,t), P_\lambda(X;q,t) \rangle_{q,t} = \prod_{s \in \lambda} \frac{(1 - q^{a_\lambda(s)+1} t^{\ell_\lambda(s)})}{(1 - q^{a_\lambda(s)} t^{\ell_\lambda(s)+1}) } 
\end{equation}
see Section \ref{sec:arms-legs} for the definition of arms and legs. This is not at all obvious from the definition which we have given. To prove it requires an alternative definition of the Macdonald polynomials as eigenvectors of a certain family of commuting operators \cite{macdonald}.

\subsubsection{Pieri formulae}\label{sec:macdonald-pieri}

\index{Macdonald polynomial!Pieri formula}
\index{Macdonald polynomial!dual Pieri formula}
Like the Schur functions, there is a Macdonald Pieri formula, and a Macdonald dual Pieri formula:
\begin{equation} \label{pieri}
 \Omega[Xz]_{q,t} \, P_\mu(X;q,t) = \sum_{\lambda \in U(\mu)}\psi_{\lambda / \mu}(q,t) \, P_\lambda(X;q,t) z^{|\lambda|-|\mu|}
\end{equation}
\begin{equation} \label{pieri2}
 \Omega^*[Xz]_{q,t} \, P_\lambda(X;q,t) = \sum_{\mu \in D(\lambda)}\varphi_{\lambda / \mu}(q,t) \, P_\mu(X;
q,t) z^{|\lambda|-|\mu|}
\end{equation}
The Macdonald Pieri coefficients are given (\cite{macdonald} page 341) by:

\begin{equation}\label{pieri-coeff}
 \varphi_{\lambda / \mu}(q,t) = \prod_{s \in C_{\lambda / \mu}}\frac{1 - q^{a_\lambda(s)} t^{\ell_\lambda(s)+1}}{1 -q^{a_\lambda(s)+1} t^{\ell_\lambda(s)}} 
\prod_{s \in {C}_{\lambda / \mu}}\frac{1 - q^{a_\mu(s)+1} t^{\ell_\mu(s)}}{1 -q^{a_\mu(s)} t^{\ell_\mu(s)+1}}
\end{equation}
\begin{equation}\label{pieri-coeff2}
 \psi_{\lambda / \mu}(q,t) = \prod_{s \not\in C_{\lambda / \mu}}\frac{1 - q^{a_\lambda(s)+1} t^{\ell_\lambda(s)}}{1 - q^{a_\lambda(s)} t^{\ell_\lambda(s)+1}} 
\prod_{s \not\in {C}_{\lambda / \mu}}\frac{1 -q^{a_\mu(s)} t^{\ell_\mu(s)+1}}{1 -q^{a_\mu(s)+1} t^{\ell_\mu(s)}}
\end{equation}

Here $C_{\lambda / \mu}$ denotes the set of columns of $\lambda$ which are longer than the corresponding columns of $\mu$.
When $q=t$ the Pieri coefficients are equal to $1$. 

\subsubsection{Key lemmas}

The following two lemmas are essentially due to Garcia, Haiman and Tesler \cite{garsia}. 
They constitute a $(q,t)$-analog of the commutation relations for ``vertex operators'' to be found in Jimbo and Miwa \cite{soliton}

\begin{lemma}
\begin{equation} \Omega^*_{q,t}[Xz] \, P_\lambda(X;q,t) = P_\lambda(X+z;q,t) \end{equation}
\begin{proof}
Let $\{ Q_\lambda(X;q,t) \}$ denote the dual basis to the $\{P_\lambda(X;q,t)\}$ with respect to the Macdonald inner product.
We have:
\begin{align*}
\Omega^*_{q,t}[Xz] \, P_\lambda(X;q,t) 
& = \langle \Omega^*_{q,t}[Yz] \, Q_\lambda(Y;q,t)\,|\, \Omega_{q,t}[XY] \rangle_{q,t} \\
& = \langle  Q_\lambda(Y;q,t)\,|\, \Omega_{q,t}[Yz] \, \Omega_{q,t}[XY] \rangle_{q,t} \\
& = \langle  Q_\lambda(Y;q,t)\,|\, \Omega_{q,t}[(X+z)Y] \rangle_{q,t} \\
& = P_\lambda(X+z;q,t)
\end{align*}
\end{proof}
\end{lemma}

\begin{lemma}\label{commutation-mac}
\begin{equation} 
\Omega_{q,t}^*[Xu] \, \Omega_{q,t}[Xv] = \frac{(tuv;q)_\infty}{(uv;q)_\infty} \, \Omega_{q,t}[Xv] \, \Omega_{q,t}^*[Xu]
\end{equation}
\begin{proof}
\begin{align*}
\Omega_{q,t}^*[Xu] \, \Omega_{q,t}[Xv] \, P_\lambda(X;q,t) 
& = \Omega_{q,t}[(X+u)z] \, P_\lambda(X + u;q,t)\\
& = \Omega_{q,t}[uz] \, \Omega_{q,t}[Xz] \, \Omega^*_{q,t}[Xz] \, P_\lambda(X;q,t) \\
& = \prod_{n \geq 0} \frac{(tuv;q)_\infty}{(uv;q)_\infty} \, \Omega_{q,t} [Xz] \, \Omega^*_{q,t}[Xz] \, P_\lambda(X;q,t)
\end{align*}
\end{proof}
\end{lemma}

\subsection{Hall--Littlewood polynomials}
\index{Hall--Littlewood polynomials}

As a final remark, when $q=0$ the Macdonald polynomials reduce to the \emph{Hall--Littlewood polynomials}.

\chapter{Bijective proof of Borodin's identity}

\section{Symmetric function proof of Borodin's identity}\label{sec:algebraic}
We shall begin by sketching an algebraic proof of Borodin's identity.
\begin{equation} \label{borodin}
\sum_{\mathfrak{c} \in \CPP(\pi)} z^{|\mathfrak{c}|} = 
\prod_{n \geq 0} \left ( \frac{1}{1-z^{(n+1)T}} \prod_{\substack{i < j \\ \pi_i > \pi_j}} \frac{1}{1-z^{j-i + nT}} 
\prod_{\substack{i > j \\ \pi_i > \pi_j }} \frac{1}{1 - z^{j-i + (n+1)T}} \right ) 
\end{equation}

\index{Borodin identity!Schur!refined}

We shall actually be proving a refined version of equation \ref{borodin}
in which on the left hand side we replace:
\begin{equation}
z^{\mathfrak{|c|}} \mapsto z_0^{|\mu_0|} z_1^{|\mu_1|} \cdots z_{T-1}^{|\mu_{T-1}|}
\end{equation}
while on the right hand side we replace:
\begin{align} 
z^{nT} & \mapsto z_0^{n} z_1^{n} \cdots z_{T-1}^{n} \\
z^{j - i + nT} & \mapsto z_0^{n} z_1^n \cdots z_i^{n} z_{i+1}^{n+1} \cdots z_{j}^{n+1} z_{j+1}^n \cdots z_{T-1}^n \quad \quad \quad \quad \text{ when } i < j \\
z^{j-i +(n+1)T} & \mapsto z_0^{n+1} z_1^{n+1} \cdots z_j^{n+1} z_{j+1}^n + \cdots z_i^{n} z_{i+1}^{n+1} \cdots z_{T-1}^{n+1} \quad \text{ when } i > j
\end{align}

The refined version of the reverse plane partition case is due to Gasner \cite{gasner-rpp}.
We would like to emphasize the fact that the structure of the bijective proof follows very closely the structure of the algebraic proof.

 
\subsection{Notation}\label{sec:formulation}

Let $D_z$ denote the ``degree'' operator:
\begin{equation} \label{commutation}
D_z S_\lambda[X] = z^{|\lambda|} S_\lambda[X] 
\end{equation}

The degree operator satisfies the following commutation relations:
\begin{lemma} \label{eq:degree}
\begin{align}
D_z \, \Omega[Xu] & = \Omega[Xuz] \, D_z \\
D_z \, \Omega^*[Xu] & = \Omega^*[Xuz^{-1}] \, D_z 
\end{align}
\begin{proof}
This fact follows immediately from equations \ref{p1} and \ref{p2}.
\end{proof}
\end{lemma}

For notational convenience we shall define:
\begin{align}
G^0(z) & = \Omega[Xz] \\
G^1(z) & = \Omega^*[Xz]
\end{align}

\subsection{Algebraic interpretation of cylindric plane partition}

\begin{lemma}  \label{lhs}
The left hand side of the refined version of equation \ref{borodin} may be expressed in the form:
\begin{equation} 
\lhs(\pi) = \sum_\mu \langle S_\mu \,|\, G^{\pi_0}(u_0) G^{\pi_1}(u_1)  \cdots G^{\pi_T}(u_T) D_w \, S_\mu  \rangle 
\end{equation}
where:
\begin{align} \label{specialization}
w & = z_0 z_1 \cdots z_{T-1} \\ 
u_k & = \begin{cases}
z_0 z_1 \cdots z_{k-1} & \text{ if $\pi_k = 1$} \\
z_0^{-1} z_1^{-1} \cdots z_{k-1}^{-1} & \text{ if $\pi_k = 0$}
\end{cases}\label{specialization2}
\end{align}

\begin{proof}
From the ``interlacing sequence'' definition of a cylindric plane partitions \ref{cpp} it is clear that a cylindric plane partition is constructed
by successively adding and removing horizontal strips. 
The degree operator $D_z$ is used to keep track of the number of cubes in the resulting cylindric plane partition.

Using the fact that the Schur functions are orthonormal with respect to the hall inner product
we may write:
\begin{equation}
\lhs(\pi) = \sum_\mu \langle S_\mu \,|\, D_{z_0} \, G^{\pi_0}(1) \, D_{z_1} \, G^{\pi_1}(1)   \cdots D_{z_{T-1}} \, G^{\pi_T}(1) \, S_\mu  \rangle
\end{equation}

It remains to commute all the shift operators to the right hand side using Lemma \ref{eq:degree}.
\end{proof}
\end{lemma}

Note that the above expression for $\lhs(\pi)$ can also be thought of as the \emph{trace} of the operator:
\[ G^{\pi_0}(u_0) G^{\pi_1}(u_1)  \cdots G^{\pi_T}(u_T) D_w \]
acting on symmetric functions.

\subsection{Some lemmas}

Let us define:
\begin{definition}\label{M-def}
\begin{equation} 
M_\pi(m)  = \sum_\mu \langle S_\mu \,|\, \prod_{\substack{k =1\\ \pi_k = 0}}^T \Omega [Xu_kw^m] 
\prod_{\substack{k =1 \\ \pi_k = 1}}^T \Omega^*[Xu_k] \, D_w \, S_\mu \rangle 
\end{equation}
\end{definition}

We have:
\begin{lemma} \label{cylindric-shift}
\begin{equation} M_\pi(m) =  \prod_{\substack{(i,j) \\ \pi_i \neq \pi_j}} \frac{1}{1-u_i u_j w^{m+1}} \, M_\pi(m+1) \end{equation}
\begin{proof}
This is a straightforward calculation. Using the fact that the Schur functions are orthogonal with respect to the Hall inner product, as well as the fact that:
\[ \tr(AB) = \tr(BA) \]
We may write:
\begin{align}
M_\pi(m) & = \sum_{\mu, \lambda} 
\langle S_\mu \,|\, \prod_{\substack{k =1\\ \pi_k = 0}}^T \Omega[Xu_kw^m] \, S_\lambda \rangle 
\langle S_\lambda \,|\,\prod_{\substack{k =1 \\ \pi_k = 1}}^T \Omega^*[Xu_k] \, D_w \, S_\mu \rangle \\
& = \sum_{\mu, \lambda} 
\langle S_\lambda \,|\,\prod_{\substack{k =1 \\ \pi_k = 1}}^T \Omega^*[Xu_k] \, D_w \, S_\mu \rangle 
\langle S_\mu \,|\, \prod_{\substack{k =1\\ \pi_k = 0}}^T \Omega[Xu_kw^m] \, S_\lambda \rangle \\ \label{trace}
& = \sum_{\lambda} 
\langle S_\lambda \,|\,\prod_{\substack{k =1 \\ \pi_k = 1}}^T \Omega^*[Xu_k] \, 
D_w \,\prod_{\substack{k =1\\ \pi_k = 0}}^T \Omega[Xu_kw^m] \, S_\lambda \rangle 
\end{align}
Next applying the commutation relations of Lemma \ref{eq:degree} and Proposition \ref{commute} we have:
\begin{align}
M_\pi(m) & = \sum_{\lambda} 
\langle S_\lambda \,|\,\prod_{\substack{k =1 \\ \pi_k = 1}}^T \Omega^*[Xu_k] \, 
D_w \,\prod_{\substack{k =1\\ \pi_k = 0}}^T \Omega[Xu_kw^m] \, S_\lambda \rangle  \\ 
& = \sum_{\lambda} 
\langle S_\lambda \,|\,\prod_{\substack{k =1 \\ \pi_k = 1}}^T \Omega^*[Xu_k] \, 
\,\prod_{\substack{k =1\\ \pi_k = 0}}^T \Omega[Xu_kw^{m+1}] \, D_w \, S_\lambda \rangle \\
& = \prod_{\substack{(i,j) \\ \pi_i \neq \pi_j}} \frac{1}{1-u_i u_j w^{m+1}} \sum_{\lambda} 
\langle S_\lambda \,| \,\prod_{\substack{k =1\\ \pi_k = 0}}^T \Omega[Xu_kw^{m+1}] \, 
\prod_{\substack{k =1 \\ \pi_k = 1}}^T \Omega^*[Xu_k] \, D_w \, S_\lambda \rangle \\
& =  \prod_{\substack{(i,j) \\ \pi_i \neq \pi_j}} \frac{1}{1-u_i u_j w^{m+1}} \, M_\pi(m+1)
\end{align}
\end{proof}
\end{lemma}
In the limit we have:

\begin{lemma} \label{infinity}
\begin{equation} 
M_\pi(\infty)  = \prod_{n \geq 1} \frac{1}{1-w^n} 
\end{equation}
\begin{proof}
In order for this limit to even make sense, we must have $|z_i| < 1$ for all $i$, in which case:
\begin{equation} \lim_{m \to \infty} \Omega[X u_k \omega^m] = 1 \end{equation}

Since $\Omega^*[X u_k]$ is a degree lowering operator, it follows that:
\begin{align}
 \lim_{m \to \infty} M_\pi(m) 
& =  \sum_\mu \langle S_\mu \,|\,  \prod_{\substack{k =1 \\ \pi_k = 1}}^T \Omega^*[Xu_k] \, D_w \, S_\mu \rangle \\
& = \sum_\mu \langle S_\mu | D_w \, S_\mu \rangle \\
& = \sum_\mu w^{|\mu|} \\
& = \prod_{n \geq 1} \frac{1}{1-w^n} 
\end{align}

\end{proof}
\end{lemma}

\subsection{The proof}\label{sec:proof}

The proof of the refined version of Theorem \ref{borodin} now proceeds as follows. We begin by applying Lemma \ref{lhs}

\begin{align}
\phantom{=} \sum_{\mathfrak{c} \in \CPP(\pi)}  z^{|\mathfrak{c}|} 
& = \sum_\mu \langle S_\mu \,|\, G^{\pi_0}(u_0) G^{\pi_1}(u_1)  \cdots G^{\pi_T}(u_T) D_w \, S_\mu  \rangle \\
\end{align}
Next we repeatedly applies the commutation relations of Lemma \ref{commutation}, followed by definition \ref{M-def}.
\begin{align*}
& = \prod_{\substack{i < j \\ \pi_i > \pi_j} }  \frac{1}{1-u_i u_j}
\sum_\mu \langle S_\mu \,|\, \prod_{\substack{k =1\\ \pi_k = 0}}^T \Omega[Xu_k] 
\prod_{\substack{k =1 \\ \pi_k = 1}}^T \Omega^*[Xu_k] \, D_w \, S_\mu \rangle \\
& = \prod_{\substack{i < j \\ \pi_i > \pi_j} } \frac{1}{1-u_i u_j} M_\pi(0) 
\end{align*}
We then repeatedly apply Lemma \ref{cylindric-shift}.
\begin{align*}
& = \prod_{\substack{i < j \\ \pi_i > \pi_j} } \frac{1}{1-u_i u_j} \prod_{m \geq 0} 
\left ( \prod_{\substack{(i,j) \\ \pi_i \neq \pi_j}} \frac{1}{1-u_i u_j w^{m+1}} \right ) M_\pi (\infty) 
\end{align*}
Splitting the second product into two, and combining it with the first we have:

\begin{align*}
& = \prod_{m \geq 1}
\left ( \prod_{\substack{i < j \\ \pi_i > \pi_j} } \frac{1}{1-u_i u_jw^{m-1}} \right )
\left ( \prod_{\substack{i > j \\ \pi_i > \pi_j} } \frac{1}{1-u_i u_jw^m}\right ) M_\pi (\infty)
\end{align*}
Finally, applying Lemma \ref{infinity} we have:
\begin{align*}
& = \prod_{m \geq 1} \frac{1}{1-w^m}
\left ( \prod_{\substack{i < j \\ \pi_i > \pi_j} } \frac{1}{1-u_i u_jw^{m-1}} \right )
\left ( \prod_{\substack{i > j \\ \pi_i > \pi_j} } \frac{1}{1-u_i u_jw^m}\right ) 
\end{align*}

To obtain the non-refined version of the Theorem, it suffices to take the following specialization of variables on both sides:
\begin{align} \label{specialization-borodin}
w & = z^{|T|} \\
u_k & = \begin{cases}
z^k & \text{ if $\pi_k = 1$} \\
z^{-k} & \text{ if $\pi_k = 0$}
\end{cases}\label{specialization2-borodin}
\end{align}

\section{Local rule as higher order function}\label{higher-order}
\subsection{Combinatorial formulation}

We saw in section \ref{alcd} that Borodin's identity (equation \ref{borodin}) may be expressed in the form:
\begin{align}
\sum_{\mathfrak{c} \in \CPP(\pi)} z^{|\mathfrak{c}|} & = \left ( \sum_{\gamma} z^{T\, |\gamma|} \right ) 
\left ( \sum_{s \in \widehat{\lambda}(\pi)} \frac{1}{1 - z^{h_{\widehat{\lambda}(\pi)}(s)}} \right ) 
\end{align}
where $\widehat{\lambda}(\pi)$ denotes the \emph{cylindric diagram} with profile $\pi$
and $h_{\widehat{\lambda}(\pi)}(s)$ denotes the \emph{cylindric hook length} of the box $s$.

The right hand side may be interpreted combinatorially as a weighted sum over pairs $(\gamma, \mathfrak{d})$
where $\gamma$ is an integer partition and $\mathfrak{d}$ is an \emph{arbitrarily labelled cylindric diagram}. In other words, we may write Borodin's identity as:
\begin{equation}
\sum_{\mathfrak{c} \in \CPP(\pi)} z^{|\mathfrak{c}|} =  \sum_{(\gamma,\mathfrak{d}) \in (\IP, \ALCD(\pi))} z^{|\mathfrak{d}| + T \, |\gamma|} 
\end{equation}

Here $\IP$ denotes the set of all integer partitions and $\ALCD(\pi)$ denotes the set of all arbitrarily labelled cylindric diagrams with profile $\pi$.

\index{Borodin identity!Bijection!strongly weight preserving}
Our goal is thus to find, for each possible profile $\pi$, a weight-preserving bijection between the
sets $\CPP(\pi)$ and the tuple $(\IP, \ALCD(\pi))$.

\begin{equation}\label{bijection-type} 
\psi_\pi : (\IP, \ALCD(\pi)) \to \CPP(\pi) 
\end{equation}

 Our bijection will be such that it actually proves the following refined identity:

\index{Cylindric diagram!refined weight}

\begin{multline} 
\sum_{\mathfrak{c} \in \CPP(\pi)} z_1^{|\mu_1|}z_2^{|\mu_2|} \cdots z_T^{|\mu_T|} = \\
 \sum_{(\gamma,\mathfrak{d}) \in (\IP, \ALCD(\pi))} z_1^{|\gamma| + |\diag(1)|_\mathfrak{d}} z_2^{|\gamma| + |\diag(2)|_\mathfrak{d}} \cdots z_T^{|\gamma| + |\diag(T)|_\mathfrak{d}} 
\end{multline}
See section \ref{sec:diag} for the definition of \emph{diag}.

\subsection{Definition of local rule}\label{sec:local-weight}

In section \ref{sec:local} we  defined the local rule to be a map with type signature:
\begin{equation}
\mathfrak{U}_{\alpha,\beta} : U(\alpha) \cap U(\beta) \to (\mathbb{Z}_{\geq 0}, D(\alpha) \cap D(\beta)) 
\end{equation}
which satisfies the following weight condition. If $\lambda = \mathfrak{D}_{\alpha,\beta}(m,\mu)$ then 
\begin{equation} 
|\lambda|  = |\alpha| + |\beta| - |\mu| + m 
\end{equation}
The inverse local rule has type signature:
\begin{equation}
\mathfrak{D}_{\alpha,\beta} : (\mathbb{Z}_{\geq 0}, D(\alpha) \cap D(\beta)) \to U(\alpha) \cap U(\beta) 
\end{equation}

An explicit map satisfying the conditions of the local rule was given in section \ref{sec:local-explicit}.

Sometimes the local rule is represented graphically as follows:
\begin{center}
\begin{tikzpicture}[scale=0.8]

\begin{scope}
\path (-2,1) node[shape=circle,inner sep=0.3mm](B0){$\alpha$};
\path (0,0) node[shape=circle,inner sep=0.3mm](C0){$\nu$};
\path (2,1) node[shape=circle,inner sep=0.3mm](D0){$\beta$};
\path (0,2) node[shape=circle,inner sep=0.3mm](C1){$\lambda$};

\draw[thick]  (B0) -- (C0) -- (D0);

\draw[thick] (B0) -- (C1);
\draw[thick] (D0) -- (C1);

\path (0,1) node[shape=circle,inner sep=0.3mm](c0){$m$};

\end{scope}

\end{tikzpicture}
\end{center}
See for example, the diagram in section \ref{growth-bijection}.

\subsubsection{Local rule as higher order function}\label{local-higher}

At slight risk of confusion, we shall also use the term ``local rule'' to refer to a certain \emph{higher order function}, in the sense of functional programming \cite{haskell}.
In functional programming, a higher order function is a function which takes as input a function, and returns as output a different function. 

Referring back to section \ref{sec:bin-order}, for any $\pi \prec \pi'$ such that $\pi'$ is obtained from $\pi$ by adding an inversion at position $i$, the input function for our local rule $\mathfrak{L}_i$ will be a weight-preserving bijection of the form:
\begin{equation} \psi_\pi : (\IP, \ALCD(\pi)) \to \CPP(\pi) \end{equation}
while the output function is a weight preserving bijection of the form:
\begin{equation} \psi_{\pi'} : (\IP, \ALCD(\pi')) \to \CPP(\pi') \end{equation}
That is to say, the local rule $\mathfrak{L}_i$ will have type signature:

\begin{equation} \mathfrak{L}_i : ((\IP, \ALCD(\pi)) \to \CPP(\pi)) \to ((\IP, \ALCD(\pi')) \to \CPP(\pi')) \end{equation}
In other words:
\begin{equation} \mathfrak{L_i}[ \psi_\pi ] = \psi_{\pi'} \end{equation}

Let $\varphi_\pi$ and $\varphi_{\pi'}$ denote the inverse of $\psi_\pi$ and $\psi_{\pi'}$ respectively.
The ``inverse local rule'' $\mathfrak{M}_i$ is the higher order function with type signature:
\begin{equation} \mathfrak{M}_i : (\CPP(\pi) \to (\IP, \ALCD(\pi)) ) \to (\CPP(\pi') \to (\IP, \ALCD(\pi'))) \end{equation}
That is:
\begin{equation} \mathfrak{M}_i[ \varphi_{\pi} ] = \varphi_{\pi'} \end{equation}
Note that $\mathfrak{M}_i$ is only ``inverse'' to $\mathfrak{L}_i$ in the sense that:
\begin{align*}
\mathfrak{L}_i [ \psi_\pi ] \circ \mathfrak{M}_i [ \varphi_{\pi} ] & = \mathfrak{1}_{\CPP(\pi')} \\
\mathfrak{M}_i [ \varphi_{\pi} ] \circ \mathfrak{L}_i [ \psi_\pi ] & = \mathfrak{1}_{(\IP, \ALCD(\pi'))} 
\end{align*}
It is not possible to compose $\mathfrak{L}_i$ and $\mathfrak{M}_i$ directly due to incompatible type signatures.

\subsubsection{Adding and removing boxes}

\index{Cylindric diagram!inside corner}
\index{Cylindric diagram!outside corner}

An \emph{inside corner} of an arbitrarily labelled cylindric diagram is 
an inversion in the profile of the form $(i,i+1)$ or $(T,1)$.
An \emph{outside corner} of an arbitrarily labelled cylindric diagram is 
a \emph{co-inversion} in the profile of the form $(i,i+1)$ or $(T,1)$.
A cylindric diagram always has the same number of inside corners as outside corners.

An inside corner at position $i$ of an arbitrarily labelled cylindric diagram $\mathfrak{d}'$ with profile $\pi'$ can always be removed to obtain an arbitrarily labelled cylindric diagram with profile $\pi$ where $\pi \prec \pi'$ (see Section \ref{sec:bin-order}). 

Keeping track of the label of the box which we have removed, we shall denote this operator by:
\begin{equation} \mathfrak{l}_i : \ALCD(\pi') \to (\mathbb{Z}_{\geq 0}, \ALCD(\pi)) \end{equation}

Conversely, if $\mathfrak{d}$ is an arbitrarily labelled cylindric diagram with profile $\pi$, then given an integer $m$ we may create a new arbitrarily labelled cylindric diagram $\mathfrak{d}'$ with profile $\pi' \succ \pi$ by adding an outside corner at position $i$ and giving it the label $m$.

We shall denote this operator by:
\begin{equation} \mathfrak{r}_i : (\mathbb{Z}_{\geq 0}, \ALCD(\pi)) \to \ALCD(\pi') \end{equation}

\subsubsection{Definition of Local Rule}\label{local-def}

Choose any $(\gamma, \mathfrak{d}') \in (\IP, \ALCD(\pi'))$ and let 
\begin{equation} (m, \mathfrak{d}) = \mathfrak{l}_i[\mathfrak{d}'] \end{equation}
Suppose that:
\begin{equation} \psi_\pi (\gamma, \mathfrak{d}) = \mathfrak{c} = (\mu^0, \mu^1, \ldots, \mu^T) \end{equation}
Let: 
\begin{align*}
\alpha & = \mu^{i-1} \\
\gamma & = \mu^i \\
\beta & = \mu^{i+1} \\
\end{align*}
and let 
\begin{equation} \boxed{\lambda = \mathfrak{U}_{\alpha,\beta} (\gamma, m)} \end{equation}
We define:
\begin{equation} \mathfrak{L}_i[\psi_\pi] (\gamma, \mathfrak{d}') = \mathfrak{c}' = (\mu^0, \ldots, \mu^{i-1}, \lambda, \mu^{i+1}, \ldots \mu^T) \end{equation} 

Note that the horizontal strip condition in the definition of $\mathfrak{U}_{\alpha,\beta}$ ensures that this definition is well-defined.

\subsubsection{Definition of inverse Local Rule}

The inverse local rule is defined similarly. Choose any cylindrical plane partition $\mathfrak{c}' = (\mu^0, \mu^1, \ldots, \mu^T)$ with profile $\pi'$.
Let us define:
\begin{align*}
\alpha & = \mu^{i-1} \\
\lambda & = \mu^i \\
\beta & = \mu^{i+1} \\
\end{align*}
Next let 
\[ \boxed{(m,\nu) = \mathfrak{D}_{\alpha,\beta}(\lambda)} \]
and let $\mathfrak{c}$ be the cylindric plane partition with profile $\pi$
given by $\mathfrak{c} = (\mu^0, \ldots, \mu^{i-1}, \nu, \mu^{i+1}, \ldots \mu^T)$.

If $\varphi_\pi(\mathfrak{c}) = (\gamma, \mathfrak{d})$ then we define 
\[ \mathfrak{M}_i[\varphi_\pi](\mathfrak{c}') = (\gamma, \mathfrak{d}')\]
where $\mathfrak{d}' = \mathfrak{r}_i[(m,\mathfrak{d})]$

\section{The bijection}\label{sec:bijection}
\subsection{Idea of bijection}

In the special case when the arbitrarily labelled cylindric diagram $\mathfrak{d}$ has depth zero (see definition \ref{depth}), the strongly weight-preserving bijection is particularly simple:
\begin{equation} \psi_{\pi}(\gamma, \emptyset) = (\gamma,\gamma, \ldots \gamma) \end{equation}

The idea is to recursively construct bijections, starting from this base case, by repeated application of the local rule (section \ref{local-higher}) and the rotation operator (definition \ref{shift}). 

Each application of the local rule corresponds to an application of the commutation relation (equation \ref{commute}) in the algebraic proof of Borodin's identity while each application of the rotation operator corresponds to an application of the cylindric invariance of the trace.

The recursion is not only over the number of inversions in the profile, but also over the \emph{depth} of the arbitrarily labelled cylindric diagram upon which the bijection is acting. 

Although the precise structure of the recursion is a little complicated to describe in words, it may be neatly encoded in a geometric object called a \emph{cylindric growth diagram}.

\subsection{Cylindric growth diagrams}\label{growth}

\subsubsection{Cylindric diagram}

The idea of a \emph{growth diagram} was first introduced by Fomin \cite{fomin-comm, fomin-growth}. Krattenthaler \cite{krattenthaler} made use of this framework to give a new bijective proof of Stanley's identity (equation \ref{stanley}). In the cylindric case we change the underlying poset, but the essential idea remains the same.

\index{Cylindric poset}

\begin{definition}
For any $n,m \geq 1$, the \emph{cylindric poset} $\mathfrak{G}(n,m)$ is the quotient of $\mathbb{Z}^2$ via the equivalence relation:
\begin{equation}
(x,y) \equiv (x+n,y-m)
\end{equation}
We shall write $v \lhd w$ to indicate that the vertex $w$ covers the vertex $v$ in the cylindric poset.
\end{definition} 
One need to check that this definition is well-defined. The equivalence classes are given by:
\[ [(x,y)] = \{ (x + kn,x-km), k \in \mathbb{Z}\} \]
Reflexivity and transitivity are obvious, we shall just prove antisymmetry. Suppose that $[(x,y)] \leq [(x',y')]$. This implies that there exists a $k' \in \mathbb{Z}$ such that:
\begin{align*}
& x \leq x' + k'n \\
& y \leq y' - k'm
\end{align*}
If we have also that $[(x',y')] \leq [(x,y)]$ then there exists $k \in \mathbb{Z}$ such that:
\begin{align*}
& x' \leq x + kn \\
& y' \leq y - km
\end{align*}
Putting these together we have $x \leq x + (k + k')n$ which implies $(k + k') \geq 0$.
Similarly we have $y \leq y - (k + k')m$ which implies that $(k + k') \leq 0$.
The only way that this is possible is if $k + k' = 0$, or $k = - k'$. 

Now:
\[x + kn \leq x' + k'n + kn = x'\]
which implies that $x \leq x' + k'n$.
But: 
\[x' + k'n \leq x + kn + k'n = x\]
So $x = x' + k'n$. Similarly $y = y + km$. It follows that $[(x,y)] = [(x',y')]$.

\subsubsection{Bijection with cylindric diagrams}

\index{Cylindric growth diagram!face}

\begin{definition}
A \emph{face} in the cylindric poset is a set of four vertices $(u ; v1, v2 ; w) $ satisfying
$u \lhd v1 \lhd w$ and $u \lhd v2 \lhd w$. 
\end{definition}

We say that the face $(u ; v1, v2; w)$ lies \emph{above} the vertex $u$ and \emph{below} the vertex $w$


\index{Cylindric growth diagram!path}
\begin{definition}\label{def:path}
For any binary string $\pi$ containing $n$ zeros and $m$ ones, a \emph{path} in a cylindric growth diagram $\mathfrak{G}(\pi)$ with profile $\pi$ is a sequence of vertices: 
\begin{equation} p = (v_0, v_1, \ldots, v_{n+m-1}, v_{n+m}) \qquad \text{ with } \qquad v_0 = v_{n+m}\end{equation} 
 satisfying $v_{k-1} \lhd v_k$ if $\pi_k = 0$, otherwise $v_{k-1} \rhd v_k$.
\end{definition}

\begin{lemma}
There is a natural bijection between the cylindric diagram with profile $\pi \in \bin(n,m)$, and the subposet of the cylindric poset $\mathfrak{G}(\pi)$ which lies below any given path with profile $\pi$.
\end{lemma}

This bijection maps the \emph{boxes} of the cylindric diagram (section \ref{sec:cylindric-diagram}) to the \emph{faces} of the cylindric poset.

\subsubsection{Cylindric growth diagram}

\index{Cylindric growth diagram!definition}
\begin{definition}\label{growth-def}
A \emph{cylindric growth diagram} with profile $\pi \in \bin(n,m)$ is the subset of the cylindric diagram $\mathfrak{G}(\pi)$ which lies below a path with profile $\pi$ (referred to as the upper boundary), whose vertices are labelled by integer partitions, and whose faces are labelled by non-negative integers, in such a way that the following three conditions are satisfied:

\begin{enumerate}
\item If $v \lhd w$, and if $\lambda$ is the integer partition labelling the vertex $v$ and $\mu$ is the integer partition labelling the vertex $w$, then $\lambda / \mu$ is a horizontal strip.
\item All but finitely many of the vertices are labelled with the same integer partition $\gamma$.
\item If $( u; v1, v2; w)$ is a face with label $m$, and if the labels of $u$, $v1$, $v2$ and $w$ are $\mu$, $\alpha$, $\beta$ and $\lambda$ respectively then:
\[ \mathfrak{U}_{\alpha,\beta}(m,\mu) = \lambda \]
\end{enumerate}
\end{definition}

The invertibility of the local rule implies that condition $(3)$ could be equivalently formulated as follows:
\begin{lemma}
If $( u; v1, v2; w)$ is a face of a cylindric diagram $\mathfrak{G}(\pi)$ with label $m$, and if the labels of $u$, $v1$, $v2$ and $w$ are $\mu$, $\alpha$, $\beta$ and $\lambda$ respectively then:
\[ \mathfrak{D}_{\alpha,\beta}(\lambda) = (m,\mu) \]
\end{lemma}

Observe that the weight condition of the local rule \ref{sec:local-weight} we have:
\begin{lemma}\label{local-help2}
$\mathfrak{D}_{\gamma,\gamma}(\gamma) = (0,\gamma)$
\end{lemma}

It follows that below a certain point in a cylindric growth diagram, all the vertices will have the same label, and all the faces will be labelled by zero.

\subsection{The bijection}\label{growth-bijection}

\index{Cylindric growth diagram!example}

Here is an example of a growth diagram with $\pi = 00101$. Note that the profile is read from right to left. Zero corresponds to a $SW$ step and one corresponds to a $NW$ step:

\begin{center}
\begin{tikzpicture}[scale=1]

\begin{scope}
\path (-4,2) node[shape=circle,inner sep=0.3mm](A0){$(3,2)$};
\path (-2,1) node[shape=circle,inner sep=0.3mm](B0){$(3,2)$};
\path (0,0) node[shape=circle,inner sep=0.3mm](C0){$(3,2)$};
\path (2,1) node[shape=circle,inner sep=0.3mm](D0){$(3,2)$};
\path (4,2) node[shape=circle,inner sep=0.3mm](E0){$(3,2)$};
\path (6,3) node[shape=circle,inner sep=0.3mm](F0){$(3,2)$};

\path (-4,4) node[shape=circle,inner sep=0.3mm](A1){$(4,3,2)$};
\path (-2,3) node[shape=circle,inner sep=0.3mm](B1){$(3,2,2)$};
\path (0,2) node[shape=circle,inner sep=0.3mm](C1){$(3,2,1)$};
\path (2,3) node[shape=circle,inner sep=0.3mm](D1){$(3,2,1)$};
\path (4,4) node[shape=circle,inner sep=0.3mm](E1){$(3,2,1)$};
\path (6,5) node[shape=circle,inner sep=0.3mm](F1){$(4,3,2)$};

\path (-4,6) node[shape=circle,inner sep=0.3mm](A2){$(6,4,3,2)$};
\path (-2,5) node[shape=circle,inner sep=0.3mm](B2){$(4,3,2)$};
\path (0,4) node[shape=circle,inner sep=0.3mm](C2){$(3,2,2)$};
\path (2,5) node[shape=circle,inner sep=0.3mm](D2){$(3,2,2)$};
\path (4,6) node[shape=circle,inner sep=0.3mm](E2){$(5,3,2)$};
\path (6,7) node[shape=circle,inner sep=0.3mm](F2){$(6,4,3,2)$};

\path (0,6) node[shape=circle,inner sep=0.3mm](C3){$(4,3,2,1)$};

\draw[thick] (A0) -- (B0) -- (C0) -- (D0) -- (E0) -- (F0);
\draw[thick] (A1) -- (B1) -- (C1) -- (D1) -- (E1) -- (F1);
\draw[thick] (A2) -- (B2) -- (C2) -- (D2) -- (E2) -- (F2);

\draw[thick] (A1) -- (B2) -- (C3);
\draw[thick] (A0) -- (B1) -- (C2);
\draw[thick] (B0) -- (C1);
\draw[thick] (D0) -- (C1);
\draw[thick] (E0) -- (D1) -- (C2);
\draw[thick] (F0) -- (E1) -- (D2) -- (C3);
\draw[thick] (F1) -- (E2);

\path (-4,3) node[shape=circle,inner sep=0.3mm](a0){$1$};
\path (-2,2) node[shape=circle,inner sep=0.3mm](b0){$1$};
\path (0,1) node[shape=circle,inner sep=0.3mm](c0){$1$};
\path (2,2) node[shape=circle,inner sep=0.3mm](dD0){$0$};
\path (4,3) node[shape=circle,inner sep=0.3mm](e0){$0$};
\path (6,4) node[shape=circle,inner sep=0.3mm](f0){$1$};

\path (-4,5) node[shape=circle,inner sep=0.3mm](a1){$5$};
\path (-2,4) node[shape=circle,inner sep=0.3mm](b1){$0$};
\path (0,3) node[shape=circle,inner sep=0.3mm](c1){$0$};
\path (2,4) node[shape=circle,inner sep=0.3mm](d1){$0$};
\path (4,5) node[shape=circle,inner sep=0.3mm](e1){$0$};
\path (6,6) node[shape=circle,inner sep=0.3mm](f1){$5$};

\path (0,5) node[shape=circle,inner sep=0.3mm](v2){$1$};
\end{scope}

\end{tikzpicture}
\end{center}

Note that we have truncated the diagram below the lower boundary, where all vertices have the same label.

\begin{lemma}
For any path $p$ in a cylindric growth diagram $\mathfrak{G}(\pi)$ with profile $\pi$, the sequence of profiles associated
to the vertices of the path form a cylindric plane partition.
\end{lemma}

In particular, the sequence of partitions labelling the vertices along the upper boundary form an element of $\CPP(\pi)$.

\begin{center}
\begin{tikzpicture}[scale=1]

\begin{scope}
\path (-4,2) node[shape=circle,inner sep=0.3mm](A0){};
\path (-2,1) node[shape=circle,inner sep=0.3mm](B0){};
\path (0,0) node[shape=circle,inner sep=0.3mm](C0){};
\path (2,1) node[shape=circle,inner sep=0.3mm](D0){};
\path (4,2) node[shape=circle,inner sep=0.3mm](E0){};
\path (6,3) node[shape=circle,inner sep=0.3mm](F0){};

\path (-4,4) node[shape=circle,inner sep=0.3mm](A1){};
\path (-2,3) node[shape=circle,inner sep=0.3mm](B1){};
\path (0,2) node[shape=circle,inner sep=0.3mm](C1){};
\path (2,3) node[shape=circle,inner sep=0.3mm](D1){};
\path (4,4) node[shape=circle,inner sep=0.3mm](E1){};
\path (6,5) node[shape=circle,inner sep=0.3mm](F1){};

\path (-4,6) node[shape=circle,inner sep=0.3mm](A2){$(6,4,3,2)$};
\path (-2,5) node[shape=circle,inner sep=0.3mm](B2){$(4,3,2)$};
\path (0,4) node[shape=circle,inner sep=0.3mm](C2){};
\path (2,5) node[shape=circle,inner sep=0.3mm](D2){$(3,2,2)$};
\path (4,6) node[shape=circle,inner sep=0.3mm](E2){$(5,3,2)$};
\path (6,7) node[shape=circle,inner sep=0.3mm](F2){$(6,4,3,2)$};

\path (0,6) node[shape=circle,inner sep=0.3mm](C3){$(4,3,2,1)$};

\draw[thick] (A0) -- (B0) -- (C0) -- (D0) -- (E0) -- (F0);
\draw[thick] (A1) -- (B1) -- (C1) -- (D1) -- (E1) -- (F1);
\draw[thick] (A2) -- (B2) -- (C2) -- (D2) -- (E2) -- (F2);

\draw[thick] (A1) -- (B2) -- (C3);
\draw[thick] (A0) -- (B1) -- (C2);
\draw[thick] (B0) -- (C1);
\draw[thick] (D0) -- (C1);
\draw[thick] (E0) -- (D1) -- (C2);
\draw[thick] (F0) -- (E1) -- (D2) -- (C3);
\draw[thick] (F1) -- (E2);

\path (-4,3) node[shape=circle,inner sep=0.3mm](a0){};
\path (-2,2) node[shape=circle,inner sep=0.3mm](b0){};
\path (0,1) node[shape=circle,inner sep=0.3mm](c0){};
\path (2,2) node[shape=circle,inner sep=0.3mm](dD0){};
\path (4,3) node[shape=circle,inner sep=0.3mm](e0){};
\path (6,4) node[shape=circle,inner sep=0.3mm](f0){};

\path (-4,5) node[shape=circle,inner sep=0.3mm](a1){};
\path (-2,4) node[shape=circle,inner sep=0.3mm](b1){};
\path (0,3) node[shape=circle,inner sep=0.3mm](c1){};
\path (2,4) node[shape=circle,inner sep=0.3mm](d1){};
\path (4,5) node[shape=circle,inner sep=0.3mm](e1){};
\path (6,6) node[shape=circle,inner sep=0.3mm](f1){};

\path (0,5) node[shape=circle,inner sep=0.3mm](v2){};
\end{scope}

\end{tikzpicture}
\end{center}

Note that this is just a rotation of the example cylindric plane partition given in section \ref{sec:cylindric-diagram}.

\begin{lemma}
If the labels of the vertices of a cylindric growth diagram $\mathfrak{G}(\pi)$ with profile $\pi$ are forgotten, then we obtain an arbitrarily labelled cylindric diagram. 
\end{lemma}

\begin{center}
\begin{tikzpicture}[scale=0.8]

\begin{scope}
\path (-4,2) node[shape=circle,inner sep=0.3mm](A0){};
\path (-2,1) node[shape=circle,inner sep=0.3mm](B0){};
\path (0,0) node[shape=circle,inner sep=0.3mm](C0){};
\path (2,1) node[shape=circle,inner sep=0.3mm](D0){};
\path (4,2) node[shape=circle,inner sep=0.3mm](E0){};
\path (6,3) node[shape=circle,inner sep=0.3mm](F0){};

\path (-4,4) node[shape=circle,inner sep=0.3mm](A1){};
\path (-2,3) node[shape=circle,inner sep=0.3mm](B1){};
\path (0,2) node[shape=circle,inner sep=0.3mm](C1){};
\path (2,3) node[shape=circle,inner sep=0.3mm](D1){};
\path (4,4) node[shape=circle,inner sep=0.3mm](E1){};
\path (6,5) node[shape=circle,inner sep=0.3mm](F1){};

\path (-4,6) node[shape=circle,inner sep=0.3mm](A2){};
\path (-2,5) node[shape=circle,inner sep=0.3mm](B2){};
\path (0,4) node[shape=circle,inner sep=0.3mm](C2){};
\path (2,5) node[shape=circle,inner sep=0.3mm](D2){};
\path (4,6) node[shape=circle,inner sep=0.3mm](E2){};
\path (6,7) node[shape=circle,inner sep=0.3mm](F2){};

\path (0,6) node[shape=circle,inner sep=0.3mm](C3){};

\draw[thick] (A0) -- (B0) -- (C0) -- (D0) -- (E0) -- (F0);
\draw[thick] (A1) -- (B1) -- (C1) -- (D1) -- (E1) -- (F1);
\draw[thick] (A2) -- (B2) -- (C2) -- (D2) -- (E2) -- (F2);

\draw[thick] (A1) -- (B2) -- (C3);
\draw[thick] (A0) -- (B1) -- (C2);
\draw[thick] (B0) -- (C1);
\draw[thick] (D0) -- (C1);
\draw[thick] (E0) -- (D1) -- (C2);
\draw[thick] (F0) -- (E1) -- (D2) -- (C3);
\draw[thick] (F1) -- (E2);

\path (-4,3) node[shape=circle,inner sep=0.3mm](a0){$1$};
\path (-2,2) node[shape=circle,inner sep=0.3mm](b0){$1$};
\path (0,1) node[shape=circle,inner sep=0.3mm](c0){$1$};
\path (2,2) node[shape=circle,inner sep=0.3mm](dD0){$0$};
\path (4,3) node[shape=circle,inner sep=0.3mm](e0){$0$};
\path (6,4) node[shape=circle,inner sep=0.3mm](f0){$1$};

\path (-4,5) node[shape=circle,inner sep=0.3mm](a1){$5$};
\path (-2,4) node[shape=circle,inner sep=0.3mm](b1){$0$};
\path (0,3) node[shape=circle,inner sep=0.3mm](c1){$0$};
\path (2,4) node[shape=circle,inner sep=0.3mm](d1){$0$};
\path (4,5) node[shape=circle,inner sep=0.3mm](e1){$0$};
\path (6,6) node[shape=circle,inner sep=0.3mm](f1){$5$};

\path (0,5) node[shape=circle,inner sep=0.3mm](v2){$1$};
\end{scope}

\end{tikzpicture}
\end{center}

Note that the arbitrarily labelled cylindric diagram is a rotation of the example given in section \ref{alcd}.

A growth diagram should be thought of as an object which interpolates between the LHS and the RHS of the bijection which we wish to establish.

\begin{proposition}
To every cylindric plane partition, there is a uniquely associated cylindric growth diagram.
\begin{proof}
Once the labels on the upper boundary have been specified, property (3) of definition \ref{growth-def} ensures that there is a unique way in which to label the remaining faces and vertices.
\end{proof}
\end{proposition}

\begin{proposition}
To every pair $(\gamma, \mathfrak{d})$ where $\gamma$ is an integer partition and $\mathfrak{d}$ is an arbitrarily labelled cylindric diagram, there exists a uniquely defined cylindric growth diagram.
\begin{proof}
Let $d$ denote the depth of $\mathfrak{d}$. To every face with depth greater than $d$ assign the label $0$. To every vertex lying below a face with depth greater than $d$, assign the label $\gamma$. Property $(3)$ of definition \ref{growth-def} ensures that there is a unique way to label the remaining faces and vertices. Lemma \ref{local-help2} ensures that the resulting growth diagram is well-defined.
\end{proof}
\end{proposition}

\subsection{Remarks}

The following lemma guarantees that, when there are multiple inversions in the profile string, the order in which the local rules are applied is of no importance.

\begin{lemma}\label{localcommute}
If $\pi$ has inversions at both positions $i$ and $j$ then:
\begin{equation} \mathfrak{L}_i \circ \mathfrak{L}_j [ \psi_\pi ] = \mathfrak{L}_j \circ \mathfrak{L}_i [ \psi_\pi ] \end{equation}
\begin{proof}
Without loss of generality we may assume that $j > i$. If $\pi$ has inversions at both positions $i$ and $j$ then $\pi_i = 0 = \pi_j$ and $\pi_{i+1} = 1 = \pi_{j+1}$
thus $j-i \geq 2$. Application of the local rule $\mathfrak{L}_i$ does not effect the $(j-1)$th diagonal. Similarly, application
of the local rule $\mathfrak{L}_j$ does not effect the $(i+1)$th diagonal, thus the two operators commute. 
\end{proof}
\end{lemma}

Given a path $p$ with profile $\pi_{\max}$ it is not possible to apply the local rule $\mathfrak{M}_i$ for any $i$. It is however possible to rotate the cylinder, and thus obtain a new path $\sigma(p)$ with profile $\pi_{\min}$ (see Section \ref{shift}).

If $\mathfrak{d}$ is an arbitrarily labelled cylindric diagram with profile $\pi$ and depth $d > 0$ then in order to construct the bijection $\psi_\pi(\gamma, \mathfrak{d})$ the cylindric shift operator will have to be applied $d-1$ times. 

It remains to prove that our bijection is weight preserving.

\section{The weight}\label{sec:weight}
\subsection{Alternative definition of weight}\label{sec:diag}

\index{Cylindric diagram!cohook}

Let us define the cohook of a box $b$ in an arbitrarily labelled cylindric diagram $\mathfrak{d}$ to be:
\begin{equation}
 \cohook(b) = \{ b' | b \in \hook(b') \} 
\end{equation}

\[
 \tableau{
\missingcell & \missingcell & \missingcell & \missingcell & \missingcell & \missingcell & \missingcell & \missingcell 1  \\
\missingcell & \missingcell & \missingcell & \missingcell & \missingcell & \missingcell & \missingcell 1&0  \\
\missingcell & \missingcell & \missingcell & \missingcell & \missingcell & \missingcell &0& \\
\missingcell & \missingcell & \missingcell & \missingcell & \missingcell & \missingcell 1&0& \\
\missingcell & \missingcell & \missingcell & \missingcell & \missingcell 1&0& \thickcell&  \\
\missingcell & \missingcell & \missingcell & \missingcell &0&&  & \thickcell\\
\missingcell & \missingcell & \missingcell & \missingcell 1&0&&  & \\
\missingcell & \missingcell & \missingcell 1&0& \thickcell &  && \\
\missingcell & \missingcell &0& \redcell & \yellowcell & \thickyellowcell& \yellowcell& \yellowcell \\
\missingcell & \missingcell 1&0& \yellowcell &  &&\thickcell& \\
\missingcell1 &  0 &  \thickcell & \yellowcell && && \thickcell\\
0&&& \thickyellowcell &  &&& \\
0&&& \yellowcell & \thickcell &&& \\
}
\]

Although the number of boxes in the cohook of a given box is always infinite, since only a finite number of the labels are zero it still makes sense to define the weight of a cohook:
\begin{equation}
 |\cohook(b)|_\mathfrak{d} = \sum_{b' \in \cohook(b)} \lab(b') 
\end{equation}

\index{Cylindric diagram!diag}

Let $\diag(k)$ denote the set of all boxes on the $k$-th diagonal.
Furthermore let us define the weight of the diagonal of an arbitrarily labelled cylindric diagram $\mathfrak{d}$ to be:
\begin{equation} 
|\diag(k)|_\mathfrak{d} = \sum_{b \in \diag(k)} |\cohook(b)|_\mathfrak{d} 
\end{equation}

Note that if $b$ is taken to be the box of $\diag(k)$ lying furthest to the ``north-west'' then $|\diag(k)|_\mathfrak{d}$ is none other than
the sum of all the labels of boxes lying ``south-east'' of $b$.

\index{Cylindric diagram!weight}

With these definitions we may give an alternative definition of the weight of $\mathfrak{d}$:
\begin{align}
|\mathfrak{d} |
& = \sum_{b' \in \lambda} \lab(b') \hook(b)  \\
& = \sum_{b' \in \lambda} \lab(b') \sum_{b \in \hook(b')} 1  \\
& =  \sum_{b \in \lambda} \, \sum_{b' \in \cohook(b)} \lab(b') \\
& = \sum_{b \in \lambda} |\cohook(b)|_\mathfrak{d} \\
& = \sum_{k=1}^T |\diag(k)|_\mathfrak{d} \\
\end{align}
The \emph{refined weight} of a cylindric plane partition is given by:
\begin{equation}  
z_1^{|\diag(1)|_\mathfrak{d}} z_2^{|\diag(2)|_\mathfrak{d}} \cdots z_T^{|\diag(T)|_\mathfrak{d}} 
\end{equation}

\index{Borodin identity!Bijection!strongly weight preserving}

\begin{definition}\label{strong-weight}
We shall say that a bijection:
\[ \psi_\pi : (\IP, \ALCD(\pi)) \to \CPP(\pi) \] is \emph{strongly weight preserving}
if whenever:
\[ \psi_\pi(\nu, \mathfrak{d}) = \mathfrak{c} = (\mu^0, \mu^1, \ldots, \mu^T) \] 
we have for all $1 \leq k \leq T$ that:
\begin{equation} |\mu^k| = |\gamma| + |\diag(k)|_\mathfrak{d} \end{equation}
\end{definition}

\begin{lemma}
Strongly weight preserving implies weight preserving.
\end{lemma}

\subsection{Proof that the bijection is strongly weight preserving}
In this section we prove that our bijection is strongly weight preserving.

\begin{proposition}\label{weight-preserving1}
If $(m, \mathfrak{d}) = \mathfrak{l}_i[\mathfrak{d}']$ then:
\begin{equation} 
|\diag(i)|_{\mathfrak{d}'} = m + |\diag(i-1)|_{\mathfrak{d}} + |\diag(i+1)|_{\mathfrak{d}} - |\diag(i)|_{\mathfrak{d}} 
\end{equation}
\begin{proof}

Let $b$ denote the box of $\mathfrak{d}'$ with cylindric inversion coordinates $(i,i+1,0)$.
The sum of all the labels in boxes in the same cylindric column as $b$ is given by: 
\begin{equation} |\diag(i+1)|_{\mathfrak{d}} - |\diag(i)|_{\mathfrak{d}} \end{equation}
while the sum of all the labels in boxes lying in the same cylindric row as $b$ is given by: 
\begin{equation} |\diag(i-1)|_{\mathfrak{d}} - |\diag(i)|_{\mathfrak{d}} \end{equation}
In other words:
\begin{align}
\cohook(i,i+1,0)_{\mathfrak{d}'} & = m + |\diag(i-1)|_{\mathfrak{d}}  + |\diag(i+1)|_{\mathfrak{d}} - 2 |\diag(i)|_{\mathfrak{d}}
\end{align}

Now:
\begin{align}
|\diag(i)|_{\mathfrak{d}'} & = |\diag(i)|_{\mathfrak{d}} + \cohook(i,i+1,0)_{\mathfrak{d}'} \\
& = m + |\diag(i-1)|_{\mathfrak{d}} + |\diag(i+1)|_{\mathfrak{d}} - |\diag(i)|_{\mathfrak{d}}
\end{align}
\end{proof}
\end{proposition}

\begin{proposition}\label{weight-preserving}
If $\psi_\pi$ is strongly weight preserving, then so is $\mathfrak{L}_i(\psi_\pi)$.
\begin{proof}

Let: 
\begin{equation} (m,\mathfrak{d}) = \mathfrak{l}_i[\mathfrak{d}'] 
\end{equation}
and let: 
\begin{equation}
\psi_\pi(\gamma, \mathfrak{d}) = \mathfrak{c} = (\mu^0, \ldots, \alpha, \nu, \beta, \ldots \mu^T)
\end{equation} 

Then: 
\begin{equation}
\mathfrak{L}_i[\psi_\pi](\gamma, \mathfrak{d}') = \mathfrak{c}' = (\mu^0, \ldots, \alpha, \lambda, \beta, \ldots \mu^T)
\end{equation}
where:
\begin{equation} \lambda = \mathfrak{U}_{\alpha,\beta}(m,\nu) \end{equation}

The weight condition for the local rule assures us that:
\begin{equation} |\lambda| = m + |\alpha| + |\beta| - |\nu| \end{equation}

For all $j \neq i$ we have 
\begin{equation}|\diag(j)|_\mathfrak{d} = |\diag(j)|_{\mathfrak{d}'} = |\mu^j| - |\gamma| \end{equation}
We must show that:
\begin{equation} |\diag(i)|_{\mathfrak{d}'} = |\lambda| - |\gamma| \end{equation}

Now, by proposition \ref{weight-preserving1} and the assumption that $\psi_\pi$ is strongly weight-preserving, we have:
\begin{align}
|\diag(i)|_{\mathfrak{d}'} & = m +  |\diag(i-1)|_{\mathfrak{d}} + |\diag(i+1)|_{\mathfrak{d}} - |\diag(i)|_{\mathfrak{d}}\\
& = m + (|\alpha| - |\gamma|) + (|\beta| - |\gamma|) - (|\nu| - |\gamma|) \\
& = |\lambda| - |\gamma|
\end{align}
The result follows.
\end{proof}
\end{proposition}

\section{Conclusion}

We have made use of Fomin's growth diagram framework \cite{fomin-comm, fomin-growth} to give a bijective proof of a refined version of Borodin's identity (equation \ref{eq:borodin}). 
Our proof generalizes known proofs for Stanley's identity (equation \ref{stanley}) and MacMahon's identity (equation \ref{macmahon}). 

Due to the fact that there are \emph{two} equally natural extensions of Fomin's local rule to the horizontal strip case (see section \ref{sec:local}) we have actually given two distinct bijective proofs of Borodin's identity, one corresponding to the RSK correspondence, and the other corresponding to the Burge correspondence:
\begin{equation}
\psi_{R},\psi_{B} : (\IP, \ALCD(\pi)) \to \CPP(\pi) 
\end{equation}

These two maps are closely related by the Sch\"utzenberger involution \cite{vanLeeuwen}. Observe that we have:
\begin{equation}\label{crystal}
\psi_{R} \circ \psi_{B}^{-1} : \CPP(\pi) \to \CPP(\pi) 
\end{equation}
What can we say about this map?

Tingley \cite{tingley} showed that cylindric plane partitions can be understood as crystal bases for $\widehat{sl}_n$.
Does the map in equation \ref{crystal} have any interpretation at the level of representation theory?

It is known \cite{soliton} that the tensor product of the highest weight representation of $\widehat{sl}_n$ associated to the profile $\pi$ by the evaluation representation of $\widehat{sl}_n$ associated to the natural representation of $sl_n$ decomposes as a direct sum of highest weight representations, where those highest weight representations occurring in the sum are precisely those associated to the profile $\pi'$ where $\pi \prec \pi'$ in the partial order described in section \ref{sec:bin-order}.
Is it possible to make use of the local rule formalism to better understand this decomposition on the crystal basis level?

\chapter{Macdonald polynomial analog}
\section{$(q,t)$-Borodin identity}\label{sec:qtborodin}
\subsection{Statement}

\begin{theorem}\label{qt-borodin}
For any binary string $\pi$ we have:
\begin{equation} \label{maintheorem}
\sum_{\mathfrak{c} \in CPP(\pi)} W_\mathfrak{c}(q,t) z^{|\mathfrak{c}|} = 
\prod_{n \geq 0} \left ( \frac{1}{1 - z^{(n+1)T}}  
\prod_{\substack{i < j \\ \pi_i > \pi_j}} \frac{(tz^{j - i + nT};q)_\infty}{(z^{j - i + nT};q)_\infty} 
\prod_{\substack{i > j \\ \pi_i > \pi_j}} \frac{(tz^{j-i + (n+1)T};q)_\infty}{(z^{j-i + (n+1)T};q)_\infty} \right )
\end{equation}
\end{theorem}

where if $\mathfrak{c} = (\mu^0,\mu^1, \ldots \mu^T)$ then the weight function is given by:

\begin{equation}  \label{weight}
W_{\mathfrak{c}}(q,t) = \prod_{\substack{k =1 \\ \pi_k = 1}}^T \varphi_{\mu^k / \mu^{k-1}}(q,t) \prod_{\substack{k = 1\\ \pi_k = 0}}^T \psi_{\mu^{k-1} / \mu^k}(q,t) 
\end{equation}
See equations \ref{pieri-coeff} and \ref{pieri-coeff2} for the definition of the Macdonald Pieri coefficients.

Note that when $q=t$ the weight function reduces to one, and equation \ref{qt-borodin} reduces to equation \ref{borodin}. 
The proof of the Hall--Littlewood case of this identity is due to Corteel and Savelief, Cyrille and Vuleti{\'c} \cite{vuletic} while the proof of the Macdonald version of the reverse plane partition case is due to Okada \cite{okada}.

\index{Borodin identity!Macdonald!refined}
As in the Schur case, the nature of the proof is such that the identity remains true if on the left hand side we replace:
\[z^{\mathfrak{|c|}} \mapsto z_0^{|\mu_0|} z_1^{|\mu_1|} \cdots z_{T-1}^{|\mu_{T-1}|}\]
while on the right hand side we replace:
\begin{align*} 
z^{nT} & \mapsto z_0^{n} z_1^{n} \cdots z_{T-1}^{n} \\
z^{j - i + nT} & \mapsto z_0^{n} z_1^n \cdots z_i^{n} z_{i+1}^{n+1} \cdots z_{j}^{n+1} z_{j+1}^n \cdots z_{T-1}^n \quad \quad \quad \quad \text{ when } i < j \\
z^{j-i +(n+1)T} & \mapsto z_0^{n+1} z_1^{n+1} \cdots z_j^{n+1} z_{j+1}^n + \cdots z_i^{n} z_{i+1}^{n+1} \cdots z_{T-1}^{n+1} \quad \text{ when } i > j
\end{align*}

\subsection{Proof}

Our proof is almost identical to that given in section \ref{sec:algebraic}. The only difference is that instead of using the Hall inner product $\langle -,-\rangle$ and the operators $\Omega[Xu]$ and $\Omega^*[Xv]$ from section \ref{sec:schur}, we use the Macdonald inner product $\langle -,-\rangle_{q,t}$ and the operators $\Omega_{q,t}[Xu]$ and $\Omega_{q,t}^*[Xv]$ from section \ref{sec:macdonald}.

\subsubsection{Notation}

Using the same ``degree'' operator as in section \ref{sec:formulation}:
\begin{equation} 
D_z P_\lambda[X] = z^{|\lambda|} P_\lambda[X] 
\end{equation}

We have:
\begin{lemma} \label{degree-mac}
\begin{align}
D_z \, \Omega_{q,t}[Xu] & = \Omega_{q,t}[Xuz] \, D_z \\
D_z \, \Omega^*_{q,t}[Xu] & = \Omega^*_{q,t}[Xuz^{-1}] \, D_z 
\end{align}
\begin{proof}
This fact follows immediately from the action of $\Omega_{q,t}[Xu]$ and $\Omega_{q,t}[Xv]$ on Macdonald polynomials 
(equations \ref{pieri} and \ref{pieri2}).
\end{proof}
\end{lemma}

For notational convenience we shall define:
\begin{align}
H^0(z) & = \Omega_{q,t}[Xz] \\
H^1(z) & = \Omega^*_{q,t}[Xz]
\end{align}
These shall play an analogous role to the operators $G^0(z)$ and $G^1(z)$ defined in section \ref{sec:formulation}.

\subsubsection{Algebraic interpretation of left hand side}

\begin{lemma}  \label{lhs-mac}
The left hand side of the refined version of equation \ref{maintheorem} may be expressed in the form:
\begin{equation} 
\lhs(\pi) = \sum_\mu \langle Q_\mu \,|\, H^{\pi_0}(u_0) H^{\pi_1}(u_1)  \cdots H^{\pi_T}(u_T) D_w \, P_\mu  \rangle_{q,t} 
\end{equation}
where:
\begin{align} \label{specialization-weight}
w & = z_0 z_1 \cdots z_{T-1} \\
u_k & = \begin{cases}
z_0 z_1 \cdots z_{k-1} & \text{ if $\pi_k = 1$} \\
z_0^{-1} z_1^{-1} \cdots z_{k-1}^{-1} & \text{ if $\pi_k = 0$}
\end{cases}\label{specialization2-weight}
\end{align}

\begin{proof}
From the ``interlacing sequence'' definition of a cylindric plane partitions \ref{cpp} it is clear that a cylindric plane partition is constructed
by successively adding and removing horizontal strips. 
As in the Schur case, the degree operator $D_z$ is used to keep track of the number of cubes in the resulting cylindric plane partition.

The new feature in the Macdonald case is presence of the $(q,t)$-Pieri coefficients in the definition of the weight function (equation \ref{weight}). This
comes directly from the action of the operators $\Omega_{q,t}(Xu)$ and $\Omega^*_{q,t}(Xu)$ on Macdonald polynomials given in equations \ref{pieri} and \ref{pieri2}.

Using the fact that the Macdonald $P$-functions are orthogonal to the Macdonald $Q$-functions with respect to the Macdonald inner product 
we may write:
\begin{equation}
\lhs(\pi) = \sum_\mu \langle Q_\mu \,|\, D_{z_0} \, H^{\pi_0}(1) \, D_{z_1} \, H^{\pi_1}(1)   \cdots D_{z_{T-1}} \, H^{\pi_T}(1) \, P_\mu  \rangle_{q,t} 
\end{equation}

It remains to commute all the shift operators to the right hand side using Lemma \ref{degree-mac}.
\end{proof}
\end{lemma}

Note that the above expression can also be understood as the \emph{trace} of the operator:
\[H^{\pi_0}(u_0) H^{\pi_1}(u_1)  \cdots H^{\pi_T}(u_T) D_w\]
acting on symmetric functions over $\mathbb{Q}(q,t)$.

\subsubsection{Some lemmas}

Let us define:
\begin{definition}\label{N-def}
\begin{equation} 
N_\pi(m)  = \sum_\mu \langle Q_\mu \,|\, \prod_{\substack{k =1\\ \pi_k = 0}}^T \Omega_{q,t} [Xu_kw^m] 
\prod_{\substack{k =1 \\ \pi_k = 1}}^T \Omega_{q,t}^*[Xu_k] \, D_w \, P_\mu \rangle_{q,t} 
\end{equation}
\end{definition}

\begin{lemma} \label{cylindric-shift-mac}
\[ N_\pi(m) =  \prod_{\substack{(i,j) \\ \pi_i \neq \pi_j}} \frac{(tu_i u_j w^{m+1};q)_\infty}{(u_i u_j w^{m+1};q)_\infty} \, N_\pi(m+1) \]
\begin{proof}
This is the same calculation as in Lemma \ref{cylindric-shift}, only using the Macdonald commutation relation (lemma \ref{commutation-mac}) rather than the Schur commutation relation (proposition \ref{commute}).
 
Since the $\{P_\lambda\}$ are orthogonal to the $\{Q_\lambda\}$ with respect to the Macdonald inner product, and since:
\[\tr(AB) = \tr(BA) \]
We may write:
\begin{align}
N_\pi(m) & = \sum_{\mu, \lambda} 
\langle Q_\mu \,|\, \prod_{\substack{k =1\\ \pi_k = 0}}^T \Omega_{q,t}[Xu_kw^m] \, P_\lambda \rangle_{q,t} 
\langle Q_\lambda \,|\,\prod_{\substack{k =1 \\ \pi_k = 1}}^T \Omega_{q,t}^*[Xu_k] \, D_w \, P_\mu \rangle_{q,t} \\
& = \sum_{\mu, \lambda} 
\langle Q_\lambda \,|\,\prod_{\substack{k =1 \\ \pi_k = 1}}^T \Omega_{q,t}^*[Xu_k] \, D_w \, P_\mu \rangle_{q,t} 
\langle Q_\mu \,|\, \prod_{\substack{k =1\\ \pi_k = 0}}^T \Omega_{q,t}[Xu_kw^m] \, P_\lambda \rangle_{q,t} \\
& = \sum_{\lambda} 
\langle Q_\lambda \,|\,\prod_{\substack{k =1 \\ \pi_k = 1}}^T \Omega_{q,t}^*[Xu_k] \, 
D_w \,\prod_{\substack{k =1\\ \pi_k = 0}}^T \Omega_{q,t}[Xu_kw^m] \, P_\lambda \rangle_{q,t} 
\end{align}
Next applying the commutation relations of Lemma \ref{degree-mac} and Lemma \ref{commutation-mac} we have:
\begin{align}
N_\pi(m) & = \sum_{\lambda} 
\langle Q_\lambda \,|\,\prod_{\substack{k =1 \\ \pi_k = 1}}^T \Omega_{q,t}^*[Xu_k] \, 
D_w \,\prod_{\substack{k =1\\ \pi_k = 0}}^T \Omega_{q,t}[Xu_kw^m] \, P_\lambda \rangle_{q,t}  \\ 
& = \sum_{\lambda} 
\langle Q_\lambda \,|\,\prod_{\substack{k =1 \\ \pi_k = 1}}^T \Omega_{q,t}^*[Xu_k] \, 
\,\prod_{\substack{k =1\\ \pi_k = 0}}^T \Omega_{q,t}[Xu_kw^{m+1}] \, D_w \, P_\lambda \rangle_{q,t} \\
& = \prod_{\substack{(i,j) \\ \pi_i \neq \pi_j}} \frac{(tu_i u_j w^{m+1};q)_\infty}{(u_i u_j w^{m+1};q)_\infty} \sum_{\lambda} 
\langle Q_\lambda \,| \,\prod_{\substack{k =1\\ \pi_k = 0}}^T \Omega_{q,t}[Xu_kw^{m+1}] \, 
\prod_{\substack{k =1 \\ \pi_k = 1}}^T \Omega_{q,t}^*[Xu_k] \, D_w \, P_\lambda \rangle_{q,t} \\
& =  \prod_{\substack{(i,j) \\ \pi_i \neq \pi_j}} \frac{(tu_i u_j w^{m+1};q)_\infty}{(u_i u_j w^{m+1};q)_\infty} \, N_\pi(m+1)
\end{align}
\end{proof}
\end{lemma}

As in the Schur case, we have in the limit:
\begin{lemma} \label{infinity-mac}
\begin{equation} 
N_\pi(\infty)  = \prod_{n \geq 1} \frac{1}{1-w^n} 
\end{equation}
\begin{proof}
In order for this limit to even make sense, we must have $|z_i| < 1$ for all $i$, in which case:
\[ \lim_{m \to \infty} \Omega_{q,t}[X u_k \omega^m] = 1 \]

Since $\Omega_{q,t}^*[X u_k]$ is a degree lowering operator, it follows that:
\begin{align*}
 \lim_{m \to \infty} N_\pi(m) 
& =  \sum_\mu \langle Q_\mu \,|\,  \prod_{\substack{k =1 \\ \pi_k = 1}}^T \Omega_{q,t}^*[Xu_k] \, D_w \, P_\mu \rangle_{q,t} \\
& = \sum_\mu \langle Q_\mu | D_w \, P_\mu \rangle_{q,t} \\
& = \sum_\mu w^{|\mu|} \\
& = \prod_{n \geq 1} \frac{1}{1-w^n} 
\end{align*}

\end{proof}
\end{lemma}


\subsubsection{The proof}

The proof of the refined version of Theorem \ref{qt-borodin} now proceeds in the same manner as in section \ref{sec:proof}. We begin by applying Lemma \ref{lhs-mac}

\begin{align*}
\phantom{=} \sum_{\mathfrak{c} \in \CPP(\pi)} W_\mathfrak{c}(q,t) z^{|\mathfrak{p}|} 
& = \sum_\mu \langle Q_\mu \,|\, H^{\pi_0}(u_0) H^{\pi_1}(u_1)  \cdots H^{\pi_T}(u_T) D_w \, P_\mu  \rangle_{q,t}  \\
\end{align*}
Next we repeatedly applies the commutation relations of Lemma \ref{commutation-mac}, followed by definition \ref{N-def}.
\begin{align*}
& = \prod_{\substack{i < j \\ \pi_i > \pi_j} }  \frac{(t u_i u_j;q)_\infty}{(u_i u_j;q)_\infty}
\sum_\mu \langle Q_\mu \,|\, \prod_{\substack{k =1\\ \pi_k = 0}}^T \Omega_{q,t} [Xu_k] 
\prod_{\substack{k =1 \\ \pi_k = 1}}^T \Omega_{q,t}^*[Xu_k] \, D_w \, P_\mu \rangle_{q,t} \\
& = \prod_{\substack{i < j \\ \pi_i > \pi_j} } \frac{(t u_i u_j;q)_\infty}{(u_i u_j;q)_\infty} N_\pi(0) 
\end{align*}
We then repeatedly apply Lemma \ref{cylindric-shift}.
\begin{align*}
& = \prod_{\substack{i < j \\ \pi_i > \pi_j} } \frac{(t u_i u_j;q)_\infty}{(u_i u_j;q)_\infty} \prod_{m \geq 0} 
\left ( \prod_{\substack{(i,j) \\ \pi_i \neq \pi_j}} \frac{(tu_i u_j w^{m+1};q)_\infty}{(u_i u_j w^{m+1};q)_\infty} \right ) N_\pi (\infty) 
\end{align*}
Splitting the second product into two, and combining it with the first we have:

\begin{align*}
& = \prod_{m \geq 1}
\left ( \prod_{\substack{i < j \\ \pi_i > \pi_j} } \frac{(t u_i u_j w^{m-1};q)_\infty}{(u_i u_jw^{m-1};q)_\infty} \right )
\left ( \prod_{\substack{i > j \\ \pi_i > \pi_j} } \frac{(t u_i u_j w^m;q)_\infty}{(u_i u_jw^m;q)_\infty}\right ) N_\pi (\infty)
\end{align*}
Finally, applying Lemma \ref{infinity-mac} we have:
\begin{align*}
& = \prod_{m \geq 1} \frac{1}{1-w^m}
\left ( \prod_{\substack{i < j \\ \pi_i > \pi_j} } \frac{(t u_i u_j w^{m-1};q)_\infty}{(u_i u_jw^{m-1};q)_\infty} \right )
\left ( \prod_{\substack{i > j \\ \pi_i > \pi_j} } \frac{(t u_i u_j w^m;q)_\infty}{(u_i u_jw^m;q)_\infty}\right ) 
\end{align*}

As in the Schur case, to obtained the non-refined version of the Theorem, it suffices to take the following specialization of variables on both sides:
\begin{align} 
w & = z^{|T|} \\ 
u_k & = \begin{cases}
z^k & \text{ if $\pi_k = 1$} \\
z^{-k} & \text{ if $\pi_k = 0$}
\end{cases}
\end{align}

\section{The Macdonald weight}\label{sec:mac-weight}
Recall that in the plethystic notation \cite{garsia,lascoux}, if 
\[ a(q,t) = \sum_{n,m} a_{n,m} \, q^n t^m \]
with $a_{n,m} \in \mathbb{Z} $ and $a_{0,0} = 0$, then we have:
\begin{equation*} \Omega \left [ a(q,t) \right ] = \prod_{n,m} \frac{1 }{(1 - q^n t^m)^{a_{n,m}}} \end{equation*}

Making use of this notation, the cylindric weight function (equation \ref{weight}) may be given an explicit combinatorial description:

\begin{theorem}\label{comb-weight}
\begin{equation}\label{cylindric-weight}
W_\mathfrak{c}(q,t) = \Omega \left [ (q-t) \mathcal{D}_\mathfrak{c}(q,t) \right ] 
\end{equation}\label{eq:peak-valley}
where the alphabet $\mathcal{D}_\mathfrak{c}(q,t)$ is given by:
\begin{equation}\mathcal{D}_\mathfrak{c}(q,t) =  \sum_{s \in peak(\mathfrak{c})} \, q^{a_\mathfrak{c}(s)} t^{\ell_\mathfrak{c}(s)} - \sum_{s \in valley(\mathfrak{c})} \, q^{a_\mathfrak{c}(s)} t^{\ell_\mathfrak{c}(s)}
\end{equation}
\end{theorem}

\index{Cube!arm length}
\index{Cube!leg length}
The arm length of a cube in a cylindric plane partition is equal to the arm length of the corresponding box with respect to the partition in which it lies. Likewise, the leg length of a cube in a cylindric plane partition is equal to the leg length of the corresponding box with respect to the partition in which it lies (see Section \ref{sec:arms-legs}).

The precise definition of ``valley'' and ``peak'' cubes will be given in subsection \ref{sec:cube-lattice}. A proof of theorem \ref{comb-weight} will be given in subsection \ref{sec:weight-proof}. 
In section \ref{sec:hall} we remark that Theorem \ref{comb-weight} reduces to the combinatorial formula for the Hall--Littlewood weight function given in \cite{vuletic2} and \cite{vuletic}.

\subsection{Lattice paths on the cylinder}\label{sec:lattice}

The goal of this section is to give a bijection between cylindric plane partitions, defined as periodic interlacing sequences, and certain families of non-intersecting lattice paths on the cylinder - 
or equivalently rhombus tilings on the cylinder. 
In particular we need to understand what the meaning of a \emph{cube} is in the non-intersecting lattice path picture (see Definition \ref{def:cube})

Our construction is ``dual'' to the construction given in Krattenthaler's original paper \cite{gessel-1997}. Before proceeding any further, here is an example:

\begin{center}
\begin{tikzpicture}[scale=0.8]

\begin{scope}
\draw[step=0.5cm,color=gray] (0,0) grid (2.5,11);

\foreach \i in {0,...,2}
{
	\foreach \j in {0,...,11}
	{
		\path (\i,\j) node[fill,shape=circle,inner sep=0.3mm](A\i\j){};
	}
}

\foreach \i in {0,...,2}
{
	\foreach \j in {0,...,10}
	{
		\path (\i+0.5,\j+0.5) node[fill,shape=circle,inner sep=0.3mm](B\i\j){};
	}
}

\draw[thick] (A01) -- (B01) -- (A11) -- (B11) -- (A21) -- (B21);
\draw[thick] (A02) -- (B02) -- (A12) -- (B12) -- (A23) -- (B22);
\draw[thick] (A05) -- (B04) -- (A14) -- (B14) -- (A25) -- (B25);
\draw[thick] (A07) -- (B06) -- (A16) -- (B16) -- (A27) -- (B27);
\draw[thick] (A08) -- (B07) -- (A18) -- (B18) -- (A29) -- (B28);
\draw[thick] (A09) -- (B09) -- (A19) -- (B19) -- (A210) -- (B29);
\draw[thick] (A010) -- (B010) -- (A111) -- (B110) -- (A211) -- (B210);

\path (0,0) node[draw,shape=circle,inner sep=1mm]{};
\path (0,3) node[draw,shape=circle,inner sep=1mm]{};
\path (0,4) node[draw,shape=circle,inner sep=1mm]{};
\path (0,6) node[draw,shape=circle,inner sep=1mm]{};

\path (0,1) node[fill,shape=circle,inner sep=0.7mm]{};
\path (0,2) node[fill,shape=circle,inner sep=0.7mm]{};
\path (0,5) node[fill,shape=circle,inner sep=0.7mm]{};
\path (0,7) node[fill,shape=circle,inner sep=0.7mm]{};
\path (0,8) node[fill,shape=circle,inner sep=0.7mm]{};
\path (0,9) node[fill,shape=circle,inner sep=0.7mm]{};
\path (0,10) node[fill,shape=circle,inner sep=0.7mm]{};

\path (0.5,0.5) node[draw,shape=circle,inner sep=1mm]{};
\path (0.5,3.5) node[draw,shape=circle,inner sep=1mm]{};
\path (0.5,5.5) node[draw,shape=circle,inner sep=1mm]{};
\path (0.5,8.5) node[draw,shape=circle,inner sep=1mm]{};

\path (0.5,1.5) node[fill,shape=circle,inner sep=0.7mm]{};
\path (0.5,2.5) node[fill,shape=circle,inner sep=0.7mm]{};
\path (0.5,4.5) node[fill,shape=circle,inner sep=0.7mm]{};
\path (0.5,6.5) node[fill,shape=circle,inner sep=0.7mm]{};
\path (0.5,7.5) node[fill,shape=circle,inner sep=0.7mm]{};
\path (0.5,9.5) node[fill,shape=circle,inner sep=0.7mm]{};
\path (0.5,10.5) node[fill,shape=circle,inner sep=0.7mm]{};

\path (1,0) node[draw,shape=circle,inner sep=1mm]{};
\path (1,3) node[draw,shape=circle,inner sep=1mm]{};
\path (1,5) node[draw,shape=circle,inner sep=1mm]{};
\path (1,7) node[draw,shape=circle,inner sep=1mm]{};
\path (1,10) node[draw,shape=circle,inner sep=1mm]{};

\path (1.5,1.5) node[fill,shape=circle,inner sep=0.7mm]{};
\path (1.5,2.5) node[fill,shape=circle,inner sep=0.7mm]{};
\path (1.5,4.5) node[fill,shape=circle,inner sep=0.7mm]{};
\path (1.5,6.5) node[fill,shape=circle,inner sep=0.7mm]{};
\path (1.5,8.5) node[fill,shape=circle,inner sep=0.7mm]{};
\path (1.5,9.5) node[fill,shape=circle,inner sep=0.7mm]{};
\path (1.5,10.5) node[fill,shape=circle,inner sep=0.7mm]{};

\path (1,1) node[fill,shape=circle,inner sep=0.7mm]{};
\path (1,2) node[fill,shape=circle,inner sep=0.7mm]{};
\path (1,4) node[fill,shape=circle,inner sep=0.7mm]{};
\path (1,6) node[fill,shape=circle,inner sep=0.7mm]{};
\path (1,8) node[fill,shape=circle,inner sep=0.7mm]{};
\path (1,9) node[fill,shape=circle,inner sep=0.7mm]{};
\path (1,11) node[fill,shape=circle,inner sep=0.7mm]{};

\path (1.5,0.5) node[draw,shape=circle,inner sep=1mm]{};
\path (1.5,3.5) node[draw,shape=circle,inner sep=1mm]{};
\path (1.5,5.5) node[draw,shape=circle,inner sep=1mm]{};
\path (1.5,7.5) node[draw,shape=circle,inner sep=1mm]{};

\path (2,0) node[draw,shape=circle,inner sep=1mm]{};
\path (2,2) node[draw,shape=circle,inner sep=1mm]{};
\path (2,4) node[draw,shape=circle,inner sep=1mm]{};
\path (2,6) node[draw,shape=circle,inner sep=1mm]{};
\path (2,8) node[draw,shape=circle,inner sep=1mm]{};

\path (2.5,0.5) node[draw,shape=circle,inner sep=1mm]{};
\path (2.5,3.5) node[draw,shape=circle,inner sep=1mm]{};
\path (2.5,4.5) node[draw,shape=circle,inner sep=1mm]{};
\path (2.5,6.5) node[draw,shape=circle,inner sep=1mm]{};

\path (2,1) node[fill,shape=circle,inner sep=0.7mm]{};
\path (2,3) node[fill,shape=circle,inner sep=0.7mm]{};
\path (2,5) node[fill,shape=circle,inner sep=0.7mm]{};
\path (2,7) node[fill,shape=circle,inner sep=0.7mm]{};
\path (2,9) node[fill,shape=circle,inner sep=0.7mm]{};
\path (2,10) node[fill,shape=circle,inner sep=0.7mm]{};
\path (2,11) node[fill,shape=circle,inner sep=0.7mm]{};

\path (2.5,1.5) node[fill,shape=circle,inner sep=0.7mm]{};
\path (2.5,2.5) node[fill,shape=circle,inner sep=0.7mm]{};
\path (2.5,5.5) node[fill,shape=circle,inner sep=0.7mm]{};
\path (2.5,7.5) node[fill,shape=circle,inner sep=0.7mm]{};
\path (2.5,8.5) node[fill,shape=circle,inner sep=0.7mm]{};
\path (2.5,9.5) node[fill,shape=circle,inner sep=0.7mm]{};
\path (2.5,10.5) node[fill,shape=circle,inner sep=0.7mm]{};

\path (-5,5) node[shape=circle,inner sep=0.3mm](valley1){valley: $q^2t$};
\draw[dashed] (A02) -- (valley1);
\draw[dashed] (A06) -- (valley1);
\path (-2,4.5) node[shape=circle,inner sep=0.3mm](peak){peak: $-qt$};
\draw[dashed] (B02) -- (peak);
\draw[dashed] (B05) -- (peak);

\end{scope}


\end{tikzpicture}
\qquad \qquad
\includegraphics[height=3.6in,width=1.4in]{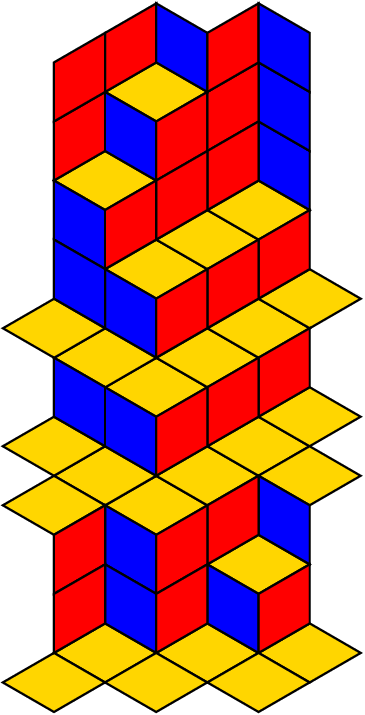}
\end{center}

\begin{equation*} 
 \mathfrak{c} = ((3,2,2), (4,3,2,1), (4,3,2), (6,4,3,2), (5,3,2), (3,2,2)) 
\end{equation*}

The bijection between the path model and the tiling model is clear. The white vertices correspond to the yellow tiles.
Each upstep of a path corresponds to a red tile. Each downstep of a path corresponds to a blue tile. Note that, due to cylindricity,
the first column corresponds to the last column, only shifted one step.

\subsection{Cubes in lattice path picture}\label{sec:cube-lattice}

\begin{definition}
We say that a vertex of the lattice is \emph{occupied} or \emph{black} if there is a path passing through that vertex,
otherwise we say that the vertex is \emph{unoccupied} or \emph{white}
\end{definition}

\index{Cube!lattice path definition}
\begin{definition}
A ``cube'' in the non-intersecting lattice path model corresponds to a pair of vertices $u = (x,y_1)$ and $v = (x,y_2)$
with $u$ coloured black, $v$ coloured white and $y_1 < y_2$. 
\end{definition}

\index{Cube!associated lattice path}
\begin{definition}
The \emph{path associated to the cube $(u,v)$} is the path which passes through the black vertex $v$.
\end{definition}

\index{Cube!valley cube}
\begin{definition}
A \emph{valley} cube is a cube $(u,v)$ for which the associated path takes a down step just before passing through $v$,
followed immediately by an upstep.
\end{definition}

\index{Cube!peak cube}
\begin{definition}
A \emph{peak} cube is a cube $(u,v)$ for which the associated path takes an up step just before passing through $v$,
followed immediately by a downstep.
\end{definition}

We have marked one peak cube and one valley cube on the diagram, together with their contribution to the alphabet $\mathcal{D}_\mathfrak{c}(q,t)$ in equation \ref{eq:peak-valley}. Since we are working on a cylinder, the first vertical is identified with the last vertical in such a way that each path forms a closed loop.

\subsection{Hall--Littlewood case}\label{sec:hall}

\index{Cube!surface cube}\label{def:cube-lattice}
\begin{definition}
 A \emph{surface cube} is a cube $(u,v)$ such that if $u = (x,y_1)$ and $v = (x,y_2)$
then for all $w = (x,y')$ with $y_1 < y' < y_2$ the vertex $w$ is coloured white. 
\end{definition}

\begin{definition}
 The \emph{level} of a surface cube $(u,v)$ is $y_2 - y_1$ where $u = (x,y_1)$ and $v = (x,y_2)$.
\end{definition}

Surface cubes are naturally in bijection with the yellow tiles in the rhombus tiling model.
In the Hall--Littlewood case we have $q=0$, thus the only boxes which contribute to the alphabet $\mathcal{D}_\mathfrak{c}(q,t)$ are those with arm-length zero. Since there is a bijection between such cubes and the ``yellow'' tiles of the rhombus tiling model, and since the leg length of the cube is precisely the \emph{level}, as defined in \cite{vuletic2} and \cite{vuletic}, it follows that Theorem \ref{comb-weight} reduces the combinatorial formula for the Hall--Littlewood weight case in \cite{vuletic2} and \cite{vuletic}.

\subsection{Bijection}

We shall now describe explicitly the bijection between cylindric plane partitions and families of
non-intersecting lattice paths on the cylinder.
The impatient reader may wish to skip the technical details in the section, and simply remark that the parts of the individual
partitions in the interlacing sequence model may be read, from right to left, from the ``heights'' of the corresponding surface cubes
in the rhombus tiling model. 

Equivalently, interpreting the black vertices as ones and the white vertices as zeros, the binary strings obtained from reading upwards along a given vertical gives the profile of the corresponding partition in the ``interlacing sequence model''.

\index{Cylinder}
\begin{definition}
The \emph{cylinder of period $T$} is the triangular lattice with vertices $(x,y)$ where either both $x$ and $y$ are even or both $x$ and $y$ are odd, and for which $0 \leq x \leq T$.
\end{definition}

In our example we have $T = 5$. We have drawn a segment of the cylinder corresponding to $0 \leq y \leq 22$.

\subsubsection{Encoding of paths}

\index{Cylinder!path}
\begin{definition}
A \emph{path} on the cylinder of period $T$ is a sequence of integers $(y_0, y_1, \ldots, y_T)$ with $y_0$ even such that for each $k$ we have either $y_{k+1} = y_k +1$ or $y_{k+1} = y_k-1$.
\end{definition}


\begin{lemma}\label{binpaths}
Each path on the cylinder of period $T$ may be uniquely
encoded by its starting position $y_0$ and the binary string $p$
given by $p^k = 1$ if $y_{k+1} = y_k + 1$ and $p^k = 0$ otherwise.
\end{lemma}

\index{Cylinder!non-intersecting lattice paths}
\begin{definition}
A \emph{family of non-intersecting lattice paths} on the cylinder of period $T$ is a collection of paths:
\begin{align*}
p^1 & = (y_0^1, y_1^1, \ldots y_T^1) \\
p^2 & = (y_0^2, y_1^2, \ldots y_T^2) \\
 & \cdots \\
p^m & = (y_0^m, y_1^m, \ldots y_T^m)
\end{align*}
satisfying $y_{k}^{i+1} > y_k^i$ for all $k$ and $i$ as well as 
$y_0^{i+1} - y_0^i = y_T^{i+1} - y_T^{i}$ for all $i$.
\end{definition}

Note that the second condition is necessary in order to ensure that it is possible to take the cylindric quotient identifying the vertices $(0,y)$ with the vertices $(T,y+d)$ for $d = m-n$ where $m$ is the number of ones in the profile, and $n$ is the number of zeros.

The paths in our example may be encoded using Lemma \ref{binpaths} as follows:
\begin{align*}
p^1 & =(1,0,1,0,1), \qquad y_0^1 = 2 \\
p^2 & =(1,0,1,1,0), \qquad y_0^2 = 4 \\
p^3 & =(0,0,1,1,1), \qquad y_0^3 = 10 \\
p^4 & =(0,0,1,1,1), \qquad y_0^4 = 14 \\
p^5 & =(0,1,1,1,0), \qquad y_0^5 = 16 \\
p^6 & =(1,0,1,1,0), \qquad y_0^6 = 18 \\
p^7 & =(1,1,0,1,0), \qquad y_0^7 = 20 \\
\end{align*}

\subsubsection{Extra conditions at boundary}

\begin{definition}
A family of non-intersecting lattice paths on the cylinder of period $T$ is said to be \emph{minimal with $m$ paths} if there is some $i$ such that $y_m^i - y_{m-1}^i > 2$.
\end{definition}

\begin{definition}
The \emph{profile} of a minimal family of $m$ non-intersecting lattice paths is the binary string associated to the $m$th path.
\end{definition}

Our example family of non-intersecting lattice paths is minimal with $7$ paths. Its profile is $\pi = 11010$. 

\subsubsection{Vertical reading of lattice paths}

The reader is referred back to section \ref{sec:inv} for the bijection between integer partitions and the binary string encoding their profiles.

\index{Cylinder!non-intersecting lattice paths!vertical reading}
\begin{definition}
The \emph{vertical reading} of a minimal family of $m$ non-intersecting
lattice paths is the sequence of binary strings $\rho^0, \rho^1, \ldots \rho^T$ obtained by reading, for each $k$, vertically upwards from the vertex $(k,y^1_k)$ to the vertex $(k,y^m_k)$, and recording a $0$ each time the vertex is occupied and a $1$ each time the vertex is unoccupied. 
\end{definition}

The vertical reading of our example family of non-intersecting lattice paths is the following:
\begin{align*}
\rho^0 & = 110010111 \quad \quad \rho^3 = 110101011 \\
\rho^1 & = 110101101 \quad \quad \rho^4 = 1010101011 \\
\rho^2 & = 110101011 \quad \quad \rho^5 = 110010111
\end{align*}
Observe that $\rho_0 = \rho_T$.

\subsubsection{Horizontal strips}

\begin{proposition}\label{prop:crucial}
Let $(\rho^0, \rho^1, \ldots, \rho^T)$ be a sequence of binary strings
arising from the vertical reading of a minimal family of non-intersecting lattice paths on a cylinder of period $T$ and profile $\pi$. For each $k \in \{1,2, \ldots T\}$ let $\mu^k$ denote the partition whose profile is 
given by $\rho^k$. If $\pi_k = 1$ then $\mu^{k} / \mu^{k-1}$ is a horizontal strip, otherwise if $\pi_k = 0$ then $\mu^{k-1} / \mu^{k}$ is a horizontal strip.
\begin{proof}
Follows immediately from the characterization of horizontal strips in terms of profiles in section \ref{sec:strip}
\end{proof}
\end{proposition}

\subsubsection{Length of $i$-th column}

Let $\mathfrak{c} = (\mu^0, \mu^1, \ldots, \mu^T)$ be an arbitrary cylindric plane partition with profile $\pi$.
For each $i$ let us define: 
\begin{equation}\label{horizontal-paths}
 \boxed{p_i(\mathfrak{c}) = ((\mu^0)'_i - (\mu^1)'_i + \pi_0, \ldots,
 (\mu^{T})'_i - (\mu^{T-1})'_i + \pi_{T-1}) }
\end{equation}

Note that $p_i(\mathfrak{c})$ encodes information about the length of the $i$-th column of the successive partitions in the interlacing sequence of $\mathfrak{c}$. This fact will be very important in section \ref{sec:switch-model}.

\begin{proposition}
For each $i$ we have that $p_i(\mathfrak{c})$ as defined in Equation (\ref{horizontal-paths}) is a binary string.
\begin{proof}
 If $\pi_k = 1$ then from definition \ref{cylindric-def} it follows that $\mu^k / \mu^{k-1}$ is a horizontal strip and thus
we have that $(\mu^{k-1})' - (\mu^k)' \in \{-1,0\}$.
Similarly, if $\pi_k = 0$ then $\mu^{k-1} / \mu^k$ is a horizontal strip and $(\mu^{k-1})' - (\mu^k)' \in \{1,0\}$
\end{proof}
\end{proposition}

\subsubsection{Bijection}

\begin{theorem}\label{path-bijection}
For any binary string $\pi$ of length $T$ there is a bijection between minimal families of non-intersecting lattice paths on the cylinder of period $T$ with profile $\pi$, and cylindric plane partitions of profile $\pi$. 
\begin{proof}
The map from families of non-intersecting lattice paths on the cylinder to interlacing sequences is given by taking vertical readings, and then translating from profiles to partitions. Conversely, the family of non-intersecting lattice paths
associated to a given interlacing sequence $\mathfrak{c}$ is given by
\[ \{(p_1(\mathfrak{c}),\tau_1), (p_2(\mathfrak{c}), \tau_2), \ldots (p_m(\mathfrak{c}), \tau_m)\}\]
where $\tau_i$ is the position of $i$-th one in the profile of $\mu^0$ and $p_i(\mathfrak{c})$ is defined in Equation \ref{horizontal-paths}.
\end{proof}
\end{theorem}

\section{Proof of Theorem \ref{comb-weight}}\label{sec:weight-proof}

Making use of the plethystic notation is just like taking the logarithm in order to turn multiplication and addition and division into subtraction. The terms which arise in the alphabet of the weight function must be regrouped appropriately before cancellations can take place between contributions of the same cube to two different Pieri co-efficients.

After the cancellation, a second regrouping of terms takes place. This regrouping involves putting together contributions of boxes which come from the came column number, but a different partition in the ``interlacing sequence'' model. 

Equivalently, we put together cubes which lie on the same vertical of the interlacing sequence model.
Proposition \ref{prop:crucial} is used extensively in this step.

\subsection{Plethystic notation}

We begin by making use of the \emph{plethystic notation} \cite{garsia, lascoux} to rewrite the $(q,t)$-Pieri coefficients (equations \ref{pieri-coeff} and \ref{pieri-coeff2}) in the following form:
\begin{align}
\varphi_{\lambda / \mu}(q,t) & = \Omega \left [ (q-t) (\mathcal{A}_{\lambda / \mu}(q,t) - \mathcal{B}_{\lambda / \mu}(q,t)) \right ] \\
\psi_{\lambda / \mu}(q,t) & = \Omega \left [ (q-t) (\mathcal{B}'_{\lambda / \mu}(q,t) - \mathcal{A}'_{\lambda / \mu}(q,t)) \right ]
\end{align}
where:
\begin{align}
\mathcal{A}_{\lambda / \mu}(q,t) & = \sum_{s \in C_{\lambda / \mu}} q^{a_\lambda(s)} t^{\ell_\lambda(s)} \\
\mathcal{B}_{\lambda / \mu}(q,t) & = \sum_{s \in \overline{C}_{\lambda / \mu}} q^{a_\mu(s)} t^{\ell_\mu(s)} \\
\mathcal{A}'_{\lambda / \mu}(q,t) & = \sum_{s \not\in C_{\lambda / \mu}} q^{a_\lambda(s)} t^{\ell_\lambda(s)} \\
\mathcal{B}'_{\lambda / \mu}(q,t) & = \sum_{s \not\in \overline{C}_{\lambda / \mu}} q^{a_\mu(s)} t^{\ell_\mu(s)} 
\end{align}

Here we have changed our notation slightly from that used in Macdonald \cite{macdonald}, so that now $C_{\lambda / \mu}$ denotes the set of boxes $s = (i,j) \in \lambda$ such that $\lambda'_j > \mu'_j$ while $\overline{C}_{\lambda / \mu}$ denotes the set of boxes $s = (i,j) \in \mu$ such that $\lambda'_j > \mu'_j$.

Making use of this notation we may rewrite equation \ref{weight} as:
\begin{equation} 
W_{\mathfrak{c}}(q,t) = \Omega \bigl [ (q-t) \mathcal{D}_{\mathfrak{c}}(q,t) \bigr ] 
\end{equation}
where:
\begin{equation}\label{unseparated}
\mathcal{D}_{\mathfrak{c}}(q,t) = \sum_{\substack{k =1 \\ \pi_k = 1}}^T  (\mathcal{A}_{k / {k-1}}- \mathcal{B}_{k / {k-1}}) 
+ \sum_{\substack{k = 1\\ \pi_k = 0}}^T (\mathcal{B}'_{k-1 / {k}} - \mathcal{A}'_{k-1 / {k}}) \\
\end{equation}
To avoid unnecessary indices, we use the convention that: 
\begin{equation} 
\mathcal{X}_{k / {k-1}} = \mathcal{X}_{\mu^k / \mu^{k-1}}(q,t)
\end{equation}
Our goal is to find a simplified expression for $\mathcal{D}_{\mathfrak{c}}(q,t)$. 

\subsection{Regrouping terms}
Recall from section \ref{cpp} that in the ``interlacing sequence'' model, a \emph{cube} of the cylindric plane partition $\mathfrak{c}$ corresponds to a \emph{box} of one of the underlying partitions $\mu^k$.

Observe now that each \emph{box} $s \in \mu^k$ contributes to at most two terms in equation \ref{unseparated}, one involving the pair of partitions $\mu^k$ and $\mu^{k-1}$, the other involving the pair of partitions $\mu^k$ and $\mu^{k+1}$.

Regrouping terms, and setting $\pi_{T+1} = \pi_1$ as well as $\mu^{T+1} = \mu^1$ we may write:
\begin{equation}
\mathcal{D}_{\mathfrak{c}}(q,t) = 
 \sum_{\substack{k =1 \\ \pi_k = 1 \\\pi_{k + 1} = 1}}^{T} \mathcal{E}^k_{11}(\mathfrak{c}) 
+ \sum_{\substack{k =1 \\ \pi_k = 0 \\ \pi_{k+1} = 0}}^{T} \mathcal{E}^k_{00}(\mathfrak{c}) 
+ \sum_{\substack{k =1 \\ \pi_k = 0 \\ \pi_{k+1}=1} }^{T} \mathcal{E}^k_{01}(\mathfrak{c}) 
+ \sum_{\substack{k =1 \\ \pi_k =1 \\ \pi_{k+1} = 0}}^{T} \mathcal{E}^k_{10}(\mathfrak{c}) 
\end{equation}

where:
\begin{align}
\mathcal{E}^k_{11}(\mathfrak{c}) & = \mathcal{A}_{k/k-1} - \mathcal{B}_{k+1/k}\\
\mathcal{E}^k_{00}(\mathfrak{c}) & = \mathcal{B}'_{k-1 / k} - \mathcal{A}'_{k/k+1}\\
\mathcal{E}^k_{01}(\mathfrak{c}) & = \mathcal{B}'_{k-1/k} - \mathcal{B}_{k+1/k} \\
\mathcal{E}^k_{10}(\mathfrak{c}) & = \mathcal{A}_{k/k-1}- \mathcal{A}'_{k/k+1}
\end{align}

For each $k$, there is only one term of the form $\mathcal{E}^k_{rs}(\mathfrak{c})$ appearing in the expression for $\mathcal{D}_{\mathfrak{c}}(q,t)$, and this term groups together all contributions from the boxes $s \in \mu^k$.

\subsection{Cancellations}

The next step is to observe that we have a large number of cancellations. For example:
\begin{align*}
\mathcal{E}^k_{11}(\mathfrak{c}) & = (\mathcal{A}_{k / {k-1}} - \mathcal{B}_{{k+1} / {k}})\\
& = \sum_{s \in C_{k / k-1}} q^{a_k(s)} t^{\ell_k(s)} - \sum_{s \in \overline{C}_{k+1/k}} q^{a_k(s)} t^{\ell_k(s)} \\
& = \sum_{s \in \mu^k} \sign_{11}(s) \, q^{a_k(s)} t^{\ell_k(s)}
\end{align*}
where:
\begin{equation} 
\sign_{11}(s) = 
\begin{cases}
1 & \text{if } s \in C_{k / k-1} \text{ and } s \not\in \overline{C}_{k+1/k} \\
0 & \text{if } s \in C_{k / k-1} \text{ and } s \in \overline{C}_{k+1/k} \\
0 & \text{if } s \not\in C_{k / k-1} \text{ and } s \not\in \overline{C}_{k+1/k}\\
-1 & \text{if } s \not\in C_{k / k-1} \text{ and } s \in \overline{C}_{k+1/k}
\end{cases}
\end{equation}
Again we are using a simplified notation:
\begin{align*}
a_k(s) & = a_{\mu^k}(s) \\ 
\ell_k(s) & = \ell_{\mu^k}(s)
\end{align*}
Similarly:
\begin{align*}
\mathcal{E}^k_{00} & = (\mathcal{B}'_{k-1 / k}- \mathcal{A}'_{k/k+1})\\
& = \sum_{s \not\in \overline{C}_{k-1 /k}} q^{a_k(s)} t^{\ell_k(s)} - \sum_{s \not\in C_{k/k+1}} q^{a_k(s)} t^{\ell_k(s)} \\
& = \sum_{s \in \mu^k} \sign_{00}(s) \, q^{a_k(s)} t^{\ell_k(s)}
\end{align*}
where:
\begin{equation} \sign_{00}(s) = 
\begin{cases}
-1 & \text{if } s \in \overline{C}_{k-1 /k} \text{ and } s \not\in C_{k/k+1}\\
0 & \text{if } s \in \overline{C}_{k-1 /k}\text{ and } s \in C_{k/k+1}\\
0 & \text{if } s \not\in \overline{C}_{k-1 /k} \text{ and } s \not\in C_{k/k+1}\\
1 & \text{if } s \not\in \overline{C}_{k-1 /k} \text{ and } s \in C_{k/k+1}
\end{cases}
\end{equation}
Next:
\begin{align*}
\mathcal{E}^k_{01} & = (\mathcal{B}'_{k-1/k} - \mathcal{B}_{k+1/k})\\
& = \sum_{s \not\in \overline{C}_{k-1 / k}} q^{a_k(s)} t^{\ell_k(s)} - \sum_{s \in \overline{C}_{k+1/k}} q^{a_k(s)} t^{\ell_k(s)} \\
& = \sum_{s \in \mu^k} \sign_{01}(s) \, q^{a_k(s)} t^{\ell_k(s)}
\end{align*}
where:
\begin{equation} \sign_{01}(s) = 
\begin{cases}
0 & \text{if } s \in \overline{C}_{k-1 / k} \text{ and } s \not\in \overline{C}_{k+1/k} \\
-1 & \text{if } s \in \overline{C}_{k-1 / k} \text{ and } s \in \overline{C}_{k+1/k} \\
1 & \text{if } s \not\in \overline{C}_{k-1 / k} \text{ and } s \not\in \overline{C}_{k+1/k} \\
0 & \text{if } s \not\in \overline{C}_{k-1 / k} \text{ and } s \in \overline{C}_{k+1/k}(
\end{cases}
\end{equation}
And:
\begin{align*}
\mathcal{E}^k_{10}(\mathfrak{c}) & = (\mathcal{A}_{k/k-1} - \mathcal{A}'_{k/k+1})\\
& = \sum_{s \in C_{k / k-1}} q^{a_k(s)} t^{\ell_k(s)} - \sum_{s \not\in C_{k/k+1}} q^{a_k(s)} t^{\ell_k(s)} \\
& = \sum_{s \in \mu^k} \sign_{10}(s) \, q^{a_k(s)} t^{\ell_k(s)}
\end{align*}
where:
\begin{equation} \sign_{10}(s) = 
\begin{cases}
0 & \text{if } s \in C_{k/ k-1}\text{ and } s \not\in C_{k/k+1} \\
1 & \text{if } s \in C_{k/ k-1} \text{ and } s \in C_{k/k+1} \\
-1 & \text{if } s \not\in C_{k/ k-1} \text{ and } s \not\in C_{k/k+1} \\
0 & \text{if } s \not\in C_{k/ k-1} \text{ and } s \in C_{k/k+1}
\end{cases}
\end{equation}

\subsection{Switching models}\label{sec:switch-model}

The final step in the proof is to switch from the ``interlacing sequence'' model of cylindric plane partitions (section \ref{cpp}) to the ``non-intersecting lattice path model'' (section \ref{sec:lattice}). This entails a grouping together of all the \emph{cubes} of the cylindric plane partition which belong to the same column of possibly different partitions in the sequence.

Recall that a \emph{cube} in the non-intersecting path model corresponds to a pair of vertices $v_1 = (x,y_1)$ and $v_2 = (x,y_2)$ with $y_1 < y_2$ where $v_1$ is coloured black and $v_2$ is coloured white.

Recall also that the $i$-th path in the non-intersecting path model encodes the length of the $i$-th column in each succeeding partition of the interlacing sequence (see equation (\ref{horizontal-paths}) in section \ref{sec:lattice}).

We shall say that the cube $c = (v_1,v_2)$ is of type $i$ if the black vertex $y_1$ lies on the $i$th path. This is equivalent to saying that the cube $c$ lies in the $i$th column of $\mu^k$ for some $k$.

\subsubsection{Four cases to check}

If $\pi_k = 0$ and $\mu^k = \mu^{k-1}$ then at the $k$th step, all the paths move downwards.
More generally if $\pi_k = 0$ then $\mu^k \preceq \mu^{k-1}$ and the $i$th path moves upwards if and only if the $i$th column of $\mu^k$ is shorter than the corresponding column of $\mu^{k-1}$.

That is to say, at the $k$th step, the $i$th path moves upwards if and only if $c \in \overline{C}_{k-1/k}$ for all $c \in \mu^k$ of type $i$.

If $\pi_k = 1$ and $\mu^k = \mu^{k-1}$ then at the $k$th step, all the paths move upwards.
More generally if $\pi_k = 1$ then $\mu^{k-1} \preceq \mu^{k}$ and the $i$th path moves downwards if and only if the $i$th column of $\mu^k$ is longer than the corresponding column of $\mu^{k-1}$.

That is to say, at the $k$th step, the $i$th path moves downwards if and only if $c \in C_{k/k-1}$ for all $c \in \mu^k$ of type $i$.

In a similar spirit, $\pi_{k+1} = 0$ and $\mu^{k+1} = \mu^{k}$ then at the $(k+1)$th step, all the paths move downwards.
More generally if $\pi_{k+1} = 0$ then $\mu^{k+1} \preceq \mu^{k}$ and the $i$th path moves up if and only if the $i$th column of $\mu^{k+1}$ is shorter than the corresponding column of $\mu^{k}$.

That is to say, at the $(k+1)$th step, the $i$th path moves upwards if and only if $c \in C_{k/k+1}$ for all $c \in \mu^k$ of type $i$.

Finally, if $\pi_{k+1} = 0$ and $\mu^{k+1} = \mu^{k}$ then at the $(k+1)$th step, all the paths move upwards.
More generally if $\pi_k = 1$ then $\mu^{k} \preceq \mu^{k+1}$ and the $i$th path moves downwards if and only if the $i$th column of $\mu^{k+1}$ is longer than the corresponding column of $\mu^{k}$.

That is to say, at the $(k+1)$th step, the $i$th path moves downwards if and only if $c \in \overline{C}_{k+1/k}$ for all $c \in \mu^k$ of type $i$.

Last but not least, one may check that the signs agree in all 16 possible cases.

\section{Conclusion}

We have proved a Macdonald polynomial analog of Borodin's identity. This simultaneously generalizes results of Okada \cite{okada} and Corteel and Savelief, Cyrille and Vuleti{\'c} \cite{vuletic}. We have also given a combinatorial interpretation of the weight function which is new even in the reverse plane partition case. Our proof relies heavily on the $(q,t)$-analog of the classical commutation relations originally given by Haiman, Garcia and Tesler \cite{garsia}. 

Is it possible to give a bijective proof of these commutation relations, and thus recover a bijective proof of the Macdonald polynomial analog of Borodin's identity? Is it possible to find a Macdonald polynomial analog of equation \ref{magic}?
\[ n! = \sum_\lambda f_\lambda^2 \]

\cleardoublepage

\part{Lambda determinants}

\chapter{Lambda determinants}

\section*{Introduction}

For those readers in a hurry, the most important definitions are \ref{def:F} and \ref{def:G}.
The key lemmas are proposition \ref{prop:inversion} and proposition \ref{prop:dual-inversion}. The later relies on proposition \ref{min-max}.
Remark \ref{rmk:dual-inversion} and Lemma \ref{reflect} are trivial but important.

The main theorem is stated in section \ref{sec:main-theorem}. The proof is by recurrence. The key step depends on a certain duality between \emph{left interlacing matrices} (section \ref{sec:left-interlacing}) and \emph{right interlacing matrices} (section \ref{sec:right-interlacing}) which allows us to factorize a certain sum. This duality is stated in propositions \ref{dual} and \ref{C-dual}.

\section{Permutations}\label{lambda:perm}

A \emph{permutation of $n$}  is simply a bijection from the set $\{1,2,\ldots n\}$ to itself. 
The number of permutations of $n$ is equal to $n!$.
It is often convenient to represent a permutation as matrix. We put a $1$ in position $(i,j)$ if $\sigma(j) = i$. In all other positions we place a $0$. 
For example, the permutation:
$\sigma = 24153$
is represented by the matrix:
\[
 \sigma = \left [ \begin{matrix}
0 & 0 & 1 & 0  & 0 \\
1 & 0 & 0 & 0  & 0 \\
0 & 0 & 0 & 0  & 1 \\
0 & 1 & 0 & 0  & 0 \\
0 & 0 & 0 & 1 & 0 \\
\end{matrix} \right ]
\]

We shall denote by $I$ the identity permutation which maps $i$ to $i$ for all $i$, and we shall denote by $J$ the \emph{maximum permutation} which maps $i$ to $n-i$ for all $i$:

\[  J =
 \left [ \begin{matrix}
0 & 0 & 0 & 0  & 1 \\
0 & 0 & 0 & 1  & 0 \\
0 & 0 & 1 & 0  & 0 \\
0 & 1 & 0 & 0  & 0 \\
1 & 0 & 0 & 0 & 0 \\ 
\end{matrix} \right ]
\]

Observe that the inversions of $\sigma$ correspond to the dual inversions of $J \circ \sigma$.

\subsection{Inversions and dual inversions}\label{sec:inversions}

\index{Permutation!inversion}

An \emph{inversion} in a permutation $\sigma$ is a pair $(i,j)$ with $i < j$ and $\sigma(i) > \sigma(j)$.

The inversions of permutation correspond to the zeros of the matrix which lie in the same row but to the left of a $1$ and the same column but above a $1$.

\[
 \left [ \begin{matrix}
* & * & 1 & 0  & 0 \\
1 & 0 & 0 & 0  & 0 \\
0 & * & 0 & *  & 1 \\
0 & 1 & 0 & 0  & 0 \\
0 & 0 & 0 & 1 & 0 \\ 
\end{matrix} \right ]
\]

\index{Permutation!dual inversion}

A \emph{dual inversion} in a permutation $\sigma$ is a pair $(i,j)$ with $i > j$ and $\sigma(i) > \sigma(j)$.
The dual inversions of permutation correspond to the zeros of the matrix which lie in the same row but to the \emph{right} of a $1$ and the same column but above a $1$.

\[
 \left [ \begin{matrix}
0 & 0 & 1 & *  & * \\
1 & * & 0 & *  & * \\
0 & 0 & 0 & 0  & 1 \\
0 & 1 & 0 & *  & 0 \\
0 & 0 & 0 & 1 & 0 \\ 
\end{matrix} \right ]
\]

\subsection{Determinants}\label{sec:determinants}

\index{Determinant!regular}

Our interest in permutations comes from the well-known formula for the determinant of a matrix:
\begin{equation} \det A =  \sum_\sigma (-1)^{\inv(\sigma)} A^\sigma  \end{equation}
where the sum is over all permutations of $n$
and we are using the notation:
\begin{equation}
A^{\sigma} = \prod_{i,j} a_{i,j}^{\sigma_{i,j}}
\end{equation}

\label{Dodgson condensation}

There is a curious method for calculating determinants, known as \emph{Dodgson condensation}, which is originally due to Charles Dodgson (aka Lewis Caroll) \cite{dodgson}. One proceeds as follows:

For each $k=0 \ldots n$ let us denote by $x_n[k]$ the doubly indexed collection of variables $x_n[k]_{i,j}$ with indices running from $i,j = 1..(n-k+1)$. 
One should think of the variables as forming a square pyramid with base $n+1$ by $n+1$. The index $k$ determines the ``height'' 
of the variable in the pyramid.

The variables are initialized as follows:
\begin{align*} 
x_n[0]_{i,j} & = 1 \mbox{ for all } i,j = 1..(n+1) \\
x_n[1]_{i,j} & = M_{i,j} \mbox{ for all } i,j = 1..n \\
\end{align*}

The value of the remaining variables is calculated via the following recurrence:
\begin{align} \label{rec1}
x_n[k+1]_{i,j} =  \frac{ x_n[k]_{i,j} x_n[k]_{i+1,j+1} - \, x_n[k]_{i,j+1} x_n[k]_{i+1,j}}{x_n[k-1]_{i+1,j+1}}
\end{align}
The end result is that:
\begin{equation} x_n[n]_{1,1} = \sum_\sigma (-1)^{\inv(\sigma)} x[1]^\sigma = \det(M) \end{equation} 

\subsection{Posets and lattices}

\index{Poset}

A \emph{partial order} is a set $P$ equipped with a binary relation $\leq$ satisfying:
\begin{itemize}
\item $a \leq a$ for all $a \in P$ (reflexivity)
\item $a \leq b$ and $b \leq c$ implies $a \leq c$ (transitivity)
\item $a \leq b$ and $b \leq a$ implies $a = b$ (antisymmetry)
\end{itemize}
We say that the element $b$ \emph{covers} the element $a \neq b$ if $a \leq b$
and if there is no $a \neq c \neq b$ such that $a \leq c \leq b$.

A poset is said to be \emph{graded} if their exists some function:
\[ w : P \to \mathbb{Z}_{\geq 0} \]
with the property that $w(b) = w(a) + 1$ whenever $a$ covers $b$.

\index{Lattice}

A \emph{lattice} is a partial order $P$ with the property that for every every $a,b \in P$
there exists unique element $a \wedge b$ (meet) and $a \vee b$ (join) with the property that $a \wedge b \leq a,b \leq a \vee b$
and for all $c$ such that $a,b \leq c$ we have $a \vee b \leq c$ and for all
$d$ such that $d \leq a,b$ we have $d \leq a \wedge b$.

Not every poset is a lattice, but very poset can be completed in a unique way to form a lattice.

\subsection{Bruhat order}

\index{Bruhat order}

A \emph{transposition} is a permutation $\sigma$ with the property that there exists some $i,j$ such that
$\sigma(i) = j$ and $\sigma(j) = i$ while for all other $k$ one has $\sigma(k) = k$.

The \emph{Bruhat order} may be defined as follows. The permutation $a$ covers $b$ if and only if there exists some transposition $\tau$ such that $b = \tau \circ a$ with $\inv(b) = \inv(a) + 1$.

The strong Bruhat order is graded by the involution number. The largest element in the Bruhat order is the maximum permutation $J$. One can show that for any permutations $\sigma$ and $\mu$ we have $\mu \leq \sigma$ if and only if $J \circ \sigma \leq J \circ \mu \circ$. 

The Bruhat order does \emph{not} form a lattice, since the permutations $132$ and $213$ have no join, while $312$ and $231$ have no meet. 

\subsection{Monotone triangles}
\index{Monotone triangle}

A monotone triangle is a triangle of integers:
\[
\begin{matrix}
a_{1,1} && a_{1,2} && a_{1,3} && a_{1,4} && a_{1,5} \\
 &a_{2,1}&  &a_{2,2}&  &a_{2,3}&  &a_{2,4}& \\
 && a_{3,1} && a_{3,2} && a_{3,3} &&  \\
 &&  &a_{4,1}&  &a_{4,2}&  && \\
 &&  && a_{5,1} &&  &&  \\
\end{matrix}
\]
with the following properties:
\begin{itemize}
\item $a_{1,k} = k$ for all $k$.
\item $a_{i,j+1} > a_{i,j}$ for all $i,j$.
\item $ a_{i,j} \leq a_{i+1,j} \leq a_{i+1,j+1}$.
\end{itemize}

A monotone triangle is a special case of a Gelfand-Tsetlin triangle \cite{schutz}. It is sometimes also referred to as a \emph{Gog triangle} \cite{zeil}.

There is a way to associate a unique monotone triangle to every permutation.
The rule is as follows. The $k$th row from the bottom of the monotone triangle associated to $\sigma$ contains a list, in increasing order, of those columns of the matrix of $\sigma$ which contain a $1$ in the final $k$ rows of the matrix. 

Here is the monotone triangle for our example permutation at the beginning of section \ref{lambda:perm}:

\[
\begin{matrix}
1 && 2 && 3 && 4 && 5 \\
 &1&  &2&  &4&  &5& \\
 && 2 && 4 && 5 &&  \\
 &&  &2&  &4&  && \\
 &&  && 4 &&  &&  \\
\end{matrix}
\]

The inversions in the monotone triangle picture corresponds to pairs of equal numbers, one immediately above and to the left of the other.
Here are the four inversions of our example:

\[
\begin{matrix}
{\color{blue} 1} && {\color{green} 2} && . && . && . \\
 &{\color{blue} 1}&  &{\color{green} 2}&  &.&  &.& \\
 && {\color{red} 2} && {\color{purple} 4} && . &&  \\
 &&  &{\color{red} 2}&  &{\color{purple} 4}&  && \\
 &&  && . &&  &&  \\
\end{matrix}
\]

The dual inversions in the monotone triangle picture correspond to pairs of equal numbers, one immediately above and to the right of the other.
Here are the six dual inversions of our example:
\[
\begin{matrix}
. && . && . && \color{green}4 && \color{red}5 \\
 &.&  &\color{purple}2&  &\color{green}4&  &\color{red}5& \\
 && \color{purple}2 && \color{green} 4 && \color{red}5 &&  \\
 &&  &.&  &\color{blue}4&  && \\
 &&  && \color{blue}4 &&  &&  \\
\end{matrix}
\]

The effect on the monotone triangle of a permutation $\sigma$ of multiplying the left by the maximum permutation $J$
is to send each $i$ to $n+1-i$, and then reflect around the main horizontal.

Not every monotone triangle is associated to a permutation. For example:
\[
\begin{matrix}
1 && 2 && 3 \\
 &1&  &3&  \\
 && 2 &&  
\end{matrix}
\]
We shall see in section \ref{sec:asm} that the above monotone triangle corresponds to the following \emph{alternating sign matrix}:
\[
\left [ \begin{matrix}
0 & 1 & 0 \\
1 & -1 & 1 \\
0 & 1 & 0
\end{matrix}
\right ]
\]

\section{Alternating Sign matrices} \label{sec:asm}

An \emph{alternating sign matrix} is a square matrix of $0$'s $1$'s and $-1$'s such that the sum of each row and column is $1$ and the non-zero entries in each row and column alternate in sign. For example:
\[A = 
\left ( \begin{matrix}
0 & 0 & 1 & 0 \\
0 & 1 & -1 & 1 \\
1 & 0 & 0 & 0 \\
0 & 0 & 1 & 0
\end{matrix} \right )
\]
A permutation matrix is an alternating sign matrix with no $(-1)$'s.


The total number of alternating sign matrices of size $n$ is given by:

\[A_n = \prod_{k=0}^{n-1} \frac{(3k+1)!}{(n+k)!} \qquad \qquad  1, 1, 2, 7, 42, 429, 7436, \ldots\]

This result was a conjecture for a long time. The first proof was given by Zeilberger \cite{zeil}. Zeilberger's proof was  long and complicated and eventually a simplified proof given by Kuperberg \cite{kuperberg}. Kuperberg's proof made use of ideas from the theory of integrable systems, the Yang-Baxter equation and the six vertex model with domain wall boundary conditions. It also made use of a recurrence relation due to Izergin and Korepin \cite{korepin}.

\subsection{Inversions and dual inversions}

\index{Alternating sign matrix!inversion}

As with permutations, an \emph{inversion} of an alternating sign matrix is a zeros in the matrix with the property that the sum of the entries in the same row lying to the right is equal to one, and the sum of the entries in the same column lying below is equal to one.
Our example alternating matrix $A$ has three $3$ inversions. 
\[
\left ( \begin{matrix}
* & * & 1 & 0 \\
* & 1 & -1 & 1 \\
1 & 0 & 0 & 0 \\
0 & 0 & 1 & 0
\end{matrix} \right )
\]
\index{Alternating sign matrix!dual inversion}

A \emph{dual inversion} of an alternating sign matrix is a zero in the matrix with the property that the sum of the entries in the same row lying to the \emph{left} is equal to one, and the sum of the entries in the same column lying below is equal to one.
Our example alternating matrix $A$ has $2$ dual inversions. 
\[
\left ( \begin{matrix}
0 & 0 & 1 & * \\
0 & 1 & -1 & 1 \\
1 & 0 & * & 0 \\
0 & 0 & 1 & 0
\end{matrix} \right )
\]

\begin{remark}\label{rmk:dual-inversion}
The alternating sign matrix $M$ has a dual inversion at position $(i,j)$ if and only if the alternating sign matrix $JM$ has an inversion at position $(i,n+1-j)$.
\end{remark}

\subsection{Monotone triangles}

\index{Monotone triangle}

Just like permutation matrices, alternating sign matrices may be associated with monotone triangles. 

The rule for constructing the monotone triangle from the alternating sign matrix is that the $k$th row from the bottom contains, in increasing order, a list of all those columns of the matrix whose final $k$ entries sum to $1$.  
Again this process is reversible.

Here is the monotone triangle of our example alternating sign matrix from section \ref{sec:asm}:
\[
\begin{matrix}
1 && 2 && 3 && 4 \\
 &1&  &2&  &4& \\
 && 1 && 3 &&  \\
 &&  &3&  && \\
\end{matrix}
\]

Every monotone triangle corresponds to an alternating sign matrix. The definition of inversions and dual inversions in the monotone triangle picture is exactly the same as for permutations.

\subsection{Completion of Bruhat order}\label{lambda:bruhat-lattice}

\index{Bruhat order}

One can show that the set of all monotone triangles of size $n$, under the partial order $a \leq b$ if and only if each element of $a$ is greater than or equal to the corresponding element of $b$, forms a lattice which completes the Bruhat order
\cite{lattice}. One monotone triangle covers another in the completion of the Bruhat order if and only if they differ by exactly $1$ at a single position.

Here is the \emph{Hasse diagram} for the Bruhat order on permutations for $n=3$:

\begin{center}
\begin{tikzpicture}[scale=0.8]
\begin{scope}
\path (5,-2) node[](min){
$
\begin{matrix}
1 && 2 && 3 \\
 &2&  &3&  \\
 && 3 && 
\end{matrix}
$
};

\path (2,2) node[](left){
$
\begin{matrix}
1 && 2 && 3 \\
 &2&  &3&  \\
 && 2 && 
\end{matrix}
$
};
\path (8,2) node[](right){
$
\begin{matrix}
1 && 2 && 3 \\
 &1&  &3&  \\
 && 3 && 
\end{matrix}
$
};
\path (2,5) node[](left2){
$
\begin{matrix}
1 && 2 && 3 \\
 &1&  &3&  \\
 && 1 && 
\end{matrix}
$
};
\path (8,5) node[](right2){
$
\begin{matrix}
1 && 2 && 3 \\
 &1&  &2&  \\
 && 2 && 
\end{matrix}
$
};

\path (5,9) node[](max){
$
\begin{matrix}
1 && 2 && 3 \\
 &1&  &2&  \\
 && 1 && 
\end{matrix}
$
};

\draw (min) -- (left) -- (left2) -- (max) -- (right2) -- (right) -- (min);
\draw (left) -- (right2);
\draw (right) -- (left2);

\end{scope}
\end{tikzpicture}
\end{center}

The completion of this diagram involves adding an extra element in the center, where the two lines cross. This element is the monotone triangle:
\[
\begin{matrix}
1 && 2 && 3 \\
 &1&  &3&  \\
 && 2 &&  
\end{matrix}
\]
corresponding to the alternating sign matrix:
\[
\left [ \begin{matrix}
0 & 1 & 0 \\
1 & -1 & 1 \\
0 & 1 & 0
\end{matrix}
\right ]
\]

\subsection{Lambda determinant}

Alternating sign matrices first appeared in the literature in the context of the so called \emph{lambda determinant} of Robbins and Rumsey \cite{robbins}. We shall now define the $\lambda$-determinant.

As in the case of Dodgson condensation (see Section \ref{sec:determinants}), for each $k=0 \ldots n$ let us denote by $x_n[k]$ the doubly indexed collection of variables $x_n[k]_{i,j}$ with indices running from $i,j = 1..(n-k+1)$. 
Again one should think of the variables as forming a square pyramid with base $n+1$ by $n+1$. The index $k$ determines the ``height'' 
of the variable in the pyramid.

The variables are initialized in the same way:
\begin{align*} 
x_n[0]_{i,j} & = 1 \mbox{ for all } i,j = 1..(n+1) \\
x_n[1]_{i,j} & = M_{i,j} \mbox{ for all } i,j = 1..n \\
\end{align*}
However the value of the remaining variables is calculated via the following modified recurrence:
\begin{align} \label{lambda-rec}
x_n[k+1]_{i,j} =  \frac{ x_n[k]_{i,j} x_n[k]_{i+1,j+1} + \lambda \, x_n[k]_{i,j+1} x_n[k]_{i+1,j}}{x_n[k-1]_{i+1,j+1}}
\end{align}
The end result \cite{robbins} is that:
\index{Determinant!lambda}
\begin{equation} \label{lambda-det} x_n[n]_{1,1} = \sum_{B \,\in \, \mathfrak{A}_n} \lambda^{\inv(B)}(1+\lambda)^{N(B)} M^B  \end{equation}
Here $\mathfrak{A}_n$ denotes the set of all alternating sign matrices of size $n$, $\inv(B)$ denotes the inversion number of $B$ and $N(B)$ denotes the number of negative ones in $B$.

The $\lambda$-determinant reduces to the regular determinant when $\lambda = -1$. 
The $\lambda$ determinant exhibits what is known as the \emph{Laurent phenomenon} \cite{laurent}. From the recursive definition we expect $x[n]_{1,1}$ to be a rational function. The fact that it turns out to be a Laurent polynomial is very surprising.

Our goal is to generalize equation \ref{lambda-det} by adding additional parameters. 

\section{Interlacing matrices}\label{sec:interlacing}

\subsection{Left corner sum matrices}\label{sec:left-corner}
\index{Corner sum matrix!left}

For each $n$ by $n$ alternating sign matrix $X$ let $\overline{X}$ be the matrix whose $(i,j)$-th entry is equal to the sum of the entries lying above and to the left of the $(i,j)$-th entry of $X$.  For example:
\[ X =
\left ( \begin{matrix}
0 & 1 & 0 & 0  \\
1 & -1 & 1 & 0  \\
0 & 1 & -1 & 1  \\
0 & 0 & 1 & 0 
\end{matrix} \right )
\quad
 \overline{X} =
\left ( \begin{matrix}
0 & 1 & 1 & 1  \\
1 & 1 & 2 & 2  \\
1 & 2 & 2 & 3  \\
1 & 2 & 3 & 4 
\end{matrix} \right )
\]

We shall refer to $\overline{X}$ as the \emph{left corner sum matrix} of $X$.
The original alternating sign matrix may be recovered by the formula:

\begin{align} 
X_{ij} & = \overline{X}_{ij} + \overline{X}_{i-1,j-1} - \overline{X}_{i,j-1} - \overline{X}_{i-1,j} \label{seconddiff}
\end{align}
If the indices are out of range, then the value of $\overline{X}_{ij}$ is taken to be zero.

The left corner sum matrices of alternating sign matrices have the following properties:
\begin{itemize}
\item The last row contains the integers from $1$ to $n$ in increasing order.
\item The last column contains the integers from $1$ to $n$ in increasing order.
\item Neighboring entries differ by at most one.
\item Each row and each column is non-decreasing.
\end{itemize}
Any non-negative integer matrix satisfying these properties is the left corner sum matrix of some alternating sign matrix \cite{robbins}.

\index{Alternating sign matrix!F(X) definition}

\begin{definition}\label{def:F}
For any alternating sign matrix $X$ we define:
\begin{equation}
F(X) = \overline{I} - \overline{X}
\end{equation}
In other words:
\begin{equation}
F(X)_{i,j} = \min(i,j) - \overline{X}_{i,j}
\end{equation}
\end{definition}

There is a natural order on the set of $n$ by $n$ left corner sum matrices given by $a \leq b$ if each entry of $a$ is \emph{greater than} or
equal to the corresponding entry of $b$. The order corresponds precisely to the lattice closure of the Bruhat order discussed in the section \ref{lambda:bruhat-lattice}.

\begin{remark}\label{remark-left}
Adding $1$ to position $(i,j)$ of $\overline{X}$ (where this is allowed) has the effect on $X$ of
adding the $2$ by $2$ matrix:
\[
\left [ \begin{matrix}
1 & -1 \\
-1 & 1
\end{matrix} \right ]
\]
to the $2$ by $2$ submatrix of $X$ whose upper left hand corner is at position $(i,j)$. 
\end{remark}

\subsection{Left interlacing matrices}\label{sec:left-interlacing}

\index{Interlacing matrices!left, down}

A pair of left corner sum matrices $(\overline{A},\overline{B})$ with $\overline{A}$ of dimensions $n$ by $n$ and $\overline{B}$ of dimension $n+1$ by $n+1$
are said to be \emph{left interlacing} if the following conditions are satisfied:

\[
\left ( \begin{matrix}
\overline{B}_{1,1} & & \overline{B}_{1,2} & & \overline{B}_{1,3}& & \overline{B}_{1,4}  \\
 &  {\color{red}\overline{A}_{1,1}}&  & \color{red}\overline{A}_{1,2} &  & \color{red}\overline{A}_{1,3} &   \\
\overline{B}_{2,1} & & \overline{B}_{2,2} & & \overline{B}_{2,3} & & \overline{B}_{2,4} \\
 &  {\color{red}\overline{A}_{2,1}}&  & \color{red}\overline{A}_{2,2} &  & \color{red}\overline{A}_{2,3} &   \\
\overline{B}_{3,1} & & \overline{B}_{3,2} & & \overline{B}_{3,3} & & \overline{B}_{3,4}  \\
 &  {\color{red}\overline{A}_{3,1}}&  & \color{red}\overline{A}_{3,2} &  & \color{red}\overline{A}_{3,3} &   \\
\overline{B}_{4,1} & & \overline{B}_{4,2} & & \overline{B}_{4,3} & & \overline{B}_{4,4} 
\end{matrix} \right )
\]

For all elements $x,y,z,w$ of $\overline{B}$ and and all elements $\color{red}a$ of $\color{red}\overline{A}$ 
which are arranged in the following configuration:

\[
\left ( \begin{matrix}
x & & y \\
 &  {\color{red}a}& \\
z & & w 
\end{matrix} \right )
\]
We must have:
\begin{equation}\label{red} 
\boxed{x,w-1 \leq {\color{red}a} \leq y,z}
\end{equation}

\hfill \break

An example:

\[
\left ( \begin{matrix}
0 & & 1 & & 1 & & 1  \\
 &  {\color{red}\{0,1\}}&  & \color{red}1 &  & \color{red}1 &   \\
1 & & 1 & & 2 & & 2  \\
 &  \color{red}1&  & \color{red}\{1,2\}&  & \color{red}2 &   \\
1 & & 2 & & 2 & & 3  \\
 & \color{red} 1&  & \color{red}2 &  & \color{red}3 &   \\
1 & & 2 & & 3 & & 4 
\end{matrix} \right )
\]

\begin{remark}\label{remark-left-A}
Above and to the left of a $-1$ in the alternating sign matrix $B$ there are two possible choices for the corresponding value of the left corner sum matrix $\overline{A}$. At all other positions there is a single choice \cite{robbins}.
\end{remark}

\index{Interlacing matrices!left, up}

Let us now consider the case of all $(n+1)$ by $(n+1)$ left corner sum matrices $\color{blue} \overline{C}$
which are left interlacing with a given left corner sum matrix $\overline{B}$:

\[
\left ( \begin{matrix} 
\color{blue} \overline{C}_{1,1} & & \color{blue} \overline{C}_{1,2}  & & \color{blue} \overline{C}_{1,3}  \\
& \overline{B}_{1,1} & & \overline{B}_{1,2} & \\
\color{blue} \overline{C}_{2,1} & & \color{blue} \overline{C}_{2,2}  & & \color{blue} \overline{C}_{2,3}  \\
& \overline{B}_{2,1} & & \overline{B}_{2,2} & \\
\color{blue} \overline{C}_{3,1} & & \color{blue} \overline{C}_{3,2}  & & \color{blue} \overline{C}_{3,3}  \\
\end{matrix} \right )
\]

The rule for constructing all possible $n+1$ by $n+1$ left corner sum matrices $\overline{C}$ which are interlacing with given $n$ by $n$ left corner sum matrix $\overline{B}$ is the last row and last column must be strictly increasing from $1$ to $n+1$
as well as that for all elements $x,y,z,w$ of $\overline{B}$ and and all remaining elements $\color{blue}c$ of $\color{blue}\overline{C}$ 
which are arranged in the following configuration:

\[
\left ( \begin{matrix}
x & & y \\
 &  {\color{blue}c}& \\
z & & w 
\end{matrix} \right )
\]
We must have:
\begin{equation}\label{blue}
\boxed{
y,z \leq {\color{blue}c} \leq w,x+1}
\end{equation}
If the element ${\color{blue}c}$ lies in the first row or column, then we take the out of range indices to be equal to zero.

\hfill \break

Here is an example:
\[
\left ( \begin{matrix} 
\color{blue} 0 & & \color{blue} \{0,1\} & & \color{blue} 1 \\
& 0 & & 1 & \\
\color{blue} \{0,1\} & & \color{blue} 1 & & \color{blue} 2 \\
& 1 & & 2 & \\
\color{blue} 1 & & \color{blue} 2 & & \color{blue} 3 \\
\end{matrix} \right )
\]

\begin{remark}\label{remark-left-C}
Above and to the left of a $1$ in the alternating sign matrix $B$ there are two possible choices for the corresponding value of $\overline{C}$. At all other positions there is a single choice \cite{robbins}. 
\end{remark}

\subsection{Inversions}

\index{Alternating sign matrix!inversion} 

\begin{lemma}\label{left-corner-sum-inversion}
Let $B$ be an alternating sign matrix.
If $B_{i,j}$ is an inversion then $\overline{B}_{i,j} = \overline{B}_{i-1,j-1}$ otherwise $\overline{B}_{i,j} > \overline{B}_{i-1,j-1}$.
\end{lemma}

The ``smallest'' alternating sign matrix $A$ which is left interlacing with a given alternating sign matrix $B$ is denoted by $A^{\min}$. 
By equation \ref{red} ts left corner sum matrix satisfies:
\begin{equation}\label{overline-min}
\boxed{\overline{A}^{\min}_{i,j} = \max (\overline{B}_{i,j}, \overline{B}_{i+1,j+1} - 1)}
\end{equation}

\index{Alternating sign matrix!inversion}
\begin{proposition}\label{prop:inversion}
If the alternating sign matrix $B$ has an inversion at position $(i,j)$ then: 
\[ F(B)_{i,j} - F(A^{\min})_{i-1,j-1} = 1\]
Otherwise: 
\[F(B)_{i,j} - F(A^{\min})_{i-1,j-1} = 0\]
\begin{proof}
\begin{equation}
F(B)_{i,j} - F(A^{\min})_{i-1,j-1} = 1 + \max(\overline{B}_{i-1,j-1}, \overline{B}_{i,j} - 1) - \overline{B}_{i,j}
\end{equation}
The result follows from Lemma \ref{left-corner-sum-inversion}.
\end{proof}
\end{proposition}

\subsection{Domino tiling of Aztec diamond}

\index{Aztec diamond}

An \emph{Aztec diamond} of size $n$ is the region of the lattice $\mathbb{Z}^2$ satisfying $|x| + |y| \leq n+1$. For example, here is an Aztec diamond of size $3$:  

\begin{center}
\includegraphics[height=2in,width=2in]{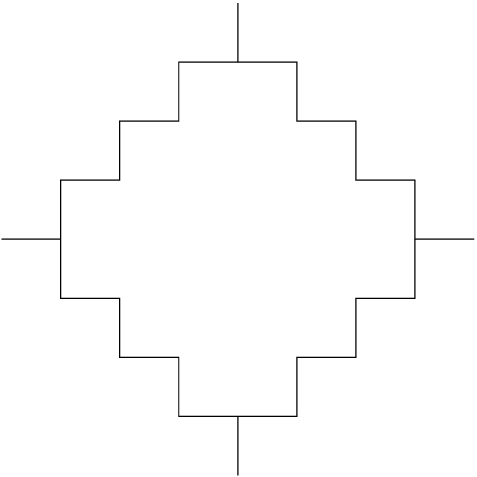}
\end{center}

A \emph{domino tiling} is a complete covering of an Aztec diamond by either $2$ by $1$ or $1$ by $2$ rectangles. For example:

\begin{center}
\includegraphics[height=2in,width=2in]{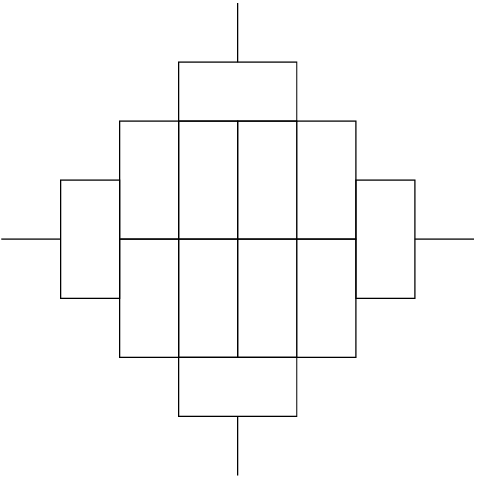}
\end{center}

The total number of domino tilings of the Aztec diamond of size $n$ is equal to $2^{n(n+1)/2}$ \cite{aztec}. There is a bijection between domino tilings of the Aztec diamond and pairs of left interlacing alternating sign matrices \cite{aztec}.

\subsubsection{An example}

We wish to find the pair of interlacing matrices $(A,B)$ which are in bijection with the domino tiling above.

If $n$ is even, then we begin by marking all the vertices $(i,j)$ on the interior of the diamond such that $i+j$ is even,
otherwise we begin by marking all the vertices $(i,j)$ on the interior of the diamond such that $i+j$ is odd.

\begin{center}
\includegraphics[height=2in,width=2in]{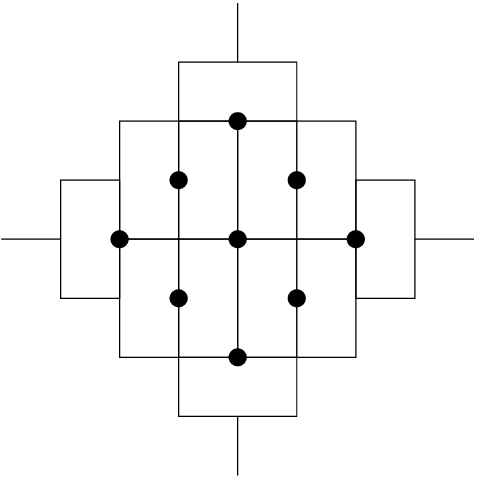}
\end{center}

Next we record the degree of each vertex. Note that the ``rows'' go up and to the right, while the ``columns'' go down and to the right:

\[
\left [ \begin{matrix}
3 & 2 & 3 \\
2 & 4 & 2 \\
3 & 2 & 3
\end{matrix} \right ]
\]

Finally we replace each $3$ with a zero, each $2$ with a one, and each $4$ with a negative one to get the matrix $A$.

\[
\left [ \begin{matrix}
0 & 1 & 0 \\
1 & -1 & 1 \\
0 & 1 & 0
\end{matrix} \right ]
\]

Next, if $n$ is even we mark all the vertices $(i,j)$ such that $i+j$ is odd, including vertices on the boundary of the diamond.
Otherwise we mark all the vertices $(i,j)$ such that $i+j$ is even, including vertices on the boundary of the diamond.
We imagine that the original tiling has been extended to infinity with horizontal tiles.

\begin{center}
\includegraphics[height=2in,width=2.1in]{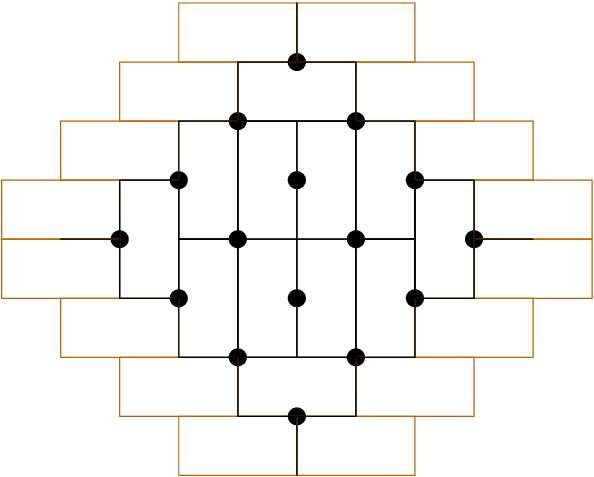}
\end{center}

Again we record the degree of each vertex:
\[
\left [ \begin{matrix}
3 & 3 & 4 & 3 \\
3 & 4 & 2 & 4 \\
4 & 2 & 4 & 3 \\
3 & 4 & 3 & 3 \\
\end{matrix} \right ]
\]

This time we replace $3$ with a zero, each $4$ with a one and each $2$ with a negative one, to obtain the matrix $B$.

\[
\left [ \begin{matrix}
0 & 0 & 1 & 0 \\
0 & 1 & -1 & 1 \\
1 & -1 & 1 & 0 \\
0 & 1 & 0 & 0 \\
\end{matrix} \right ]
\]

One may verify that in our example the two matrices are interlacing.

\subsubsection{Flips}

For the general case there are several things to check:
\begin{itemize}
\item That the two matrices obtained are always alternating sign matrices
\item That the two matrices obtained are always interlacing.
\item That the process is invertible.
\end{itemize}
The reader is referred to the paper \cite{aztec} for a complete proof. Nevertheless we shall make several remarks.

An \emph{elementary flip} involves changing two adjacent horizontal dominos into two adjacent vertical dominos, or vice versa.

\begin{center}
\includegraphics[height=0.6in,width=2in]{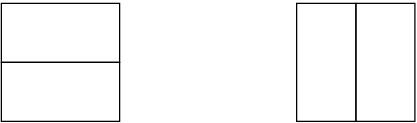}
\end{center}

One can show \cite{aztec} that every domino tiling of the Aztec diamond may be obtained from the ``minimal'' tiling, with only vertical tiles, by some sequence of elementary flips. For example, for $n=3$ we have:
\begin{center}
\includegraphics[height=3.6in,width=2in]{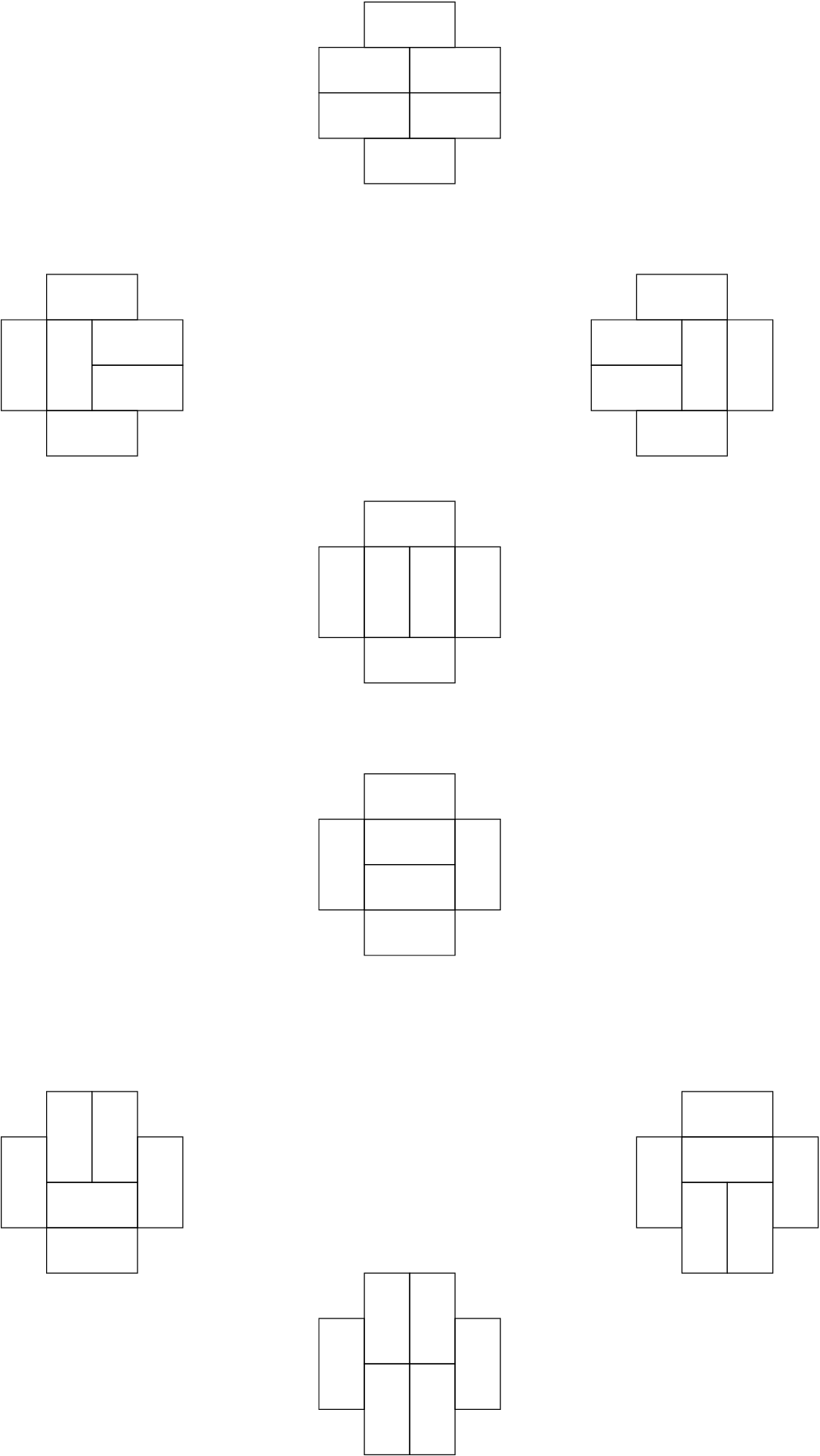}
\end{center}

The pair of interlacing matrices $(A,B)$ for the ``minimal'' tiling are a pair of maximal permutations, 
while for the ``maximal'' tiling, with all vertical tiles, they are a pair of identity permutations.

The midpoint $(i,j)$ of big square involved in an elementary flip is always of degree $2$ both before and after the flip. 

\begin{center}
\includegraphics[height=0.5in,width=2in]{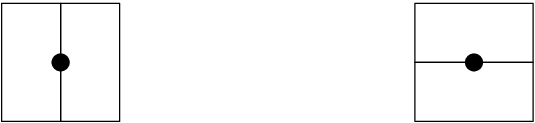}
\end{center}

There are two cases to consider. Either the midpoint belongs to the matrix $A$ and corresponds to a $1$. In this case the four vertices around the edge of the square being flipped all belong to $B$.

\begin{center}
\includegraphics[height=0.6in,width=2in]{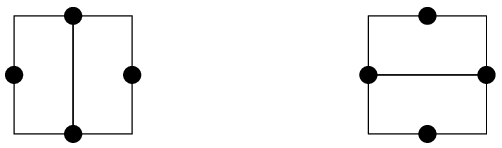}
\end{center}

The effect of of flipping from the vertical configuration to the horizontal configuration is to add the matrix:
\[\left [ \begin{matrix} -1 & 1 \\ 1 & -1 \end{matrix} \right ] \]
to these four vertices of $B$, or equivalently move up the Bruhat order by one step (see remark \ref{remark-left})

The second case to consider is when the midpoint of the square being flipped belongs to $B$ 
and corresponds to a $-1$. The four vertices around the edge of the square all belong to $A$.
In this case the effect of flipping \emph{from the horizontal configuration to the vertical configuration} is to add the matrix:
\[\left [ \begin{matrix} -1 & 1 \\ 1 & -1 \end{matrix} \right ]\]
 to these four vertices of $A$, or equivalently move up one step in the Bruhat order (see remark \ref{remark-right}).

One may conclude that since the bijection works in the base case, and since the elementary flips preserve both the alternating sign matrix condition and the interlacing condition, the bijection works in all cases. 

\section{Duality}\label{sec:duality2}

\subsection{Right corner sum matrices}\label{sec:right-cornersum}
\index{Corner sum matrix!right}

For each $n$ by $n$ alternating sign matrix $X$ let $\underline{X}
$ be the matrix whose $(i,j)$-th entry is equal to the sum of the entries lying above and to the \emph{right} of the $(i,j)$-th entry of $X$. For example:
\[ X =
\left ( \begin{matrix}
0 & 1 & 0 & 0  \\
1 & -1 & 1 & 0  \\
0 & 1 & -1 & 1  \\
0 & 0 & 1 & 0 
\end{matrix} \right )
\quad
 \underline{X} =
\left ( \begin{matrix}
1 & 1 & 0 & 0  \\
2 & 1 & 1 & 0  \\
3 & 2 & 1 & 1  \\
4 & 3 & 2 & 1 
\end{matrix} \right )
\]

We shall refer to $\underline{X}$ as the \emph{right corner matrix} of $X$.
The original alternating sign matrix may be recovered by the formula:

\begin{align} 
X_{ij} & = \underline{X}_{ij} + \underline{X}_{i-1,j+1} - \underline{X}_{i,j+1} - \underline{X}_{i-1,j}\label{seconddiff2}
\end{align}
If the indices are out of range, then the value of $\underline{X}_{ij}$ is taken to be zero. 

\begin{remark}\label{remark-right}
Adding one at position $(i,j)$ in the right corner sum matrix $\underline{X}$ is equivalent to adding the matrix 
$ \left (\begin{matrix} -1 & 1 \\ 1 & -1 \end{matrix} \right ) $ with upper \emph{right} hand corner at position $(i,j)$
to the alternating sign matrix $X$.
\end{remark}

\index{Alternating sign matrix!G(X), definition}
\begin{definition}\label{def:G}
For any alternating sign matrix $X$ we define:
\begin{equation}
G(X) = \underline{J} - \underline{X}
\end{equation}
In other words:
\begin{equation}
G(X)_{i,j} = \min(i,n+1-j) - \underline{X}_{i,j}
\end{equation}
\end{definition}

\begin{lemma}\label{reflect}
If $B' = BJ$ then $\underline{B}' = \overline{B}J$
\end{lemma}

\index{Alternating sign matrix! Bruhat order}

\begin{remark}\label{bruhat-right}
The Bruhat order on right corner sum matrices is given by $a \leq b$ if each entry of $a$ is \emph{less than} or
equal to the corresponding entry of $b$. 
\end{remark}

\subsection{Right interlacing matrices}\label{sec:right-interlacing}

\index{Interlacing matrices!right, down}

A pair of matrices $(A,B)$ of dimensions $n$ by $n$ and $n+1$ by $n+1$ respectively are said to be \emph{right interlacing}
if their right corner sum matrices  satisfy the following conditions:

\[
\left ( \begin{matrix}
\underline{B}_{1,1} & & \underline{B}_{1,2} & & \underline{B}_{1,3}& & \underline{B}_{1,4}  \\
 &  {\color{red}\underline{A}_{1,1}}&  & \color{red}\underline{A}_{1,2} &  & \color{red}\underline{A}_{1,3} &   \\
\underline{B}_{2,1} & & \underline{B}_{2,2} & & \underline{B}_{2,3} & & \underline{B}_{2,4} \\
 &  {\color{red}\underline{A}_{2,1}}&  & \color{red}\underline{A}_{2,2} &  & \color{red}\underline{A}_{2,3} &   \\
\underline{B}_{3,1} & & \underline{B}_{3,2} & & \underline{B}_{3,3} & & \underline{B}_{3,4}  \\
 &  {\color{red}\underline{A}_{3,1}}&  & \color{red}\underline{A}_{3,2} &  & \color{red}\underline{A}_{3,3} &   \\
\underline{B}_{4,1} & & \underline{B}_{4,2} & & \underline{B}_{4,3} & & \underline{B}_{4,4} 
\end{matrix} \right )
\]

For all elements $x,y,z,w$ of $\underline{B}$ and and all elements $\color{red}a$ of $\color{red}\underline{A}$ 
which are arranged in the following configuration:

\[
\left ( \begin{matrix}
x & & y \\
 &  {\color{red}a}& \\
z & & w 
\end{matrix} \right )
\]
We must have:
\begin{equation}\label{red2}
\boxed{
y,z-1 \leq {\color{red}a} \leq x,w}
\end{equation}

\hfill \break

Continuing with our example:

\[
\left ( \begin{matrix}
1 & & 1 & & 0 & & 0  \\
 &  {\color{red}1}&  & \color{red}\{0,1\} &  & \color{red}0 &   \\
2 & & 1 & & 1 & & 0  \\
 &  \color{red}2&  & \color{red}1&  & \color{red}\{0,1\} &   \\
3 & & 2 & & 1 & & 1  \\
 & \color{red} 3&  & \color{red}2 &  & \color{red}1 &   \\
4 & & 3 & & 2 & & 1 
\end{matrix} \right )
\]

\begin{remark}\label{above-right-minus}
Above and to the \emph{right} of a $-1$ in the alternating sign matrix $B$ there are two possible choices for the corresponding value of the right corner sum matrix $\underline{A}$. At all other positions there is a unique choice.
\end{remark}

\index{Interlacing matrices!right, up}

Let us now consider the case of all $(n+1)$ by $(n+1)$ left corner sum matrices $\color{blue} \underline{C}$
which are right interlacing with a given right corner sum matrix $\underline{B}$:

\[
\left ( \begin{matrix} 
\color{blue} \underline{C}_{1,1} & & \color{blue} \underline{C}_{1,2}  & & \color{blue} \underline{C}_{1,3}  \\
& \underline{B}_{1,1} & & \underline{B}_{1,2} & \\
\color{blue} \underline{C}_{2,1} & & \color{blue} \underline{C}_{2,2}  & & \color{blue} \underline{C}_{2,3}  \\
& \underline{B}_{2,1} & & \underline{B}_{2,2} & \\
\color{blue} \underline{C}_{3,1} & & \color{blue} \underline{C}_{3,2}  & & \color{blue} \underline{C}_{3,3}  \\
\end{matrix} \right )
\]

The rule for constructing all possible $n+1$ by $n+1$ right corner matrices $\underline{C}$ which are interlacing with a given $n$ by $n$ right corner sum matrix $\underline{B}$ is the first column must be strictly increasing from $1$ to $n+1$, the last row must be strictly decreasing from $n+1$ to $1$,
and for all elements $x,y,z,w$ of $\underline{B}$ and and all elements $\color{blue}c$ of $\color{blue}\underline{C}$ 
which are arranged in the following configuration:

\[
\left ( \begin{matrix}
x & & y \\
 &  {\color{blue}c}& \\
z & & w 
\end{matrix} \right )
\]
we must have:
\begin{equation}\label{blue2}
\boxed{
w,x \leq {\color{blue}c} \leq y+1,z}
\end{equation}

\hfill \break

Here is an example:
\[
\left ( \begin{matrix} 
\color{blue} 1 & & \color{blue} 1 & & \color{blue} \{0,1\} \\
& 1 & & 1 & \\
\color{blue} 2 & & \color{blue} \{1,2\} & & \color{blue} 1 \\
& 2 & & 1 & \\
\color{blue} 3 & & \color{blue} 2 & & \color{blue} 1 \\
\end{matrix} \right )
\]

\begin{remark}\label{above-right}
Above and to the \emph{right} of a $1$ in the alternating sign matrix $B$ there are two possible choices for the corresponding value of $\underline{C}$. At all other positions there is a single choice. 
\end{remark}





\subsection{Duality between left and right interlacing pairs}\label{sec:duality}
Left corner matrices
and right corner matrices are related by the following lemma:

\begin{lemma}\label{lemma}
For all $i,j$ we have:
\[\overline{B}_{i,j} + \underline{B}_{i,j+1} = i \]
\begin{proof}
The left hand side is equal to the sum of all the entries of the alternating sign matrix $B$ in the first $i$ rows. Since the sum of entries in each row of $B$ is equal to $1$, the final result is equal to $i$ as claimed. 
\end{proof}
\end{lemma}

For any binary string $\pi$, let $\tilde{\pi}$ denote the binary string in which all the zeros have been replaced by ones, and all the ones have been replaced by zeros. For example:
\begin{align*}
\pi & = 0010110100 \\
\tilde{\pi} & = 1101001011
\end{align*}

Suppose that we fix an $n+1$ by $n+1$ alternating sign matrix. In keeping with our previous examples, let us choose:
\[ B =
\left ( \begin{matrix}
0 & 1 & 0 & 0  \\
1 & -1 & 1 & 0  \\
0 & 1 & -1 & 1  \\
0 & 0 & 1 & 0 
\end{matrix} \right )
\]

If our alternating sign matrix $B$ has exactly $k$ negative ones then it followed that the number of $n$ by $n$ corner sum matrices which are interlacing with $\underline{B}$ is $2^k$.

Let us fix now an order on the negative ones, say top to bottom, left to right. We can now index all the $n$ by $n$ matrices which are right interlacing with $\underline{B}$ by a binary string. Continuing with the example from section \ref{sec:right-interlacing}, we have:

\begin{center}
\begin{tikzpicture}
\begin{scope}
\path (5,0) node[](min){
$
\underline{A}_{00} = \left ( \begin{matrix}
 1 & 0 & 0 \\
2 & 1 & 0 \\
3 & 2 & 1 \\
\end{matrix} \right ) 

$
};

\path (2,2) node[](left){
$
\underline{A}_{01} = \left ( \begin{matrix}
 1 & 0 & 0 \\
2 & 1 & 1 \\
3 & 2 & 1 \\
\end{matrix} \right ) 
$
};
\path (8,2) node[](right){
$
\underline{A}_{10} = \left ( \begin{matrix}
 1 & 1 & 0 \\
2 & 1 & 0 \\
3 & 2 & 1 \\
\end{matrix} \right ) 
$
};

\path (5,4) node[](max){
$
\underline{A}_{11} = \left ( \begin{matrix}
 1 & 1 & 0\\
2 & 1 & 1 \\
3 & 2 & 1 \\
\end{matrix} \right ) 
$
};

\draw (min) -- (left) -- (max) -- (right) -- (min);

\end{scope}
\end{tikzpicture}
\end{center}

A zero means that we chose the smaller of the two possibilities for $A$. A one means that we chose the larger.
Here are the corresponding alternating sign matrices:

\begin{center}
\begin{tikzpicture}
\begin{scope}
\path (5,0) node[](min){
$
A_{00} = \left ( \begin{matrix}
 1 & 0 & 0 \\
0 & 1 & 0 \\
0 & 0 & 1 \\
\end{matrix} \right )

$
};

\path (2,2) node[](left){
$
A_{01} = \left ( \begin{matrix}
 1 & 0 & 0 \\
0 & 0 & 1 \\
0 & 1 & 0 \\
\end{matrix} \right )

$
};
\path (8,2) node[](right){
$
A_{10} = \left ( \begin{matrix}
 0 & 1 & 0 \\
1 & 0 & 0 \\
0 & 0 & 1 \\
\end{matrix} \right ) 
$
};

\path (5,4) node[](max){
$
A_{11} = \left ( \begin{matrix}
 0 & 1 & 0 \\
1 & -1 & 1 \\
0 & 1 & 0 \\
\end{matrix} \right ) 
$
};

\draw (min) -- (left) -- (max) -- (right) -- (min);

\end{scope}
\end{tikzpicture}
\end{center}

Note that we use subscripts to indicate alternating sign matrices which are \emph{right} interlacing with $B$.
No underscore means the interlacing matrix itself, rather than the associated right corner sum matrix.
We shall denote by $A_{\min}$ the alternating sign matrix right interlacing with $B$ which is indexed by the binary string with all zeros and $A_{\max}$ the alternating sign matrix right interlacing with $B$ which is indexed by the binary string with all ones.

Let us now do the same with the \emph{left} corner sum matrix $\overline{X}$ from section \ref{sec:left-interlacing}.

\begin{center}
\begin{tikzpicture}
\begin{scope}
\path (5,0) node[](min){
$
\overline{A}^{00} = \left ( \begin{matrix}
 0 & 1 & 1 \\
1 & 1 & 2 \\
1 & 2 & 3 \\
\end{matrix} \right ) 
$
};

\path (2,2) node[](left){
$
\overline{A}^{01} = \left ( \begin{matrix}
 0 & 1 & 1 \\
1 & 2 & 2 \\
1 & 2 & 3 \\
\end{matrix} \right ) 
$
};
\path (8,2) node[](right){
$
\overline{A}^{10} = \left ( \begin{matrix}
 1 & 1 & 1 \\
1 & 1 & 2 \\
1 & 2 & 3 \\
\end{matrix} \right ) 
$
};

\path (5,4) node[](max){
$
\overline{A}^{11} = \left ( \begin{matrix}
 1 & 1 & 1 \\
1 & 2 & 2 \\
1 & 2 & 3 \\
\end{matrix} \right ) 
$
};

\draw (min) -- (left) -- (max) -- (right) -- (min);

\end{scope}
\end{tikzpicture}
\end{center}

Here are the corresponding alternating sign matrices:

\begin{center}
\begin{tikzpicture}
\begin{scope}
\path (5,0) node[](min){
$
A^{00} = \left ( \begin{matrix}
 0 & 1 & 0 \\
1 & -1 & 1 \\
0 & 1 & 0 \\
\end{matrix} \right )
$
};

\path (2,2) node[](left){
$
A^{01} = \left ( \begin{matrix}
 0 & 1 & 0 \\
1 & 0 & 0 \\
0 & 0 & 1 \\
\end{matrix} \right )
$
};
\path (8,2) node[](right){
$
A^{10} = \left ( \begin{matrix}
 1 & 0 & 0 \\
0 & 0 & 1 \\
0 & 1 & 0 \\
\end{matrix} \right )
$
};

\path (5,4) node[](max){
$
A^{11} = \left ( \begin{matrix}
 1 & 0 & 0 \\
0 & 1 & 0 \\
0 & 0 & 1 \\
\end{matrix} \right )
$
};

\draw (min) -- (left) -- (max) -- (right) -- (min);

\end{scope}
\end{tikzpicture}
\end{center}

Note that we use \emph{superscripts} to indication alternating sign matrices which are \emph{left} interlacing with $B$.
No overline indicates the alternating sign matrix itself, rather than the associated left corner sum matrix.
We shall denote by $A^{\min}$ the alternating sign matrix left interlacing with $B$ which indexed by the binary string with all zeros and $A^{\max}$ the alternating sign matrix left interlacing with $B$ which is indexed by the binary string with all ones.

\begin{proposition}\label{min-max}
\begin{equation} {A}_{\min} = A^{\max} \end{equation}
\begin{proof}
Consider the following segments of left and right interlacing matrices respectively:
\[
\left ( \begin{matrix}
a & & b & & c & & d  \\
 &  \color{red}x&  & \color{red}y &  & \color{red}z &   \\
e & & f & & g & & h \\
 &  \color{red}w&  & \color{red}u &  & \color{red}v &   \\
i & & j & & k & & \ell  \\
 &  \color{red}t&  & \color{red}s &  & \color{red}r &   \\
m & & n & & o & & p 
\end{matrix} \right )
\qquad
\left ( \begin{matrix}
a^* & & b^* & & c^* & & d^*  \\
 &  \color{red}x^*&  & \color{red}y^* &  & \color{red}z^* &   \\
e^* & & f^* & & g^* & & h^* \\
 &  \color{red}w^*&  & \color{red}u^* &  & \color{red}v^* &   \\
i^* & & j^* & & k^* & & \ell^*  \\
 &  \color{red}t^*&  & \color{red}s^* &  & \color{red}r^* &   \\
m^* & & n^* & & o^* & & p^* 
\end{matrix} \right )
\]

The elements $a$,$b$,$c$, etc... belong to the left corner sum matrix $\overline{B}$ while the elements
$x$,$y$,$z$ etc... belong to the left-interlacing matrix $\overline{A}^{\max}$.

Similarly the elements $a^*$,$b^*$,$c^*$, etc... belong to the right corner sum matrix $\underline{B}$ while the elements
$x^*$,$y^*$,$z^*$ etc... belong to the right-interlacing matrix $\underline{A}_{\min}$.

We wish to show that the value of the alternating sign matrix $A^{\max}$ at position $u$ is equal to the value of the alternating sign matrix $A_{\min}$ at position
$u^*$. That is, by equations (\ref{seconddiff}) and (\ref{seconddiff2}) we want to show that:

\[ x + u - w - y = u^* + z^* - y^* - v^* \]

As a consequence of lemma \ref{lemma} there is some $\gamma$ such that:
\begin{align*}
a + b^* & = b + c^* = c + d^* = \gamma \\
e + f^* & = f + g^* = g + h^* = \gamma+1 \\
i + j^* & = j + k^* = k + \ell^* = \gamma+2 \\
m + n^* & = n + o^* = o + p^* = \gamma+3 \\
\end{align*}

Now, by equations \ref{red} and \ref{red2}, we have:
\begin{align*}
& x + u - w - y \\
& = \min(a,f-1) + \min(f,k-1) - \min(e,j-1) - \min(b,g-1) \\
& = \min(\gamma - b^*,\gamma - g^*) + \min(\gamma + 1 - g^*, \gamma + 1 - \ell^*) \\
& \phantom{*******} - \min(\gamma + 1 - f^*,\gamma+1- k^*) - \min(\gamma - c^*,\gamma - h^*) \\
& = - \max(b^*,g^*) - \max( g^*, \ell^*) + \max( f^*, k^*) + \max( c^*, h^*) \\
& = -y^* - v^* + u^* + z^*
\end{align*}

The result follows.
\end{proof}
\end{proposition}

\begin{proposition}\label{dual}
If $\pi$ is a binary string, and $\tilde{\pi}$ is its complement, then $A^\pi = A_{\tilde{\pi}}$

\begin{proof}
This follows immediately from remarks \ref{remark-left} and \ref{remark-right}. It is possible to remove a $1$ from position $(i,j)$ of $\overline{A}^\pi$
if and only if it is possible to add a $1$ at position $(i,j+1)$ of $\underline{A}_{\tilde{\pi}}$.

\end{proof}
\end{proposition}

\begin{corollary}
\begin{equation}
A^{\min} = A_{\max}
\end{equation}
\end{corollary}

Suppose now that we fix an  $n$ by $n$ alternating sign matrix $B$ containing exactly $s$ ones. It follows that the number of $n+1$ by $n+1$ right corner sum matrices which are right interlacing with $\underline{B}$ is equal to $2^s$. 

If we fix an order on the ones, then each of these right corner sum matrices may be indexed by a binary string.
Let $C_\pi$ denote the $n+1$ by $n+1$ alternating sign matrix associated with $\underline{C}_\pi$

Similarly the number of $n+1$ by $n+1$ left corner sum matrices which are left interlacing with $\overline{B}$ is equal to $2^s$.
Again we may index these by binary strings. Let $C^\pi$ denote the $n+1$ by $n+1$ alternating sign matrix associated to $\overline{C}^\pi$.
\begin{proposition}\label{C-dual}
We have:
\[ C^{\max} = C_{\min}\] 
and more generally, if $\tilde{\pi}$ is the complement of $\pi$ then 
\begin{equation}
C^\pi = C_{\tilde{\pi}}\end{equation}
\begin{proof}
The proof is essentially identical to proposition \ref{min-max} and proposition \ref{dual} 
\end{proof}
\end{proposition}

\index{Alternating sign matrix!dual inversion}
\begin{proposition}\label{prop:dual-inversion}
If the $n+1$ by $n+1$ alternating sign matrix $B$ has a dual inversion at position $(i,n+1-j)$ then: 
\[ G(B)_{i,n+1-j} - G(A_{\max})_{i-1,n+1-j} = 1\]
Otherwise: 
\[G(B)_{i,n_1-j} - G(A_{\max})_{i-1,n+1-j} = 0\]
\begin{proof}
By Lemma \ref{reflect} we have:
\begin{align}G(B)_{i,n+1-j} - G(A_{\max})_{i-1,n+1-j} & = F(B')_{i,j} - F((A_{\max})')_{i-1,j-1} \\
& = F(B')_{i,j} - F((A^{\min})')_{i-1,j-1} 
\end{align}
But by remark \ref{reflect} $B$ has a dual inversion at position $(i,n+1-j)$ if and only if $B'$ has an inversion at position $(i,j)$.
The result now follows from proposition \ref{prop:inversion}.
\end{proof}
\end{proposition}

\section{Main Theorem}\label{sec:main-theorem}
\subsection{Notation}

\index{Alternating sign matrix!lambda weight}

Before we can state our formula, we need a few more definitions.
Let us define the \emph{lambda weight} of a $k$ by $k$ alternating sign matrix $B$ to be:
\begin{equation}
\boxed{F_\lambda(B) = \lambda^{F(B)} = \prod_{i,j=1}^k \lambda_{i,j}^{\min(i,j) - \overline{B}_{i,j}} }
\end{equation}
(see Section \ref{sec:left-corner})

\index{Alternating sign matrix!mu weight}

Similarly, let us define the \emph{mu weight} of an $k$ by $k$ alternating sign matrix $A$ to be:
\begin{equation} \boxed{G^n_\mu(B) = \prod_{i,j=1}^k \mu_{i,n+1-j}^{\min(i,k+1-j) - \underline{B}_{i,j}}} \end{equation}
(see Section \ref{sec:right-cornersum}).

\index{Alternating sign matrix!inversion}
\index{Alternating sign matrix!dual inversion}
Let $\Inv(B)$ denote the matrix with the property that $\Inv(B)_{i,j}=1$ if $B_{i,j}$ is an inversion, otherwise $\Inv(B)_{i,j}=0$.
Similarly, let $\Dinv(B)$ denote the matrix with the property that $\Dinv(B)_{i,j}=1$ if $B_{i,j}$ is a dual inversion, otherwise
$\Dinv(B)_{i,j} = 0$.

\index{Alternating sign matrix!shifts}

Finally, two more pieces of notation. We shall define the following two operators on laurent polynomials in the variables $\{\lambda_{i,j}\}$ and 
$\{\mu_{i,j}\}$:
\begin{align}
s(\lambda_{i,j}) & = \lambda_{i+1,j+1} \\
t(\mu_{i,j}) & = \mu_{i+1,j}
\end{align}

\subsection{Statement of theorem}
For each $k=0 \ldots n$ let us denote by $x_n[k]$ the doubly indexed collection of variables $x_n[k]_{i,j}$ with indices running from $i,j = 1..(n-k+1)$. 
One should think of the variables as forming a square pyramid with base $n+1$ by $n+1$. The index $k$ determines the ``height'' 
of the variable in the pyramid.

The variables are initialized as follows::
\begin{align*} 
x_n[0]_{i,j} & = Y_{i,j} \mbox{ for all } i,j = 1..(n+1) \\
x_n[1]_{i,j} & = X_{i,j} \mbox{ for all } i,j = 1..n \\
\end{align*}

The remaining variables are defined by the following recurrence:
\begin{align} \label{rec2}
x_n[k+1]_{i,j} =  \frac{ \mu_{i,n-k+1-j} x_n[k]_{i,j} x_n[k]_{i+1,j+1} + \lambda_{i,j} x_n[k]_{i,j+1} x_n[k]_{i+1,j}}{x_n[k-1]_{i+1,j+1}}
\end{align}

\index{Determinant!lambda, multi-parameter}

\begin{theorem}\label{theorem:lambda}
\begin{equation} \label{double}
x_n[k+1]_{1,1} = \sum_{\substack{(A,B) \\ |B| = k, |A| = k-1}}
\frac{F_\lambda(B)}{s(F_\lambda(A))} \frac{G^{n}_\mu(B)}{t(G^{n}_\mu(A))} x_n[1]^B s(x_n[0])^{-A} 
\end{equation}
The sum is over all pairs of left interlacing matrices.
\end{theorem}

If instead we initialize the variables as:
\begin{align*} 
x_n[0]_{i,j} & = 1 \mbox{ for all } i,j = 1..(n+1) \\
x_n[1]_{i,j} & = M_{i,j} \mbox{ for all } i,j = 1..n \\
\end{align*}
then we obtain the following corollary:

\begin{corollary}\label{lambda:corollary}
\begin{equation}
x_n[n]_{1,1} = \sum_{|B| = n} M^B \left ( 
\lambda^{\Inv(B)}
\mu^{\Dinv(B)'}
\prod_{B_{i,j} = -1} 
\left (\mu_{i,n+1-j} + \lambda_{i,j} \right )
\right )
\end{equation}
\end{corollary}

\subsection{Proof of Theorem \ref{theorem:lambda}}
This proof follows closely the original proof of Robbins and Rumsey \cite{robbins} for the $\lambda$-determinant \cite{robbins}.
\subsubsection{Base case}
In the case $k=1$ we have from the recurrence (equation \ref{rec2}) that:
\begin{equation}\label{base-case}
x_n[2]_{1,1} = \frac{\mu_{1,n-1} X_{1,1} X_{2,2} - \lambda_{1,1}X_{1,2}X_{2,1}}{Y_{2,2}} 
\end{equation}
while the closed form expression in equation \ref{double} tells us that we must take the sum over all interlacing pairs $(A,B)$ where $A$ is $1$ by $1$ and 
$B$ is $2$ by $2$.
There are exactly two such pairs of interlacing matrices, corresponding to the two terms in equation \ref{base-case}.

\subsubsection{Induction hypothesis}
Suppose now that the proposition is true for all $m \leq k$. 
We wish to show that:
\begin{equation}\label{toprove} x_n[k+2]_{1,1} = \sum_{\substack{(B,C) \\ |C| = k+1, |B| = k}} \frac{F_\lambda(C)} {s(F_\lambda(B))}
\frac{G^{n}_\mu(C)}{ t(G^{n}_\mu(B))} x[1]^C s(x[0])^{-B} 
\end{equation}
By the induction hypothesis, we know that (using a different set of variables):
\[ y_{n-1}[k+1]_{1,1} = \sum_{\substack{(A,B) \\ |B| = k, |A| = k-1}} 
\frac{F_\lambda(B)} {s(F_\lambda(A))}
\frac{G^{n-1}_\mu(B)}{ t(G^{n-1}_\mu(A))} y_{n-1}[1]^B s(y_{n-1}[0])^{-A} \]
If we now set $x_n[k+1] = y_{n-1}[k]$ then we have:
\begin{equation}\label{shifted}
x_n[k+2]_{1,1} = \sum_{\substack{(A,B) \\ |B| = k, |A| = k-1}} \frac{F_\lambda(B)} {s(F_\lambda(A))}
\frac{G^{n-1}_\mu(B)}{ t(G^{n-1}_\mu(A))} x_n[2]^B s(x_n[1])^{-A}
\end{equation}
But by the recurrence we know that:
\[ x_n[2]_{i,j} = \frac{ \mu_{i,n-j} x_n[1]_{i,j} x_n[1]_{i+1,j+1} + \lambda_{i,j} x_n[1]_{i,j+1} x_n[1]_{i+1,j}}{x_n[0]_{i+1,j+1}} \]
To simplify notation we will write:
\[
x_n[2]^B = \frac{D(x_n[1])^B}{s(x_n[0])^B}
\]
where:
\begin{equation}
D(x_n[k+1])_{i,j} = \mu_{i,n-k+1-j} x_n[k]_{ij} x_n[k]_{i+1,j+1} + \lambda_{ij} \, x_n[k]_{i,j+1} x_n[k]_{i+1,j} 
\end{equation}

\subsubsection{Term by term comparison of coefficients}
Let us now fix an arbitrary alternating sign matrix $B$ of dimensions $k$ by $k$. The coefficient of $s(x_n[0])^{-B}$ on the right hand side of equation \ref{toprove} is given by:
\begin{equation} \label{rhs}
\sum_{|C| = k+1} \frac{F_\lambda(C)} {s(F_\lambda(B))}
\frac{G^{n}_\mu(C)}{ t(G^{n}_\mu(B))} x_n[1]^C 
\end{equation}
where the sum is over all $C$ which are interlacing with $B$.
Similarly, the coefficient of $s(x_n[0])^{-B}$ on the right hand side of equation \ref{shifted} is given by:
\begin{equation}\label{lhs-lambda}
\sum_{ |A| = k-1} \frac{F_\lambda(B)} {s(F_\lambda(A))}
\frac{G^{n-1}_\mu(B)}{ t(G^{n-1}_\mu(A))} D(x_n[1])^B s(x_n[1])^{-A} 
\end{equation}
where the sum is over all $A$ which are interlacing with $B$.

We shall now use the results of section \ref{sec:duality} to simplify these two expressions.

\subsubsection{Duality: $C^{\max} = C_{\min}$}

Beginning with equation \ref{rhs} we use the fact that the set of $k+1$ matrices interlacing with $B$ decomposes as a boolean lattice, as well as remark \ref{remark-left} to factorize equation \ref{rhs} into:
\[
\frac{F_\lambda(C^{\max})} {s(F_\lambda(B))}
\frac{G^{n}_\mu(C_{{\max}})}{ t(G^{n}_\mu(B))} x_n[1]^{C^{\max}} 
\prod_{B_{ij} = 1}\left (\mu_{i,n+1-j} + \lambda_{ij} \frac{x_n[1]_{i+1,j} x_n[1]_{i,j+1}}{x_n[1]_{ij} x_n[1]_{i+1,j+1}}\right )
\]
Re-arranging slightly, this gives:
\begin{equation}\label{rhs2}
\frac{F_\lambda(C^{\max})} {s(F_\lambda(B))}
\frac{G^{n}_\mu(C_{\max})}{ t(G^{n}_\mu(B))} x_n[1]^{C^{\max}} 
\prod_{B_{ij} = 1}
\frac{D(x_n[1]_{i,j})}{x_n[1]_{ij} s(x_n[1]_{i,j})}
\end{equation}

\subsubsection{Duality: $A^{\max} = A_{\min}$}

Now for equation \ref{lhs-lambda}. We use the fact that the set of $k$ matrices interlacing with $B$ also decomposes as a boolean lattice, as well as remark \ref{remark-right} to factorize equation \ref{lhs-lambda} into:
\begin{align*}
& \frac{F_\lambda(B)} {s(F_\lambda(A^{\min}))}
\frac{G^{n-1}_\mu(B)}{ t(G^{n-1}_\mu(A_{{\min}}))} D(x_n[1])^B s(x_n[1])^{-A^{\min}} \\
& \phantom{=======----}\prod_{B_{i,j} = -1} \left (\mu_{i,n+1-j} + \lambda_{i,j} \frac{x_n[1]_{i+1,j} x_n[1]_{i,j+1}}{x_n[1]_{ij} x_n[1]_{i+1,j+1}}\right )
\end{align*}
Re-arranging slightly, this gives:
\[  \frac{F_\lambda(B)} {s(F_\lambda(A^{\min}))}
\frac{G^{n-1}_\mu(B)}{ t(G^{n-1}_\mu(A_{{\min}}))}  s(x_n[1])^{-A^{\min}} 
D(x_n[1])^B\prod_{B_{i,j} = -1} \frac{D(x_n[1]_{i,j})}{x[1]_{ij} s(x_n[1]_{i,j})} \]
which by combining the terms of the form $D(x[1])$ may be simplified even further to give:
\begin{equation}\label{lhs2}
 \frac{F_\lambda(B)} {s(F_\lambda(A^{\min}))}
\frac{G^{n-1}_\mu(B)}{ t(G^{n-1}_\mu(A_{{min}}))}  s(x_n[1])^{-A^{\min}} 
\prod_{B_{i,j} = 1}D(x_n[1]_{i,j}) \prod_{B_{i,j} = -1} \frac{1}{x_n[1]_{ij} s(x_n[1]_{i,j})}
\end{equation}

\subsubsection{Summarizing what we need to prove}

We wish to show that equation \ref{rhs2} is equal to equation \ref{lhs2}. Equivalently, we wish to show that:

\begin{align}\label{fish} 
&{s(F_\lambda(A^{\min}))}{F_\lambda(C^{\max})} 
{ t(G^{n-1}_\mu(A_{{\min}}))}{G^{n}_\mu(C_{{\max}})} x_n[1]^{C^{\max}} s(x_n[1])^{A^{\min}} \\
= & \qquad{F_\lambda(B)} {s(F_\lambda(B))}
{G^{{n-1}}_\mu(B)}{ t(G^{
{n}}_\mu(B))} (x_n[1] s(x_n[1]))^B \notag
\end{align}

\subsubsection{Comparing power of $\lambda_{i,j}$ on both sides}

In this section we show that:
\begin{equation}\label{cat}
{s(F_\lambda(A^{\min}))}{F_\lambda(C^{\max})} = {F_\lambda(B)} {s(F_\lambda(B))}
\end{equation}
Observe firstly that:
\[ \min(x+1, y) + \max(x,y-1) = x+y \]
Now, by equations \ref{red} and \ref{blue} we have:
\begin{align*}
\overline{A}^{\min}_{i,j} & = \max(\overline{B}_{i,j}, \overline{B}_{i+1,j+1} - 1) \\
\overline{C}^{\max}_{i,j} & = \min(\overline{B}_{i,j} , \overline{B}_{i-1,j-1} +1) 
\end{align*}
and so:
\begin{align*} 
\overline{C}^{\max}_{i,j} + \overline{A}^{\min}_{i-1,j-1} 
& = \min(\overline{B}_{i-1,j-1} + 1, \overline{B}_{i,j}) + \max(\overline{B}_{i-1,j-1}, \overline{B}_{i,j} - 1) \\
& = \overline{B}_{i,j} + \overline{B}_{i-1,j-1}
\end{align*}
This gives us the same power of $\lambda_{i,j}$ on both sides of equation \ref{fish} 

\subsubsection{Comparing power of $\mu_{i,n+1-j}$ on both sides}
In this section we show that:
\begin{equation}
t(G^{n-1}_\mu(A_{{\min}}))G^{n}_\mu(C_{{\max}}) = G^{{n-1}}_\mu(B) t(G^{{n}}_\mu(B))
\end{equation}
By equations \ref{red2} and \ref{blue2}
we have:
\begin{align*}
\underline{A}_{\min}^{i,j} & = \max(\overline{B}_{i,j+1}, \overline{B}_{i+1,j}-1) \\
\underline{C}_{\max}^{i,j} & = \min(\overline{B}_{i-1,j} +1, \overline{B}_{i,j-1})
\end{align*}
and so:
\begin{align*} 
\underline{A}_{\min}^{i-1,j-1} + \underline{C}_{\max}^{i,j}
& = \max(\overline{B}_{i-1,j}, \overline{B}_{i,j-1}-1) + \min(\overline{B}_{i-1,j}+1 , \overline{B}_{i,j-1}) \\
& = \overline{B}_{i-1,j} + \overline{B}_{i,j-1}
\end{align*}
This gives us the same power of $\mu_{i,n+1-j}$ on both sides of equation \ref{fish}. 

\subsubsection{Comparing power of $x_n[k]_{i,j}$ on both sides}

In this section we show that:
\begin{equation}
x_n[1]^{C^{\max}} s(x_n[1])^{A^{\min}} = (x_n[1] s(x_n[1]))^B
\end{equation}

The result follows from equation \ref{cat} together with equations \ref{seconddiff} and \ref{seconddiff2} for expressing the original alternating sign matrix in terms of the corner sum matrices.

\subsection{Special case}

In this section we shall prove corollary \ref{lambda:corollary}. Begin by observing that if:
\begin{align*} 
x_n[0]_{i,j} & = 1 \mbox{ for all } i,j = 1..(n+1) \\
x_n[1]_{i,j} & = M_{i,j} \mbox{ for all } i,j = 1..n \\
\end{align*}
then:
\begin{equation} 
x_n[n]_{1,1} = \sum_{|B| = n} M^B \left ( \sum_{ |A| = n-1} 
\frac{F_\lambda(B)}{s(F_\lambda(A))} \frac{G^{n}_\mu(B)}{t(G^{n}_\mu(A))} \right )
\end{equation}
The second sum is over all matrices $A$ which are interlacing with $B$. 
The above equation may be rewritten in the form:
\begin{equation}
x_n[n]_{1,1} = \sum_{|B| = n} M^B \left ( 
 \frac{F_\lambda(B)} {s(F_\lambda(A^{\min}))}
\frac{G^{n}_\mu(B)}{ t(G^{n}_\mu(A_{\max}))} 
\prod_{B_{i,j} = -1} 
\left (\mu_{i,n+1-j} + \lambda_{i,j} \right )
\right )
\end{equation}
By \ref{prop:inversion} we have:
\begin{equation}
 \frac{F_\lambda(B)} {s(F_\lambda(A^{\min}))} = \lambda^{\Inv(B)}
\end{equation}
while by proposition \ref{prop:dual-inversion} we have:
\begin{equation}
\frac{G^{n}_\mu(B)}{ t(G^{n}_\mu(A_{\max}))} = \mu^{\Dinv(B)'}
\end{equation}
The result follows.

\section{Conclusion}

We have given a multi-parameter generalization of the $\lambda$-determinant of Robbins and Rumsey \cite{robbins}.
Our result exhibits the Laurent phenomenon \cite{laurent} and generalizes a previous result given by diFrancesco \cite{cluster}.

\cleardoublepage

\part{Commutants in semicircular systems}
\chapter{Commutators in semicircular systems}

\section{Hilbert spaces}

In this section we recall some basic definitions from functional analysis \cite{kreyszig} which will perhaps be helpful to combinatorialists.
Professional analysts may skip this section.

\subsubsection{Inner product spaces}

An \emph{inner product space} is a complex vector space $V$ equipped with a map:
\[ \langle - | - \rangle : V \times V \to \mathbb{C} \]
which is linear in the second variable, and which satisfies:
\begin{itemize}
\item $\langle w | v \rangle = \overline{\langle v | w \rangle} \qquad \qquad \mbox{ (conjugate symmetry) }$ 
\item $\langle v | v \rangle > 0 \mbox{ for } v \neq 0 \qquad \mbox{ (positive definiteness)}$
\end{itemize}

A basis $\{e_1, e_2, \ldots \}$ of $V$ is said to be \emph{orthonormal} if
$\langle e_i | e_j \rangle = \delta_{ij}$.
A pair of bases $\{ w_1, w_2, \ldots \}$ and $\{r_1, r_2, \ldots \}$ are said to be \emph{dual} if
$\langle w_i | r_j \rangle = \delta_{i,j}$.

\subsubsection{Projections}\label{semisircular:gramm}
For $S$ a vector subspace of $V$, a \emph{projection operator} with image $S$ is a linear map $\pi_S : V \to V$
such that:
\begin{itemize}
\item $\img(\pi_S) = S$
\item $\pi_S[v] = v \mbox{ for all } v \in S$
\end{itemize}
If $\{w_i\}$ is a basis for $S$ and $\{r_i\}$ is the corresponding dual basis, then:
\[ \pi_S[v] = \sum_i \langle r_i | v \rangle w_i \]

In the finite dimensional case, let $P$ be the rectangular matrix $P_{ij} = \langle e_i | w_i \rangle$ and let $Q$ be the rectangular matrix $Q_{i,j} = \langle r_i | e_i \rangle$. One can show that:
\[ Q = (PP^*)^{-1} P\] 
This is known as the \emph{left pseudo inverse}. The projection operator may be expressed in the form:
\[ \pi_S = P^* Q \]

\subsubsection{Gramm matrices}
 
The matrices $PP^T$ and $QQ^T$ are known as the \emph{Gramm matrices}. We have:
\begin{align}
 (PP^*)_{i,j} & = \langle w_i \,|\, w_j \rangle \\
 (QQ^*)_{i,j} & = \langle r_i \,|\, r_j \rangle \\
\end{align} 
One can show that:
\begin{equation} 
(PP^*)^{-1} = QQ^* 
\end{equation}

\subsubsection{Topology}

A topology is a pair $(X,\Sigma)$ consisting of a set $X$ and a collection $\Sigma$ of subsets of $X$, called open sets, satisfying the following three axioms:
\begin{itemize}
\item The union of open sets is an open set.
\item The finite intersection of open sets is an open set.
\item $X$ and the empty set $\emptyset$ are open sets.
\end{itemize}
A function $f: X \to Y$ between topological spaces is said to be \emph{continuous} if the inverse image of every open set is open.

\subsubsection{Metric spaces}
A \emph{metric space} is a pair $(X,d)$ where $X$ is a set and $d$ is a function:
\[ d : X \times X \to \mathbb{R} \]
which satisfies:
\begin{itemize}
\item $d(x,y) \geq 0$
\item $d(x,y) = 0$ if and only if $x=y$
\item $d(x,y) = d(y,x)$
\item $d(x,z) \leq d(x,y) + d(y,x)$
\end{itemize}
For any inner product space we may define the \emph{norm} of a vector to be:
\[ ||v|| = \langle v | v \rangle \]
One can show that every inner product space is also a metric space with distance function:
\[ d(v,w) = ||v - w|| \]
Every metric space carries a natural topology generated by open sets of the form:
\[ B_r(v) = \{w \in X \,|\, d(w,v) < r\} \]

\subsubsection{Cauchy sequences}\label{semicircular:cauchy}

For any metric space $(X,d)$, a \emph{Cauchy sequence} is a sequence $x_1,x_2,x_3, \ldots $ of elements of $X$ such that for every $\epsilon$ there exists an $N$ such that: 
\[ d(x_n,x_m) \leq \epsilon \mbox{ for all }n,m > N\] 
An element $x$ of $X$ is said to be the \emph{limit} of the Cauchy sequence $x_1,x_2,x_3, \ldots$ if for every $\epsilon$ there exists an $N$ such that such that: 
\[d(x_m,x) < \epsilon \mbox{ for all } n > N\]

A \emph{complete} metric space $(X,d)$ is one for which every Cauchy sequence of $X$ converges to some limit in $X$. 
Every metric subspace admits a unique \emph{metric completion}

A subset $S$ of a metric space $(X,d)$ is \emph{closed} if for every $x \in X$ there exists a Cauchy sequence in $S$ whose limit is also in $S$. Every subspace $S$ of a complete metric space admits a unique \emph{closure} which which is denoted by $\overline{S}$. 

A subspace of a metric space is said to be \emph{dense} if its closure is the whole space.

\subsubsection{Hilbert spaces}\label{semicircular:hilbert}

A \emph{Hilbert space} is an inner product space which is also a complete metric space with respect to the norm induced by the inner product.
All finite dimensional inner product spaces are automatically Hilbert spaces.
All Hilbert spaces of the same dimension are isomorphic. 

We shall refer to the topology on a Hilbert space induced by the norm associated to the inner product as the \emph{natural topology} of the space.

\subsubsection{Hilbert basis}
Let $\{e_1,e_2,e_3, \ldots \}$ denote a sequence of linearly independent vectors in some infinite dimensional Hilbert space $H$. The \emph{span} of these vectors, consisting of all finite linear combinations of the form:
\[ \sum_{i=1}^n \alpha_i e_i \]
is an infinite dimensional vector space, but it is \emph{not} a Hilbert space, because it is \emph{not} closed under the natural topology on $H$.

The vector space spanned by the linearly independent vectors
$\{e_1, e_2, e_2 \ldots \}$ is said to be \emph{dense} in $H$ if for every $v \in H$ there exists
a sequence of vectors $\{v_n\}$ of the form:
\[v_n = a_{1,n} e_1 + a_{2,n} e_2 + \cdots a_{n,n} e_n\] 
such that:
\[ \lim_{n \to \infty} ||v - v_n|| = 0 \]
If this is the case we say that $\{e_1, e_2, \ldots \}$ is a \emph{Hamel basis} for $V$.
If in addition we have $\langle e_i, e_j \rangle = \delta_{i,j}$ then we say that $\{e_1, e_2, \ldots \}$ is a \emph{Hilbert basis}.

\subsubsection{Example of a dense subspace}\label{dense}

It is possible for a proper vector subspace of a Hilbert space to be dense in the whole space. We shall now give an example

Let $\{e_1, e_2, e_3, \ldots \}$ be a Hilbert basis for some infinite dimensional Hilbert space $H$. For each $k \geq 1$ let $V_k$ denote the finite dimensional subspace spanned by Hilbert $\{e_1, \ldots e_k\}$, let $d_k = e_k - e_{k+1}$ and let $W_k$ denote the finite dimensional subspace spanned by the vectors $\{ d_1, \ldots, d_k \}$. We shall denote by $W$ the infinite dimensional vector subspace spanned by \emph{all} the $d_k$. It is clear that $e_1 \not\in W$. We shall see however that $e_1$, lies in the \emph{closure} of $W$.

The subspace $W_k$ may be characterized as the orthogonal complement of the vector $e_1 + e_2 + \cdots e_{k+1}$ in the subspace $V_{k+1}$. 
For each $k$ let us define the vector:
\[ f_k = \sqrt{\frac{1}{k(k+1)}} (e_1 + e_2 + \cdots + e_k) - \sqrt{\frac{k}{k+1}} e_{k+1} \]
We have $||f_k|| = 1$.
We also have:
\[ \langle f_k, f_\ell \rangle = \delta_{k,\ell} \]
It follows that the $\{f_k\}$ form an orthonormal basis for the subspace $W$.

Let us define the following sequence of vectors:
\begin{align*}
q_n & = \sum_{k=1}^n \langle e_1, f_k \rangle f_k \\
& = \sum_{k=1}^n \sqrt{\frac{1}{k(k+1)}} f_k
\end{align*} 
The vector $q_n$ is the projection of the vector $e_1$ onto the subspace $W_n$. We have:
\begin{align*}
|| q_n || & = \sum_{k=1}^n \frac{1}{k(k+1)} = 1 - \frac{1}{n+1}
\end{align*}
In other words $\lim_{n \to \infty} ||q_n|| = 1$.
Since $||e_1 - q_n|| \leq ||e_1 || - ||q_n||$ it follows that $q_1, q_2, q_3, \ldots$ forms a Cauchy sequence in $W$ which converges towards $e_1$. That is $e_1$ lies in the closure of $W$. This implies that $W$ is a dense subspace of $H$. 

\subsubsection{Operator norm}

Let $U$ and $V$ be Hilbert spaces. A linear operator:
\[ L : V \to U \] 
is said to be \emph{bounded} if there exists some $M > 0$ such that for all $v \in V$ we have:
\[ ||L(v)||_U \leq M ||v||_V \]

The set of bounded linear operators on a Hilbert space $H$ is an algebra which is denoted by $B(H)$.
Note that $B(H)$ is \emph{not} a Hilbert space since it does not admit an inner product. We can nevertheless define a \emph{norm} on $B(H)$ via:
\begin{equation}
||L||_{B(H)} = \sup_{v \in V} \frac{||L(v)||_U}{||v||_V}
\end{equation}

This norm induces a topology on $B(H)$ which we shall refer to as the \emph{operator norm topology} of $B(H)$

\subsubsection{Weak operator topology}\label{wot}

The operator norm topology is not the only topology which can be defined on $B(H)$.

Let $B(H)$ denote the set of bounded linear operators on some Hilbert space $H$. For each $x,y \in H$ we can define a map:
\[ \psi(x,y) : B(H) \to \mathbb{R} \]
which is given by:
\[ \psi(x,y)[T] = \langle x \,|\,  T y \rangle \] 

The \emph{weak operator topology} is the weakest topology on $B(H)$, such that the map $\psi(x,y)$ is continuous for all $x,y \in H$.

\subsubsection{von Neumann algebras}
 
The \emph{adjoint} of a linear map $A : V \to V$ is the unique map $A^* : V \to V$ satisfying:
\[ \langle A^* v | u \rangle = \langle v | A u \rangle \]

A linear operator is said to be \emph{self-adjoint} if $A^* = A$. 
We say that a \emph{subalgebra} of $\mathfrak{A}$ of $B(H)$ is self adjoint if: 
\[a \in \mathfrak{A} \Longleftrightarrow a^* \in \mathfrak{A}\]
A \emph{von Neumann algebra} a self-adjoint algebra $\mathfrak{A}$ of bounded operators on a Hilbert space which contains the identity and which is closed under the weak operator topology.

For any algebra $\mathfrak{A}$ and any subset $S \subseteq \mathfrak{A}$ the \emph{commutant} of $S$ in $\mathfrak{A}$
is defined to be the subalgebra of $\mathfrak{A}$ given by:
\[ S' = \{ a \in A \,|\, sa = as \mbox{ for all } s \in S\} \]
The \emph{double commutant} of $S$ is defined to be the subalgebra of $\mathfrak{A}$ given by:
\[ S'' = \{ a \in A \,|\, ta = at \mbox{ for all } t \in S'\} \]
By von Neumann's \emph{double commutant theorem} a self-adjoint subalgebra $\mathfrak{A}$ of $B(H)$ is a von Neumann algebra if and only if:
\[ \mathfrak{A}'' = \mathfrak{A}\]

A von Neumann algebra is said to be a \emph{factor} if $\mathfrak{A}' \cap \mathfrak{A} = \mathbb{C}.1$. That is to say, a factor is a von Neumann algebra with trivial center.

\subsubsection{Random variables}

For our purposes a \emph{random variable} is simply a self-adjoint operator $X$ acting on some Hilbert space space $H$ with inner product $\langle - | - \rangle $ and distinguished vector $\Omega$ satisfying:
\[ \langle \Omega \,|\, \Omega \rangle = 1 \]
In the language of quantum mechanics, random variables are referred to as \emph{observables}.
The \emph{moments} of a random variable $X$ are defined to be:
\[ \mu_n = \mathbb{E}(X^n) = \langle \Omega \,|\, X^n \Omega \rangle \] 
Two operators are are said to be \emph{identically distributed} if they have the same moments.
The distinguished vector $\Omega$ is often referred to as a \emph{state}. 

When working with an \emph{algebra of random variables} $\mathfrak{A}$ with the property that $\mathbb{E}(XY) = \mathbb{E}(YX)$ for all $X,Y \in \mathbb{A}$ the state $\Omega$ is said to be \emph{tracial}.

\section{Chebyschev polynomials}\label{chebyshev}
In what follows we shall make extensive use of the operator ``multiplication by $x$'' acting on the Hilbert space $\mathfrak{h}$.

\subsubsection{Measurable functions}

A \emph{measure space} is a pair $(X,M)$ where $X$ is a set and $M$ is a set of subsets of $X$ the property that if $A \in M$ then $X \backslash A \in M$, and if $A_1, A_2, A_3, \ldots$ is a sequence of elements of $M$ then $\cup_{i} A_i$ is an element of $M$. We also require that $\emptyset \in M$.

The \emph{Borel $\sigma$-algebra} of $\mathbb{R}$ is the smallest $\sigma$ algebra on $\mathbb{R}$ which is generated by the open sets of $\mathbb{R}$. Similarly, the \emph{Borel $\sigma$-algebra} of $\mathbb{C}$ is the smallest $\sigma$ algebra on $\mathbb{C}$ which contains all the open sets of $\mathbb{C}$.

A function between measure spaces is said to be \emph{measurable} if the pre-image of every measurable set is measurable.

\subsubsection{Square integrable functions}

Let $\mathfrak{H}$ denote the set of all measurable functions. 
$f : [-2,2] \to \mathbb{C}$, modulo the equivalence relation in which two functions are equal if they differ only on a set of measure zero,
which satisfy:
\[ \int_{-2}^2 |f(x)|^2 \sqrt{4 - x^2}\,dx < \infty \]
The inner product is given by
\[ \langle f \,|\, g \rangle = \int_{-2}^2 f(x) \overline{g(x)} \sqrt{4-x^2}\,dx \]
The space of polynomials $\mathbb{C}[x]$ sits inside of $\mathfrak{H}$ as a dense subspace.
The closure of $\mathbb{C}[x]$ with respect to the natural topology on $\mathfrak{H}$ is the whole space $\mathfrak{H}$.
In other words, the monomials $\{1, x, x^2, x^3, \cdots \}$ form a Hamel basis for $\mathfrak{h}$. They are not orthogonal however. Let us define:
\[ P(x,z) = \frac{1}{1-2tx + t^2} = \sum_n P_n(x) t^n \]
One may check that:
\[ \langle P_n(x), P_m(x) \rangle = \pi \delta_{n,m} \]
The polynomials $\{P_n(x)\}$ are known as the \emph{Chebyshev polynomials of the second kind}.
They form a Hilbert basis for the space $\mathfrak{H}$.


\subsubsection{Three term recurrence}\label{multiplication-x}

The Chebyshev polynomials satisfy the following three term recurrence:
\begin{equation} \label{three-term-recurrence}
x P_n(x) = P_{n+1}(x) + P_{n-1}(x) 
\end{equation}
The operator ``multiplication by $x$'' is bounded and may be thought of as living in the space $B(\mathfrak{H})$.

\subsubsection{Embedding}

Let us consider the algebra $\Poly(x)$ of polynomials in the operator $x$. This is a subalgebra
of $B(H)$ of bounded linear operators on $\mathfrak{h}$. 
There is a natural embedding of $\Poly(x)$ into $\mathfrak{H}$ which sends the operator $p(x) \in B(H)$ to the vector $v \in \mathfrak{H}$ given by:
\begin{equation}\label{embedding-chebyshev}
v = p(x)[1] 
\end{equation}

\subsubsection{Topological considerations}

By the \emph{Stone--Weierstrass theorem} the closure of $\Poly(x)$ with respect to the \emph{operator norm topology} on $B(H)$
is isomorphic to the space of continuous functions in $x$.
Note that since not every square integrable function is continuous, the closure of 
$\Poly(x)$ with respect to the operator norm is \emph{not} the entire space $\mathfrak{H}$.

A measurable function:
\[ f :\mathbb{R} \to \mathbb{C} \] 
is said to be \emph{bounded} if there exists some $K > 0$ such that $|f(x)| < K$ for all $x \in \mathbb{C}$.
One can show that the closure of $\Poly(x)$ under the \emph{weak operator topology} is isomorphic to the space of all \emph{bounded} measurable functions in $x$. Again, since not every square integrable function is bounded, this is not the entire Hilbert space $\mathfrak{H}$.

In general, the weaker the topology, the less open sets and the more continuous functions. The closure of a subspace with respect to a weaker topology is larger than or equal to the closure of a subspace with respect to a stronger topology.
As the name suggests, the weak operator topology is weaker than the operator norm topology.

\subsubsection{Moments}\label{moments}

We remark that the operator ``multiplication by $x$'' is self-adjoint, and may be thus thought of as a random variable. We will now study its moments. One may calculate:
\begin{equation}\label{catalan3}
\int_{-2}^2 x^{2n} \sqrt{4 - x^2} dx = \frac{1}{n+1} \left ( \begin{matrix} 2n \\ n\end{matrix}\right )
\end{equation}
while:
\begin{equation}\label{catalan4}
\int_{-2}^2 x^{2n+1} \sqrt{4 - x^2} dx = 0
\end{equation}
The numbers:
\begin{equation}
C_n = \frac{1}{n+1} \left ( \begin{matrix} 2n \\ n\end{matrix}\right )
\end{equation}
are known as the \emph{Catalan numbers}.

\section{Semicircular Systems}\label{sec:semicircular}

The \emph{semicircular system} was first introduced by Voiculescu to study free group von Neumann algebras \cite{semicircular}.
The ``multiplication by $x$'' operator in the Chebyshev basis arises naturally in the the context of Wigner's semi-circular law for random matrices. In free probability theory Wigner's semicircular distribution plays a role analogous to the \emph{Gaussian distribution} in classical probability theory.

\subsubsection{Fock space}

Fix a finite dimensional inner product space $H$ with inner product $\langle - | - \rangle$ and orthonormal basis:
\[\{e_1 , e_2 , \cdots , e_k \}\] 
Let us define:
\begin{equation}
T(H) = \mathbb{C} \Omega \oplus_{n\geq1} H^{\otimes n}
\end{equation}

As a vector space $T(H)$ is naturally graded. The dimension of the graded component of degree $d$ is equal to $n^d$.
A basis for the graded component of degree $d$ is given by:
\begin{equation}\label{fock:basis}
 \{e_{i_1} \otimes e_{i_2} \otimes \cdots \otimes e_{i_d}\} \,\,|\,\, i_1,i_2, \ldots i_d \in \{1,2,\ldots k\}\} 
\end{equation}

The inner product of $H$ lifts to an inner product of $T(H)$ in a natural way. If $v$ and $w$ are homogeneous, but not of the same degree, then $\langle v | w \rangle = 0$. Otherwise, if $v = e_{i_1} \otimes e_{i_2} \otimes \cdots \otimes e_{i_d}$
and $w = e_{j_1} \otimes e_{j_2} \otimes \cdots \otimes e_{j_d}$ then:
\[ \langle v | w \rangle = \delta_{i_1,j_1} \delta_{i_2,j_2} \cdots \delta_{i_d,j_d} \]
In otherwords, the inner product on $T(V)$ is defined in such a way that the basis described in equation \ref{fock:basis} is orthonormal.
The \emph{Fock space} $\mathfrak{F}$ of $H$ is defined to be the \emph{metric completion} of $T(H)$.

\subsubsection{Special subspaces}

Note that since all Hilbert spaces of the same dimension are isomorphic, for any $i$ the closure with respect to the natural topology on $\mathfrak{F}$ of the subspace generated by 
the  $\{ e_i^{\otimes n}\}_{n \geq 0}$ for any given $i$ is isomorphic to the vector space $\mathfrak{h}$ of square integrable functions discussed in section \ref{chebyshev}.

\subsubsection{Creation and annihilation operators}

For each $i$ let $c_i$ be the \emph{creation operator}:
\begin{equation}
c_i [v ] = e_i \otimes v
\end{equation}
The adjoint operator $c_i^*$ is known as the \emph{annihilation operator}. It acts via:
\begin{align}
c_i^* [\Omega] & = 0 \\
c_i^* [(v_1 \otimes v_2 \otimes \cdots \otimes v_n )] & = 
\langle e_i |v_1 \rangle (v_2 \otimes \cdots \otimes v_n)
\end{align}
The operators $A_i = c_i + c_i^*$ are both self-adjoint and bounded.

\subsubsection{Moments}

One can easily convince
oneself that:
\begin{equation}\label{catalan1}
\langle \Omega \,|\, A_i^{2n} \Omega \rangle = \frac{1}{n+1} \left ( \begin{matrix} n \\ 2n\end{matrix} \right )
\end{equation}
while:
\begin{equation}\label{catalan2}
\langle \Omega \,|\, A_i^{2n+1} \Omega \rangle = 0 
\end{equation}
Comparing equations \ref{catalan1} and \ref{catalan2} with equations \ref{catalan3} and \ref{catalan4} one can conclude that the operator $A_i$ acts on the natural closure of the vector space spanned by $\{e_i^{\otimes n}\}_{n \geq 0}$ in the same way as the operator ``multiplication by $x$'' acts on the space $\mathfrak{H}$ of section \ref{chebyshev}.
 
In particular, if $P_n(x)$ is the $n$-th Chebyshev polynomial then we have:
\begin{equation}\label{one-variable}
P_n(A_i)[\Omega] = e_i^{\otimes n}
\end{equation}

\subsubsection{von Neumann algebra}
We shall use the notation $\mathfrak{A}_n$ to denote the \emph{Von Neumann algebra} generated by the semi-circular elements $A_1, A_2, \ldots A_n$. Thinking of each of the operators $A_i$ as a \emph{random variable} we may define the following \emph{state}
\begin{equation}
\tau(X) = \langle \Omega \,|\, X \Omega \rangle
\end{equation}
It is not hard to show that $\tau$ is in fact a \emph{tracial state}:
\[
\tau(XY) = \tau(YX)
\]

\subsubsection{Embedding}

It is not hard to convince oneself that for all $\{i_1 , i_2 , \ldots , i_r \} \in \mathbb{N}_{\geq 1}$ and for all  $j_1 , j_2 , \ldots j_r \in \{1, 2, \ldots n\}$ 
such that $j_\ell \neq j_{\ell + 1}$ we have: 
\begin{equation}
P_{j_1}(A_{i_1})P_{j_2}(A_{i_2})\cdots P_{j_k}(A_{i_k})[\Omega] = e_{i_1}^{\otimes j_1}e_{i_2}^{\otimes j_2} \cdots e_{i_k}^{\otimes j_k}
\end{equation} 
In other words, we have an embedding of $\mathfrak{A}_n \subseteq B(\mathfrak{F})$ into $\mathfrak{F}$ which is given by:
\begin{equation}
X \mapsto X[\Omega]
\end{equation}

\subsubsection{Topological considerations}

Let us remark that there is no natural inner product on $B(\mathfrak{F})$. Nevertheless, the embedding of $\mathfrak{A}_n$ into $\mathfrak{F}$ allows us to \emph{induce} an inner product on $\mathfrak{A}_n$.
Even with this induced inner product, the vector space $\mathfrak{A}_n$ is \emph{not} a Hilbert space, since it is \emph{not} closed under the topology induced by the inner product.

\hfill \break

\section{Main result}\label{sec:semicircular:main}

\subsection{Statement}

For each $i = 1,2, \ldots n$, let $\Poly(A_i)$ denote the subalgebra of $B(\mathfrak{F})$ generated by $A_i$ and let $\alg(A_i)$ denote the closure of $\Poly(A_i)$ in the weak operator topology of $B(\mathfrak{F})$.

Similarly let
$\comm(A_i)$ the sub-algebra of $\mathfrak{A}_n$ consisting of all elements which commute with $A_i$. One can show that $\comm(A_i)$ is already closed under the weak operator topology. 

The goal of the next two sections is to prove the following result:
\begin{theorem}\label{semicircular:main}
\begin{equation}
\boxed{\alg(A_i) = \comm(A_i)}
\end{equation}
\end{theorem}

As a consequence of Theorem \ref{semicircular:main} we have:

\begin{theorem} \label{semicircular:factor}
For $n \geq 2$ the center of $\mathfrak{A}_n$ is trivial.
\begin{proof}
Suppose that $x$ lies in the center of $\mathfrak{A}_n$. It follows that $x$ commutes with $A_1$, and $x$ commutes with $A_2$. Now, by Theorem \ref{semicircular:main} we have $x \in \alg(A_1)$ and $x \in \alg(A_2)$. 
But:
\[ \alg(A_i)\cap \alg(A_j) = \mathbb{C}.1 \]
The result follows.
\end{proof}
\end{theorem}

Note that in the language of Von Neumann algebras Theorem \ref{semicircular:factor} implies that $\mathfrak{A}_n$ is a \emph{factor}. This result is already known \cite{semicircular}. Our proof is still interesting however as it yields a number of explicit projection formulae (see section \ref{semicircular:matrix}). 


\subsection{Topological considerations}

We shall use the notation $\mathbb{A}_n$ to denote the closure of $\mathfrak{A}_n$ under the natural topology of $\mathfrak{F}$.
Furthermore we shall use the notation $\mathbb{A}_n^d$ to denote the intersection of $\mathbb{A}^n$ with the finite dimensional subspace of $\mathfrak{F}$ consisting of elements of degree  $d$.

\begin{lemma}
For any vector subspace $V$ of some Hilbert space $H$, the orthogonal complement $V^\perp$ is closed.
\end{lemma}

Let $\overline{\alg(A_i)}$ denote the closure of $\alg(A_i)$ with respect to the topology of $\mathfrak{F}$.

\begin{lemma}
\begin{equation}
\overline{\alg(A_i)} \cap \mathfrak{A}_n = \alg(A_i)
\end{equation}
\begin{proof}
Let $V_i$ denote the vector space spanned by $\{e_i^{\otimes n}\}_{n \geq 0}$.
Under the embedding of $\mathfrak{A}_n$ into $\mathfrak{F}$ we have:
\begin{equation}
\Poly(A_i) = V_i
\end{equation}
Let  $\mathfrak{h}_i$ denote the closure of $V_i$ under the topology of $\mathfrak{F}$. We have
\begin{equation}
\overline{\alg(A_i)} = \mathfrak{h}_i 
\end{equation}
But we also have that
$\mathfrak{A}_n \cap \mathfrak{h}_i  = \alg(A_i) \cap \mathfrak{h}_i$.
Referring back to remarks made in section \ref{chebyshev}, every bounded measurable function is also square integrable, that is to say:
\[  \alg(A_i) \subset \mathfrak{h}_i\]
 and so:
\[  \alg(A_i) \cap \mathfrak{h}_i = \alg(A_i)\]
The result follows.
\end{proof}
\end{lemma}

Let $\overline{\comm(A_i)}$ denote the closure of $\comm(A_i)$ with respect to the topology of $\mathfrak{F}$. 

\begin{lemma}
\begin{equation}
\overline{\comm(A_i)} \cap \mathfrak{A}_n = \comm(A_i)
\end{equation}
\begin{proof}
We have:
\begin{align}
\overline{\comm(A_i)} 
& =  \overline{\{ x \in \mathfrak{A}_n \,|\, x A_i = A_i x \}} \\
\end{align}
and so:
\begin{align}
\overline{\comm(A_i)} \cap \mathfrak{A}_n
& =  \overline{\{ x \in \mathfrak{A}_n \,|\, x A_i = A_i x \}} \cap \mathfrak{A}_n \\
& = \{ x \in \mathfrak{A}_n \,|\, x A_i = A_i x \} \\
& = \comm(A_i)
\end{align}
The result follows.
\end{proof}
\end{lemma}

\subsubsection{Strategy}

We shall begin by proving the following:
\begin{proposition}\label{fucked-topology}
\begin{equation}
\boxed{\comm(A_i)^\perp = \alg(A_i)^\perp}
\end{equation}
\end{proposition}
This implies in particular that:
\begin{equation}
\overline{\alg(A_i)} = \overline{\comm(A_i)}
\end{equation}
Taking the intersection of both sides with $\mathfrak{A}_n$ we obtain Theorem \ref{semicircular:main}

Our goal is now to prove proposition \ref{fucked-topology}. It is clear that we have:
\begin{equation}\label{semicircular:perp}
\comm(A_i)^\perp \subseteq \alg(A_i)^\perp
\end{equation}
We need only to prove the reverse inclusion of equation \ref{semicircular:perp}. 
\begin{equation}\label{semicircular:toprove}
\boxed{\alg(A_i)^\perp \subseteq \comm(A_i)^\perp}
\end{equation}

\subsection{Commutators}

Recall that the \emph{commutator} of two operators $X$ and $Y$ is given by:
\begin{equation}
[X,Y] = XY - YX
\end{equation}

\begin{lemma}\label{lemma:trace}
For all $x,y \in \mathfrak{A}$ we have:
\begin{equation}
\tau(x[A_i,y]) = \tau(y[x,A_i])
\end{equation}
\begin{proof}
By the linearity and the traciality of the state $\tau$, we have:
\begin{align*}
\tau(x[A_i,y]) & = \tau(xA_iy) - \tau(xyA_i) \\
& = \tau(yxA_i) - \tau(yA_ix) \\
& = \tau(y[x,A_i])
\end{align*}
\end{proof}
\end{lemma}
Our plan is to study the following family of finite dimensional vector subspaces:
\begin{equation}
\boxed{\{ [A_i,y], y \in \mathbb{A}_n^d \}}
\end{equation}

\begin{lemma}\label{semicircular:not-closed}
For all $d \geq 1$ we have the following inclusion:
\begin{equation} 
 \{ [A_i,y], y \in \mathbb{A}_n^d \} \subseteq \comm(A_i)^\perp
\end{equation}
\begin{proof}
Choose any $y \in \mathbb{A}_n^d$. By lemma \ref{lemma:trace} we have that if $x \in \comm(A_i)$ then: 
\[\tau(x[A_i,y]) = \tau(y[x,A_i]) = 0 \]
In other words $[A_i,y]$ lies in the orthogonal complement of $\comm(A_i)$.
\end{proof}
\end{lemma}

\begin{lemma}
For all $d \geq 1$
\begin{equation*}
\alg(A_i)^\perp \not\subseteq \{ [A_i,y], y \in \mathbb{A}_n^d \} 
\end{equation*}
\begin{proof}
For any $j \neq i$ we have $A_j \in \alg(A_i)^\perp$ but we do not have $A_j \in \{ [A_i,y], y \in \mathbb{A}_n^d \}$ for any $d$.
\end{proof}
\end{lemma}

Our goal in the next section will be to show that:
\begin{proposition}\label{prop:center}
\begin{equation}
\alg(A_i)^\perp \subseteq \overline{\bigcup_{d \geq 1}\{ [A_i,y], y \in \mathbb{A}^d_n \}} 
\end{equation}
\end{proposition}
That is to say, for any $B \in \alg(A_i)^\perp$ there exists a sequence of operators $\{y_1,y_2,y_3, \ldots\}$
with: $y_m \in  \mathbb{A}^d_n $ for some $d$ which depends on $m$.
such that, as a vector:
\begin{equation}
B = \lim_{m \to \infty}[A_i,y_m] 
\end{equation}
Note that this is to be understood as convergence in the metric space induced by the inner product of the Fock space $\mathfrak{F}$.

Proposition \ref{prop:center} and Lemma \ref{semicircular:not-closed} will together give us equation \ref{semicircular:toprove}, thus proving Theorem \ref{semicircular:main}.

\subsection{Projection formula}

In this section we prove proposition \ref{prop:center} for $i=1$. The general result follows by symmetry.

\subsubsection{Notation}

Let $S$ be any non-trivial element of $\mathfrak{A}$ of the form:
\[ S = P_{j_1}(A_{i_1})P_{j_2}(A_{i_2})\cdots P_{j_r}(A_{i_r}) \] 
such that $j_1,j_2,\ldots j_r \geq 1$ and $j_\ell \neq j_{\ell + 1}$ as well as $i_1 \neq 1 \neq i_r$.
We have:
\[ S \in \alg(A_1)^\perp \]
Next, for each $n,k \geq 1$ let:
\[ S_{n,k} = P_n(A_1) S P_k(A_1) \]
Again we have:
\[ S_{n,k} \in \alg(A_1)^\perp \]
In fact, the collection of all operators of the form $S_{n,k}$ for all possible choices of $S$ form a Hilbert basis for $\alg(A_1)^\perp$ 

\subsubsection{Simplification}

\begin{lemma}\label{semicircular:commutator-simplification}
If $a$ and $c$ both commute with $x$ then:
\[ [x,abc] = a[x,b]c \]
\begin{proof}
\begin{align*}
[x,abc] & = xabc - abcx \\
& = axbc - abxc \\
& = a[x,b]c
\end{align*}
\end{proof}
\end{lemma}

\begin{proposition}\label{semicircular:first-row}

If:
\[ S \in \overline{\bigcup_{d \geq 1}\{ [A_i,y], y \in \mathbb{A}^d_n \}} \]
Then, for all $n,k \geq 0$ we have:
\[ S_{n,k} \in \overline{\bigcup_{d \geq 1}\{ [A_i,y], y \in \mathbb{A}^d_n \}} \]
\begin{proof}
Suppose that we already have a sequence of operators $\{y_1,y_2,y_3, \ldots \}$ with: 
$y_m \in \mathbb{A}_n $
such that:
\begin{equation*}
S = \lim_{m \to \infty}[A_1,y_m] 
\end{equation*}
It follows that:
\begin{align*}
S_{n,k} & = \lim_{m \to \infty}P_n(A_1)[A_1,y_m]P_k(A_1) \\
& = \lim_{m \to \infty}[A_1,P_n(A_1)\,y_m\,P_k(A_1)] \\
& = \lim_{m \to \infty}[A_1,z_m]
\end{align*}
where $z_m = P_n(A_1)\,y_m\,P_k(A_1)$. The result follows
\end{proof}
\end{proposition}

\subsubsection{Projection formula for arbitrary $S$:}
\begin{equation}\label{semicircular:projection}
\boxed{S = \lim_{m \to \infty} \frac{1}{2m+3} \sum_{k=0}^{m-1} (m-k) ([S_{k,k+1},A_1] - [S_{k+1,k},A_1]) }
\end{equation}

In order to prove this formula we will need some more notation.
For fixed $S$, let $W$ denote the vector subspace spanned by all vectors of the form:
\begin{align*}
v_{n,k} & = [ S_{n,k}, A_1] \\
& = S_{n,k+1} - S_{n+1,k} + S_{n,k-1} - S_{n-1,k}
\end{align*}
The second equality above is a consequence of the three-term recurrence for Chebyshev polynomials (see equation \ref{three-term-recurrence}).
It will be convenient to make the following change of notation:
\begin{align*}
e_{n,k} & = \frac{1}{\sqrt{2}} (S_{n,k} + S_{k,n}) \mbox { for } n < k \\
e_{k,k} & = S_{k,k} \\
w_{n,k} & = \frac{1}{\sqrt{2}}(v_{n,k} - v_{k,n}) \mbox { for } n < k
\end{align*}
We may now rewrite equation \ref{semicircular:projection} as follows:
\begin{align*}
s_m & = \frac{1}{2m+3} \sum_{k=0}^{m-1} (m-k) ([S_{k,k+1},A_1] - [S_{k+1,k},A_1]) \\
& = \frac{\sqrt{2}}{2m+3} \sum_{k=0}^{m-1} (m-k) w_{k,k+1} \\
& = \frac{1}{2m+3} \sum_{k=0}^{m-1} (m-k) (\sqrt{2}e_{k,k+2} - 2 e_{k+1,k+1} +  2e_{k,k} -\sqrt{2} e_{k-1,k+1}) \\
& = \frac{1}{2m+3} \left (  2m e_{0,0} + 
\sqrt{2} \sum_{k=0}^{m-1} e_{k,k+2} - 2 \sum_{k=1}^m e_{k,k} \right )
\end{align*}

\newpage

Suppose that $S$ is an element of degree $t = j_1 + j_2 + \cdots + j_r$. For each $m \geq 0$ let us define:
\begin{equation}
E_S^{2m} \subset \bigcup_{d=0}^m \mathbb{A}_n^{2d+t}
\end{equation}
to be the vector space spanned by the $e_{n,k}$ for $n+k$ even and $n+k \leq 2m$.
Similarly, let us define, for each $m \geq 0$ the vector subspace:
\begin{equation}
W_S^{2m-1} \subset E_S^{2m}
\end{equation}
to be the vector space spanned by the $w_{n,k}$ for $n+k$ odd, $n<k$ and $n+k \leq 2m-1$.

\begin{proposition}
The vector $s_m$ is the orthogonal projection of the vector $S$ onto the finite dimensional subspace $W_S^{2m-1}$.
\begin{proof}
We shall show that the vector $(S - s_m)$ lies in the orthogonal complement of $W_S^{2m-1}$ with respect to $E_S^{2m}$.
We begin by observing that if $r > 2$ then: 
\[ \langle S, e_{k,k+r} \rangle  = \langle s_m , e_{k,k+r} \rangle= 0 \] Thus we may restrict ourselves even further to the subspace: 

\begin{equation}
F \subseteq E_S^{2m}
\end{equation}
spanned by vectors of the form $e_{k,k}$ or $e_{k,k+2}$ with degree less than or equal to $2m$.
Next, if $j-i \geq 4$ then: 
\[\langle e_{k,k}, w_{i,j}  \rangle = \langle e_{k,k+2}, w_{i,j} \rangle = 0\]

Let $U$ denote the vector subspace of $W_S^{2m-1}$ spanned by vectors of the form $w_{k,k+1}$ for $k < m$ 
and $w_{k,k+3}$ for $k < m-1$. We must show that for all $u \in U$ we have:
\[ \langle S - S_m | u \rangle = 0 \]

Now:
\[ (2m+3) \langle S | w_{0,1} \rangle = \langle e_{0,0} | \sqrt{2} e_{0,0} \rangle = -(2m+3)\sqrt{2} \]
While:
\begin{align*}
(2m+3) \langle s_m | w_{0,1} \rangle  &= 
\langle 2m e_{0,0} | \sqrt{2} e_{0,0} \rangle
+ \langle \sqrt{2} e_{0,2}, e_{0,2} \rangle 
+ \langle -2 e_{1,1} | -\sqrt{2} e_{1,1} \rangle \\
& = (2m+3)\sqrt{2}
\end{align*}
Next, for all $i > 0$ we have:
\[ \langle S | w_{i,i+1} \rangle = 0 \]
and:
\begin{align*}
(2m+3) \langle s_m | w_{i,i+1} \rangle  & =
\langle \sqrt{2} e_{i,i+2}, e_{i,i+2} \rangle 
- \langle 2 e_{i,i}, \sqrt{2} e_{i,i} \rangle \\
& \phantom{====} 
+ \langle 2 e_{i+1,i+1}, \sqrt{2} e_{i+1,i+1} \rangle 
- \langle \sqrt{2} e_{i-1,i+1} | e_{i+1,i+1} \rangle \\
& = \sqrt{2} - 2 \sqrt{2} + 2 \sqrt{2} - \sqrt{2} \\
& = 0
\end{align*}
Finally, for all $i \geq 0$ we have:
\[ \langle S | w_{i,i+3} \rangle = 0 \]
While:
\begin{align*}
(2m+3) \langle s_m | w_{i,i+1} \rangle & =
\langle \sqrt{2} e_{i,i+2}, e_{i,i+2} \rangle 
- \langle 2 e_{i,i}, \sqrt{2} e_{i,i} \rangle \\
& \phantom{====} 
+ \langle 2 e_{i+1,i+1}, \sqrt{2} e_{i+1,i+1} \rangle 
- \langle \sqrt{2} e_{i-1,i+1} | e_{i+1,i+1} \rangle \\
& = \sqrt{2} - 2 \sqrt{2} + 2 \sqrt{2} - \sqrt{2} \\
& = 0
\end{align*}
The result follows.

\end{proof}
\end{proposition}

It remains to verify that:
\begin{align*}
\lim_{m \to \infty}|| s_m  || & = \frac{1}{(2m+3)^2} \left ( 4m^2 + 2m + 4m \right ) \\
& = \lim_{m \to \infty}\frac{2m}{2m+3}\\
& = 1
\end{align*}
We have proven equation \ref{semicircular:projection} and shown that:
\[ 
S \in  \overline{\bigcup_{d \geq 1}\{ [A_i,y], y \in \mathbb{A}^d_n \}}
\]
Theorems \ref{semicircular:main} and \ref{semicircular:factor} follow as a consequence of proposition \ref{semicircular:first-row}.

\subsection{An interesting matrix}\label{semicircular:matrix}

Let us fix the following order on the wet of vectors $w_{n,k}$ for $n+k$ odd and $n < k$ and $n+k \leq 2m-1$:
\[
w_{0,1}, w_{0,3}, w_{1,2}, w_{0,5}, w_{1,4}, w_{2,3},w_{0,7}, w_{1,6}, w_{2,5}, w_{3,4}, w_{0,9}, w_{1,8},w_{2,7}, w_{3,6}, w_{4,5}
\]

For each $m \geq 1$ let $A_m$ denote the \emph{Gramm matrix} of inner products $\langle w_{n,k}, w_{n',k'} \rangle$ for $n+k$ odd and $n+k \leq 2m-1$.
Let $B_m = (2m+3) A_m^{-1}$. Here is the matrix $B_5$:

\[
B_5 = 
\left [
\begin{array}{c|cc|ccc|cccc|ccccc}
5   & \color{gray}0 & 4   & \color{gray}0 & \color{gray}0 & 3   & \color{gray}0 & \color{gray}0 & \color{gray}0 & 2   & \color{gray}0 & \color{gray}0 & \color{gray}0 & \color{gray}0 & 1 \\

\hline

\color{gray}0   & \color{red}12& \color{red}4   & \color{gray}0 & \color{red}9 & \color{red}3   
& \color{gray}0 & \color{gray}0 & \color{red}6 & \color{red}2   
& \color{gray}0 & \color{gray}0 & \color{gray}0 & \color{red}3 & \color{red}1 \\

4   & \color{red}4 & 8   & \color{gray}0 & \color{red}3 & 6   
& \color{gray}0 & \color{gray}0 & \color{red}2 & 4   
& \color{gray}0 & \color{gray}0 & \color{gray}0 & \color{red}1 & 2 \\

\hline

\color{gray}0   & \color{gray}0 & \color{gray}0   & \color{blue}15& \color{Emerald}9 & \color{Emerald}3   
& \color{gray}0 & \color{blue} 10& \color{Emerald}6 & \color{Emerald}2   
& \color{gray}0 & \color{gray}0 & \color{blue}5 & \color{Emerald}3 & \color{Emerald}1 \\

\color{gray}0   & \color{red}9 & \color{red}3   & \color{Emerald}9 & \color{red}18& \color{red}6   
& \color{gray}0 & \color{Emerald}6 & \color{red}12& \color{red}4   
& \color{gray}0 & \color{gray}0 & \color{Emerald}3 & \color{red}6 & \color{red}2\\

3   & \color{red}3 & 6   & \color{Emerald}3 & \color{red}6 & 9   
& \color{gray}0 & \color{Emerald}2 & \color{red}4 & 6   
& \color{gray}0 & \color{gray}0 & \color{Emerald}1 & \color{red}2 & 3\\

\hline

\color{gray} 0   & \color{gray}0 & \color{gray}0   & \color{gray}0 & \color{gray}0 & \color{gray}0   
& \color{purple}14& \color{green}10& \color{green}6 & \color{green}2   
& \color{gray}0 & \color{purple}7 & \color{green}5 & \color{green}3 & \color{green}1 \\

\color{gray}0   & \color{gray}0 & \color{gray}0   & \color{blue}10& \color{Emerald}6 & \color{Emerald}2   
& \color{green}10& 
\color{blue} 20& \color{Emerald}12& \color{Emerald}4   
& \color{gray}0 & \color{green}5 & \color{blue}10& \color{Emerald}6 & \color{Emerald}2 \\

\color{gray}0   & \color{red}6 & \color{red}2   & \color{Emerald}6 & \color{red}12& \color{red}4  
 & \color{green}6 & \color{Emerald}12& \color{red}18& \color{red}6   
& \color{gray}0 & \color{green}3 & \color{Emerald}6 & \color{red}9 & \color{red}3 \\

2   & \color{red}2 & 4   & \color{Emerald}2 & \color{red}4 & 6   
& \color{green}2 & \color{Emerald}4 & \color{red}6 & 8   
& \color{gray}0 & \color{green}1 & \color{Emerald}2 & \color{red}3 & 4 \\

\hline

\color{gray}0   & \color{gray}0 & \color{gray}0   & \color{gray}0 & \color{gray}0 & \color{gray}0   
& \color{gray}0 & \color{gray}0 & \color{gray}0 & \color{gray}0   
& \color{Melon}9 & \color{yellow}7 & \color{yellow}5 & \color{yellow}3 & \color{yellow}1 \\

\color{gray}0   & \color{gray}0 & \color{gray}0   & \color{gray}0 & \color{gray}0 & \color{gray}0   
& \color{purple}7 & \color{green}5 & \color{green}3 & \color{green}1   
& \color{yellow}7 & \color{purple}14& \color{green}10& \color{green}6 & \color{green}2 \\

\color{gray}0   & \color{gray}0 & \color{gray}0   & \color{blue} 5 & \color{Emerald}3 & \color{Emerald}1   
& \color{green}5 & \color{blue}10& \color{Emerald}6 & \color{Emerald}2   
& \color{yellow}5 & \color{green}10& \color{blue} 15& \color{Emerald}9 & \color{Emerald}3 \\

\color{gray}0   & \color{red}3 & \color{red}1   & \color{Emerald}3 & \color{red}6 & \color{red}2   
& \color{green}3 & \color{Emerald}6 & \color{red}9 & \color{red}3   
& \color{yellow}3 & \color{green}6 & \color{Emerald}9 & \color{red}12& \color{red}4 \\

1   & \color{red}1 & 2   & \color{Emerald}1 & \color{red}2 & 3   
& \color{green}1 & \color{Emerald}2 & \color{red}3 & 4   
& \color{yellow}1 & \color{green}2 & \color{Emerald}3 & \color{red}4 &  5 \\

\end{array}
\right ]
\]

We have used some colour to highlight the self-similar structure of this matrix. 

We may block diagonalize the matrix $B_5$ by conjugating by the following matrix:
\[
U = \left [ \begin{array}{c|cc|ccc|cccc|ccccc}
1   & \color{gray}0 & \color{gray}0   & \color{gray}0 & \color{gray}0 & \color{gray}0   & \color{gray}0 & \color{gray}0 & \color{gray}0 & \color{gray}0   & \color{gray}0 & \color{gray}0 & \color{gray}0 & \color{gray}0 & \color{gray}0 \\
\hline
\color{gray}0   & 1 & \color{gray}0   & \color{gray}0 & \color{gray}0 & \color{gray}0   & \color{gray}0 & \color{gray}0 & \color{gray}0 & \color{gray}0   & \color{gray}0 & \color{gray}0 & \color{gray}0 & \color{gray}0 & \color{gray}0 \\
-\frac{4}{5}   & \color{gray}0 & 1   & \color{gray}0 & \color{gray}0 & \color{gray}0   & \color{gray}0 & \color{gray}0 & \color{gray}0 & \color{gray}0   & \color{gray}0 & \color{gray}0 & \color{gray}0 & \color{gray}0 & \color{gray}0 \\
\hline
\color{gray}0   & \color{gray}0 & \color{gray}0   & 1 & \color{gray}0 & \color{gray}0   & \color{gray}0 & \color{gray}0 & \color{gray}0 & \color{gray}0   & \color{gray}0 & \color{gray}0 & \color{gray}0 & \color{gray}0 & \color{gray}0 \\
\color{gray}0   & -\frac{3}{4} & \color{gray}0   & \color{gray}0 & 1 & \color{gray}0   & \color{gray}0 & \color{gray}0 & \color{gray}0 & \color{gray}0   & \color{gray}0 & \color{gray}0 & \color{gray}0 & \color{gray}0 & \color{gray}0 \\
\color{gray}0   & \color{gray}0 & -\frac{3}{4}   & \color{gray}0 & \color{gray}0 & 1   & \color{gray}0 & \color{gray}0 & \color{gray}0 & \color{gray}0   & \color{gray}0 & \color{gray}0 & \color{gray}0 & \color{gray}0 & \color{gray}0 \\
\hline
\color{gray}0   & \color{gray}0 & \color{gray}0   & \color{gray}0 & \color{gray}0 & \color{gray}0   & 1 & \color{gray}0 & \color{gray}0 & \color{gray}0   & \color{gray}0 & \color{gray}0 & \color{gray}0 & \color{gray}0 & \color{gray}0 \\
\color{gray}0   & \color{gray}0 & \color{gray}0   & -\frac{2}{3} & \color{gray}0 & \color{gray}0   & \color{gray}0 & 1 & \color{gray}0 & \color{gray}0   & \color{gray}0 & \color{gray}0 & \color{gray}0 & \color{gray}0 & \color{gray}0 \\
\color{gray}0   & \color{gray}0 & \color{gray}0   & \color{gray}0 & -\frac{2}{3} & \color{gray}0   & \color{gray}0 & \color{gray}0 & 1 & \color{gray}0   & \color{gray}0 & \color{gray}0 & \color{gray}0 & \color{gray}0 & \color{gray}0 \\
\color{gray}0   & \color{gray}0 & \color{gray}0   & \color{gray}0 & \color{gray}0 & -\frac{2}{3}   & \color{gray}0 & \color{gray}0 & \color{gray}0 & 1   & \color{gray}0 & \color{gray}0 & \color{gray}0 & \color{gray}0 & \color{gray}0 \\
\hline
\color{gray}0   & \color{gray}0 & \color{gray}0   & \color{gray}0 & \color{gray}0 & \color{gray}0   & \color{gray}0 & \color{gray}0 & \color{gray}0 & \color{gray}0   & 1 & \color{gray}0 & \color{gray}0 & \color{gray}0 & \color{gray}0 \\
\color{gray}0   & \color{gray}0 & \color{gray}0   & \color{gray}0 & \color{gray}0 & \color{gray}0   & -\frac{1}{2} & \color{gray}0 & \color{gray}0 & \color{gray}0   & \color{gray}0 & 1 & \color{gray}0 & \color{gray}0 & \color{gray}0 \\
\color{gray}0   & \color{gray}0 & \color{gray}0   & \color{gray}0 & \color{gray}0 & \color{gray}0   & \color{gray}0 & -\frac{1}{2} & \color{gray}0 & \color{gray}0   & \color{gray}0 & \color{gray}0 & 1 & \color{gray}0 & \color{gray}0 \\
\color{gray}0   & \color{gray}0 & \color{gray}0   & \color{gray}0 & \color{gray}0 & \color{gray}0   & \color{gray}0 & \color{gray}0 & -\frac{1}{2} & \color{gray}0   & \color{gray}0 & \color{gray}0 & \color{gray}0 & 1 & \color{gray}0 \\
\color{gray}0   & \color{gray}0 & \color{gray}0   & \color{gray}0 & \color{gray}0 & \color{gray}0   & \color{gray}0 & \color{gray}0 & \color{gray}0 & -\frac{1}{2}   & \color{gray}0 & \color{gray}0 & \color{gray}0 & \color{gray}0 & 1 \\
\end{array}
\right ]
\]

The matrix $U B_5 U^{T}$ is given by:

\[ [5] \oplus \frac{1}{5} \left [ \begin{matrix} \color{green}60 & \color{green}20 \\ \color{green}20 & \color{purple}24\end{matrix}\right ]
\oplus \frac{1}{4} \left [ \begin{matrix} \color{blue}60 & \color{blue}36 & \color{blue}12 \\ \color{blue}36 & \color{green}45 & \color{green}15 \\ \color{blue}12 & \color{green}15 & \color{purple}18 \end{matrix}\right ]
\oplus \frac{1}{3} \left [ \begin{matrix} \color{red}42 & \color{red}30 & \color{red}18 & \color{red}6 \\ \color{red}30 & \color{blue}40 & \color{blue}24 & \color{blue}8 \\ \color{red}18 & \color{blue}24 & \color{green}30 & \color{green}10 \\ \color{red}6 & \color{blue}8 & \color{green}10 & \color{purple}12\end{matrix}\right ]
\oplus \frac{1}{2} \left [ \begin{matrix} 18 & 14 & 10 & 6 & 2 \\ 14 & \color{red}21 & \color{red}15 & \color{red}9 & \color{red}3 \\ 10 & \color{red}15 & \color{blue}20 & \color{blue}12 & \color{blue}4 \\ 6 & \color{red}9 & \color{blue}12 & \color{green}15  & \color{green}5\\ 2 & \color{red}3 & \color{blue}4 & \color{green}5& \color{purple}6\end{matrix}\right ]
\] 

The matrix $A_5$ is given by:
\[
\left[\begin{array}{r|rr|rrr|rrrr|rrrrr}
5 & 1 & -3 & 0 & 0 & 0 & 0 & 0 & 0 &
0 & 0 & 0 & 0 & 0 & 0 \\
\hline
1 & 3 & -2 & 1 & -2 & 1 & 0 & 0 & 0
& 0 & 0 & 0 & 0 & 0 & 0 \\
-3 & -2 & 6 & 0 & 1 & -3 & 0 & 0 & 0
& 0 & 0 & 0 & 0 & 0 & 0 \\
\hline
0 & 1 & 0 & 3 & -2 & 0 & 1 & -2 & 1
& 0 & 0 & 0 & 0 & 0 & 0 \\
0 & -2 & 1 & -2 & 4 & -2 & 0 & 1 & -2
& 1 & 0 & 0 & 0 & 0 & 0 \\
0 & 1 & -3 & 0 & -2 & 6 & 0 & 0 & 1
& -3 & 0 & 0 & 0 & 0 & 0 \\
\hline
0 & 0 & 0 & 1 & 0 & 0 & 3 & -2 & 0 &
0 & 1 & -2 & 1 & 0 & 0 \\
0 & 0 & 0 & -2 & 1 & 0 & -2 & 4 & -2
& 0 & 0 & 1 & -2 & 1 & 0 \\
0 & 0 & 0 & 1 & -2 & 1 & 0 & -2 & 4
& -2 & 0 & 0 & 1 & -2 & 1 \\
0 & 0 & 0 & 0 & 1 & -3 & 0 & 0 & -2
& 6 & 0 & 0 & 0 & 1 & -3 \\
\hline
0 & 0 & 0 & 0 & 0 & 0 & 1 & 0 & 0 &
0 & 3 & -2 & 0 & 0 & 0 \\
0 & 0 & 0 & 0 & 0 & 0 & -2 & 1 & 0 &
0 & -2 & 4 & -2 & 0 & 0 \\
0 & 0 & 0 & 0 & 0 & 0 & 1 & -2 & 1 &
0 & 0 & -2 & 4 & -2 & 0 \\
0 & 0 & 0 & 0 & 0 & 0 & 0 & 1 & -2 &
1 & 0 & 0 & -2 & 4 & -2 \\
0 & 0 & 0 & 0 & 0 & 0 & 0 & 0 & 1 &
-3 & 0 & 0 & 0 & -2 & 6
\end{array}\right]
\]
The matrix $V = U^{-1}$ is given by:
\[
U^{-1} = 
\left[\begin{array}{r|rr|rrr|rrrr|rrrrr}
1 & 0 & 0 & 0 & 0 & 0 & 0 & 0 & 0 &
0 & 0 & 0 & 0 & 0 & 0 \\
\hline
0 & 1 & 0 & 0 & 0 & 0 & 0 & 0 & 0 &
0 & 0 & 0 & 0 & 0 & 0 \\
\frac{4}{5} & 0 & 1 & 0 & 0 & 0 & 0 & 0
& 0 & 0 & 0 & 0 & 0 & 0 & 0 \\
\hline
0 & 0 & 0 & 1 & 0 & 0 & 0 & 0 & 0 &
0 & 0 & 0 & 0 & 0 & 0 \\
0 & \frac{3}{4} & 0 & 0 & 1 & 0 & 0 & 0
& 0 & 0 & 0 & 0 & 0 & 0 & 0 \\
\frac{3}{5} & 0 & \frac{3}{4} & 0 & 0 & 1 & 0
& 0 & 0 & 0 & 0 & 0 & 0 & 0 & 0 \\
\hline
0 & 0 & 0 & 0 & 0 & 0 & 1 & 0 & 0 &
0 & 0 & 0 & 0 & 0 & 0 \\
0 & 0 & 0 & \frac{2}{3} & 0 & 0 & 0 & 1
& 0 & 0 & 0 & 0 & 0 & 0 & 0 \\
0 & \frac{1}{2} & 0 & 0 & \frac{2}{3} & 0 & 0
& 0 & 1 & 0 & 0 & 0 & 0 & 0 & 0 \\
\frac{2}{5} & 0 & \frac{1}{2} & 0 & 0 & \frac{2}{3}
& 0 & 0 & 0 & 1 & 0 & 0 & 0 & 0 & 0
\\
\hline
0 & 0 & 0 & 0 & 0 & 0 & 0 & 0 & 0 &
0 & 1 & 0 & 0 & 0 & 0 \\
0 & 0 & 0 & 0 & 0 & 0 & \frac{1}{2} & 0
& 0 & 0 & 0 & 1 & 0 & 0 & 0 \\
0 & 0 & 0 & \frac{1}{3} & 0 & 0 & 0 &
\frac{1}{2} & 0 & 0 & 0 & 0 & 1 & 0 & 0 \\
0 & \frac{1}{4} & 0 & 0 & \frac{1}{3} & 0 & 0
& 0 & \frac{1}{2} & 0 & 0 & 0 & 0 & 1 &
0 \\
\frac{1}{5} & 0 & \frac{1}{4} & 0 & 0 & \frac{1}{3}
& 0 & 0 & 0 & \frac{1}{2} & 0 & 0 & 0 &
0 & 1
\end{array}\right]
\]

The matrix $V^T A_5 V$ is given by:
\[
\frac{1}{5} [{\color{Emerald}13}] \oplus \frac{1}{4} \left [ \begin{matrix} \color{green}6 & \color{blue}-5 \\ \color{blue}-5 & \color{purple}15 \end{matrix}\right ]
\oplus \frac{1}{3} \left [ \begin{matrix} \color{green}5 & \color{blue}-4 & 0 \\ \color{blue}-4 & \color{red}8 & \color{blue}-4 \\ 0 & \color{blue}-4 & \color{purple}12 \end{matrix}\right ]
\oplus \frac{1}{2} \left [ \begin{matrix} \color{green}4 & \color{blue}-3 & 0 & 0 \\ \color{blue}-3 & \color{red}6 & \color{blue}-3 & 0 \\ 0 & \color{blue}-3 & \color{red}6 & \color{blue}-3 \\ 0 & 0 & \color{blue}-3 & \color{purple}9\end{matrix}\right ]
\oplus \left [ \begin{matrix} \color{green}3 & \color{blue}-2 & 0 & 0 & 0 \\ \color{blue}-2 & \color{red}4 & \color{blue}-2 & 0 & 0 \\ 0 & \color{blue}-2 & \color{red}4 & \color{blue}-2 & 0 \\ 0 & 0 & \color{blue}-2 & \color{red}4 & \color{blue}-2  \\ 0 & 0 & 0 & \color{blue}-2&\color{purple}6 \end{matrix}\right ]
\]

Note that:
\[ {\color{Emerald}13} = {\color{green}7} + {\color{purple}18}- {\color{blue} 6} -{\color{blue}6} \]

These examples suggest that there is some rich underlying structure which is worth further pursuit.

\section{Conclusion}

We have given a elementary proof that semicircular algebra is a factor in the sense of Von Neumann algebras. Although this result is already known, our proof is particularly simple. We have also given an explicit formula for the projection of any element onto the dense subspace of commutators. 
It would be interesting to see if this approach could be extended to the case of the $q$-semicircular algebra \cite{speicher}, or even the $(q,t)$-semicircular algebra \cite{natasha}. 

Examining the Gramm matrix of the most natural basis for the space of commutators, we have found it to have an intriguing self-similar structure which may be interesting to investigate further.

\cleardoublepage

\backmatter
\lhead[\oldstylenums \thepage]{\rightmark}
\rhead[\leftmark]{\oldstylenums \thepage}


\lhead[\oldstylenums \thepage]{Bibliographie}
\rhead[Bibliographie]{\oldstylenums \thepage}
\bibliographystyle{alpha}
\bibliography{references}
\cleardoublepage

\lhead[\oldstylenums \thepage]{Index}
\rhead[Index]{\oldstylenums \thepage}
\printindex
\cleardoublepage



\end{document}